\documentclass[10pt]{article}

\usepackage{amssymb}
\usepackage{amsmath}
\usepackage{amsthm}
\usepackage{amsfonts}
\usepackage{xcolor}
\usepackage{mathbbol}
\usepackage{bm}
\usepackage{hyperref}
\usepackage{upgreek}
\usepackage{cleveref}
\usepackage{dsfont}
\usepackage{mathrsfs}
\usepackage{times}
\usepackage{graphicx}
\usepackage{bm}
\usepackage{mathtools}

\usepackage{dsfont}
\RequirePackage[algo2e]{algorithm2e}
\RequirePackage{enumitem}
\RequirePackage{algorithm}
\usepackage{titling}
\usepackage{tabularx}
\usepackage{accents}
\usepackage[round]{natbib}
\newcommand{\E}[1]{\mathbb{E}\left[#1\right]}
\newcommand{\Enb}[1]{\mathbb{E}#1}

\newcommand{\N}{\mathbb{N}}
\newcommand{\R}{\mathbb{R}}
\newcommand{\prob}[1]{\mathbb{P}(#1)}

\newcommand{\abs}[1]{\left | #1 \right |}

\newcommand{\norm}[1]{\|#1\|}

\newcommand{\cind}{\overset{\mathrm{d}}{\rightarrow}}
\newcommand{\cinp}{\xrightarrow[]{P}}

 \renewcommand{\bm}[1]{\boldsymbol{#1}}

\newcommand{\longconditional}{\hspace{1mm}\middle \vert \hspace{1mm}}

\newcommand{\idummy}{\ell}
\newcommand{\RNum}[1]{\uppercase\expandafter{\romannumeral #1\relax}}

\newcommand{\QexcludingVj}{}
\newcommand{\QexcludingVi}{}

\newcommand{\QexcludingVdummy}{}
\newcommand{\QexcludingVone}{}
\newcommand{\QexcludingVperturbi}{^{i}}
\newcommand{\QexcludingVoneperturb}{^{i}}
\newcommand{\QexcludingVdummyperturbj}{^{j}}
\newcommand{\QexcludingViperturbdummy}{^{\ell}}
\newcommand{\QexcludingV}{}

\newcommand{\dexcludingVi}{}

\newcommand{\DexcludingV}{}
\newcommand{\DincludingV}{^{[v]}}
\newcommand{\DexcludingVi}{}
\newcommand{\DexcludingVone}{}
\newcommand{\DexcludingVdummy}{}
\newcommand{\DexcludingVdummyperturbj}{^{j}}
\newcommand{\DexcludingViperturbdummy}{^{\ell}}

\newcommand{\KexcludingV}{}
\newcommand{\KexcludingW}{}
\newcommand{\KexcludingVone}{}
\newcommand{\KincludingV}{^{[v]}}
\newcommand{\KexcludingViperturbdummy}{^{\ell}}

\newcommand{\intermediatevarianceone}{\hat{\sigma}^2_{r, 1}}
\newcommand{\intermediatevariancetwo}{\hat{\sigma}^2_{r, 2}}

\newcommand{\wscenterr}{_{r,s}}

\usepackage[utf8]{inputenc}
\usepackage{setspace}
\usepackage[margin=1in]{geometry}

\DeclareMathOperator*{\argmin}{argmin}
\title{Winners with Confidence: \\
Discrete Argmin Inference with an Application to Model Selection}
\author{
    Tianyu Zhang\textsuperscript{1}, Hao Lee\textsuperscript{2}, and Jing Lei\textsuperscript{2} \\
    \textsuperscript{1}Department of Statistics and Applied Probability, University of California, Santa Barbara \\
    \textsuperscript{2}Department of Statistics \& Data Science, Carnegie Mellon University
}

 \date{}
\usepackage[page]{appendix}
\newtheorem{theorem}{Theorem}[section]
\newtheorem{lemma}[theorem]{Lemma}
\newtheorem{corollary}[theorem]{Corollary}
\newtheorem{proposition}[theorem]{Proposition}
\newtheorem{definition}[theorem]{Definition}

\newtheorem{remark}[theorem]{Remark}
\begin{document}

\maketitle
\vspace{-7mm}
\textbf{Abstract} We study the problem of finding the index of the minimum value of a vector from noisy observations. This problem is relevant in population/policy comparison, discrete maximum likelihood, and model selection. We develop an asymptotically normal test statistic, even in high-dimensional settings and with potentially many ties in the population mean vector, by integrating concepts and tools from cross-validation and differential privacy. The key technical ingredient is a central limit theorem for globally dependent data.  We also propose practical ways to select the tuning parameter that adapts to the signal landscape. Numerical experiments and data examples demonstrate the ability of the proposed method to achieve a favorable bias-variance trade-off in practical scenarios.

\section{Introduction} \label{intro}

Let $X_1,..., X_n$ be independent and identically distributed (IID) random vectors in $\mathbb R^p$ with the mean vector $\mu=\mathbb E X_1$.
We are interested in finding the index set of the minimum entries of $\mu$:
\begin{equation}\label{eq: min-mean set}
\Theta=\{r\in [p]\mid \mu_r \leq \mu_s \text{ for all }s\in [p]\},
.\end{equation}
where $[p] = \{1,...,p\}$.

Uncertainty quantification in the inference of $\Theta$ is naturally motivated by various real-world problems regarding best or optimal choices. One example is the prediction of election outcomes and the analysis of polling data. When multiple candidates are competing for a single position, we model each voter's preference as a binary vector $X_i$ where $X_{i,r} = 1$ indicates a vote for candidate $r$ from voter $i$. Constructing a confidence set for the candidate(s) with the highest population support rate allows us to directly forecast the most likely winners, accounting for voter randomness and variability. As discussed in \cite{xie2009confidence, Hung2019rank,mogstad2024inference}, similar data types appear in social science, institution evaluation, clinical trial and education (see our real-data example in \Cref{section: real_data}). Confidence sets that acknowledge the limitations of data are considered highly important in practice \citep{goldstein1996league}.  

One may also consider comparing the performance of $p$ agents at a task of interest. Given an environment random variable $W\sim P_W$, the performance of the agents under this specific environment is quantified as $\ell(f_1, W),...,$ $\ell(f_p, W)$, where $f_1,...,f_p$ are the agents and $\ell$ is a pre-specified loss function. For regression tasks, the environment variable is a pair of predictors and an outcome of interest. In this case, $W = (Z, Y)$, $f_r(Z)$ is an estimate of $\mathbb{E}[Y\mid Z]$ and $\ell(f_r,W) = (f_r(Z) - Y)^2$. It is often of interest to identify the agent that performs the best on average (over the randomness from the environment), that is, to identify $\argmin_{r\in [p]} \mathbb{E}[\ell(f_r, W)]$. In the notation of \eqref{eq: min-mean set}, $X_{i,r} = \ell(f_r, W_i)$, given some sampled environments $W_1,..., W_n$. Methods that offer users a confident set of best-performing agents can help examine the robustness of the decision. Moreover, rather than a single estimated best performer subject to the insufficiency of the data, the users are theoretically justified to choose agents in the confidence set that offer better computational properties, enhanced interpretability, or greater financial feasibility.

While in some applications $\Theta$ is a singleton set $ \{r^\star\}$, in this work we consider the general situation where there may be arbitrarily many tied values in $\mu$, at the minimum and/or other values. For example, we allow for the vector $\mu$ to have constant entries: $\mu_1=\mu_2= \cdots =\mu_p$, in which case $\Theta=[p]$. We also allow the number of coordinates, $p$, to be comparable or larger than the sample size $n$.

The index of minimal \emph{empirical} mean is a natural choice as a point estimate of $r^*$. However, it implicitly assumes a unique minimum index in $\mu$. When all entries of $\mu$ are the same and $X$ has a continuous distribution, the empirical argmin will only return a single coordinate, missing all other $p-1$ coordinates. Therefore, to quantify the uncertainty in estimating the location of the minimum, we aim to construct a confidence set $\widehat C\subset [p]$ that accounts for the variability in the data. Formally, we require the confidence set $\widehat C$ to satisfy
\begin{equation}\label{eq: marginal validity}
     \lim_{n\rightarrow\infty}\mathbb{P}(r \in \widehat C) \geq 1-\alpha\,,~~\forall~r\in \Theta\,,
\end{equation}
where $\alpha$ is a given significance level (commonly $=0.05$).

In this work, we develop a novel method to construct confidence sets of the argmin index set $\Theta$ that asymptotically satisfy \eqref{eq: marginal validity}.  The idea is to compare each index with ``the best of others".  Intuitively, in order to decide whether $r\in \widehat{C}$ for a specific $r\in[p]$, we only need to test $\mu_r\le\mu_{s_r}$ for some $s_r\in \{s\neq r \mid \mu_s = \min_{t\neq r}\mu_t\}$. However, the best index excluding $r$, $s_r$, is not available and must be adaptively estimated from the data.  If we use an empirical version of such an $s_r$, then there will be a double-dipping issue (which is also widely discussed in the post-selection inference literature). As a main methodological contribution, we employ a combination of cross-validation and exponential mechanism, a technique originated from the differential privacy literature  \citep{dwork2014algorithmic}, which is known to limit the dependence between nuisance parameter estimates and the final inference. Our theoretical analysis relies on a central limit theorem for globally dependent data, which may be of general interest.  To the best of our knowledge, this is the only method that uses asymptotic normality for the argmin inference of a discrete random vector. More importantly, our method comes with an intuitive and simple way to choose the tuning parameter in a data-driven manner.

The coverage guarantee described in Equation~\ref{eq: marginal validity} is \textit{marginal}. It ensures that each $r \in \Theta$ is included in $\widehat{C}$ with high probability, individually. Two alternatives are (asymptotic) \textit{non-empty} coverage, defined as $\lim_{n\rightarrow\infty}\mathbb{P}(\Theta \cap \widehat{C} \neq \emptyset) \ge 1 - \alpha$, and (asymptotic) \textit{simultaneous} coverage, defined as $\lim_{n\rightarrow\infty}\mathbb{P}(\Theta \subseteq \widehat{C}) \ge 1 - \alpha$. Non-empty coverage is implied by marginal coverage~\eqref{eq: marginal validity}, which is in turn implied by simultaneous coverage. The three notions become equivalent when $|\Theta|=1$. All coverage modes have been discussed in the literature, and the appropriate choice depends on the application. In general, stronger modes of coverage increase the likelihood that $\widehat{C}$ includes dimensions outside $\Theta$. We will elaborate on these trade-offs and provide supporting evidence as the paper develops.

\paragraph{Related work.} 

Inference of argmin indices has a long history in the statistical literature,  dating back to the early works of \cite{gibbons1977selecting,gupta1979multiple}. A refinement was proposed in \cite{futschik1995confidence}, assuming known marginal distributions of $X_{i,r}$ and independence between dimensions of $X_i$. Methods that are strongly dependent on these conditions are theoretically valid but may have restricted applicability.
\cite{mogstad2024inference} developed a general framework—valid under broad distributional assumptions—for constructing both marginal and simultaneous confidence sets, based on pairwise comparisons between the entries of $\mu$. 
Variants of bootstrap methods for the argmin inference are also available in the model selection setting \citep{hansen2011model}. Their proposed Model Confidence Set (MCS) is guaranteed to achieve simultaneous coverage asymptotically with fixed $p$. 
However, the standard implementation \cite{bernardi2018model} of it is computationally demanding and may not yield satisfactory power as a trade-off. 
\cite{dey2024anytime} proposes a martingale- and e-value-based method for constructing argmin confidence sets, achieving finite sample non-empty coverage $\mathbb P(\Theta\cap\widehat C \neq \emptyset)\ge 1-\alpha$ for any sample size $n$ (the current work focuses on asymptotic marginal coverage). In a recent work \citep{kim2025locally}, the authors developed a sample-splitting procedure for argmin inference with asymptotic marginal coverage and included a comparison with the proposal of this work.

Although the main focus of this work is on the stationary setting with an IID learning environment, sequential inferential frameworks of the best performer have also been studied in the literature. For example, a recent work \cite{arnold2024sequential} builds on the principle in  \cite{hansen2011model} and develops an online argmin inference procedure. In \cite{chen2023inference}, the authors investigate uncertainty quantification for optimal policy estimation problems. This earlier work employed a similar (softmax) exponential mechanism as in the current work, leading to asymptotically normal statistics as well. Ideas related to softmax and stable statistics are also studied in the recent work \cite{adrian2024stabilizing}.

The argmin inference problem can be treated as a dual problem of \emph{rank inference/verification}. In rank verification, the parameter of interest is the rank of an index $r$: $R_r=1+\sum_{s\neq r}\mathds{1}(\mu_r<\mu_s)$, and the inference task is to establish confidence set $\widehat C_r$ such that $\mathbb P(R_r\in \widehat C_r)\ge 1-\alpha$ \citep{goldstein1996league, hall2009, xie2009confidence, Hung2019rank, mogstad2024inference, fan2024ranking}. Although it is conceptually possible to construct argmin index confidence sets from corresponding rank confidence sets, many rank verification methods (e.g. \cite{Hung2019rank}) would perform poorly or degenerate when the cardinality of $\Theta$ is greater than $1$, making it hard to transfer them to the argmin inference setting where ties or near-ties are prevalent. 

The study of the argmin index is also central to discrete stochastic optimization \citep{kleywegt2002sample}, in which discrete Maximum Likelihood Estimation (MLE) is a subbranch most relevant to statistics \citep{Choirat2012estimation, seri2021model}. Unlike standard MLE where the parameter of interest is allowed to take values in a continuous open set, some applications only permit integer-valued parameters due to natural constraints. For example, in astrophysics the parameter of interest may be the number of planets in a star system. Unlike in the continuous case, results on confidence sets for discrete MLE are scarce due to the irregularity of the problem \citep[see][and references therein]{Choirat2012estimation}.

The problem of argmin confidence sets can also be approached using methods in post-selection inference \citep[PoSI][]{taylor2015statistical}, or selective inference (SI) due to its multiple comparison nature. PoSI/SI methods usually require known and easy-to-compute noise distributions (such as isotropic Gaussian), which are impractical in most natural argmin inference scenarios.  In practice, we also find PoSI/SI-based methods less powerful compared to other alternatives.

\paragraph{Notation} Denote the integer set $\{1,...,p\}$ as $[p]$. We will use $V$ to denote the number of folds in $V$-fold cross-validation and assume $n/V$ is an integer. Without loss of generality, we will split the samples into folds sequentially and use the index-set notations $I_v = \{(v-1)n/V + 1, ...,vn/V\}$ and $I_v^c = [n]\symbol{92} I_v$. Given a sample index $i$, the notation $v_i$ maps it to the fold-index that sample $i$ belongs to: $i\in I_{v_i}$. The symbol $\bm{X}$ denotes the whole data set $\{X_i: i \in [n]\}$, and $\bm{X}^{(-v)} = \{X_i: i\notin I_v\}$. The notation $\bm{X}^j$ denotes the sample $\bm{X}$ but replaces $X_j\in\bm{X}$ by an IID copy $X_j^\prime$ while keeping everything else intact. Similarly, $\bm{X}^{j, k}$ replaces $X_j,X_k$ by the same IID copies $X_j^\prime, X_k^\prime$ as in $\bm{X}^{j}, \bm{X}^{k}$. That is, $\bm{X}^{j,k}$ differ from $\bm{X}^k$ or $\bm{X}^j$ by only one sample. 
For two positive sequences $a_n,b_n$, $a_n = o(b_n)$ means $\lim_n a_n/b_n = 0$, and $a_n = \omega(b_n)$ means $b_n=o(a_n)$. For a set $\mathbb{C}$, $|\mathbb{C}|$ returns its cardinality. We will use $\Phi(x)$ to denote the cumulative distribution function of the standard normal.

\section{Methods}
We propose a cross-validated exponential weighting scheme to construct a confidence set for the argmin index. The procedure is formally presented in Algorithm~\ref{algorithm: exp weighting}. Additionally, this section provides some intuition behind the approach, along with other simpler proposals one might consider (but do not share the favorable properties).

\subsection{Reduction to a Selective Mean-testing Problem}
The coverage requirement \eqref{eq: marginal validity} is marginal for each individual index $r \in \Theta$. Therefore, we can focus on each $r\in[p]$ to decide whether it is in the argmin set.
We start from the simple observation that $r\in \Theta$ if and only if $\mu_r\le \min_{s\neq r}\mu_s$.  Let $s_r$ be an index in $[p]\backslash\{r\}$ such that $\mu_{s_r}=\min_{s\neq r}\mu_s$. Therefore, we have $r\in\Theta$ if and only if $\mu_r\le \mu_{s_r}$. If we know the value of $\mu_{s_r}$ (or $s_r$, resp.), then the decision of whether $r \in \widehat{C}$ can be made by a simple one-sample (or two-sample, resp.) one-sided $t$-test.


However, we have access to neither $\mu_{s_r}$ nor $s_r$ in practice, and any one of them has to be inferred from the noisy data. Suppose we instead use $ \min_{s \neq r} \hat \mu_s$ as an estimate of $\mu_{s_r}$ with $\hat\mu \in \mathbb{R}^p$ being the empirical mean vector, the constructed $\widehat C$ would not have the desired coverage due to the well-known ``double-dipping'' or ``selective inference'' issue \citep{taylor2015statistical}. To illustrate it, consider the case when all the dimensions of $X$ are independent standard normal so that $\Theta = [p]$. The minimal sample mean, $\min_{s\neq r}\hat \mu_s$, is related to the Gumbel distribution with expectation
$\asymp -\sqrt{\log p/n}$. Combined with recent finite-sample concentration inequalities \citep[Theorem 3]{tanguy2015some}, we know with high probability, $ \min_{s\neq r} \hat \mu_s$ is less than $-c\sqrt{\log p/n}$ with a constant $c > 0$. When $p$ is large, it is much smaller than the sample mean $\hat \mu_r =O_P (1/\sqrt{n})$. Naively comparing $\min_{s\neq r} \hat \mu_s$ with $\hat \mu_r$ using standard two-sample tests would falsely reject the hypothesis $\mu_r \leq \mu_{s_r}$ and almost always exclude $r$ from the confidence set even when $r\in \Theta$.

\subsection{Initial Fix: Removing Dependence by Cross-validation}\label{section: simple split}

To avoid the ``double-dipping" bias, one may consider a cross-validation type of scheme using a part of the data to obtain an estimate $\hat s_r$ and compare $\mu_r$ with $\mu_{\hat s_r}$ on the left-out sample point(s), aggregating the two-sample comparison by rotating the left-out set. 
More concretely, consider a leave-one-out (LOO) version of this idea. For each sample point $i\in[n]$, define the $i$-th LOO argmin index (without $r$) 
\begin{equation}\label{eq:loo-arg-min}
\hat s_r^{(-i)} := \argmin_{s\neq r} \hat{\mu}_s^{(-i)}, \hspace{5mm} \hat{\mu}^{(-i)}_s = \frac{1}{n 
- 1}\sum_{j\neq i} X_{j,s}\,,
\end{equation}
with any arbitrary tie-breaking rules.
Now by construction $\hat s_r^{(-i)}$ and $X_i$ are independent, and one would expect the cross-validation-type statistic
\begin{equation}\label{eq:simple-cv-sum}
\frac{1}{\sqrt{n}}\sum_{i=1}^n \left(X_{i,r}-X_{i,\hat s_r^{(-i)}}\right)
\end{equation}
may be asymptotically normal (after being properly centered).

\begin{figure}[!htbp]
    \centering
    \includegraphics[width =0.8\linewidth]{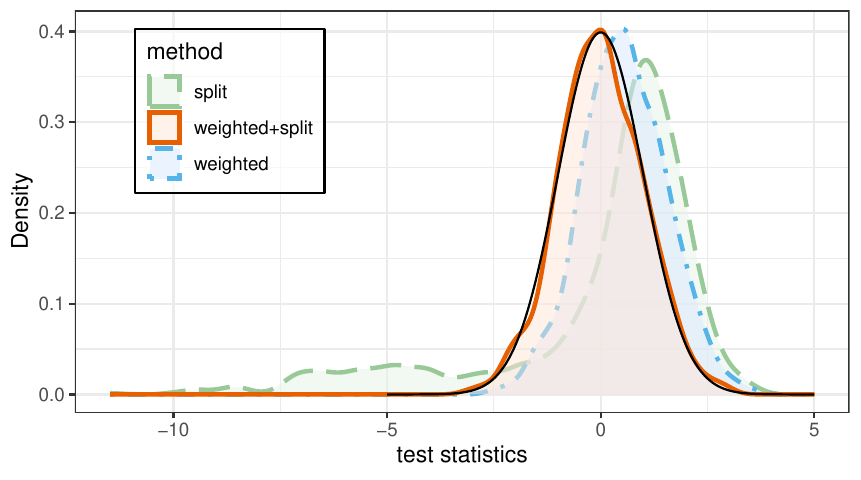}
    \caption{Sample splitting and exponential weighting are both crucial for normality. Smoothed histograms of the normalized $T_r$ in \Cref{algorithm: exp weighting} and its related variants. We take $r = 1$. {\tt weighted+split} is the normalized $T_r$ presented in \Cref{algorithm: exp weighting}; {\tt split} is described in Section~\ref{section: simple split}, and {\tt weighted} is the non-split version of {\tt weighted+split}, discussed in Remark~\ref{remark: softmin nosplit}. The solid black line is the density curve of the standard normal. A LOO ($V = n$) sample-splitting scheme is employed in {\tt split} and {\tt weighted+split}.}
    \label{fig: motivation}
\end{figure}

Unfortunately, this is not the case in general.
In Figure~\ref{fig: motivation}, we demonstrate this phenomenon with dimension $p = 50$ and a sample size $ n= 100$. It is the all-tie case $\Theta = [p]$ with IID standard normal coordinates of $X$. The empirical distribution is obtained from $10^3$ repeats. We observe the simple LOO method {\tt split} is left-skewed with a visible irregular tail. Note that the test statistics produced by the LOO method are more dispersed than normal on both tails, which may negatively impact both the type I error control and power.

\begin{algorithm}[!htbp]
\SetAlgoLined
\SetKwInOut{Input}{Input}
\SetKwInOut{Output}{Output}
\SetKwInOut{Initialize}{Initialize}
\vspace{3pt}

\Input{
Data $\bm{X}$, the number of folds $V$, a significance level $\alpha$, and a weighting parameter $\lambda$.
}
\Output{
Confidence set $\widehat C$.
}
\Initialize{
The confidence set $\widehat C = \varnothing$
}
\vspace{3pt}
\For{dimension index $r$ in $[p]$}{
\For{fold index $v$ in $[V]$}{
Compute the sample mean from the out-sample data: $\hat{\mu}^{(-v)}:=\left(n-n/V\right)^{-1} \sum_{j \in I_v^c} X_j$\\
\For{sample index $i$ in $I_v$}{
Calculate weighted competitor $Q_{i, r}\QexcludingV =\sum_{s \neq r}  \hat w_{r,s}^{(-v)} X_{i, s}$ with weights $\left( \hat w_{r,s}^{(-v)}, s \neq r\right)$ satisfying
\begin{equation}\label{eq:exp_weight}
\sum_{s \neq r} \hat w_{r,s}^{(-v)}=1 \text { and }  \hat w_{r,s}^{(-v)} \propto \exp \left(-\lambda\hat{\mu}_{s}^{(-v)}\right).
\end{equation}
}
}
Calculate test statistics $T_r$ with an estimated standard deviation $\hat\sigma_r >0 $, e.g. \Cref{eq: variance out}, 
\begin{equation*}
T_r := \frac{1}{\sqrt{n} \hat\sigma_r} \sum_{i=1}^n\left(X_{i, r}-Q_{i, r}\QexcludingVi \right)\,.
\end{equation*}

Add $r$ to $\widehat{C}$ if $T_r<z_{1-\alpha}$, where $z_{1-\alpha}$ is the $(1-\alpha)$-quantile of standard normal.
}
\caption{Proposal: Exponentially Weighted Argmin Confidence Set}
\label{algorithm: exp weighting}
\end{algorithm}

\subsection{Final Fix: Cross-validated Exponential Mechanism}\label{section: final fix}
The failure of asymptotic normality for the statistic in \eqref{eq:simple-cv-sum} is related to the perplexed distribution of general cross-validation-type statistics. Some recent work in the cross-validation literature imposed various ``stability conditions'' that are crucial to establishing asymptotic normality for cross-validated risks \citep{austern2020,kissel2023blackbox}. Roughly speaking, these sufficient conditions require the quantity $X_{i,\hat s_r^{(-i)}}$ to have a distortion in $L_2$ norm being much smaller than $1/\sqrt{n}$ when one sample in $\{X_j:j\neq i\}$ is replaced by an IID copy. This is not the case since a change of a single sample point may result in $\hat s_r^{(-i)}$ changing to a completely different value. Particularly, this can occur with a probability proportional to $1/\sqrt{n}$, and lead to a constant level change in $X_{i, \hat s_r^{(-i)}}$.

As a concrete example, we consider a simple case when $p = 3$, $\mu_s = 0$ for all $s \in [p]$ and the underlying samples are drawn from $X_1 \sim \mathcal{N}(0, I_3)$. Suppose that we are testing the first dimension $r=1$. When $\hat{\mu}_2^{(-i)} > \hat{\mu}_3^{(-i)}$, there is a $\sim 1/\sqrt{n}$ probability this ordering reverses if we replace a single $X_j, j\neq i$ by an IID copy of it (\Cref{lemma: sign_flipping}). Equivalently, when the selected coordinate is $\hat{s}_r^{(-i)}=2$, there is a $\sim$ $1 / \sqrt{n}$ chance that modifying just one observation causes $\hat{s}_r^{(-i)}$ to switch to 3, so that the selected variable $X_{i, \hat{s}_r^{(-i)}}$ changes from $X_{i, 2}$ to $X_{i, 3}$. Although this event vanishes with $n$, its impact is large enough to disrupt asymptotic normality.

Our solution to this lack of stability is inspired by the differential privacy literature \citep{dwork2014algorithmic}, where the distortion of a statistic under the perturbation of a single data entry is referred to as sensitivity. Many techniques have been developed to produce insensitive counterparts of standard statistics. For the argmin index, a differentially private version can be obtained by the Exponential Mechanism \citep{mcsherry2007mechanism}. The original exponential mechanism will randomly sample a single coordinate as the insensitive argmin. In our problem, it is more convenient to use a weighted average with the weights corresponding to the sampling probabilities in the exponential mechanism.  The resulting algorithm replaces
$X_{i,\hat s_r^{(-i)}}$ by a weighted average:
$$
Q_{i,r}\QexcludingVi := \sum_{s\neq r} \hat w_{r,s}^{(-i)} X_{i,s}
$$
with weights $\{\hat w_{r,s}^{(-i)} \geq 0, s\neq r\}$ satisfying
$$
\sum_{s \neq r} \hat{w}_{r, s}^{(-i)}=1 \text{ and } \hat w_{r,s}^{(-i)}\propto \exp\left(-\lambda \hat\mu_{s}^{(-i)}\right).
$$
Here $\hat\mu^{(-i)}=(n-1)^{-1}\sum_{j\neq i}X_{j}$ is the $i$-th LOO empirical mean and $\lambda\geq 0$ is a tuning parameter to be chosen by users. Our final algorithm for constructing a confidence set of the argmin indices is formally presented in Algorithm~\ref{algorithm: exp weighting}, which allows both LOO version ($V = n$) and finite fold scheme ($V = O(1)$). 

Instead of identifying one single dimension $X_{i, \hat{s}_r^{(-i)}}$ as the quantity to compare with, the competitor statistic $Q_{i,r}\QexcludingVi$ is a weighted sum of multiple competitive dimensions of $X_i$. 
The quantity $Q_{i,r}\QexcludingVi$ can be viewed as a cross-validated soft-min of the vector $(X_{i,s}:s\neq r)$, and is more stable than $X_{i,\hat s_r^{(-i)}}$ in the sense that perturbing any one sample point (other than the $i$-th one) can only perturb $Q_{i,r}\QexcludingVi$ by a small amount. 

A smaller $\lambda$ implies a smaller perturbation and stronger stability.
In one extreme case, $\lambda=0$ implies perfect stability as the weights do not depend on $\{X_j:j\neq i\}$. By contrast, a larger value of $\lambda$ will more effectively eliminate the contribution from dimensions with ``obviously" larger sample means, leading to smaller confidence sets. In the other extreme case, $\lambda=\infty$ results in $Q_{i,r}\QexcludingVi=X_{i,\hat s_r^{(-i)}}$. To achieve a worst-case valid inference procedure, our theoretical results suggest $\lambda$ needs to grow slower than $\sqrt{n}$. In practice, the best choice of $\lambda$ would be just small enough to ensure the stability required by asymptotic normality. We provide detailed theoretical and practical guidance in choosing the tuning parameter $\lambda$ in Sections~\ref{section: asymptotic normality} - \ref{section: data-driven}. 

\begin{remark}
{(Weighted differences and standardization)}
An equivalent way to write the test statistic in \eqref{eq:exp_weight} is
$$
T_r := \frac{1}{\sqrt{n} \hat{\sigma}_r} \sum_{i=1}^n\left(X_{i, r}-Q_{i, r}\QexcludingVi \right) = \frac{1}{\sqrt{n}\hat{\sigma}_r} \sum_{i=1}^n \sum_{s \ne r} \hat{w}^{(-v_i)}_{r,s}(X_{i,r} - X_{i,s}),
$$
with weights $\hat w_{r,s}^{(-v_i)}\propto \exp(-\lambda (\hat\mu_s^{(-v_i)}-\hat\mu_r^{(-v_i)}))$. In practice, it is beneficial to pre-standardize the differences $\{X_{i,r}-X_{i,s}: i \in [n]\}$ so that the algorithm remains invariant under constant scaling of the input data, and the comparison of mean values is easier to interpret when the variances are the same.
All simulation studies and real-data applications are conducted with this pre-processing.
\end{remark}

\begin{remark}\label{remark: softmin nosplit}
    The importance of exponential weighting is emphasized in Sections~\ref{section: simple split} and~\ref{section: final fix}. 
    A natural follow-up question is whether a method utilizing exponential weighting alone without sample splitting suffices for the argmin inference problem. From the stability perspective, such an approach would still achieve asymptotic normality, but the double-dipping issue would sustain, resulting in a positive mean shift and violation of coverage guarantee. We demonstrate the failure of this choice in Figure~\ref{fig: motivation}. 
    It shows both sample splitting and weighting are crucial ingredients to achieving valid argmin inference. 
\end{remark}

\section{Asymptotic Normality and Coverage}\label{section: asymptotic normality}
In this section, we show that for each $r$, the $T_r$ statistics are asymptotically normal under proper choices of $\lambda$, which further implies asymptotic coverage of the confidence set $\widehat C$ constructed in \Cref{algorithm: exp weighting}. The result is formally stated as Theorem~\ref{th: random center} below.  In this section, we present our results for general $V$-fold cross-validation schemes, where the number of folds $V$ can either grow with $n$ or remain constant.

\begin{theorem}\label{th: random center}
    Let $X_i \in \R^p, i \in [n]$, be IID samples with uniformly bounded entries: $\sup_{s\in [p]} \abs{X_{i,s}} \leq M$ almost surely for a constant $M$. The dimension $p$ can depend on $n$ so long as the assumptions below are satisfied. We further assume
    \begin{enumerate}
        \item The smallest eigenvalue of covariance matrix $\operatorname{Cov}(X_1)$ is bounded away from zero by a positive constant.
        \item The weighting parameter in Algorithm~\ref{algorithm: exp weighting} satisfies $\lambda=\lambda_n = o(\sqrt{n})$.
    \end{enumerate}
    Define the centered version of $T_r$:
    \begin{equation}\label{eq: centered test stat}
    \tilde{T}_r := \frac{1}{\sqrt{n} \sigma_r}\sum_{v=1}^V  \sum_{i \in I_v}\left(X_{i, r}-Q_{i, r}\QexcludingVi- d_{i,r}\dexcludingVi \right),
    \end{equation}
    where $\sigma_r^2 = \operatorname{Var}\left[X_{1, r}-Q_{1, r}\QexcludingVone \right]$, and 
    \begin{equation*}
        d_{i,r} = \mathbb{E}\left[X_{i, r}-Q_{i, r}\QexcludingVi \mid \boldsymbol{X}^{(-v_i)}\right] \,.
    \end{equation*}
Then for any $x\in\mathbb{R}$:
    \begin{equation}
       \lim _{n \rightarrow \infty} \max _{r \in [p]}\left|\mathbb{P}\left(\tilde{T}_r \leq x\right)-\Phi(x)\right|=0,
    \end{equation}
    where $\Phi(x)$ is the cumulative distribution function of the standard normal.
\end{theorem}

Key ingredients in the proof of Theorem~\ref{th: random center} are presented in Section~\ref{section: proof sketch}. The formal argument is provided in \Cref{app:thm3.1_cor3.4}. The argument is dissected into two steps: a general weakly dependent Central Limit Theorem (CLT) based on stability and verifying the stability conditions for the exponential weighting mechanism. 

\begin{remark}
    (Boundedness) In Theorem~\ref{th: random center} the critical growth rate of $\lambda=o(\sqrt{n})$ required for asymptotic normality does not depend on the dimensionality $p$. This property is a consequence of the entry-wise $L^\infty$-norm bound on the random vector $X_1$. We chose to present the bounded-case result for simplicity, and it can be directly generalized to the unbounded case. In \Cref{app: remove boundedness} we present the technical details with a formal result stated in \Cref{th: unbounded random center}. Specifically, when $\lambda= o(\sqrt{n}/\log^{A} p)$ for a fixed $A>0$ depending on the distribution of $X_1$, we can still guarantee $\tilde{T}_r$ is approximately normal. (For sub-Gaussian data, $A = 1/2$.) The key is deriving sufficient bounds for $\Delta_1$ and $\Delta_2$ without $L^\infty$ conditions on $X_1$: these results are presented in \Cref{lemma: moment bound K} and \Cref{lemma: moment bound second order K}.
\end{remark}

\begin{remark}\label{rem:worst-case-lambda}
    (Worst-case stability guarantee) We require $\lambda = o(\sqrt{n})$ in Theorem~\ref{th: random center} to guarantee the weights are stable enough to establish the asymptotic normality. This requirement applies to the worst case where many of the $\mu_r$'s are equal and can be relaxed if there is enough separation between leading coordinates. In Lemma~\ref{lemma: ultra stable}, we present the case when there is an ``obvious" winner within $[p]\symbol{92}\{r\}$ (that is, one dimension is significantly better than the others as a competitor of the $r$-th dimension). In this case, the condition of $\lambda$ can be substantially relaxed without compromising the asymptotic normality of $\tilde T_r$. In \Cref{section: data-driven}, we present an automatic data-driven procedure to detect such situations and choose $\lambda$ that best balances the validity and power.
\end{remark}

\begin{remark}\label{rem:random_centering}
    (Random centering) The center $d_{i,r}$ in Theorem~\ref{th: random center} is a random quantity depending on the data $\bm{X}^{(-v_i)}$.
It is straightforward to verify that
$$d_{i,r}\dexcludingVi =\sum_{s\neq r}\hat w_{r,s}^{(-v_i)}(\mu_r-\mu_s)\,,$$
 where $\hat w_{r,s}^{(-v_i)}$ is the exponential weight defined in \eqref{eq:exp_weight}.
As a result, when $r\in \Theta$ we have $d_{i,r}\le 0$ for all $i$. This simple but crucial fact bridges the gap between Theorem~\ref{th: random center} and the inference validity guarantee in Corollary~\ref{cor: coverage} below. The quantities $d_{i,r}$ can be viewed as the signal strength as it reflects the gap between $\mu_r$ and other $\mu_s$.  The coverage will be close to $1-\alpha$ if $n^{-1/2}\sigma_r^{-1}\sum_i d_{i,r}$ is close to $0$. 
\end{remark}

The coverage guarantee/validity result follows directly from the asymptotic normality in \Cref{th: random center} and \Cref{rem:random_centering}. The proof of the following result is in \Cref{app:thm3.1_cor3.4}.
\begin{corollary}\label{cor: coverage}
    Under the same assumptions as in Theorem~\ref{th: random center}, for each $r\in \Theta$, we have 
    \begin{equation*}
\lim_{n\rightarrow\infty}\mathbb{P}(r \in \widehat{C}) \ge 1-\alpha,
    \end{equation*}
    for the confidence set $\widehat C$ constructed in Algorithm~\ref{algorithm: exp weighting} with a consistent standard deviation estimate $\hat \sigma_r$. 
\end{corollary}

\begin{remark}
(Marginal coverage versus simultaneous coverage)\label{remark: marginal coverate} 
Our coverage guarantee is marginal for each $r\in\Theta$.
As discussed in \Cref{intro}, this validity is stronger than the non-empty coverage--- $\widehat C \cap \Theta \neq \varnothing$ with high probability, but does not imply the simultaneous coverage $\mathbb{P}(\Theta \subset \widehat{C})\ge 1-\alpha+o(1)$ \citep{hansen2011model, mogstad2024inference}. When $|\Theta| > 1$, these three modes of coverage are distinct, and the alternative guarantees can be of practical interest depending on the context. We chose to focus on the marginal coverage as it strikes a balance between power and validity, which is also a common choice in the literature \cite{futschik1995confidence, mogstad2024inference}.

Despite different validity goals, there is broad consensus that $\widehat{C}$ should avoid including irrelevant dimensions $s \notin \Theta$. Methods aiming for weaker modes yield higher rejection power. Some parts of our real-data/numerical analysis demonstrate this feature (\Cref{fig: data.application} and \Cref{fig: LOO vs MCS}). To provide a comprehensive comparison, we also report the finite-sample simultaneous coverage of our method in \Cref{fig: fail simultaneous coverage}, which is typically below $1 - \alpha$ when $|\Theta| > 1$. In contrast, MCS can achieve its promised guarantee in more settings, presented in \Cref{fig: simultaneous power comparison}, right panel.
\end{remark}


\subsection{Variance Estimation} \label{sec: variance}
Now we consider estimating the variance  $\sigma_r^2$ as required by Algorithm~\ref{algorithm: exp weighting}. Motivated by the literature of cross-validation~\cite{bayle2020}, a natural estimator is
\begin{equation} \label{eq: variance out}
 \hat{\sigma}_r^2 = \frac{1}{n} \sum_{i=1}^n\left(X_{i, r}-Q_{i, r}\QexcludingVi -\frac{1}{n} \sum_{j=1}^n\left(X_{j, r}-Q_{j, r}\QexcludingVj \right)\right)^2 .
\end{equation}
Although the estimator takes a simple form, showing such consistency is not a trivial task due to the dependence among $\{Q_{i,r}\QexcludingVi, i\in [n]\}$. In fact, we need to leverage the same stability property that is also a critical component in the proof of Theorem~\ref{th: random center}. We present the formal statement below. Its proof can be found in Appendix~\ref{app: variance estimation}. 

\begin{theorem}\label{th: consistency of variance estimator}
    Under the same assumptions as in Theorem~\ref{th: random center}, we have, for each $r\in[p]$,
    \begin{equation}
\hat\sigma_{r}^2/ \sigma_r^2 \stackrel{P}{\longrightarrow} 1. 
    \end{equation}
This result holds whether the number of folds, $V$, is a fixed integer or equal to 
$n$ (LOO).
\end{theorem}

One may alternatively consider the estimator by averaging the within-fold empirical variances (see \eqref{eq: variance in} in Appendix~\ref{app: variance estimation}), but it can only handle the finite $V$ case where the cardinalities of $I_v$ and $I_v^c$ both diverge to infinite, which is not the case for LOO ($|I_v| = 1$). The estimate $\hat \sigma_r^2$ in \eqref{eq: variance out} can cover both cases, which makes it more suitable for this work. Its analysis may be of independent interest to readers concerning LOO procedures. The proof of Theorem~\ref{th: consistency of variance estimator} involves several ``variance varieties" closely related to $\sigma_r^2$, which are also discussed in Appendix~\ref{app: variance estimation}.

\subsection{Key Ingredients to Theorem~\ref{th: random center}}
\label{section: proof sketch}
There are two main steps to prove Theorem~\ref{th: random center}: 1) establish a general CLT for sums of nearly IID random variables under general stability conditions, and 2) prove our statistics satisfy such stability conditions.  In this subsection we present the main ingredients. The proof of \Cref{th: random center} given these tools is very short and is provided in \Cref{app:thm3.1_cor3.4}.

The CLT related to ``weakly dependent" data transforms often, but not exclusively, appear in cross-validation-type methods. The most standard CLT is for normalized sums of independent random variables, but the quantities of interest in modern statistics do not necessarily take such a simple form. Intuitively, suppose each summand term is mostly driven by a random variable that is independent of others, the total sum should behave like the sum of independent random variables. All we need is to control the residuals induced by weak dependence, and a CLT is expected to hold.

In the dependent CLT below, we will consider $K_{i} = \mathcal{K}(i;\bm{X})$ with a pre-specified mapping $\mathcal K:[n]\times \mathcal{X}^{n}\rightarrow \mathbb{R}$. To derive a CLT for $n^{-1/2}\sum_{i}K_i$, a sufficient condition is that $K_i$ is essentially a function of the $i$-th sample $X_i$. We need the following difference operator to quantify this notion:

\begin{definition}
    For distinct $i,j,l\in [n]$, we define the (stability) operator as follows:
\begin{equation}\label{eq: stability terms}
\begin{aligned}
\nabla_j K_i & :=\mathcal{K}\left(i, \bm{X}\right)-\mathcal{K}\left(i, \bm{X}^{j}\right)\,, \\
\nabla_l \nabla_j K_i & :=\mathcal{K}\left(i, \bm{X}\right)-\mathcal{K}\left(i, \bm{X}^{j}\right) 
 -\left\{\mathcal{K}\left(i, \bm{X}^{l}\right)-\mathcal{K}\left(i, \bm{X}^{j, l}\right)\right\},
\end{aligned}
\end{equation}
where the perturbed data sets are defined in \Cref{intro} ($\bm{X}^j$ replaces the sample $X_j$ in $\bm{X}$ by an IID copy). 
\end{definition}

\begin{remark}
    The quantity $\nabla_j K_i$ is one way to measure how much the statistic $K_i$ depends on the sample $X_j$. When it is small for all $j\neq i$, then $K_i$ can be viewed as an approximately deterministic function of $X_i$ only. For example, when $K_i = f(X_i)$ for some given $f$, we have $\nabla_j K_i=0$ for all $j \neq i$, and in this case, standard CLT holds directly. If $K_i = X_i + X_{i+1}$, perturbing $j = i+1$ induces an $O(1)$ change in $K_i$; thus we do not view $K_i$ as approximately a function of $X_i$. 
\end{remark}

In general, CLT requires $\nabla_j K_i$ having a negligible moment; formally, we require
\begin{equation*}
        \Delta_1^2 =\max _{i \neq j \in[n]} \mathbb{E}\left[\left(\nabla_j K_i\right)^2\right]\qquad 
        \Delta_2^2 = \max _{i \neq j \neq l \in[n]} \mathbb{E}\left[\left(\nabla_l \nabla_i K_j\right)^2\right].
\end{equation*}
to vanish as $n\rightarrow\infty$ in the order of
\begin{equation}\label{eq: stability}
    \Delta_1 = o(n^{-1/2}) \text{ and } \Delta_2 = o(n^{-1}).
\end{equation}
The following dependent CLT is essentially a re-statement of Theorem 1 of \cite{austern2020}, adapted to a more general setting with slight technical modifications on tail and boundedness conditions.

\begin{theorem}\label{th: weakly dependent CLT}
Let $\bm{X} = \{X_i\in\mathcal{X}, i \in [n]\}$ be a collection of IID random vectors. We consider $K_i=\mathcal{K}(i ; \boldsymbol{X}) \in [-M,M]\subset\mathbb{R}$ and assume $K_i$ has the same distribution for all $i\in[n]$. We further assume $\mathbb{E}[K_i\mid \bm{X}^{(-i)}] = 0$ where $\bm{X}^{(-i)}\coloneqq \bm{X}\symbol{92}\{X_i\}$ and the variance $\omega_n^2 \coloneqq\operatorname{Var}\left(K_{1}\right)$ satisfies $\liminf_n \omega_n > 0$.

Then for any $x\in\mathbb{R}$ and $\epsilon > 0$, there exists a constant $C_\epsilon >0$ depending on $\epsilon$ such that
\begin{equation}\label{eq: explicit bound by stability}
\begin{aligned}
& \left|\mathbb{P}\left(\omega_n^{-1} n^{-1 / 2}\left(\sum_{i=1}^n K_i\right) \leq x\right)-\Phi(x)\right| \\
\leq & C_\epsilon\left(n^{1 / 2} \Delta_1+n \Delta_2+n^{-1 / 2}+n^{1/2} \Delta_1^2+n^{3 / 2} \Delta_2^2\right)+2 \epsilon
\end{aligned}
\end{equation}
where $\Phi(x)$ is the CDF of $\mathcal{N}(0,1)$. Specifically, if 
\eqref{eq: stability} holds, we have
\begin{equation*}
\omega_n^{-1} n^{-1 / 2}\bigg(\sum_{i=1}^n K_i\bigg) \stackrel{d}{\longrightarrow} \mathcal{N}(0,1).
\end{equation*}
\end{theorem}

We present the proof of Theorem~\ref{th: weakly dependent CLT} in Appendix~\ref{app: weakly dependent CLT}, which streamlines the proof in \cite{austern2020} using a modified Slepian's interpolation. Other CLT results for general functionals of $\bm{X}$ are available in the literature. For example, Chatterjee's
generalized perturbative approach \citep{chatterjee2008new} for  Wasserstein distance bounds is very generic and fostered many applications. However, its result does not directly imply \Cref{th: weakly dependent CLT}. Although both our framework and \cite{chatterjee2008new} are related to Stein's method, 
their resulting bound includes an absolute third-moment term that is difficult to control. The recent work of \cite{shao2025berry} built on \cite{chatterjee2008new} and introduced a new version of the perturbative approach for deriving Berry–Esseen bounds. 
Leveraging their Theorem 2.1, we obtain a qualitatively similar bound for our weak dependence CLT statistic. For completeness, we provide the proof—without the boundedness assumption—in Appendix~\ref{sec: BE bound from Shao and Zhang}.

The mapping $\mathcal{K}$ is determined by the statistic of interest. The requirement of $\mathcal K(i;\bm{X})$ having the same distribution for all $i\in[n]$ can be dropped by setting $\omega_n$ to the average of ${\rm Var}(K_i)$ over $i\in[n]$. When applying \Cref{th: weakly dependent CLT} in the proof of \Cref{th: random center}, for each $r$, we will analyze
\begin{equation}\label{eq:K_i}
K_i=K_{i,r}\coloneqq X_{i, r}-Q_{i, r}\QexcludingVi-d_{i,r}\dexcludingVi,
\end{equation}
and verify that the stability conditions \eqref{eq: stability} indeed hold for this specific $K_i$ (Appendix~\ref{app: stability}). Establishing the second-order stability condition is the most technically involved step. Our proof is inspired by the differential privacy literature using the connection between the stabilized soft-min and the exponential mechanism \citep[Section~3.4]{dwork2014algorithmic}. 

\begin{figure}[!t]
    \centering    \includegraphics[width =0.6\linewidth]{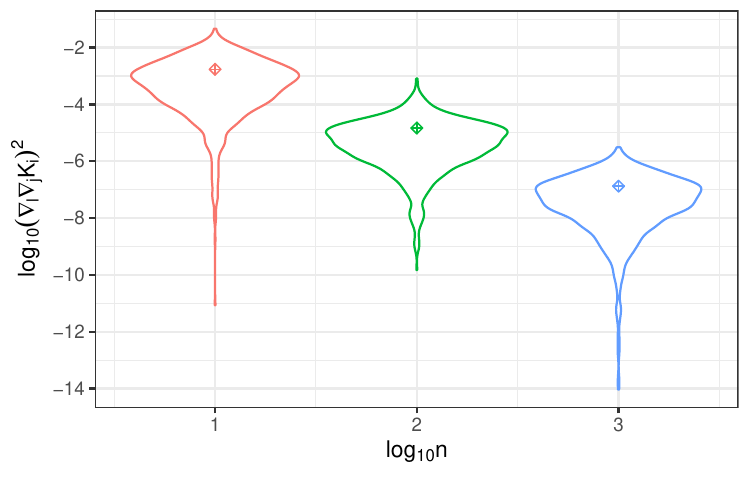}
    \caption{The second-order stability term vanishes at the rate predicted by our theoretical analysis. The violin plots illustrate the distribution of $\log_{10} \left(\nabla_l \nabla_j K_i\right)^2$ stratified by sample size. The points are the estimated $\log_{10}\mathbb{E}[\left(\nabla_l \nabla_j K_i\right)^2]$ over $10^3$ simulations repeats.}
    \label{fig: second order stability}
\end{figure}

To provide some empirical evidence on the scaling of $\left(\nabla_l \nabla_j K_i\right)^2$, 
in Figure~\ref{fig: second order stability} we plot the distribution of $\left(\nabla_l \nabla_j K_i\right)^2$ at various sample sizes. 
Our theoretical results predict $\Delta_2^2$ is of order $n^{-2}$ when $\lambda \sim \sqrt{n}$. This is reflected numerically as the mean points in Figure~\ref{fig: second order stability} roughly lie on a line of slope $-2$. The simulation settings are identical to those in Figure~\ref{fig: motivation}, but the sample size varies in $\{10^1,10^2,10^3\}$. Since we require $\Delta_2^2$ to be a smaller order than $n^{-2}$ to achieve asymptotic normality, in practice, one should implement $\lambda$ slightly smaller than $\sqrt{n}$.

\section{Bias Analysis and Power Guarantees}\label{section:power}

The asymptotically normal test statistic developed in the previous section relies on the stability of the soft-min $Q_{i,r}$, which can be viewed as a bias-variance trade-off. In this section, we analyze the bias of our test statistic as well as the power of the confidence set.

\subsection{Bias Analysis}\label{sec:4.1-bias}
By construction, the original target of our test statistic is $\theta_r\coloneqq \mu_r-\mu_{s_r}$, where $s_r=\arg\min_{s\neq r} \mu_s$, and the inference is based on a confidence lower bound of $\theta_r$.  In order to ensure the stability and hence asymptotic normality of the test statistic, we use a softmin estimate for $\mu_{s_r}$, and the test statistic $T_r$ centers at 
$\hat\theta_r = \frac{1}{n}\sum_{i=1}^n \sum_{s\neq r}\hat w_{r,s}^{(-i)}(\mu_r-\mu_s)$,
which is always less than or equal to $\theta_r$. The difference between $\hat\theta_r$ and $\theta_r$ can be interpreted as the bias due to the use of softmin, which equals
$\theta_r-\hat\theta_r=n^{-1}\sum_{i=1}^n\sum_{t\neq r}\hat w_{r,t}^{(-i)}(\mu_t-\mu_{s_r})$.

We have the following general, worst-case bound of this bias term.
\begin{theorem}
    \label{thm:bias}
     Let $\{ X_i \in \mathbb{R}^p, i \in [n]\}$ be IID samples from a distribution $P_{X,n}$ and $P_{X,n}$ is allowed to vary with $n$ (the triangular array setting). We assume $X_i$ has uniformly bounded entries: $\sup_{s\in [p]} |X_{i,s}| \leq M$. For $\lambda = o(\sqrt{n})$, we have for some constant $C >0$,
    $$
    \theta_r-\hat{\theta}_r \leq C \log (p+n) / \lambda_n+O_P\left((p+n)^{-1}\right).
    $$
\end{theorem}

The proof is presented in \Cref{app: proof of bias}. Theorem \ref{thm:bias} ensures that, in the worst case, the bias incurred by the softmin method is at most $\log(p+n)/\lambda_n$, which can be arbitrarily close to $\log(p+n)/\sqrt{n}$, but it is still larger than $n^{-1/2}$. The center of the test statistics is $\hat\theta_r$ rather than the more interpretable $\theta_r$.

However, the bias result can be significantly improved when considering IID samples from a fixed distribution $P_X$.

\begin{corollary}\label{cor: improved bias}
Let $\{X_i \in \mathbb{R}^p, i \in[n]$\} be IID samples from a fixed distribution $P_X$ with uniformly bounded entries: $\sup _{s \in[p]}\left|X_{i, s}\right| \leq M$. We have 
\begin{equation*}
  \theta_r - \hat \theta_r \leq p M \exp (-\lambda_n m / 2)+o_P\left(n^{-1 / 2}\right)  
\end{equation*}
for any fixed $c >\sqrt{6}M$ and $n \geq 16c^2 \log(p+n)/m^2$. Here $m$ is the smallest gap against $\mu_{s_r}$:
\begin{equation*}
m:=\min_t \left\{\mu_t-\mu_{s_r}: t \neq r, \mu_t>\mu_{s_r}\right\} .
\end{equation*}
When $\lambda_n = \omega(\log n)$, $\theta_r - \hat \theta_r = o_P(n^{-1/2})$.
\end{corollary}

This result implies, similar to Corollary~5 in \cite{chen2023inference}, that the test statistic will center at $\min_{s\neq r} \mu_s$ and the bias does not impact its asymptotic distribution under the listed conditions (so long as \Cref{th: random center} holds). The key difference between \Cref{thm:bias} and \Cref{cor: improved bias} is that when considering the triangular array setting, one can construct adversarial landscapes depending on $\lambda_n$ to amplify the difference between $\theta_r$ and $\hat \theta_r$, while under fixed distribution asymptotics, such challenging landscapes are no longer permitted. Any $\lambda_n = \omega(\log n)$ suffices to assign small enough weights to all suboptimal dimensions.



\subsection{Power Guarantees}\label{subsec:4.2-power}
The power of the confidence set construction refers to its ability to exclude dimension indices whose population means are not minimal. 
Let $\tilde \mu = \min_{s \in [p]}\mu_s$. For each $r\in[p]$, let $\alpha_n(r)=(\mu_r-\tilde\mu)\lambda/2$ be the scaled gap.
Let $\beta_n=4(\log p+3\sqrt{\log V})$. Define 
\begin{equation}
\mathbb C(r)=|\{s:\alpha_n(r)/\lambda<\mu_s-\tilde\mu \le \beta_n/\lambda\}|\,.\label{eq:mathbb_C}    
\end{equation}

The magnitude of $\mathbb C(r)$ measures the cardinality of the ``confusing set'' and reflects the hardness of rejecting the hypothesis $r\in\Theta$. Intuitively, if $\mu_s-\tilde\mu\le \beta_n/\lambda$, then the coordinate $s$ will possibly receive some non-trivial weight in the soft-min calculation of $Q_{i,r}$. If $\mu_s-\tilde\mu>\alpha_n(r)/\lambda$, then $\mu_s$ is somewhat close to $\mu_r$---In order to detect the sub-optimality of coordinate $r$, the exponential mechanism cannot assign too much weight to coordinates whose $\mu_s$ value is close to $\mu_r$.


In the statement and the proof of Theorem \ref{thm:power_new}, the 
limits are all taken for $(n,p)\rightarrow \infty$.
Here we consider a triangular array type of 
asymptotic setting, where $\mu$ may change as $(n,p)$ increases. The proof of this result is presented in \Cref{app: new adaptive}.

\begin{theorem}\label{thm:power_new}
    Under the same assumptions as in \Cref{th: random center}, 
    \begin{enumerate}
    \item If $\mathbb C(r)>0$, then $\lim_{n\rightarrow\infty}\mathbb P(r\in \widehat{C}) = 0$
    if \begin{equation}\label{eq:adaptive_gap_condition}
        \mu_r-\tilde\mu\ge c\lambda^{-1}\left[\log \mathbb C(r)+ \log\log p+\log\log V\right]
    \end{equation}
    for a large enough constant $c>0$;
        \item   If $\mathbb C(r)=0$, then $\lim_{n\rightarrow\infty}\mathbb P(r\in \widehat{C})=0$ when $\mu_r-\tilde\mu = \omega(1/\sqrt{n})$.
    \end{enumerate}
\end{theorem}

\begin{remark}
    This result reflects the adaptivity of the method.  When there is no confusion set ($\mathbb C(r)=0$) the power guarantee almost achieves the parametric rate.  When there is a non-empty confusion set, it is the logarithm of the size of the confusion set that matters in the power guarantee.  The worst case value of $V$ is $n$, and the dependence on $p$ and $n$ are in the iterated logarithm.  The value $\mathbb C(r)$ depends on $r$, the coordinate of interest. It can be upper bounded by
$|\{s:0< \mu_s- \tilde\mu \le \beta_n/\lambda\}|$.
\end{remark}

\section{Data-driven Selection of the Weighting Parameter}
\label{section: data-driven}

As mentioned in \Cref{rem:worst-case-lambda}, the choice of $o(\sqrt{n})$-order $\lambda$ is to cover the worst-case scenario, which does not always hold in real-world applications. In practice, it is important to choose $\lambda$ in a data-driven fashion to achieve better statistical power.

Moreover, the optimal choice of $\lambda$ can depend on $r$, the dimension under comparison, because the variance of each dimension $X_{i,s}, s\in [p]$ enters the procedure differently when switching $r$. Algorithm~\ref{algorithm: exp weighting} describes the procedure with an $r$-agnostic choice of $\lambda$ for ease of presentation. In practice, $r$-dependent choices of $\lambda=\lambda_{n,r}$ often perform better.



The goal of $\lambda$ tuning is to maximize power while attaining the desired coverage. Recall that a large $\lambda$ tends to give better power but an overly large $\lambda$ may violate the stability conditions required by the asymptotic normality. 
With these considerations in mind, we pick the largest $\lambda$ that satisfies the first-order stability condition:
\begin{equation}\label{eq: stability term need to be small}
    n\cdot \mathbb{E}\left(\nabla_i K_{1,r}\right)^2 
    = o (\mathrm{Var}[K_{1,r}])\,,
\end{equation}
with $i\neq 1$ for LOO implementation (or $i\notin I_{v_1}$ for $V$-fold implementation), $K_{1,r}$ defined as in \eqref{eq:K_i}, and operator $\nabla_i$ defined in \eqref{eq: stability terms}. This criterion is motivated by \Cref{eq: stability bound} in the proof of Theorem~\ref{th: weakly dependent CLT}. Having it satisfied, one could expect the test statistic $\tilde{T}_r$ to be nearly standard normal for valid inference. 

\subsection{Iterative Data-driven Selection}
We apply the following iterative algorithm to conduct data-driven parameter selection. For each $r$:
\begin{itemize}
    \item[(i)] Set $\lambda$ to be $\lambda_0$ as a small initial candidate (details on how $\lambda_0$
is determined are presented in Appendix~\ref{appendix: initial lambda_0});
\item[(ii)] Run the initial parts of Algorithm~\ref{algorithm: exp weighting} with $2\lambda$ until Step~\eqref{eq:exp_weight} and examine if Condition \eqref{eq: stability term need to be small} 
is approximately satisfied. Specifically, we use the sample-version criterion
\begin{equation} \label{eq: iterative condition}
n \cdot \widehat{\mathbb{E}}\left(\nabla_i K_{1,r}\right)^2 \le \varepsilon \widehat{\operatorname{Var}}[K_{1,r}]. 
\end{equation}
Here the small number $\varepsilon$ is set to be $0.08$ for the LOO procedure. 
\item[(iii)] If the criterion in step (ii) is not satisfied or $2\lambda \geq \phi$, return $\lambda$ as the selected parameter. Otherwise, set $\lambda \leftarrow 2\lambda$ and repeat step (ii). The threshold $\phi$ is a sufficiently large number (e.g. $n^5$). 
\end{itemize}

\noindent In Section~\ref{sec: sensitivity of data-driven lambda}, we investigate the sensitivity of our LOO method’s performance to the choice of the weighting parameter $\lambda$. The results show that the proposed iterative data-driven selection contributes to boosting power while achieving sharp control over coverage guarantee.  


\subsection{Estimation of Relevant Quantities}\label{section: estimate relevant quantities}

To compute the sample expectation $\widehat{\mathbb{E}}\left(\nabla_i K_{1,r}\right)^2$ in \eqref{eq: iterative condition}, we leverage the identity 
\begin{equation}\label{eq: LTO ki}
    \nabla_i K_{1,r}=\sum_{s \neq r}\left(\hat{w}_{r,s}^{\left(-v_1\right)}-\hat{w}_{r,s}^{\left(-v_1\right),i}\right)\left(X_{1,r} - X_{1, s} -(\mu_r - \mu_s)\right)\,,
\end{equation}
where $\hat w_{r,s}^{(-v_1),i}$ is the exponential soft-min weight $\hat w_{r,s}^{(-v_1)}$ computed with $X_i$ replaced by an IID copy $X_i'$. 

To estimate $\nabla_i K_{1,r}$, we approximate $\hat w_{r,s}^{(-v_1),i}$ and $\hat w_{r,s}^{(-v_1)}$ using a leave-two-out technique which was also employed in~\cite{austern2020} and~\cite{kissel2023blackbox} for quantities related to the $\nabla_i$ operator. For $i\neq k \in I^c_{v_1}$, we first approximate $\hat w_{r,s}^{(-v_1)}$ by $\hat w_{r,s}^{(-v_1),-k}$, i.e., using sample means from $\{X_j, j\in I_{v_1}^c \setminus\{k\}\}$. Similarly, $\hat w_{r,s}^{(-v_1),i}$ is approximated by $\hat w_{r,s}^{(-v_1),-i}$, which uses the sample means from $\{X_j, j\in I_{v_1}^c \setminus\{i\}\}$.  Therefore these two weights differ in one sample replacement. 
We obtain 
\begin{align*}
    \widehat{\mathrm{K}}_{1,r}(i,k) &\coloneqq 
    \sum_{s \neq r}\left(\hat w_{r,s}^{(-v_1),-k}-\hat w_{r,s}^{(-v_1),-i}\right)\left(X_{1,r} - X_{1, s} -(\hat\mu_r - \hat\mu_s)\right),
\end{align*}
as an approximated version of $\nabla_i K_{1, r}$.


The estimator $\widehat{\mathbb{E}}\left(\nabla_i K_{1,r}\right)^2$ is the sample average $\abs{\mathcal{B}}^{-1} \sum_{(j,i,k) \in \mathcal{B}} \left ( \widehat{\mathrm{K}}_{j,r} (i,k)\right)^2$ over the set $\mathcal{B} := \{(j,i,k) \in [n] \times I^c_{v_j} \times I^c_{v_j} \mid i \ne k\}$. 
For each triplet $(j,i,k) \in \mathcal{B}$, the quantity $\widehat{\mathrm{K}}_{j,r}(i,k)$ is analogously estimated as described above. For large values of $n$, we further approximate the estimate using a random subset of $\mathcal B$ of size $100$ drawn uniformly at random. 

As for the variance, we use $\hat{\sigma}^2_{r}$ as an estimate of $\widehat{\mathrm{Var}}(K_{1,r})$ in \eqref{eq: iterative condition}, inspired by Proposition~\ref{prop: equivalency among variances} and Theorem~\ref{th: consistency of variance estimator}.

\section{Simulation Study}\label{section: simulation}
\subsection{Method Comparison} \label{subsection: method comparison}

To evaluate the performance of the proposed procedure, we compare it with three methods that are either proposed in existing literature or readily adaptable to our argmin inference problem. In particular, our investigation will focus on how the methods respond to data dependencies and characteristics of mean landscapes. All of the methods under comparison are guaranteed to achieve marginal coverage for each $r\in\Theta$ but not simultaneous coverage. 

\subsubsection{Compared Methods}
The first method is the \textit{Bonferroni correction} which may be loosely viewed as a benchmark for the class of multiple testing procedures. In our context, a dimension $r \in[p]$ is included in the confidence set $\widehat C$ if and only if all the pairwise null $H_0:\mu_r\leq \mu_s$, $s\neq r \in [p]$ are not rejected at the Bonferroni-adjusted significance level $\alpha/(p-1)$. We implemented paired two-sample $t$-tests for the mean comparisons.

The second method is the two-step procedure in~\cite{futschik1995confidence} built upon the selection rule developed by~\cite{gupta1965some}. Such a selection constructs an $(1 - \alpha)$ argmin confidence set by collecting all the dimensions $r \in [p]$ satisfying the inequality
\begin{equation} \label{eq: selection rule Gupta}
    \sqrt{n}\left( \frac{\hat{\mu}_r}{\sigma_r} - \min_{s \ne r}\frac{\hat{\mu}_s}{\sigma_s} \right) \le q_{(1 - \alpha), p}, 
\end{equation}
where $\sigma^2_s$ is the true variance of $X_{1,s}$ for all $s \in [p]$, and the threshold $q_{1-\alpha, p}$ is the $(1 - \alpha)$-quantile of the random variable $\varepsilon_r - \min_{s \ne r} \varepsilon_s$ with $\varepsilon_s \overset{\mathrm{IID}}{\sim} \mathcal{N}(0, 1)$. The variant introduced by~\cite{futschik1995confidence} applies a two-step selection rule to enhance power.
Given the proper choices of $\alpha_1, \alpha_2 \in (0,1)$ such that $(1 - \alpha_1) (1 - \alpha_2) = 1 - \alpha$, the first selection is performed to generate an $(1 - \alpha_1)$ argmin confidence set $\widehat{C}_1$ using the threshold $q_{(1 - \alpha_1), p}$. Then, the following selection adapts the cardinality of $\widehat{C}_1$ for the generation of an $(1 - \alpha_2)$ argmin confidence set $\widehat{C}_2$ using the threshold $q_{(1 - \alpha_2), \abs{\widehat{C}_1} - 1}$. The final argmin confidence set is $\widehat{C}=\widehat{C}_1 \cap \widehat{C}_2$.
Here we follow \cite{futschik1995confidence}'s choice of $\alpha_1 = \alpha/10$. 
The validity of this method relies on (i) the variance $\sigma^2_s$ being known and the same across all dimensions, and (ii) the coordinates are mutually independent. Both of these restrict the method's applicability in practice.
Even replacing the true $\sigma_s$ by its estimate leads to validity violations (see Appendix~\ref{app:futschik.validity.violation}).

The third method is included to illustrate how one can construct a valid argmin confidence set from rank confidence intervals. \cite{mogstad2024inference} and~\cite{fan2024ranking} reduce confidence intervals for ranks to that of the pairwise mean differences, which can be further reformulated as a problem of \textit{testing many moment inequalities} (see~\cite{romano2014practical, chernozhukov2016testing}). In particular, the final component is handled via Gaussian comparison that accounts for dependence across dimensions, with inference for the associated max-statistic carried out through bootstrap methods. In this way, the approach of \cite{mogstad2024inference} can be regarded as an extension of \cite{futschik1995confidence} that relaxes the assumption of dimension-wise independence.

Here we use the R package {\tt csranks} by \cite{R-csranks} to construct a confidence lower bound $\hat{L}_r$ for the population rank $R_r$ of the mean $\mu_r$ for each $r \in [p]$.  We then include a dimension $r \in [p]$ in the confidence set $\widehat{C}$ if and only if $\hat{L}_r = 1$. This confidence set yields the validity~\eqref{eq: marginal validity} because for each $r \in \Theta$, we have $\mathbb{P}(r \in \widehat{C}) = \mathbb{P}(\hat{L}_r = 1) = \mathbb{P}(\hat{L}_r \le 1) = \mathbb{P}(\hat{L}_r \le R_r) \ge 1 - \alpha + o(1)$. The package implements 
$t$-bootstrap procedures (see \cite{romano1988bootstrapping, romano2014practical}) for inference under fixed
$p$, following the main discussion in~\cite{mogstad2024inference}. As noted in that work, however, the theoretical framework naturally extends to high-dimensional settings by incorporating the high-dimensional Gaussian comparison results of \cite{chernozhukov2013gaussian, chernozhukov2014gaussian, chernozhukov2017central}.


\subsubsection{Setups and Results} \label{sec: method comparison setup}
Samples are drawn from multivariate normal distributions with Toeplitz covariance matrices. We take type I error size $\alpha = 0.05$. All the methods achieve $95\%$ coverage for the true argmin index in all settings (see Appendix~\ref{app: validity simulation}). 

Two types of mean landscapes---denoted as ``increasing" and ``3-tier" --- are explored. For each type of landscape, we vary the signal strength (size of the difference in true means) as well as the dependency strength across dimensions of $X$ and investigate their impact on the statistical power. 

\begin{figure}[!t]
    \centering
    \includegraphics[width=\textwidth]{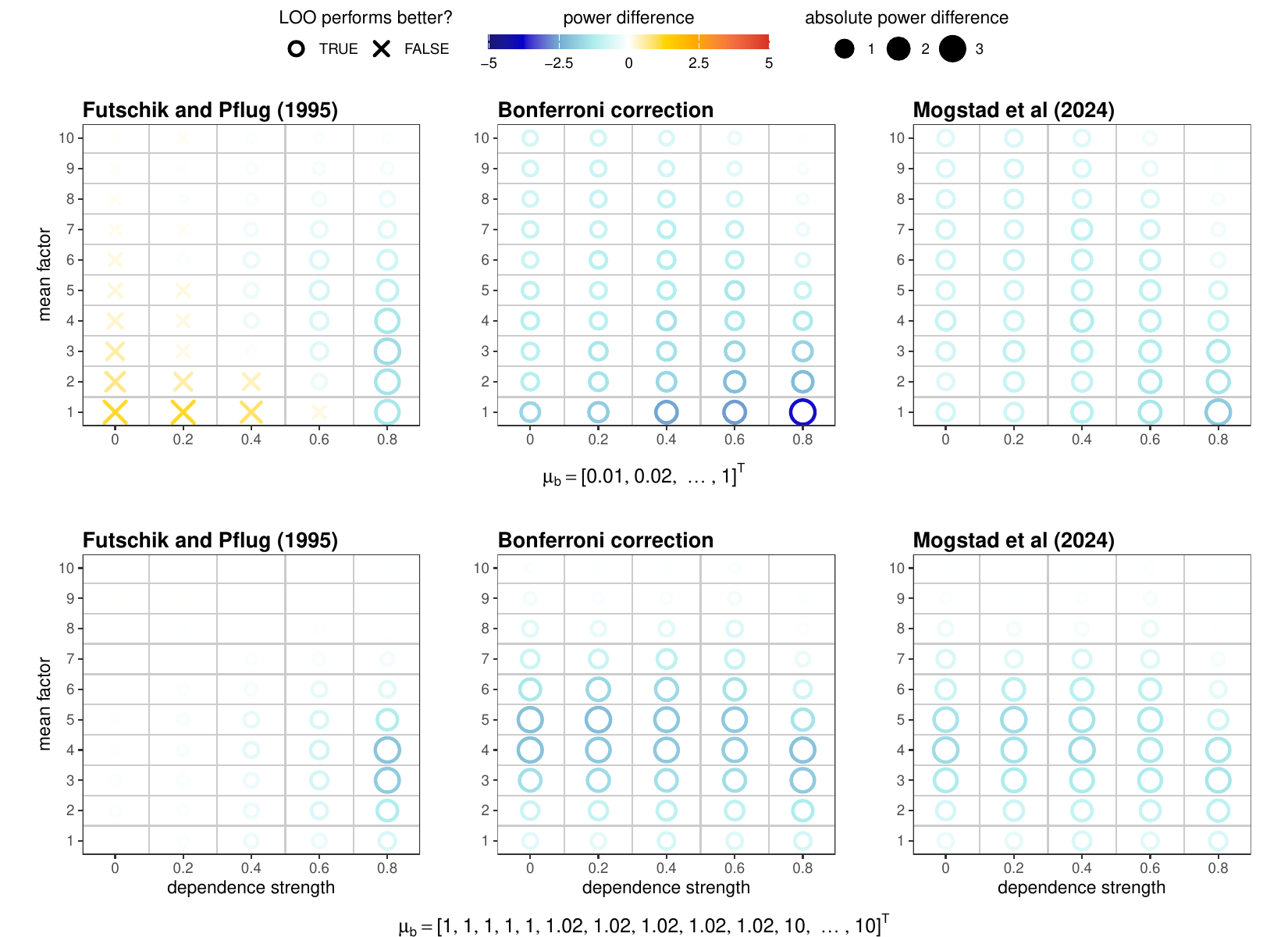}
    \caption{Method comparison, ``increasing" (top) and ``3-tier" (bottom) landscapes. Comparison between the proposed LOO method and three other methods. Each cell in the heatmaps corresponds to a different simulation setting. The x-axis corresponds to different dependency strength $\varrho$, and in the y-axis, signal strength $f$ is varied.  The color (and shape size) in each cell illustrates the difference in the average number of false negatives between the proposed LOO method and one literature method. A more negative value indicates a greater advantage of the proposed method over its competitor in rejecting sub-optimal dimensions. 
    }
    \label{fig: method comparison}
\end{figure}

Formally, the true means are of the form $\mu = f \times \mu_b \in \R^p$ for the mean factor $f \in \{1, 2, \ldots, 10\}$ with the base mean vector $\mu_b$ specified in Figure~\ref{fig: method comparison}. As $f$ increases, the difference between different coordinates of $\mu$ will be enlarged, making it easier to exclude the sub-optimal dimensions from the confidence set. The covariance matrices are Toeplitz with the $(r,s)$-entry
$\sigma_{rs}^2 = \varrho^{\abs{r - s}}$ for $r, s \in [p]$.
We consider the dependency strengths $\varrho \in \{0, 0.2, \ldots, 0.8\}$, where $\varrho = 0$ leads to an identity covariance matrix and $\varrho = 0.8$ signifies a highly correlated scheme. In total, we have $2\times 10 \times 5= 100$ settings. 

In Figure~\ref{fig: method comparison}
we present the \emph{difference} in the number of \textit{false negatives} within confidence sets, which is computed as the number of false negatives produced in the proposed method minus that produced by a competing method. 
Here the number of false negatives is the cardinality of $\widehat{C}\backslash\Theta$. We set the dimensionality $p = 10^2$ and a sample size $n = 10^3$. A more negative value indicates a greater advantage of the proposed method over its competitor in rejecting sub-optimal dimensions. 
The number of repetitions for each simulation setting is $100$. 

The first type of mean landscape---``increasing"---implements a base vector $\mu_b = [0.01, 0.02, \ldots, 1]^\top$. It depicts a situation where the true means gradually increase across the entire landscape. In Figure~\ref{fig: method comparison} top row, we observe the proposed LOO procedure outperforms the three compared methods in a majority of the experiment settings.
The advantage of the proposed procedure becomes more apparent with an increased dependency strength, which corresponds to moving from left to right in each size-encoded heatmap row. The method by~\cite{futschik1995confidence} performs better when the underlying dimensions are nearly independent thanks to its screening-like step (which cannot be justified in dependent cases). 
By contrast, the other three methods factor in the dependence structures in different ways, which helped improve their power. Specifically, for the proposed LOO procedure, the weights concentrate better on the dimensions with lower true population means as the correlation increases.

The other type of mean landscape (``3-tier") 
concerns the case when there are several close competitors having tied and near-tied means, along with many clearly inferior ones. Such a scenario often unfolds in commercial markets, where a handful of dominant brands share a similar market reputation due to competitive product qualities, while many budget brands cater to niche consumer segments. As a market researcher, one might aim to identify the most highly regarded companies based on the quantitative feedback provided in customer surveys. In Figure~\ref{fig: method comparison} bottom row, we see that the proposed LOO method typically results in finer confidence sets than the other three methods in this case. Compared to the Bonferroni correction and the procedure by~\cite{mogstad2024inference}, the proposed method initially exhibits increasingly higher power when moving from bottom to top in each column of the size-encoded heatmap. Yet the advantage eventually diminishes after the mean factor $f$ passes a point where the problem becomes too easy for all methods. 

\subsection{Sensitivity of the Data-driven Weighting Parameter} 
\label{sec: sensitivity of data-driven lambda}

\begin{figure}[!t]
\centering
\includegraphics[width=\textwidth]{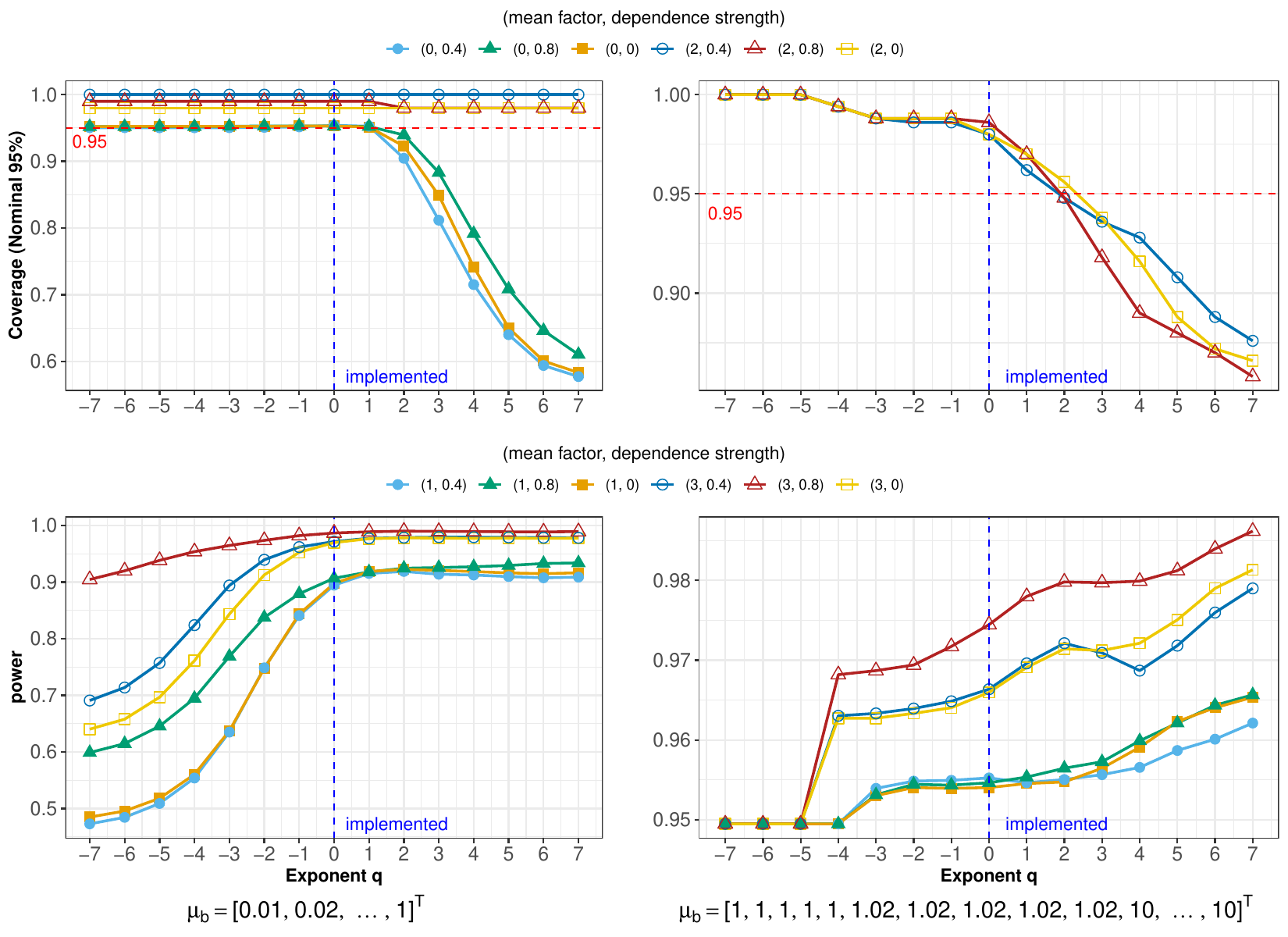}
\caption{Sensitivity analysis of the weighting parameter $\lambda$ in terms of average coverage $\overline{\nu}$ and average power $\overline{\kappa}$.
Here $q$ is the distortion exponent in $\lambda = 2^q \hat{\lambda}$. Notably, \(\hat{\lambda}\) may take different values across the settings.
For each configuration $(\mu_b, f, \varrho, \lambda)$, we perform $100$
simulation repetitions with a sample size of $1000$. The curves corresponding to the setting \((\text{mean factor}, \text{dependence strength}) = (0, \varrho)\), for \(\varrho \in \{0, 0.4, 0.8\}\), are omitted from the top-right plot as they coincide with the flat mean cases illustrated in the top-left plot.} \label{fig:sensitivity}
\end{figure}

We explore how sensitive our LOO method’s performance is to the choice of the weighting parameter $\lambda$, focusing on its impact on average coverage guarantee $\overline{\nu}$ and average power $\overline{\kappa}$. The two metrics are respectively defined as the averages
$$
    \overline{\nu} = \frac{1}{\abs{\Theta}} \sum_{r \in \Theta} \nu_r \hspace{5mm} \text{ and } \hspace{5mm} \overline{\kappa} = \frac{1}{\abs{\Theta^c}} \sum_{s \in \Theta^c} \kappa_s, 
$$
where $\nu_r, r \in \Theta$ (and $\kappa_s, s \in \Theta^c$) stands for the coverage guarantee for the optimal dimension $r$ (and the power for the suboptimal dimension $s$). We evaluate the performance metrics---coverage guarantee $\nu_r$ and power $\kappa_s$---for each optimal ($r$) and suboptimal ($s$) dimension and then take the average. We use the same simulation setting as in \Cref{sec: method comparison setup}. However, instead of adhering to the data-driven estimate $\hat{\lambda}$ from Section~\ref{section: data-driven}, we evaluate a range of scaled values defined as $\lambda = 2^q \hat{\lambda}$, where the exponent $q$ varies from $-7$ to $7$. This range typically spans from very small $\lambda$---which promotes uniform weights and high stability---to very large values that concentrates the weightings on the empirical minimum.

The top-left plot in  Figure~\ref{fig:sensitivity} displays the average coverage under the ``increasing'' mean landscape, evaluated across different mean factors $f$ (signal strength) and dependence strengths $\varrho$. When $f = 0$, the setting corresponds to a flat mean landscape, $\mu_b = [0, 0, \ldots, 0]$, in which case only the average coverage is relevant. The corresponding average coverage in the plot remains close to the nominal level up to the point where the weighting parameter exceeds roughly twice the default data-driven value $\hat{\lambda}$. Any such exceeding one would break the stability assumptions, leading to violations of asymptotic normality. This validates the effectiveness of our tuning algorithm in maintaining sharp control under the worst case. In contrast, the mean factor $f = 2$ has made the signal strength sufficiently large, so we have flexibility in choosing the weighting parameter $\lambda$, as pointed out in 
Remark~\ref{rem:worst-case-lambda} and Lemma~\ref{lemma: ultra stable}.

While any value less than or equal to $\hat{\lambda}$ provides coverage guarantees, the default data-driven $\hat{\lambda}$ is calibrated to enhance power performance. The bottom-left plot in Figure~\ref{fig:sensitivity} shows that smaller values of $\lambda$ result in reduced power across varying signal strengths $f$ and dependence strengths $\varrho$. Notably, the default data-driven choice $\hat{\lambda}$ achieves power levels comparable to those attained using an empirical min (realized by $\lambda = 2^7 \hat{\lambda}$) under this mean landscape.
  
The right panel of Figure~\ref{fig:sensitivity} presents the average coverage and power under the ``3-tier'' mean landscape. When $\lambda$ is small (corresponding to $q \in \{-7, -6, -5\}$), a substantial portion of the weights is allocated to dimensions within the second tier (e.g., $\mu_6 = 1.02$). This misallocation introduces a negative mean shift in the test statistic when testing the truly optimal dimensions in the first tier (e.g., $\mu_1 = 1$), rendering our procedure conservative. As $\lambda$ increases, the influence of the third tier (e.g., the inferior dimension $\mu_{11} = 10$) diminishes and is eventually removed, reducing the number of effective dimensions. This enables the algorithm to begin disregarding the second tier as well, gradually approaching the behavior expected under a flat mean landscape. Consequently, average coverage regresses towards the nominal level. However, when $\lambda$ exceeds approximately four times the default data-driven value $\hat{\lambda}$, the nominal coverage is no longer maintained. 

The power trend is similar to the ``increasing" mean landscape. 
When $\lambda$ is too small, the rejections are limited to the third tier, yielding an average rejection rate of exactly $90/95 \approx 0.947$. As $\lambda$ increases and the dimensions in the second tier are progressively filtered out, the power improves.


\subsection{An Application to LASSO Model Selection}\label{section: lasso}

\begin{figure}[!t]
    \centering
    \includegraphics[width =0.9\linewidth]{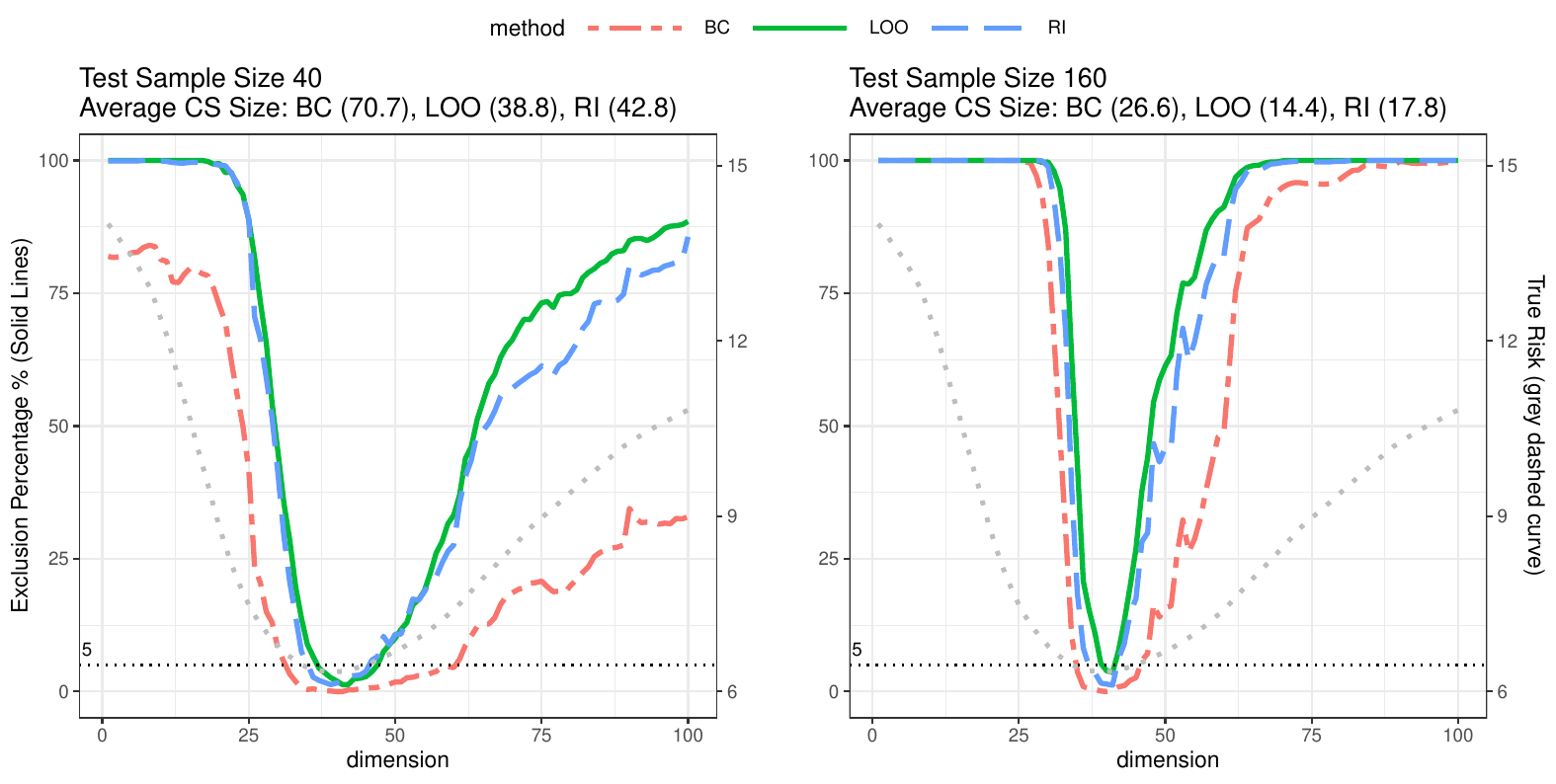}
    \caption{Average exclusion percentage, LASSO model selection. The numerical experiments are conducted over two different test sample sizes, $n = 40$ (left) and $n = 160$ (right). The gray dotted curve represents the true population risks of the $\beta_r$'s, with the risk values shown on the right y-axis. Comparison among the proposed LOO method, Bonferroni correction (BC) and rank inference approach by~\cite{mogstad2024inference} (RI). Each solid curve documents the proportion of the $100$ models---each corresponds to a $\eta_r$ parameter---being excluded from the confidence sets. The exclusion percentage is calculated over $10^3$ repeats. 
    }
    \label{fig: lasso}
\end{figure}
As discussed in the Introduction, one important application of the proposed procedure is model and tuning parameter selection. We take the LASSO \citep{tibshirani1996regression} in high-dimensional regression as an example. The goal is to relate a collection of predictors $Z\in \mathbb{R}^d$ to an outcome of interest $Y\in\mathbb{R}$. Given collected samples, each LASSO predictor is constructed by minimizing 
\begin{equation}\label{eq: lasso}
    \sum_{i=1}^{n_{\rm tr}}(Y_i - Z_i^T\beta)^2 + \eta \|\beta\|_1
\end{equation}
over $\beta\in\mathbb{R}^d$. 
The penalty parameter $\eta > 0$ controls the sparsity of the estimated regression coefficients and significantly impacts the generalization capacity of the fitted model. To simplify the discussion we consider a sample-splitting scenario, where for each $\eta_r$, $r\in[p]$ we have estimated a $\beta_r\in\mathbb R^d$ from an external training sample of size $n_{\rm tr}$, and we want to identify the $r\in[p]$ such that the estimate $\beta_r$ minimizes the future prediction error $\mathbb{E}_{Y,Z}[(Y - Z^T\beta_r)^2]$, using a testing sample $\{(Z_i, Y_i),i\in [n]\}$. This can be formulated as an argmin inference task: $X_{i,r} = (Y_i - Z_i^T\beta_r)^2$. Applying the proposed procedure, we can construct a confidence set $\widehat C$ for the argmin index of the prediction risk $\mu_r=\mathbb E_{Y,Z}(Y-Z^T\beta_r)^2$. 

In the example presented in Figure~\ref{fig: lasso}, the distribution of $Z$ is multivariate normal $\mathcal{N}(0_{300}, I_{300})$, where $0_{d}$ is an all-zero vector of length $300$ and $I_d$ is the $d \times d$ identity matrix. We define the true $\beta$ as $(1_{10}, 0 _{290})^T$ and 
\begin{equation*}
    Y = Z^T \beta + \mathcal{N}(0,2^2).
\end{equation*}
The training sample size $n_{\rm tr}=160$. The candidate $\eta$'s are automatically generated by R package {\tt glmnet} (\cite{glmnet}). Here the best index is $r^* = 39$ and the population ``true risk" $\mathbb{E}_{Y,Z}[(Y - Z^T\beta_r)^2]$ is illustrated in dashed grey for all $r \in [p] = [100]$. 

In Figure~\ref{fig: lasso}, we plot the frequency of each index being excluded from the confidence set by the proposed LOO algorithm, Bonferroni correction (BC) and Rank Inference (RI) by~\cite{mogstad2024inference}. Specifically, a point $(r,\omega)$ on a solid curve indicates that the corresponding $\beta_r$ is excluded in $\omega\%$ of the confidence sets constructed over numerical repeats. Two different test sample sizes $n \in \{40, 160\}$ are investigated. For each compared method, the observed exclusion frequency of each $r$ is positively associated with the true risk of the corresponding $\beta_r$. Increasing the sample size from $n = 40$ to $160$, the confidence set $\widehat C$ rejects sub-optimal dimensions more frequently.  Both the proposed LOO method and the RI method consistently outperform the BC method---yielding higher exclusion percentages for all sub-optimal indices---while controlling the type I errors for $r^\ast = 39$ at the level $\alpha = 0.05$.

For the LASSO problem, the dimensions of the risk vector $X$ are highly correlated, due to the similarity between $\beta_r$'s when trained under similar $\eta_r$'s. Many of them share a Pearson's $\rho$ greater than $0.8$, with the highest being close to $0.99$. This distinguishes the LASSO application from the experiments in Section~\ref{subsection: method comparison}, where only a few dimensions can have $\rho = 0.8$ even in the most highly correlated scheme.  



\section{Real Data Applications}\label{section: real_data}
In this section, we apply the proposed procedure to two real data sets. 
In machine learning study, model competitions are frequently used as pedagogical practices to allow practitioners or students to explore the strengths and weaknesses of different machine learning methods. It is essential in this context to acknowledge the merit of all competitive models while screening out the inferior ones. In a course \textit{Methods of Statistical Learning} instructed by one of the authors, students were asked to train classification algorithms over a given data set. Then, the student-trained classifiers were submitted for evaluation over a held-out testing data set. The data sets are sourced from \textit{Kaggle.com}. Here we implement the proposed LOO algorithm with a data-driven tuning parameter $\lambda$ to identify the best performers.
The identities of students and group names are anonymized.

\subsection{2023 Classification Competition}
In Spring 2023, a total number of $44$ submitted models were evaluated upon a test data set of sample size $183$. In our notation, this corresponds to an independent sample $X_1, \ldots, X_{183}$ with $X_i \in \{0, 1\}^{44}$ encoding the binary classification error ($0$ for correct, $1$ for error). Models of lower expected error rates are preferred. Within these $44$ models, there are $7$ pairs of student models---models $(3, 34), (7, 24), (9, 42), (10, 23), (21, 37), (25, 35)$ and $(36, 39)$---with identical evaluation results. We excluded one model in each pair to avoid degeneracy.

We implemented the proposed LOO algorithm, Bonferroni correction (BC), Rank Inference by~\cite{mogstad2024inference} (RI), and a method by~\cite{hansen2011model} that we refer to as the MCS procedure. Most of the methods mentioned above involve some randomized steps such as sample splitting and bootstrap. 
In the proposed LOO method, the randomness is from the sub-sampling in the leave-two-out estimate of $\nabla_j K_i$ (detailed in Section~\ref{section: data-driven}), 
while the RI and MCS methods involve bootstrap samples.
We constructed $100$ confidence sets using different random seeds, all with the same testing data.

Over the 100 repetitions, the average size of the LOO confidence sets is $31.6$, with models 4, 10, 12, 13, 14, 15, 23, 29, 33, 36, 38 and 39 frequently excluded.
The Bonferroni correction results in an average size of $41$, with models, 12, 29, and 33 excluded in most repeats. The RI approach yields an average size of $38.7$, frequently excluding models $4, 12, 13, 29$ and $33$. As for the MCS procedure, we follow the recommended implementation in~\cite{bernardi2018model}. It shows an average size of 43.

\begin{figure}[!t]
    \centering
    \includegraphics[width=\linewidth]{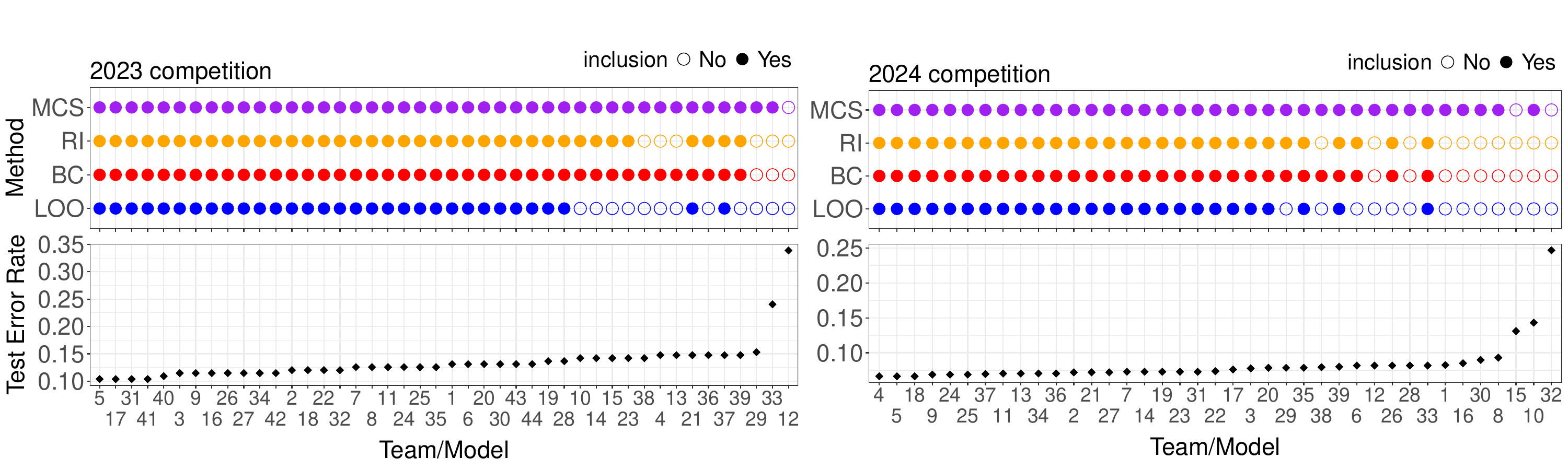}
    \vspace{-20pt}
    \caption{Confidence sets with real data. We compare the proposed LOO algorithm (LOO), the Bonferroni correction (BC), the rank inference method by~\cite{mogstad2024inference} (RI), and the model confidence set (MCS) \citep{hansen2011model} over the test results in 2023 (left) and 2024 (right) classification competitions.}
    \label{fig: data.application}
\end{figure}

In Figure~\ref{fig: data.application}, the left panel presents a comparison of the confidence set from one realization. BC and MCS can only exclude the most obvious inferior models. The RI approach is more powerful than BC and MCS, but still less so compared to the LOO method. Note that MCS can achieve the stronger simultaneous coverage (\Cref{remark: marginal coverate}). As a trade-off, it is expected to have less rejection power than methods aiming at marginal coverage, including the one in this work. The comparison among simultaneous coverage methods is given in Appendix~\ref{sec: heuristic simultaneous inference}. 

There is no strict monotonic relationship between the test error rate and exclusion from the confidence set. For example, models $21$ and $37$ possessed a higher empirical test error but were not excluded by the proposed LOO method. A similar pattern also held for model $21,36,37,39$ when implementing the RI method. This is due to the association between the compared submissions. 
When an inferior model $r$ (e.g., model 36) was strongly positively correlated with a superior model $s$, the variance of $X_{i,r}-X_{i,s}$ will be small. Compared with models with the same average performance but less associated with any superior candidates, rejecting the correlated inferior models from the confidence set is easier. In other words, the difficulty of excluding an inferior model depends on both the test-error difference and the covariance. 

\subsection{2024 Classification Competition}

In Spring 2024, a total number of $39$  models were evaluated upon a test data set of sample size $1236$, leading to a prediction discrepancy sample $\{X_i \in \{0, 1\}^{39}: i \in [1236]\}$.

Again, we constructed confidence sets 100 times over the same real testing data to account for the randomness in the algorithms. For the proposed LOO method, it yields an average size of 25.3 with models 1, 6, 8, 10, 12, 15, 16, 26, 28, 29, 30, 32, 38, and 39 frequently excluded. In comparison, the MCS procedure performs least favorably with an average size of $37$. The BC method results in an average size of $30$, while the RI approach achieves an average size of $28.7$. The former mostly excludes models 1, 8, 10, 12, 15, 16, 28, 30 and 32, while the latter excludes models 1, 8, 10, 12, 15, 16, 28, 30 and 32.  
Figure~\ref{fig: data.application}, right panel, provides an example of one realization. We can see that most approaches succeed in excluding the obviously inferior models 15 and 32, but the proposed method rejects more competitive ones, suggesting a better finite-sample statistical power.

\section{Discussion}
For a fixed component index $r$, our approach to testing $r\in\Theta$ is to test whether any of the differences $\delta_{r,s}\coloneqq \mu_r-\mu_s, s \ne r,$ is strictly positive.
If we take the differences in each sample point: $\tilde X_{i,s}=X_{i,r}-X_{i,s}$ such that $\mathbb E\tilde X = \delta_{r,s}$, then the problem is equivalent to testing whether $\max_s\delta_{r,s}>0$. The comparison between the proposed method and the rank inference-based method \citep{mogstad2024inference} reflects two sides in the bias-variance trade-off spectrum when constructing a confidence lower bound for the parameter of interest $\max_s \delta_{r,s}$. On the one hand, the exponentially weighted soft-max in our proposed method targets a biased version which is a weighted average of all $\delta_{r,s}, s \ne r$ with weights concentrating more on the large entries. Such a bias is traded in favor of a reduced variance: the resulting test statistic has variance on the scale of $1/n$ and enjoys asymptotic normality. Moreover, the bias has a known sign under the null hypothesis and hence does not affect the validity of the resulting confidence lower bound.  On the other hand, the method in \cite{mogstad2024inference} has no bias as it looks at the largest sample mean. Yet this comes at a cost of wider confidence intervals since it must account for the fluctuation of the maximum of a random vector.  Intuitively, a favorable scenario for the proposed method is when there are many ``irrelevant'' components and quite a few slightly sub-optimal components, as exemplified in the ``3-tier'' setting in \Cref{subsection: method comparison}.  A scenario that would favor the low-bias-high-variance method (e.g., \cite{mogstad2024inference}) is when the correlations among the components are very high, so the variability of the maximum of the random vector is not much larger than that of a single component.  This is exemplified in the LASSO tuning example in \Cref{section: lasso}.  However, even in this setting, the proposed method is still competitive, demonstrating its good adaptivity.

In this vein, an important and intriguing research question is to investigate optimal confidence intervals for the maximum value of a high-dimensional mean vector from IID noisy observations.  The results in the current paper suggest that the hardness of such confidence intervals must depend on the shape of the mean vector---especially the gap among the leading components---as well as the correlation between the components.  A special case is to find the most probable entry in a multinomial distribution, which has been well-studied in a sequential setting \citep{panchapakesan2006note}.

One practical motivation for argmin inference is its use in  model selection tasks.  Broadly speaking, the term ``model selection'' can refer to either identifying the underlying true model (assuming there is one), or finding the best approximation of the data generating mechanism. These two tasks are also known as \emph{model identification} and \emph{model estimation}, respectively \citep{arlot2010survey}.  
While it may be more straightforward to interpret model confidence sets in the context of model identification, a similar interpretation can also be extended to the model approximation case.

Take the example of hyperparameter tuning. Assume we are given $p$ hyperparameter values in a regression algorithm, each leading to an estimated regression function $\hat f_r$ that predicts a response $Y$ using covariate $Z$. Let $\ell(\cdot,\cdot)$ be a loss function. Then 
the best hyperparameter can be defined as the one that minimizes $\mathbb E\ell(\hat f_r(Z),Y)$, where the expectation is taken over both the randomness of $\hat f_r$ and the evaluation sample point $(Y,Z)$.
The model confidence set problem in this context becomes finding all the candidate hyperparameter values that produce (nearly) the best predictive risk.
The current proposal can be applied to the sample-splitting scheme where the fitted models come independently of the argmin inference data (the setting considered in Section~\ref{section: lasso}). However, as in the standard cross-validation, the training and evaluation samples may be (partially) swapped for selecting the best-performing hyperparameters. A collection of associated confidence sets $\{\widehat C_{m}, m \in [\mathcal{V}]\}$ can be obtained for a $\mathcal{V}$-fold cross-validation and one may consider combining them to get an overall confidence set for hyperparameter selection. Rigorously integrating these confidence sets is a challenging task due to the irregular distribution of the random center $d_{i,r}\dexcludingVi$ defined in \eqref{eq: centered test stat}, which will be pursued in future work. 

A natural future direction is to extend the current framework to achieve simultaneous coverage. We propose a heuristic solution in \Cref{sec: heuristic simultaneous inference} that shows promising empirical results, which is also direct to implement given the current software. In \cite{mogstad2024inference}, the authors shift focus from $\max_s \sigma_{r,s}$ to a ``double max” formulation $\max_{s \in J} \max_s \sigma_{r,s}$ for some $J \subseteq [p]$, and apply a stepdown procedure to construct simultaneous confidence sets. It remains an open question how to adapt the softmax mechanism in this paper in a similar fashion. As an alternative route for this extension, one could seek to establish a multivariate CLT in the case of fixed $p$---or even a high-dimensional Gaussian comparison result as $p$ diverges---for the vector $[T_1, \ldots, T_p]^\top$ with the same random centerings as in this work. By including all indices whose $T_r$ lies below the appropriate upper quantile of the corresponding max-statistic, one would obtain a confidence set with simultaneous coverage. We view this line of investigation as a particularly exciting and intriguing direction for future research.

Another promising direction, briefly noted in our literature review, is to systematically extend the asymptotic coverage to finite-sample settings. One of the few works addressing this challenge is \cite{dey2024anytime}, which develops an anytime-valid approach to constructing argmin confidence sets. 

\section*{Code Availability}

The main proposal is
implemented in R package {\tt argminCS}. Source code and vignette at: \url{https://github.com/xu3cl4/argminCS}. Differentially private versions of the data in Section~\ref{section: real_data} are also provided within the package. For reproducibility, we include the inference results based on the differentially private data in
Appendix~\ref{sec: differential private data}. 

\section*{Acknowledgments}

We thank the Reviewers and Associate Editors for their constructive comments, which helped to improve this work. Part of this work was completed while Tianyu Zhang was at Carnegie Mellon University. Tianyu Zhang's work at Carnegie Mellon University was partially supported by National Institute of Mental Health grant R01MH123184. Hao Lee and Jing Lei are partially supported by National Science Foundation grants DMS-2310764 and DMS-2515687.
\bibliographystyle{plainnat}
\bibliography{main}

@article{shao2025berry,
  title={Berry-Esseen bounds for functionals of independent random variables},
  author={Shao, Qi-Man and Zhang, Zhuo-Song},
  journal={Stochastic Processes and their Applications},
  volume={183},
  pages={104574},
  year={2025},
  publisher={Elsevier}
}

@article{chatterjee2008new,
author = {Sourav Chatterjee},
title = {{A new method of normal approximation}},
volume = {36},
journal = {The Annals of Probability},
number = {4},
publisher = {Institute of Mathematical Statistics},
pages = {1584 -- 1610},
year = {2008},
doi = {10.1214/07-AOP370},
URL = {https://doi.org/10.1214/07-AOP370}
}

@article{adrian2024stabilizing,
  title={Stabilizing black-box model selection with the inflated argmax},
  author={Adrian, Melissa and Soloff, Jake A and Willett, Rebecca},
  journal={arXiv preprint arXiv:2410.18268},
  year={2024}
}

@article{chen2023inference,
  title={Inference on Optimal Dynamic Policies via Softmax Approximation},
  author={Chen, Qizhao and Austern, Morgane and Syrgkanis, Vasilis},
  journal={arXiv preprint arXiv:2303.04416},
  year={2023}
}

@article{arnold2024sequential,
  title={Sequential model confidence sets},
  author={Arnold, Sebastian and Gavrilopoulos, Georgios and Schulz, Benedikt and Ziegel, Johanna},
  journal={arXiv preprint arXiv:2404.18678},
  year={2024}
}

@article{arlot2010survey,
  title={A survey of cross-validation procedures for model selection},
  author={Arlot, Sylvain and Celisse, Alain},
  journal={Statistics surveys},
  volume={4},
  pages={40--79},
  year={2010},
  publisher={Amer. Statist. Assoc., the Bernoulli Soc., the Inst. Math. Statist., and the~…}
}

@article{kim2025locally,
  title={Locally minimax optimal and dimension-agnostic discrete argmin inference},
  author={Kim, Ilmun and Ramdas, Aaditya},
  journal={arXiv preprint arXiv:2503.21639},
  year={2025}
}

@article{panchapakesan2006note,
  title={A note on a subset selection procedure for the most probable multinomial event},
  author={Panchapakesan, Subramanian},
  journal={Sequential Analysis},
  volume={25},
  number={2},
  pages={145--149},
  year={2006},
  publisher={Taylor \& Francis}
}

@article{taylor2015statistical,
  title={Statistical learning and selective inference},
  author={Taylor, Jonathan and Tibshirani, Robert J},
  journal={Proceedings of the National Academy of Sciences},
  volume={112},
  number={25},
  pages={7629--7634},
  year={2015},
  publisher={National Acad Sciences}
}

@inproceedings{mcsherry2007mechanism,
  title={Mechanism design via differential privacy},
  author={McSherry, Frank and Talwar, Kunal},
  booktitle={48th Annual IEEE Symposium on Foundations of Computer Science (FOCS'07)},
  pages={94--103},
  year={2007},
  organization={IEEE}
}

@article{tibshirani1996regression,
  title={Regression shrinkage and selection via the lasso},
  author={Tibshirani, Robert},
  journal={Journal of the Royal Statistical Society Series B: Statistical Methodology},
  volume={58},
  number={1},
  pages={267--288},
  year={1996},
  publisher={Oxford University Press}
}

@article{dwork2014algorithmic,
  title={The algorithmic foundations of differential privacy},
  author={Dwork, Cynthia and Roth, Aaron and others},
  journal={Foundations and Trends{\textregistered} in Theoretical Computer Science},
  volume={9},
  number={3--4},
  pages={211--407},
  year={2014},
  publisher={Now Publishers, Inc.}
}

@book{gibbons1977selecting,
  title={Selecting and ordering populations: A New Statistical Methodology},
  author={Gibbons, Jean Dickinson and Olkin, Ingram and Sobel, Milton},
  year={1977},
  publisher={Wiley},
  city={New York}
}

@book{gupta1979multiple,
  title={Multiple decision procedures: theory and methodology of selecting and ranking populations},
  author={Gupta, Shanti S and Panchapakesan, Subramanian},
  year={1979},
  publisher={Wiley},
  city={New York}
}

@article{mogstad2024inference,
  title={Inference for ranks with applications to mobility across neighbourhoods and academic achievement across countries},
  author={Mogstad, Magne and Romano, Joseph P and Shaikh, Azeem M and Wilhelm, Daniel},
  journal={Review of Economic Studies},
  volume={91},
  number={1},
  pages={476--518},
  year={2024},
  publisher={Oxford University Press US}
}

@article{tanguy2015some,
  title={Some superconcentration inequalities for extrema of stationary Gaussian processes},
  author={Tanguy, Kevin},
  journal={Statistics \& Probability Letters},
  volume={106},
  pages={239--246},
  year={2015},
  publisher={Elsevier}
}

@article{kamath2015bounds,
  title={Bounds on the expectation of the maximum of samples from a gaussian},
  author={Kamath, Gautam},
  journal={Technical Report},
  year={2015},
URL = {http://www. gautamkamath. com/writings/gaussian max. pdf}
}

@article{goldstein1996league,
  title={League tables and their limitations: statistical issues in comparisons of institutional performance},
  author={Goldstein, Harvey and Spiegelhalter, David J},
  journal={Journal of the Royal Statistical Society Series A: Statistics in Society},
  volume={159},
  number={3},
  pages={385--409},
  year={1996},
  publisher={Oxford University Press}
}

@misc{dey2024anytime,
      title={Generalized Universal Inference on Risk Minimizers}, 
      author={Neil Dey and Ryan Martin and Jonathan P. Williams},
      year={2024},
      eprint={2402.00202},
      archivePrefix={arXiv},
      primaryClass={stat.ME},
      url={https://arxiv.org/abs/2402.00202}, 
}

@article{bernardi2018model,
  title={The model confidence set package for R},
  author={Bernardi, Mauro and Catania, Leopoldo},
  journal={International Journal of Computational Economics and Econometrics},
  volume={8},
  number={2},
  pages={144--158},
  year={2018},
  publisher={Inderscience Publishers (IEL)}
}

@article{kleywegt2002sample,
  title={The sample average approximation method for stochastic discrete optimization},
  author={Kleywegt, Anton J and Shapiro, Alexander and Homem-de-Mello, Tito},
  journal={SIAM Journal on optimization},
  volume={12},
  number={2},
  pages={479--502},
  year={2002},
  publisher={SIAM}
}

@article{Hung2019rank,
author = {Kenneth Hung and William Fithian},
title = {{Rank verification for exponential families}},
volume = {47},
journal = {The Annals of Statistics},
number = {2},
publisher = {Institute of Mathematical Statistics},
pages = {758 -- 782},
year = {2019},
doi = {10.1214/17-AOS1634},
URL = {https://doi.org/10.1214/17-AOS1634}
}

@article{xie2009confidence,
  title={Confidence intervals for population ranks in the presence of ties and near ties},
  author={Xie, Minge and Singh, Kesar and Zhang, Cun-Hui},
  journal={Journal of the American Statistical Association},
  volume={104},
  number={486},
  pages={775--788},
  year={2009},
  publisher={Taylor \& Francis}
}

@article{hall2009,
  author = {Hall, Peter and Miller, Henry},
  year = {2009},
  title = {Using the bootstrap to quantify the authority of an empirical ranking},
  journal = {The Annals of Statistics},
  volume = {37},
  number = {6B},
  pages = {3929--3959},
  doi = {10.1214/09-AOS699}
}

@article{Choirat2012estimation,
author = {Christine Choirat and Raffaello Seri},
title = {{Estimation in Discrete Parameter Models}},
volume = {27},
journal = {Statistical Science},
number = {2},
publisher = {Institute of Mathematical Statistics},
pages = {278 -- 293},
year = {2012},
doi = {10.1214/11-STS371},
URL = {https://doi.org/10.1214/11-STS371}
}

@article{seri2021model,
  title={Model calibration and validation via confidence sets},
  author={Seri, Raffaello and Martinoli, Mario and Secchi, Davide and Centorrino, Samuele},
  journal={Econometrics and Statistics},
  volume={20},
  pages={62--86},
  year={2021},
  publisher={Elsevier}
}

@article{futschik1995confidence,
  title={Confidence sets for discrete stochastic optimization},
  author={Futschik, Andreas and Pflug, Georg},
  journal={Annals of Operations Research},
  volume={56},
  pages={95--108},
  year={1995},
  publisher={Springer}
}

@article{wasserman2014stein,
  title={Stein’s method and the bootstrap in low and high dimensions: A tutorial},
  author={Wasserman, Larry},
  year={2014},
  publisher={Working paper}
}

@article{chin2022short,
  title={A short and elementary proof of the central limit theorem by individual swapping},
  author={Chin, Calvin Wooyoung},
  journal={The American Mathematical Monthly},
  volume={129},
  number={4},
  pages={374--380},
  year={2022},
  publisher={Taylor \& Francis}
}

@book{vershynin2018high,
  title={High-dimensional probability: An introduction with applications in data science},
  author={Vershynin, Roman},
  volume={47},
  year={2018},
  publisher={Cambridge university press}
}

@article{hansen2011model,
  title={The model confidence set},
  author={Hansen, Peter R and Lunde, Asger and Nason, James M},
  journal={Econometrica},
  volume={79},
  number={2},
  pages={453--497},
  year={2011},
  publisher={Wiley Online Library}

}

@article{romano1988bootstrapping,
  title={Bootstrapping the mode},
  author={Romano, Joseph P},
  journal={Annals of the Institute of Statistical Mathematics},
  volume={40},
  pages={565--586},
  year={1988},
  publisher={Springer}
}

@inproceedings{bayle2020, author = {Bayle, Pierre and Bayle, Alexandre and Janson, Lucas and Mackey, Lester}, title = {Cross-validation confidence intervals for test error}, year = {2020}, isbn = {9781713829546}, publisher = {Curran Associates Inc.}, address = {Red Hook, NY, USA}, booktitle = {Proceedings of the 34th International Conference on Neural Information Processing Systems}, articleno = {1371}, numpages = {12}, location = {Vancouver, BC, Canada}, series = {NIPS'20} }

@article{austern2020,
author = {Morgane Austern and Wenda Zhou},
title = {{Asymptotics of cross-validation}},
volume = {published online},
journal = {Annales de l'Institut Henri Poincaré, Probabilités et Statistiques},
number = {},
publisher = {Institut Henri Poincaré},
year = {2024},
doi = {},
URL = {https://www.e-publications.org/ims/submission/AIHP/user/submissionFile/62392?confirm=f28e3906}
}

@misc{kissel2023blackbox,
      title={Black-Box Model Confidence Sets Using Cross-Validation with High-Dimensional Gaussian Comparison}, 
      author={Nicholas Kissel and Jing Lei},
      year={2023},
      eprint={2211.04958},
      archivePrefix={arXiv},
      primaryClass={math.ST},
      url={https://arxiv.org/abs/2211.04958}, 
}

@article{fan2024ranking,
author = {Fan, Jianqing and Lou, Zhipeng and Wang, Weichen and Yu, Mengxin},
title = {Ranking Inferences Based on the Top Choice of Multiway Comparisons},
journal = {Journal of the American Statistical Association},
volume = {},
number = {},
pages = {1--14},
year = {2024},
publisher = {ASA Website},
doi = {10.1080/01621459.2024.2316364},
URL = {https://doi.org/10.1080/01621459.2024.2316364},
eprint ={https://doi.org/10.1080/01621459.2024.2316364}}

@article{chernozhukov2013gaussian,
  author = {Chernozhukov, V. and Chetverikov, D. and Kato, K.},
  title = {GAUSSIAN APPROXIMATIONS AND MULTIPLIER BOOTSTRAP FOR MAXIMA OF SUMS OF HIGH-DIMENSIONAL RANDOM VECTORS},
  journal = {The Annals of Statistics},
  volume = {41},
  number = {6},
  pages = {2786--2819},
  year = {2013},
  doi = {10.1214/13-aos1161},
  url = {https://doi.org/10.1214/13-aos1161}
}

@article{chernozhukov2014gaussian,
  author = {Chernozhukov, V. and Chetverikov, D. and Kato, K.},
  title = {GAUSSIAN APPROXIMATION OF SUPREMA OF EMPIRICAL PROCESSES},
  journal = {The Annals of Statistics},
  volume = {42},
  number = {4},
  pages = {1564--1597},
  year = {2014},
  doi = {10.1214/14-aos1230},
  url = {https://doi.org/10.1214/14-aos1230}
}

@article{chernozhukov2017central,
  author = {Chernozhukov, V. and Chetverikov, D. and Kato, K.},
  title = {CENTRAL LIMIT THEOREMS AND BOOTSTRAP IN HIGH DIMENSIONS},
  journal = {The Annals of Probability},
  volume = {45},
  number = {4},
  pages = {2309--2352},
  year = {2017},
  doi = {10.1214/16-aop1113},
  url = {https://doi.org/10.1214/16-aop1113}
}

@article{gupta1965some,
  title={On some multiple decision (selection and ranking) rules},
  author={Gupta, Shanti S},
  journal={Technometrics},
  volume={7},
  number={2},
  pages={225--245},
  year={1965},
  publisher={Taylor \& Francis}
}

@article{romano2014practical,
  title={A practical two-step method for testing moment inequalities},
  author={Romano, Joseph P and Shaikh, Azeem M and Wolf, Michael},
  journal={Econometrica},
  volume={82},
  number={5},
  pages={1979--2002},
  year={2014},
  publisher={Wiley Online Library}
}

@techreport{chernozhukov2016testing,
  title={Testing many moment inequalities},
  author={Chernozhukov, Victor and Chetverikov, Denis and Kato, Kengo},
  year={2016},
  institution={cemmap working paper}
}

@Manual{R-csranks,
    title = {csranks: Statistical Tools for Ranks},
    author = {Daniel Wilhelm and Pawel Morgen},
    url = {https://danielwilhelm.github.io/R-CS-ranks/},
    year = {2023},
}

@article{Bentkus2007,
author = {Bentkus, Vidmantas and Jing, Bing-Yi and Shao, Qi-Man and Zhou, Wang},
address = {VOORBURG},
copyright = {Copyright 2007 International Statistical Institute/Bernoulli Society},
issn = {1350-7265},
journal={Bernoulli},
language = {eng},
number = {2},
pages = {346-364},
publisher = {International Statistics Institute / Bernoulli Society},
title = {Limiting Distributions of the Non-Central T-Statistic and Their Applications to the Power of T-Tests under Non-Normality},
volume = {13},
year = {2007},
}

@article{glmnet,
  title={Package ‘glmnet’},
  author={Friedman, Jerome and Hastie, Trevor and Tibshirani, Robert and Narasimhan, Balasubramanian and Tay, Kenneth and Simon, Noah and Qian, Junyang},
  journal={CRAN R Repositary},
  volume={595},
  year={2021}
}

@book{cramer1999mathematical,
  title={Mathematical methods of statistics},
  author={Cram{\'e}r, Harald},
  volume={9},
  year={1999},
  publisher={Princeton university press}
}

@article{korolev2017bounds,
  title={Bounds of the accuracy of the normal approximation to the distributions of random sums under relaxed moment conditions},
  author={Korolev, Victor and Dorofeeva, Alexandra},
  journal={Lithuanian Mathematical Journal},
  volume={57},
  number={1},
  pages={38--58},
  year={2017},
  publisher={Springer}
}

\newpage
\appendix



\section{A Sign-flipping Lemma}

\begin{lemma}\label{lemma: sign_flipping}
Let $Z_1, \dots, Z_n \overset{\text{iid}}{\sim} \mathcal{N}(0, 1)$. Let $j \in [n]$ a fixed index, define
\begin{equation*}
\mu = \frac{1}{n} \sum_{i=1}^n Z_i, \quad v = \mu - \frac{1}{n} Z_j, \quad \mu^{(j)} = v + \frac{1}{n} Z_j',
\end{equation*}
where $Z_j' \sim \mathcal{N}(0, 1)$ is an independent copy of $Z_j$. Then the probability of a sign flip satisfies
\begin{equation*}
cn^{-1/2}\leq \mathbb{P}(\mu > 0, \mu^{(j)} < 0) + \mathbb{P}(\mu < 0, \mu^{(j)} > 0)\leq Cn^{-1/2},
\end{equation*}
for some constants $c,C > 0$.
\end{lemma}

We provide two proofs of this lemma. 
\begin{proof}[First proof]
Note that $\mu = v + n^{-1} Z_j$ and $\mu^{(j)} = v + n^{-1} Z_j'$. Therefore,
\[
\mathbb{P}(\mu > 0, \mu^{(j)} < 0) = \mathbb{P}\left(Z_j > -n v, \, Z_j' < -n v\right).
\]
We analyze this probability by partitioning the support of $v$. Since $v \sim \mathcal{N}(0, \sigma^2)$ with $\sigma^2 = (n-1)/n^2$, it holds that
\[
\mathbb{P}(-n^{-1} < v < n^{-1}) = \mathbb{P}\left(|N(0,1)| < (n-1)^{-1 / 2} \right) \asymp n^{-1/2}.
\]

\textbf{Lower bound:} On the event $|v| < n^{-1}$, we have
\begin{align*}
&\mathbb{P}(Z_j > -n v, \, Z_j' < -n v, \, |v| < n^{-1}) \\
&= \mathbb{P}(Z_j > -n v, \, Z_j' < -n v \mid |v| < n^{-1}) \cdot \mathbb{P}(|v| < n^{-1}) \\
&\geq \mathbb{P}(Z_j > 1, \, Z_j' < -1) \cdot \mathbb{P}(|v| < n^{-1}) \\
&\gtrsim n^{-1/2}.
\end{align*}

\textbf{Upper bound:} For $k \geq 1$, define intervals $I_k = \left(\frac{k}{n}, \frac{k+1}{n}\right)$. Then
\begin{align*}
&\sum_{k=1}^{\infty} \mathbb{P}\left(Z_j > -n v, \, Z_j' < -n v, \, v \in I_k\right)\\
&\leq \sum_{k=1}^\infty \mathbb{P}\left(Z_j' < -n v, \, v \in I_k\right) \\
&\leq \sum_{k=1}^\infty \mathbb{P}(Z_j' < -k) \cdot \mathbb{P}\left(v \in I_k\right) \\
&\leq \sum_{k=1}^\infty \exp\left(-\frac{k^2}{2}\right) \cdot \mathbb{P}\left(k(n-1)^{-1/2} < N(0,1) < (k+1)(n-1)^{-1/2}\right) \\
&\leq \sum_{k=1}^\infty \exp\left(-\frac{k^2}{2}\right) \cdot \frac{1}{\sqrt{2\pi(n-1)}} \exp\left(-\frac{k^2}{2(n-1)}\right) \\
&\lesssim \frac{1}{\sqrt{n}}.
\end{align*}
Since 
\begin{equation*}
\begin{aligned}
&\mathbb{P}\left(Z_j>-n v, Z_j^{\prime}<-n v\right) \\
=&\mathbb{P}\left(Z_j>-n v, Z_j^{\prime}<-n v,-n^{-1}<v<n^{-1}\right) \\
&\hspace{2mm} +\sum_{k=1}^{\infty} \mathbb{P}\left(Z_j>-n v, Z_j^{\prime}<-n v, k / n<v<(k+1) / n\right) \\
&\hspace{2mm} +\sum_{k=-\infty}^{-1} \mathbb{P}\left(Z_j>-n v, Z_j^{\prime}<-n v, k / n<v<(k+1) / n\right),
\end{aligned}
\end{equation*}
we know $\mathbb{P}\left(Z_j>-n v, Z_j^{\prime}<-n v\right)\lesssim n^{-1/2}$.

By symmetry, $\mathbb{P}(\mu < 0, \mu^{(j)} > 0) = \mathbb{P}(\mu > 0, \mu^{(j)} < 0)$. Therefore,
\[
\mathbb{P}(\mu > 0, \mu^{(j)} < 0) + \mathbb{P}(\mu < 0, \mu^{(j)} > 0) \asymp n^{-1/2}.
\]
\end{proof}

\begin{proof}[Second proof]

We note that
\begin{equation}
\label{eq: probabilties sum to 1}
\begin{aligned}
    & \prob{\mu > 0, \mu^{(j)} < 0} + \prob{\mu > 0, \mu^{(j)} > 0} \\
    & = \prob{\mu > 0, \mu^{(j)} < 0} + \prob{\mu > 0, \mu^{(j)} \ge 0} \\
    &= \prob{\mu > 0 } = \frac{1}{2}. 
\end{aligned}
\end{equation}
Because $Z_1, \ldots, Z_n, Z_j'$ are jointly normal, we know that $(\mu, \mu^{(j)})$, as their linear combination, is a bivariate normal random vector. In particular, we have $\mathrm{Var}(\mu) = 1/n$ and $\mathrm{Var} (\mu^{(j)}) = 1/n$. Moreover, letting $A = \frac{1}{n}\sum_{i \ne j} Z_i, B = \frac{1}{n} Z_j$ and $C = \frac{1}{n} Z_j'$, we have 
    \begin{align*}
    &\mathrm{cov}(\mu, \mu^{(j)})  \\
    &= \mathrm{cov} (A + B, A + C - B) \\
    &= \mathrm{cov} (A + B, A + B+  C - 2B) \\
    &= \mathrm{Var}(A + B)- 2 \mathrm{Var}(B) \\
    &= \frac{1}{n} - \frac{2}{n^2}. 
    \end{align*}
It follows that $(\mu, \mu^{(j)})$ has the correlation coefficient 
$$
    \rho := \frac{\mathrm{cov}(\mu, \mu^{(j)})}{\sqrt{\mathrm{Var}(\mu)} \sqrt{\mathrm{Var}(\mu^{(j)})}} = \frac{1/n - 2/n^2}{1/n} = 1 - \frac{2}{n}. 
$$
The orthant probability $\prob{\mu > 0, \mu^{(j)} > 0}$ has the closed form (see page 290 in~\cite{cramer1999mathematical} for example):
\begin{equation}
\label{eq: orthant prob}
\begin{aligned}
    &\frac{1}{4} + \frac{1}{2 \pi} \arcsin(1 - \frac{2}{n}) \\
    &= \frac{1}{4} + \frac{1}{2 \pi} \arctan{\frac{1 - 2/n}{2\sqrt{1/n - 1/n^2}}} \\
    &= \frac{1}{4} + \frac{1}{2\pi} \arctan{\frac{\sqrt{n} - 2/\sqrt{n}}{2 \sqrt{1 - 1/n}}}. 
\end{aligned}
\end{equation}
Define $\beta_n := (\sqrt{n} - 2/\sqrt{n})/(2 \sqrt{1 - 1/n}) \asymp \sqrt{n}$ for simplicity. Combining~\eqref{eq: probabilties sum to 1} and~\eqref{eq: orthant prob}, we have 
\begin{align*}
    &\prob{\mu > 0, \mu^{(j)} < 0} = \frac{1}{4} - \frac{1}{2\pi}\arctan{\beta_n}= \frac{1}{2\pi} (\frac{\pi}{2} - \arctan{\beta_n})\\
    &= \frac{1}{2\pi} \mathrm{arccot}(\beta_n) \asymp \frac{1}{2\pi \beta_n} \asymp \frac{1}{\sqrt{n}}, 
\end{align*}
where one can verify the first asymptotic equivalency by applying L'Hopital's rule. By symmetry, $\prob{\mu < 0, \mu^{(j)} > 0} = \prob{\mu > 0, \mu^{(j)} < 0}$, from which we can conclude the proof. 
\end{proof}

\section{Variance Estimation}\label{app: variance estimation}

This section concerns the proofs for the consistency of the variance estimator $\hat \sigma^2_{r}$ in~\eqref{eq: variance out} and its relevant results. In fact, we show its consistency to a ``variety'' of $\sigma^2_r$, and justify the asymptotic equivalence between the two by stability. 

\begin{proposition} \label{prop: equivalency among variances}
Distinguish the following three quantities:
\begin{align*}
    &\tilde{\sigma}_r^2 = \Enb\mathrm{Var}[X_{1, r} - Q_{1, r} | \bm{X}^{(-v_1)}]; \\ 
    &\undertilde{\sigma}_r^2 = \mathrm{Var}( \E{X_{1,r} - Q_{1,r}\QexcludingVone | X_1}); \\ 
    &\sigma_r^2 = \mathrm{Var}(X_{1, r} - Q_{1, r}\QexcludingVi).
\end{align*}
Under the same assumptions as Theorem~\ref{th: random center}, the three quantities are asymptotically equivalent, i.e., $\abs{\tilde{\sigma}_r^2 - \sigma_r^2} \rightarrow 0$ and $\abs{\undertilde{\sigma}_r^2 - \sigma_r^2} \rightarrow 0$. Moreover, suppose that $\mathrm{cov}(X_1)$ is positive definite so that the variance $\sigma_r^2$ is bounded away from $0$ for all $n \in \N$. Then, these results further imply that $\tilde{\sigma}_r^2/\sigma_r^2 \rightarrow 1$ and $\undertilde{\sigma}_r^2/\sigma_r^2 \rightarrow 1$. 
\end{proposition}

\begin{remark} \label{rem: variances relation}
    The quantity $\sigma^2_r$ is the population variance of the statistic $X_{1,r} - Q_{1,r}\QexcludingVone$ which contains two critical components. One is the `center' $X_1$ which plays a role in determining the difference, and the other is the exponential weightings, derived from the out-of-fold data, which helps determine how the weighted average $Q_{1,r}\QexcludingVone$ gets computed. Intuitively, one can imagine when $n$ is sufficiently large, the dependence via exponential weightings would be weak enough so that the variance across $X_i$'s has contributed to a large source of variance in $\sigma^2_r$. This is essentially because the exponential weightings are computed from the out-of-fold mean which converges to a fixed vector. This intuition is justified by Proposition~\ref{prop: equivalency among variances}. Indeed, we know $\sigma^2_r = \tilde{\sigma}_r^2 + \delta$ with $\delta = \mathrm{Var}(\E{X_{1, r} - Q_{1, r}\QexcludingVone| \bm{X}^{(-v_1)}})$. The quantity $\delta$ captures the variance contributed by $\bm{X}^{(-v_1)}$ and the proposition shows that $\delta$ is asymptotically negligible. A similar conclusion can be made when we condition on $X_1$. 
    
\end{remark}

Now recall the definition of $\hat{\sigma}^2_r$:
$$
    \hat{\sigma}^2_r = \frac{1}{n} \sum_{i = 1}^n \left (X_{i, r} - Q_{i, r}\QexcludingVi - \frac{1}{n}\sum_{\idummy = 1}^n(X_{\idummy, r} - Q_{\idummy,r}\QexcludingVdummy) \right )^2.  
$$
This estimator has recently garnered attention for its role in exploring uncertainty quantification in cross-validation. Particularly,~\cite{bayle2020} has studied the consistency of its variant under different assumptions and notions of stability, such as mean-square stability and loss stability. The estimator is simply the sample variance of all the differences $X_{i,r} - Q_{i,r}\QexcludingVi$, which makes the definition intuitive on its own. However, we should emphasize the dependency among the differences in contrast to the classical sample variance of IID data. To clarify, their dependency are present not only within the differences centered on samples in the same fold, but also across all folds due to the overlap in out-of-fold data used for exponential weightings. 
As illustrated in Remark~\ref{rem: variances relation}, one can infer that the weak dependency aims $\hat{\sigma}^2_r$ to behave similarly as the sample variance for IID data (both yield consistency) although its existence might make a proof non-trivial.  

\begin{theorem}[Theorem~\ref{th: consistency of variance estimator}] \label{cor: consistency of all-pair estimator}
    Under the same assumptions as Theorem~\ref{th: random center}, we have that $\abs{\hat{\sigma}^2_r - \undertilde{\sigma}_{r}^2} \cinp 0$. In particular, this implies that $\hat{\sigma}^2_r/\sigma_r^2$ converges to $1$ in probability.  
\end{theorem}

One may, in turn, consider the estimator: 
\begin{equation} \label{eq: variance in}
    \hat{\varphi}_{r}^{2} = \frac{1}{V} \sum_{v =1}^V \frac{1}{(n/V) - 1} \sum_{i \in I_{v}} \left ( X_{i,r} - Q_{i, r}\QexcludingV - \frac{V}{n} \sum_{\idummy \in I_v} (X_{\idummy, r} - Q_{\idummy, r}\QexcludingV ) \right )^2.
\end{equation}
This estimator's variants have been studied in the literature concerning cross-validation (see~\cite{austern2020, bayle2020,  kissel2023blackbox}). It is simply a sample average of sample variances, where each sample variance is computed from the differences $X_{i,r} - Q_{i,r}\QexcludingVi$ centered on the samples within the same fold. 
Its definition is motivated by $\tilde{\sigma}_r^2$ in Proposition~\ref{prop: equivalency among variances}. Nonetheless, we adapt the proof in~\cite{bayle2020} to show its consistency to $\undertilde{\sigma}_r^2$ with an eye to highlighting the asymptotic equivalence between the two population quantities $\tilde{\sigma}_r^2$ and $\undertilde{\sigma}_r^2$. As a side note, this estimator $\hat{\varphi}_{r}^{2}$ cannot be applied to the LOO setting as it would be exactly $0$ otherwise.  

\begin{theorem} \label{prop: consistency of within-fold sample variance}
    Under the same assumptions as Theorem~\ref{th: random center}, we have that $\abs{\hat{\varphi}_{r}^{2} - \undertilde{\sigma}_r^2} \cinp 0$. 
    In particular, this implies that $\hat{\varphi}_{r}^{2}/\sigma_r^2$ converges to $1$ in probability.  
\end{theorem}

\subsection{Proof of Proposition~\ref{prop: equivalency among variances}}
\begin{proof}
    \textbf{Prove $\abs{\tilde{\sigma}_r^2 - \sigma_r^2} \rightarrow 0.$}
    
    The difference between $\tilde{\sigma}_r^2$ and $\sigma^2_r$ is $\delta := \mathrm{Var}(\E{X_{1, r} - Q_{1, r}\QexcludingVone| \bm{X}^{(-v_1)}})$. Let $i \notin I_{v_1}$ be arbitrary. We have 
    \begin{align*}
        &\Enb \left ( \E{X_{1,r} - Q_{1, r}\QexcludingVone | \bm{X}^{(-v_1)}} - \E{ X_{1,r}-Q_{1,r}\QexcludingVoneperturb|\bm{X}^{(-v_1), i}}\right)^2 \\
        &= \Enb \left ( \E{X_{1,r} - Q_{1, r}\QexcludingVone - (X_{1,r} -Q_{1,r}\QexcludingVoneperturb) | \bm{X}^{(-v_1)}, X'_i}\right )^2 \\
        &\le \Enb (Q_{1, r}\QexcludingVone - Q_{1,r}\QexcludingVoneperturb)^2
    \end{align*}
    by the Jensen's inequality. Modifying the stability result~\eqref{eq: bound Q difference} (we only bound the difference $\abs{Q_{i,r}\QexcludingVi - Q_{i, r}\QexcludingViperturbdummy}$), one can conclude from the Efron-Stein inequality that 
    $\delta = o(1)$. 

\textbf{Prove $\abs{\undertilde{\sigma}_r^2 - \sigma_r^2} \rightarrow 0.$}

The proof of this statement is essentially the same as the one above. The difference between $\undertilde{\sigma}^2_r$ and $\sigma_r^2$ is $\Enb \mathrm{Var} (X_{1,r} - Q_{1,r} \mid X_1)$. By the (conditional) Efron-Stein inequality (Lemma 1 in~\cite{bayle2020}), one can obtain 
\begin{equation*}
    \begin{aligned}
    \abs{\undertilde{\sigma}_r^2 - \sigma_r^2}
        & = \abs{\Enb \mathrm{Var}(X_{1,r} - Q_{1, r}\QexcludingVone | X_1)}\\
        &\le \E{ \frac{1}{2}\sum_{j \notin I_{v_1}} \E{(Q_{1,r} - Q_{1,r}^j)^2 | X_1}} \\
        &= \frac{1}{2} \sum_{j \notin I_{v_1}} \Enb(Q_{1,r} - Q_{1,r}^j)^2 \\
        &= o(1),
    \end{aligned}
\end{equation*}
where one can achieve the last equality by modifying the stability result~\eqref{eq: bound Q difference} (we only bound the difference $\abs{Q_{i,r}\QexcludingVi - Q_{i, r}\QexcludingViperturbdummy}$). 
\end{proof}

\subsection{Proof of Theorem~\ref{cor: consistency of all-pair estimator}/Theorem~\ref{th: consistency of variance estimator}}
\begin{proof}

For any $i \in [n]$, define $D_{i,r}\DexcludingVi = X_{i,r} - Q_{i,r}$ and $\overline{D}_{r} = \frac{1}{n} \sum_{i = 1}^n D_{i,r}\DexcludingVi$. Under this notation, we can rewrite $\hat \sigma_{r}^2$ as follows:
\begin{equation}
    \begin{aligned}
    \hat{\sigma}_{r}^2 &= \frac{1}{n} \sum_{i=1}^n \left (D_{i,r}\DexcludingVi - \overline{D}_{r} \right )^2\\
    & = \frac{1}{2n} \sum_{i =1}^n \left (D_{i,r}\DexcludingVi - \overline{D}_{r} \right)^2 + \frac{1}{2n} \sum_{\idummy=1}^n \left (D_{\idummy,r}\DexcludingVdummy - \overline{D}_{r} \right)^2 \\
    &= \frac{1}{2n} \sum_{i =1}^n \left (D_{i,r}\DexcludingVi - \overline{D}_{r} \right)^2 + \frac{1}{2n} \sum_{\idummy=1}^n \left (D_{\idummy,r}\DexcludingVdummy - \overline{D}_{r} \right)^2 \\
    & \quad - \frac{1}{n^2} \sum_{i = 1}^n \left (D_{i,r}\DexcludingVi - \overline{D}_r \right) \cdot \sum_{\idummy=1}^n \left(D_{\idummy,r}\DexcludingVdummy - \overline{D}_r \right) \\
    &= \frac{1}{2n^2} \sum_{i, \idummy = 1}^n \left( D_{i,r}\DexcludingVi - D_{\idummy,r}\DexcludingVdummy\right)^2.
    \end{aligned}
\end{equation}

To prove the desired result, it suffices to show $\Enb \abs{\hat{\sigma}_{r}^2- \undertilde{\sigma}_r^2} \rightarrow 0$ thanks to Proposition~\ref{prop: equivalency among variances}. 

\textbf{Part 1: split the difference  $\left|\hat{\sigma}_{r}^2-\undertilde{\sigma}_{r}^2\right| $ into three parts, and bound each separately.}

We split the difference into three parts:
\begin{equation} \label{eq: out-sample variance bound}
    \Enb \abs{\hat{\sigma}_{r}^2 - \undertilde{\sigma}_r^2} \le \Enb \abs{\hat{\sigma}_{r}^2 - \intermediatevarianceone} + \Enb \abs{\intermediatevarianceone - \intermediatevariancetwo} + \Enb \abs{\intermediatevariancetwo - \undertilde{\sigma}_r^2}. 
\end{equation}
The quantity $\intermediatevarianceone$ is similar to the reformulated $\hat \sigma_{r}^2$:
\begin{equation}
    \intermediatevarianceone = \frac{1}{2n^2} \sum_{i ,\idummy = 1}^n \left (D_{i,r}\DexcludingViperturbdummy - D^{(-v_i,\idummy)}_{\idummy,r} \right)^2,
\end{equation}
where the variable $D_{i,r}\DexcludingViperturbdummy = X_{i,r} - \sum_{s \ne r} \hat{w}^{(-v_i), \idummy}\wscenterr X_{i,s}$ in the summand is $D_{i,r}\DexcludingVi$ with $X_{\idummy}$ replaced by an IID copy $X'_{\idummy}$ if $l\notin I_{v_i}$. Otherwise, if $i, \idummy$ are within the same fold, the calculation of $D_{i,r}\DexcludingViperturbdummy$ would not involve $X_\idummy$, and therefore $D_{i,r}\DexcludingViperturbdummy$ is simply identical to $D_{i,r}\DexcludingVi$. 

The second quantity $D^{(-v_i,\idummy)}_{\idummy,r}$ in the summand is defined by $D^{(-v_i,\idummy)}_{\idummy,r} = X_{\idummy,r} - \sum_{s \ne r} \hat{w}^{(-v_i), \idummy}\wscenterr X_{\idummy,s}$. We stress that its exponential weightings are computed from $\bm{X}^{(-v_i), \idummy}$ rather than the out-of-fold data $\bm{X}^{(-v_{\idummy})}$, but both of them are independent of $X_\idummy$ and are identically distributed. Here the sample perturbation for $X_\idummy$ only occurs for exponential weightings, reflected in our choice of notation. In particular, if $i,\idummy$ and within the same fold, we have $\bm{X}^{(-v_i), \idummy} = \bm{X}^{(-v_i)}$ as $X_\idummy \notin \bm{X}^{(-v_i)} = \bm{X}^{(-v_\idummy)}$, and thereby $D^{(-v_i,\idummy)}_{\idummy,r} = D_{\idummy,r}$.  
This construction is for the purpose of making $D_{i,r}\DexcludingViperturbdummy$ and $D^{(-v_i,\idummy)}_{\idummy,r}$ share the same exponential weightings for every $(i, \idummy) \in [n]^2$ so that given the shared exponential weightings, $D_{i,r}\DexcludingViperturbdummy$ and $D^{(-v_i,\idummy)}_{\idummy,r}$ are 
identically distributed. In particular, this implies 
\begin{equation} \label{eq: same conditional expectation}
\E{D_{i,r}\DexcludingViperturbdummy \longconditional \bm{X}^{(-v_i), \idummy}} = \E{D^{(-v_i,\idummy)}_{\idummy,r} \longconditional \bm{X}^{(-v_i), \idummy}}.
\end{equation}

The other quantity $\intermediatevariancetwo$ in \eqref{eq: out-sample variance bound} is
\begin{equation}
\intermediatevariancetwo = \frac{1}{2n^2} \sum_{i,\idummy=1}^n \left ( \E{K_{i,r} \longconditional  X_i}  - \E{K_{\idummy,r} \longconditional X_{\idummy}} \right)^2,
\end{equation}
where $K_{i,r}$ is defined as in Lemma~\ref{lemma: first order stability} for all $i \in [n], r \in [p]$. 
By definition, $\Enb[K_{i,r}] = 0$. As the last note, the uniform boundedness of $X_{1}$ ensures that there exists $C > 0$ such that $\intermediatevarianceone < C$ and $\intermediatevariancetwo < C$. 

\textbf{Part 2: bound $\mathbb{E}\left|\hat{\sigma}_{r}^2-\intermediatevarianceone\right|.$}

It follows from simple algebra and the Cauchy-Schwartz inequality that 
\begin{equation} \label{eq: approximate out-sample variance 1}
\begin{aligned}
    &\Enb | \sigma_{r}^2 - \intermediatevarianceone | \\
    &= \Enb \hspace{0.6mm}\Bigg | \frac{1}{2n^2} \sum_{i,\idummy=1}^n \left( D_{i,r}\DexcludingVi - D_{\idummy,r}\DexcludingVdummy -  D_{i,r}\DexcludingViperturbdummy + D^{(-v_i,\idummy)}_{\idummy,r} \right) \\
    & \hspace{16mm} \times \left( D_{i,r}\DexcludingVi - D_{\idummy,r}\DexcludingVdummy + D_{i,r}\DexcludingViperturbdummy - D^{(-v_i,\idummy)}_{\idummy,r} \right) \Bigg | \\
    &= \Enb \hspace{0.6mm}\Bigg | \frac{1}{2n^2} \sum_{i,\idummy=1}^n \left( D_{i,r}\DexcludingVi - D_{\idummy,r}\DexcludingVdummy -  D_{i,r}\DexcludingViperturbdummy + D^{(-v_i,\idummy)}_{\idummy,r} \right)^2 \\
    &\hspace{8mm}+ \frac{1}{2n^2} \sum_{i,\idummy=1}^n 2 \left( D_{i,r}\DexcludingVi - D_{\idummy,r}\DexcludingVdummy -  D_{i,r}\DexcludingViperturbdummy + D^{(-v_i,\idummy)}_{\idummy,r} \right) \\
    & \hspace{32mm} \times \left (D_{i,r}\DexcludingViperturbdummy - D^{(-v_i,\idummy)}_{\idummy,r}\right) \Bigg |  \\
    & \le \Enb \mathfrak{D}_{r,1} + 2 \hat{\sigma}_{r,1} \Enb \sqrt{ \mathfrak{D}_{r,1}} \\
    &\le \Enb \mathfrak{D}_{r,1} + 2 \sqrt{C} \cdot \sqrt{\Enb \mathfrak{D}_{r,1}}, 
\end{aligned}
\end{equation}
where $\mathfrak{D}_{r,1} = \frac{1}{2n^2} \sum_{i,\idummy=1}^n \left( D_{i,r}\DexcludingVi - D_{\idummy,r}\DexcludingVdummy -  D_{i,r}\DexcludingViperturbdummy + D^{(-v_i,\idummy)}_{\idummy,r} \right)^2$. We can further bound $\Enb \mathfrak{D}_{r,1}$ by
\begin{equation} \label{eq: bound for Delta 1}
    \Enb \mathfrak{D}_{r,1} \le \frac{1}{n^2} \sum_{i,\idummy=1}^n \Enb \left( D_{i,r}\DexcludingVi - D_{i, r}\DexcludingViperturbdummy \right)^2 + \frac{1}{n^2} \sum_{i,\idummy=1}^n \Enb \left ( D^{(-v_i,\idummy)}_{\idummy,r}  - D_{\idummy,r}\DexcludingVdummy\right)^2. 
\end{equation}
The first summation in~\eqref{eq: bound for Delta 1} is 
$$
     \frac{1}{n^2} \sum_{i,\idummy=1}^n \Enb \left( D_{i,r}\DexcludingVi - D_{i, r}\DexcludingViperturbdummy \right)^2 = \frac{1}{n^2} \sum_{i = 1}^n \sum_{\idummy \notin I_{v_i}} \Enb \left (Q_{i,r}\QexcludingVi - Q_{i, r}\QexcludingViperturbdummy \right)^2 = o(n^{-1})
$$
by modifying the stability result~\eqref{eq: bound Q difference} (we only bound the difference $\abs{Q_{i,r}\QexcludingVi - Q_{i, r}\QexcludingViperturbdummy}$ this time). As for the second summation in~\eqref{eq: bound for Delta 1}, one can obtain 
\begin{equation}
\begin{aligned}
    &\frac{1}{n^2} \sum_{i,\idummy=1}^n \Enb \left ( D^{(-v_i,\idummy)}_{\idummy,r}  - D_{\idummy,r}\DexcludingVdummy\right)^2 \\
    &\stackrel{(I)}{=} \frac{1}{n^2} \sum_{i=1}^n \sum_{\ell \notin I_{v_i}}\E{ \E{\left ( D^{(-v_i,\idummy)}_{\idummy,r}  - D_{\idummy,r}\DexcludingVdummy\right)^2 \longconditional \bm{X}^{(-v_i), (-v_\idummy)}, X_{\idummy}}} \\
    & \stackrel{(II)}{=}\frac{2}{n^2} \sum_{i=1}^n \sum_{\ell \notin I_{v_i}} \mathbb{E}\left[\operatorname{Var}\left(D_{\idummy, r}\DexcludingVdummy \mid \bm{X}^{(-v_i), (-v_\idummy)}, X_{\ell}\right)\right]\\
    & \stackrel{(III)}{\le} \frac{1}{n^2} \sum_{i=1}^n \sum_{\idummy \notin I_{v_i}} \sum_{j \in I_{v_i}} \E{ \E{\left ( D_{\idummy,r}\DexcludingVdummy  - D_{\idummy,r}\DexcludingVdummyperturbj\right)^2 \longconditional \bm{X}^{(-v_i), (-v_\idummy)}, X_{\idummy}}} \\
    &= \frac{1}{n^2} \sum_{i=1}^n \sum_{\idummy \notin I_{v_i}} \sum_{j \in I_{v_i}} \E{ \E{\left ( Q_{\ell, r}\QexcludingVdummy  - Q_{\idummy,r}\QexcludingVdummyperturbj\right)^2 \longconditional \bm{X}^{(-v_i), (-v_\idummy)}, X_{\idummy}}} \\
    & = \frac{1}{n^2} \sum_{i=1}^n \sum_{\ell \notin I_{v_i}} \sum_{j \in I_{v_i}} \mathbb{E}\left[\left(Q_{\ell, r}\QexcludingVdummy-Q_{\idummy,r}\QexcludingVdummyperturbj\right)^2\right]\\
    &\stackrel{(\RNum{4})}{=} o(1).
\end{aligned}
\end{equation}

The step (I) follows from the definition of $D_{\idummy, r}\DexcludingVdummy$. When $\idummy\in I_{v_i}$, we know $v_i = v_\idummy$. Because $X_\idummy$ does not involve in the calculation of the exponential weightings in $D_{\idummy, r}\DexcludingVdummy$, replacing it by $X_\idummy^\prime$ does not change the value. Namely, we simply have $D_{\ell,r}^{(-v_i),\idummy} = D_{\idummy, r}\DexcludingVdummy$ in this case. 

For step (II), we used a simple identity $\mathbb{E}(X - X')^2 = 2\operatorname{Var}(X)$ with $X,X'$ being IID variables. In our case, conditioning on the presented variables, $D_{\idummy, r}^{\left(-v_i\right), \idummy}$ is a function of $\boldsymbol{X}^{\left(v_\idummy\right),\idummy}$ (samples belonging to fold $v_\idummy$ with $X_\idummy$ perturbed) and $D_{\idummy, r}\DexcludingVdummy$ is a function of $\boldsymbol{X}^{\left(v_{i}\right)}$. One can observe that they are conditionally independent and identically distributed. 

The step (III) employs the (conditional) Efron-Stein's inequality (see Lemma 1 in~\cite{bayle2020}), where the variability only takes place in $X^{(v_i)}$ since we have conditioned on the other variables. 

As for the last step (IV), it holds true again by modifying the stability result~\eqref{eq: bound Q difference}. 

Overall, we have $\Enb \mathfrak{D}_{r,1} = o(1)$ and therefore $\Enb\abs{\hat{\sigma}_{r}^2 - \intermediatevarianceone } = o(1)$.

\textbf{Part 3: Bound $\mathbb{E}\left|\intermediatevarianceone-\intermediatevariancetwo\right|.$}

To analyze the second expectation in \eqref{eq: out-sample variance bound}, we rewrite $\intermediatevarianceone$ as 
\begin{align*}
    \intermediatevarianceone &= \frac{1}{2n^2} \sum_{i,\idummy=1}^n \Enb \left (D_{i,r}\DexcludingViperturbdummy - D^{(-v_{i}),\idummy}_{\idummy,r} \right)^2 \\
    &= \frac{1}{2n^2} \sum_{i,\idummy=1}^n \Enb \Bigg (D_{i,r}\DexcludingViperturbdummy - \E{D_{i,r}\DexcludingViperturbdummy \longconditional \bm{X}^{(-v_i), \idummy}} \\
    &\hspace{20mm}- D^{(-v_i,\idummy)}_{\idummy,r} + \E{D^{(-v_i,\idummy)}_{\idummy,r} \longconditional \bm{X}^{(-v_i), \idummy}}\Bigg)^2\\
    &:= \frac{1}{2n^2} \sum_{i,\idummy=1}^n \Enb \left (K_{i,r}\KexcludingViperturbdummy - K^{(-v_i,\idummy)}_{\idummy,r}\right)^2, 
\end{align*}
where the second equality holds true as discussed in~\eqref{eq: same conditional expectation}, and we define 
$$
K_{i,r}\KexcludingViperturbdummy := D_{i,r}\DexcludingViperturbdummy - \E{D_{i,r}\DexcludingViperturbdummy \mid \bm{X}^{(-v_i), \idummy}}, \hspace{3mm} K^{(-v_i,\idummy)}_{\idummy,r} := D^{(-v_i,\idummy)}_{\idummy,r} - \E{D^{(-v_i,\idummy)}_{\idummy,r} \mid \bm{X}^{(-v_i), \idummy}}.
$$
Based on the definitions of $D_{i,r}^\idummy$ and $D_{\idummy,r}^{(-v_i,\idummy)}$, we know that $K_{i,r}^\idummy = K_{i,r}$ and $K_{\idummy,r}^{(-v_i,\idummy)} = K_{\idummy,r}$ if $i,\idummy$ are within the same fold. 

Applying a similar argument as in \textbf{Part 2}, we have $\Enb \abs{\intermediatevarianceone - \intermediatevariancetwo} \le \Enb \mathfrak{D}_{r,2} + 2 \sqrt{C \cdot \Enb \mathfrak{D}_{r,2}}$ with
\begin{equation} \label{eq: bound for Delta 2}
\begin{aligned}
    \Enb \mathfrak{D}_{r,2} &= \frac{1}{2n^2} \sum_{i,\idummy=1}^n \Enb \left (K_{i,r}\KexcludingViperturbdummy - K^{(-v_i,\idummy)}_{\idummy,r} -  \E{K_{i,r} \longconditional X_i} + \E{K_{\idummy, r}\longconditional X_{\idummy}}\right)^2 \\
    & \le \frac{1}{n^2} \sum_{i,\idummy=1}^n \Enb \left (K_{i,r}\KexcludingViperturbdummy -  \E{K_{i,r}\longconditional X_i}\right)^2  + \frac{1}{n^2} \sum_{i,\idummy=1}^n \Enb \left (K^{(-v_i,\idummy)}_{\idummy,r} -  \E{K_{\idummy,r}\longconditional X_{\idummy}}\right)^2. 
\end{aligned}
\end{equation}
By the conditional Efron-Stein inequality, the first summation in~\eqref{eq: bound for Delta 2} can be bounded by
\begin{equation} \label{eq: ES for expectation of variance}
\begin{aligned}
    &\frac{1}{n^2} \sum_{i,\idummy=1}^n \Enb \left (K_{i,r}\KexcludingViperturbdummy -  \E{K_{i,r}\longconditional X_i}\right)^2 \\ 
    &=  \frac{1}{n^2} \sum_{i, \idummy = 1}^n\E{\E{ \left (K_{i,r}\KexcludingViperturbdummy -  \E{K_{i,r}\longconditional X_i}\right)^2 \longconditional X_i}} \\
    &=  \frac{1}{n^2} \sum_{i, \idummy = 1}^n\E{\E{ \left (K_{i,r}\KexcludingViperturbdummy -  \E{K_{i,r}\KexcludingViperturbdummy \longconditional X_i}\right)^2 \longconditional X_i}} \\
    &=  \frac{1}{n^2} \sum_{i,\idummy = 1}^n \E{\mathrm{Var} \left (K_{i,r}\KexcludingViperturbdummy \longconditional X_i \right ) } \\
    &\le \frac{1}{2n^2} \sum_{i, \idummy=1}^n \sum_{j \notin I_{v_i}} \E{ \E{ \left (\nabla_j K_{i,r}\KexcludingViperturbdummy \right)^2 \longconditional X_i }} \\
    &= \frac{1}{2n^2} \sum_{i, \idummy=1}^n \sum_{j \notin I_{v_i}} \E{ \left (\nabla_j K_{i,r}\KexcludingViperturbdummy \right)^2} \\
    &= o(1),
\end{aligned}
\end{equation}
where the last equality follows from Lemma~\ref{lemma: first order stability}. Also, the second equality holds true because (1) when $i, \idummy$ are within the same fold, we directly have $K_{i,r} = K_{i,r}^\idummy$ and (2) when $\ell \notin I_{v_i}$, the variables $K_{i,r}$ and $K_{i,r}^\idummy$ are identically distributed, conditioning on $X_i$. 

Similarly, one can show that the second summation in~\eqref{eq: bound for Delta 2} is $o(1)$, using that $\E{K_{\idummy,r}\mid X_\idummy} = \E{K^{(-v_i,\idummy)}_{\idummy,r}\mid X_\idummy}$. Indeed, note that (1) when $i,\idummy$ are within the same fold, we have $K_{\idummy,r} = K^{(-v_i, \idummy)}_{\idummy,r}$; (2) when $\idummy \notin I_{v_i}$, the variables $K_{\idummy,r}$ and $K_{\idummy,r}^{(-v_i,\idummy)}$ are identically distributed, conditioning on $X_\idummy$. Overall, $\Enb \mathfrak{D}_{r,2} = o(1)$ and therefore $\Enb \abs{\intermediatevarianceone - \intermediatevariancetwo} = o(1)$. 

\textbf{Part 4: bound $\Enb \left|\intermediatevariancetwo-\undertilde{\sigma}_r^2\right|.$}

To prove the convergence of the third expectation~\eqref{eq: out-sample variance bound}, it suffices to show that $\abs{\intermediatevariancetwo - \undertilde{\sigma}_r^2} \cinp 0$, given the boundedness of $\intermediatevariancetwo$ and $\undertilde{\sigma}_r^2$. Observe that 
\begin{align*}
    \intermediatevariancetwo &= \frac{1}{2n^2} \sum_{i,\idummy=1}^n \left ( \E{K_{i,r} \longconditional  X_i}  - \E{K_{\idummy,r} \longconditional X_{\idummy}} \right)^2 \\ 
    &= \frac{1}{2n^2} \sum_{i,\idummy=1}^n \Bigg( \Enb^2 \left [ K_{i,r} \longconditional X_i \right] + \Enb^2 \left [ K_{\idummy,r} \longconditional X_{\idummy} \right] \\
    &\hspace{24mm}- 2 \Enb \left [ K_{i,r} \longconditional X_i \right] \Enb \left [ K_{\idummy,r} \longconditional X_{\idummy} \right] \Bigg) \\
    &= \frac{1}{n} \sum_{i=1}^n  \Enb^2 \left [ K_{i,r} \longconditional X_i \right ] - \left (\frac{1}{n} \sum_{i=1}^n \E{K_{i,r} \longconditional X_i} \right)^2. 
\end{align*}
and that by the independence between $X_1$ and $\bm{X}^{(-v_1)}$, 
\begin{align}
\label{eq: variance of expectation equivalency}
\mathrm{Var} \left ( \E{K_{1,r} \longconditional X_1}\right) &= \mathrm{Var}(\E{D_{1,r} - \E{D_{1,r} \mid \bm{X}^{(-v_1)}} \mid X_1})\\
&= \mathrm{Var}(\E{D_{1,r}\mid X_1} - \E{D_{1,r}})\\
&= \mathrm{Var} \left( \E{D_{1,r}\DexcludingVone \longconditional X_1} \right) = \undertilde{\sigma}_r^2. 
\end{align}
By the uniform boundedness of $X_1$, we know that $\E{K_{i,r} \longconditional X_i}$ is bounded for any $i \in [n]$, and therefore the sufficient condition for the weak law concerning its triangular array must be satisfied. 

We have thus established $\abs{\hat{\sigma}_{r}^2 - \undertilde{\sigma}_r^2} \cinp 0$. Together with Proposition~\ref{prop: equivalency among variances}, we know $\abs{\hat{\sigma}_{r}^2 - \sigma_r^2} \cinp 0$. Because the entries of $X_1$ are uniformly bounded and assuming $\mathrm{cov}(X_1)$ has strictly positive eigenvalues assures that $\sigma^2_r$ is bounded away from $0$ for all $n \in \N$, it can be concluded that $\hat{\sigma}_{r}^2/\sigma_r^2 \cinp 1$.

\end{proof}

\subsection{Proof of Theorem~\ref{prop: consistency of within-fold sample variance}}
\begin{proof} 
    To prove the desired result, it suffices to show $\Enb \abs{\hat{\varphi}_{r}^{2}- \undertilde{\sigma}_r^2} \rightarrow 0$ with $\undertilde{\sigma}_r^2 = \mathrm{Var}( \E{X_{1,r} - Q_{1,r} | X_1})$ by Proposition~\ref{prop: equivalency among variances}. 
    
\textbf{Part 1: split the difference  $\left|\hat{\varphi}_{r}^2-{\undertilde{\sigma}}_{r}^2\right| $ into two parts, and bound each separately.}

To prove the convergence, we consider the bound 
\begin{equation} \label{eq: in-sample variance bound}
    \Enb \abs{\hat{\varphi}_{r}^{2} - \undertilde{\sigma}_r^2} \le \Enb \abs{\hat{\varphi}_{r}^{2} - \hat{\varphi}^2_{r,1}} + \Enb \abs{\hat{\varphi}^2_{r,1} - \undertilde{\sigma}_r^2}, 
\end{equation}
where the estimator $\hat{\varphi}^2_{r,1}$ is defined by 
\begin{align*}
\hat{\varphi}^2_{r,1} &= \frac{1}{V} \sum_{v = 1}^V \frac{1}{(V/n) - 1} \sum_{i \in I_v} \left (\E{K_{i,r}\KexcludingV \longconditional X_i} - \frac{V}{n} \sum_{\idummy \in I_{v}} \E{K_{\idummy,r}\KexcludingV \longconditional X_\idummy} \right)^2 \\
&=: \frac{1}{V} \sum_{v=1}^V \frac{1}{(V/n) - 1} \sum_{i \in I_v} \left (\E{K_{i,r}\KexcludingV \longconditional X_i} - \tilde{K}_r\KincludingV \right)^2,
\end{align*}
where the variable $K_{i,r}$ is defined as in Lemma~\ref{lemma: first order stability}, and for any $v \in [V]$, the variable $\tilde{K}_r\KincludingV$ is defined by $\tilde{K}_r\KincludingV = \frac{V}{n} \sum_{\idummy \in I_{v}} \E{K_{\idummy,r}\KexcludingV \longconditional X_\idummy}$. The uniform boundedness of $X_{1}$ ensures that there exists $C > 0$ such that $\hat{\varphi}^2_{r,1} < C$ and $\undertilde{\sigma}_r^2 < C$.

\textbf{Part 2: bound $\mathbb{E}\left|\hat{\varphi}_{r}^2-\hat{\varphi}^2_{r,1}\right|.$}

For any $v \in [V]$ and $i \in I_v$, define $D_{i,r}\DexcludingV = X_{i,r} - \sum_{s \ne r} \hat{w}^{(-v)}\wscenterr  X_{i, s}$ and $\overline{D}_{r}\DincludingV = \frac{V}{n} \sum_{\idummy \in I_
    {v}}^n D_{\idummy,r}\DexcludingV$. To analyze the first expectation in~\eqref{eq: in-sample variance bound}, we first rewrite the sample variance $\hat{\varphi}_{r}^{2}$ by 
\begin{align*}
\hat{\varphi}_{r}^{2} 
    &= \frac{1}{V} \sum_{v = 1}^V \frac{1}{(n/V) - 1} \sum_{i \in I_v} \left (D_{i,r}\DexcludingV - \overline{D}_{i,r}\DincludingV\right)^2 \\
    & \stackrel{(\RNum{1})}{=} \frac{1}{V} \sum_{v=1}^V \frac{1}{(n/V) - 1} \sum_{i \in I_v} \Bigg (D_{i,r}\DexcludingV-\E{D_{i,r}\DexcludingV\longconditional \bm{X}^{(-v)}} \\
    & \hspace{42mm}- \overline{D}_{i,r}\DincludingV + \frac{V}{n} \sum_{\idummy \in I_{v}} \E{D_{\idummy,r}\DexcludingV \longconditional \bm{X}^{(-v)}}\Bigg )^2 \\
    & := \frac{1}{V} \sum_{v=1}^V \frac{1}{(n/V) - 1} \sum_{i \in I_v} \left (K_{i,r}\KexcludingV - \overline{K}_r\KincludingV\right)^2, 
\end{align*}
with the variable $\overline{K}_r\KincludingV$ defined by
\begin{align*}
\overline{K}_{r}\KincludingV &= \overline{D}_{i,r}\DincludingV - \frac{V}{n}\sum_{\idummy \in I_{v}} \E{D_{\idummy,r}\DexcludingV\longconditional\bm{X}^{(-v)}}\\
& = \frac{V}{n} \sum_{\idummy \in I_{v}} \left(D_{\idummy,r}\DexcludingV - \E{D_{\idummy,r}\DexcludingV \longconditional \bm{X}^{(-v)}} \right) \\
&= \frac{V}{n} \sum_{\idummy \in I_v} K_{\idummy,r}\KexcludingV
\end{align*}
for any $v \in  [V]$. The step (I) holds true because for any $\idummy \in I_{v}$, we know that conditioning on $\bm{X}^{(-v)}$, $D_{\idummy,r}\DexcludingV$ and $D_{i,r}\DexcludingV$ are IID random variables for all $i, \idummy \in I_v$. This particularly implies that $\E{D_{\idummy,r}\DexcludingV\longconditional \bm{X}^{(-v)}} = \E{D_{i,r}\DexcludingV\longconditional \bm{X}^{(-v)}}$ for all $i, \idummy \in I_v$.

Following a similar argument as \textbf{Part 2} in the proof of Proposition~\ref{cor: consistency of all-pair estimator}, one can obtain that $\Enb \abs{\hat{\varphi}_{r}^{2} - \hat{\varphi}^2_{r,1}} \le \Enb \mathcal{D}_{r,1} + 2 \sqrt{ C\cdot \Enb \mathcal{D}_{r,1}}$ with 
\begin{align*}
    \Enb \mathcal{D}_{r,1} &= \frac{1}{V} \sum_{v=1}^V \frac{1}{(n/V) - 1} \sum_{i \in I_v} \Enb \left( K_{i,r}\KexcludingV - \overline{K}_r\KincludingV - \E{K_{i,r}\KexcludingV \longconditional X_i} + \tilde{K}_r\KincludingV\right)^2 \\
    &\lesssim \frac{1}{V} \sum_{v=1}^V \frac{1}{(n/V) -1}\sum_{i \in I_v} \Enb \left (K_{i,r}\KexcludingV- \E{K_{i,r}\KexcludingV\longconditional X_i} \right)^2 \\
    &\hspace{6mm}+ \frac{1}{V} \sum_{v=1}^V \frac{1}{(n/V) -1} \sum_{i \in I_v} \Enb \left (\overline{K}_r\KincludingV - \tilde{K}_r\KincludingV \right)^2 \\
    &=\frac{1}{V} \sum_{v=1}^V \frac{1}{(n/V) -1}\sum_{i \in I_v} \Enb \left (K_{i,r}\KexcludingV-\E{K_{i,r}\KexcludingV\longconditional X_i} \right)^2 \\
    &\hspace{6mm}+ \frac{1}{V} \sum_{v=1}^V \frac{(n/V)}{(n/V) -1}  \Enb \left (\overline{K}_r\KincludingV - \tilde{K}_r\KincludingV \right)^2 \\
    &\le \frac{1}{V} \sum_{v=1}^V \frac{1}{(n/V) - 1} \sum_{i \in I_v} \Enb \left (K_{i,r}\KexcludingV- \E{K_{i,r}\KexcludingV\longconditional X_i} \right)^2 \\
    &\hspace{6mm}+ \frac{1}{V} \sum_{v=1}^V \frac{(n/V)}{(n/V) -1}  \left [ \frac{V}{n} \sum_{\idummy \in I_v}  \Enb \left (K_{\idummy,r}\KexcludingV- \E{K_{\idummy,r}\KexcludingV\longconditional X_\idummy} \right)^2 \right ], 
\end{align*}
where we employ the Jensen's inequality for the last step. From now on, we denote the fraction $(n/V)/[(n/V) - 1]$ by $C_n$ to ease our notation. Remark that $C_n = O(1)$. 

Using the conditional Efron-Stein inequality and Lemma~\ref{lemma: first order stability}, we have 
\begin{align*}
    &\Enb \left (K_{\idummy,r}\KexcludingV- \E{K_{\idummy,r}\KexcludingV\longconditional X_\idummy} \right)^2 \\
    &= \E{\E{ \left (K_{\idummy,r}\KexcludingV- \E{K_{\idummy,r}\KexcludingV\longconditional X_\idummy} \right)^2\longconditional X_\idummy}} \\
    &= \E{\mathrm{Var} \left (K_{\idummy,r}\KexcludingV \longconditional X_\idummy \right)} \\
    & \le \frac{1}{2}\sum_{j \notin I_{v_\idummy}} \E{\E{\left( \nabla_j K_{\idummy,r}\KexcludingV \right)^2\longconditional X_\idummy}} \\
    &= \frac{1}{2}\sum_{j \notin I_{v_\idummy}} \Enb \left (\nabla_j K_{\idummy,r}\KexcludingV\right )^2 \\
    &= o(1),
\end{align*}
which in turn gives $\Enb \mathcal{D}_{r,1} = o(1)$ and therefore $\Enb \abs{\hat{\varphi}_{r}^{2} - \hat{\varphi}^2_{r,1}} = o(1)$. 

\textbf{Part 3: show $\mathbb{E}\left|\hat{\varphi}^2_{r,1}-\undertilde{\sigma}_r^2\right| \le \sqrt{\mathrm{Var} ( \hat{\varphi}^2_{r,1}) }.$}

To prove the desired inequality, it reduces to justify $\Enb \hat{\varphi}^2_{r,1} = \undertilde{\sigma}_r^2$ according to the Jensen's inequality. We first rewrite the estimator $\hat{\varphi}^2_{r,1}$ by 
\begin{align*}
    \hat{\varphi}^2_{r,1} 
    &= \frac{1}{V} \sum_{v=1}^V \frac{1}{(n/V) - 1} \sum_{i \in I_v} \left (\E{K_{i,r}\KexcludingV \longconditional X_i} - \tilde{K}_r\KincludingV \right)^2 \\
    &= \frac{1}{V} \sum_{v=1}^V \frac{1}{(n/V) - 1} \sum_{i \in I_{v}} \Bigg( \Enb^2 \left [K_{i,r}\KexcludingV \longconditional X_i \right] + \left\{ \tilde{K}_r\KincludingV \right\}^2 \\
    & \hspace{40mm}- 2 \tilde{K}_r\KincludingV \cdot \E{K_{i,r}\KexcludingV \longconditional X_i} \Bigg) \\
    &= \frac{1}{V} \sum_{v=1}^V \frac{1}{(n/V) - 1}\sum_{i \in I_v} \Enb^2 \left [K_{i,r}\KexcludingV \longconditional X_i \right] -\frac{C_n}{V} \sum_{v=1}^V \left \{\tilde{K}_r\KincludingV \right\}^2. 
\end{align*}
Expanding the second summation as
\begin{align*} 
    &\frac{C_n V}{n^2} \sum_{v=1}^V \sum_{i \in I_v} \Enb^2 \left[K_{i,r}\KexcludingV \longconditional X_i \right] + \frac{C_n V}{n^2} \sum_{v=1}^V \sum_{i,\idummy \in I_v, i \ne \idummy} \E{K_{i,r}\KexcludingV \longconditional X_i} \E{K_{\idummy,r}\KexcludingV \longconditional X_\idummy} \\
    &= \frac{C_n V}{n^2} \sum_{v=1}^V \sum_{i \in I_v} \Enb^2 \left[K_{i,r}\KexcludingV \longconditional X_i \right] + \frac{2C_n V}{n^2} \sum_{v=1}^V \sum_{i,\idummy \in I_v, i  < \idummy} \E{K_{i,r}\KexcludingV \longconditional X_i} \E{K_{\idummy,r}\KexcludingV \longconditional X_\idummy},
\end{align*}
we can achieve the expression
\begin{equation} \label{eq: dummay variance expression}
\begin{aligned}
     \hat{\varphi}^2_{r,1} &= \frac{C_n}{n} \left(1 - \frac{V}{n} \right) \sum_{v=1}^V \sum_{i \in I_v} \Enb^2 \left[K_{i,r}\KexcludingV \longconditional X_i \right] \\
     &\hspace{3mm}-\frac{2 C_n V}{n^2} \sum_{v=1}^V \sum_{i,\idummy \in I_v, i < \idummy} \E{K_{i,r}\KexcludingV \longconditional X_i} \E{K_{\idummy,r}\KexcludingV \longconditional X_\idummy}.
\end{aligned}
\end{equation}
The estimator $\hat{\varphi}^2_{r,1}$ therefore has the expectation 
\begin{equation} \label{eq: expectation of sigma in 1}
\begin{aligned}
    \Enb\hat{\varphi}^2_{r,1} &= \frac{C_n}{n} \left(1 - \frac{V}{n} \right) \sum_{v=1}^V \sum_{i \in I_v} \E{\Enb^2 \left[K_{i,r}\KexcludingV \longconditional X_i \right]} \\
     &\hspace{3mm}-\frac{2 C_n V}{n^2} \sum_{v=1}^V \sum_{i,\idummy \in I_v, i < \idummy} \E{\E{K_{i,r}\KexcludingV \longconditional X_i} \E{K_{\idummy,r}\KexcludingV \longconditional X_\idummy}} \\
     &\stackrel{(\RNum{1})}{=} \frac{C_n}{n} \left ( 1 - \frac{V}{n}\right) \sum_{v=1}^V \sum_{i \in I_v} \undertilde{\sigma}_r^2\\
     &\hspace{3mm}-\frac{2 C_n V}{n^2} \sum_{v=1}^V \sum_{i,\idummy \in I_v, i < \idummy} \E{\E{K_{i,r}\KexcludingV \longconditional X_i} \E{K_{\idummy,r}\KexcludingV \longconditional X_\idummy}} \\
     &\stackrel{(\RNum{2})}{=} \frac{C_n}{n} \left ( 1 - \frac{V}{n}\right) \sum_{v=1}^V \sum_{i \in I_v} \undertilde{\sigma}_r^2\\
     &\hspace{1mm}=\undertilde{\sigma}_r^2.
\end{aligned}
\end{equation}
For the step (I), we used the fact that 
\begin{equation} \label{eq: undertilde sigma}
\begin{aligned}
    &\E{\Enb^2 \left[K_{i,r}\KexcludingV \longconditional X_i \right]} \\
    &=\E{\Enb^2 \left[K_{1,r}\KexcludingV \longconditional X_1 \right]} \\
    &= \mathrm{Var}\left (\E{K_{1,r}\KexcludingVone \longconditional X_1}\right) + \left (\E{\E{K_{1,r}\KexcludingVone \longconditional X_1}} \right)^2\\
    &= \mathrm{Var}\left (\E{D_{1,r}\DexcludingVone \longconditional X_1}\right) + \left (\E{\E{K_{1,r}\KexcludingVone \longconditional X_1}} \right)^2\\
    &= \undertilde{\sigma}_r^2 + \left ( \E{K_{1,r}\KexcludingVone} \right )^2\\
    &= \undertilde{\sigma}_r^2,
\end{aligned}
\end{equation}
where the third equality holds true by~\eqref{eq: variance of expectation equivalency}. 

The second step (II) holds true because for any $i, \idummy \in [n]$ such that $i \ne \idummy$, the variable $\E{K_{i,r}\KexcludingV \longconditional X_i}$ is independent of the counterpart $\E{K_{\idummy,r}\KexcludingV \longconditional X_\idummy}$ and we have 
\begin{equation} \label{eq: expection of product is zero}
\E{\E{K_{i,r}\KexcludingV \longconditional X_i}} = \E{K_{i,r}\KexcludingV} = 0.
\end{equation}

\textbf{Part 4: bound $\sqrt{\mathrm{Var} ( \hat{\varphi}^2_{r,1})}.$}

Based on the expression~\eqref{eq: dummay variance expression}, we have $\mathrm{Var} \left( \hat{\varphi}^2_{r,1} \right) = \mathrm{Var} \left ( S_1 \right ) + \mathrm{Var} \left (S_2\right) - 2\mathrm{cov} (S_1, S_2)$ with 
\begin{align*}
    &S_1 = \frac{C_n}{n} \left(1 - \frac{V}{n} \right) \sum_{v=1}^V \sum_{i \in I_v} \Enb^2 \left[K_{i,r}\KexcludingV \longconditional X_i \right]; \\
    &S_2 = \frac{2C_n V}{n^2} \sum_{v=1}^V \sum_{i,\idummy \in I_v, i < \idummy} \E{K_{i,r}\KexcludingV \longconditional X_i} \E{K_{\idummy,r}\KexcludingV \longconditional X_\idummy}.
\end{align*}
Because $\{X_i : i \in [n]\}$ is a set of independent random variables and $K_{i,r}\KexcludingV$ is uniformly bounded by the uniform boundness of $X_1$, the variance $\mathrm{Var}(S_1)$ can be bounded by
\begin{align*}
    \mathrm{Var}(S_1) 
    &= \frac{C_n^2}{n^2} \left ( 1 - \frac{V}{n} \right)^2 \sum_{v=1}^V \sum_{i \in I_v}\mathrm{Var} \left ( \Enb^2 \left [K_{i,r} \longconditional X_i \right ]\right ) \\
    &= \frac{C_n^2}{n^2} \left ( 1 - \frac{V}{n} \right)^2 \sum_{v=1}^V \sum_{i \in I_v} \left \{ \E{\Enb^4 \left [K_{i,r} \longconditional X_i \right ]} - \undertilde{\sigma}^4_r \right \} \\
    &\le \frac{C^2_n}{n} \left (1 - \frac{V}{n}\right)^2 \tilde{M} +  \frac{C^2_n}{n} \left (1 - \frac{V}{n}\right)^2 \undertilde{\sigma}^4_r
\end{align*}
for some $\tilde{M} > 0$, where the second equality holds true by the fact~\eqref{eq: undertilde sigma}. Hence, $\mathrm{Var}(S_1) = O(n^{-1})$.

We also have 
\begin{align*}
    \mathrm{Var}(S_2) 
    &\stackrel{(\RNum{1})}{=}\Enb \left ( \frac{2C_n V}{n^2} \sum_{v=1}^V \sum_{i,\idummy \in I_v, i  < \idummy} \E{K_{i,r}\KexcludingV \longconditional X_i} \E{K_{\idummy,r}\KexcludingV \longconditional X_\idummy} \right)^2 \\
    &\hspace{-1mm}\stackrel{(\RNum{2})}{=} \frac{4C_n^2 V^2}{n^4} \sum_{v = 1}^V \sum_{i,\idummy \in I_v, i < \idummy} \E{\Enb^2 \left [ K_{i,r}\KexcludingV\longconditional X_i \right]} \E{\Enb^2 \left [ K_{\idummy,r}\KexcludingV\longconditional X_{\idummy} \right]} \\
    &\hspace{-2mm}\stackrel{(\RNum{3})}{=} \frac{2C_n^2 V^2}{n^4} \cdot V \cdot \frac{n}{V} \cdot \left(\frac{n}{V} -1\right) \cdot \undertilde{\sigma}^4_r \\
    &\hspace{-1mm}\stackrel{(\RNum{4})}{\le} \frac{C_n^2 \undertilde{\sigma}^4_r}{2n} =O(n^{-1}).
\end{align*}
In the step (I), we used the argument for the step (II) in~\eqref{eq: expectation of sigma in 1} that essentially gives $\Enb S_2 = 0$, and the step (III) is ensured by~\eqref{eq: undertilde sigma}. As for the above step (II), it holds because for any $v, w \in [V]$ and any quadruplet $(i, \idummy, j, \iota)$ such that $i, \idummy \in I_v$ with $i < \idummy$ and $j, \iota \in I_w$ with $j < \iota$, the expectation
$$
    \E{\Enb \left [ K_{i,r}\KexcludingV\longconditional X_i \right] \cdot \Enb \left [ K_{\idummy,r}\KexcludingV\longconditional X_{\idummy} \right] \cdot \Enb \left [ K_{j,r}\KexcludingW\longconditional X_j \right] \cdot \Enb \left [ K_{\iota,r}\KexcludingW\longconditional X_{\iota} \right]}
$$
would be exactly zero, by employing~\eqref{eq: expection of product is zero}, whenever $v \ne w$ (in the case, the four variables $X_i, X_\idummy, X_j$ and $X_\iota$ are independent), or $v = w$ yet $(i, \idummy) \ne (j, \iota)$ (in the case, at least one of $X_i, X_\idummy, X_j$ and $X_\iota$ is independent of the others). 
The step (IV) takes advantage of the inequality $V (n - V) \le n^2/4$.

Finally, let us examine $\mathrm{cov}(S_1, S_2)$. Observe that for any $v, w\in [V]$ and any triplet $(j, i, \idummy)$ such that $j \in I_v$ and $i, \idummy \in I_w$ with $i < \idummy$, we must have 
$$
    \E{\Enb^2 \left [K_{j,r}\KexcludingV \longconditional X_j \right]} \cdot \E{ \E{K_{i,r}\KexcludingW \longconditional X_i } \cdot \E{K_{\idummy,r}\KexcludingW \longconditional X_{\idummy}}} = 0.
$$
Also, we have
\begin{align*}
    &\E{\Enb^2 \left [K_{j,r}\KexcludingV \longconditional X_j \right] \cdot \E{K_{i,r}\KexcludingW \longconditional X_i } \cdot \E{K_{\idummy,r}\KexcludingW \longconditional X_{\idummy}}} = 0. 
\end{align*}
Indeed, the expectation is, by independence, equal to either 
$$
\E{\Enb^2 \left [K_{j,r}\KexcludingV \longconditional X_j \right] \cdot \E{K_{i,r}\KexcludingW \longconditional X_i }} \cdot \E{\E{K_{\idummy,r}\KexcludingW \longconditional X_{\idummy}}}
$$
or 
$$
\E{\Enb^2 \left [K_{j,r}\KexcludingV \longconditional X_j \right] \cdot \E{K_{\idummy,r}\KexcludingW \longconditional X_\idummy}} \cdot \E{\E{K_{i,r}\KexcludingW \longconditional X_{i}}},
$$
(at most one of $i$ or $\idummy$ can be $j$ when $v = w$) and each of them is zero by~\eqref{eq: expection of product is zero}. Together with this observation, the argument for the step (II) in~\eqref{eq: expectation of sigma in 1} ensures $\mathrm{cov}(S_1, S_2) = 0$.

We thus have shown $\sqrt{\mathrm{Var}(\hat{\varphi}^2_{r,1})} = O(n^{-1/2}) = o(1)$ and therefore $\Enb \abs{\hat{\varphi}^2_{r,1} - \undertilde{\sigma}_r^2} = o(1)$. Together with Proposition~\ref{prop: equivalency among variances}, we have $\abs{\hat{\varphi}_{r}^{2} - \sigma_r^2} \cinp 0$. Because the entries of $X_1$ are uniformly bounded and assuming $\mathrm{cov}(X_1)$ has strictly positive eigenvalues assures that $\sigma^2_r$ is bounded away from $0$ for all $n \in \N$, it can be concluded that $\hat{\varphi}_{r}^2/\sigma_r^2 \cinp 1$.  
 \end{proof}

\section{Proof of Theorem~\ref{th: weakly dependent CLT}}\label{app: weakly dependent CLT}


\begin{proof}
We will first show the statistic of interest converges in distribution to some random variable $Y$. Afterwards we will establish $Y\rightarrow\mathcal{N}(0,1)$. 

\textbf{Part 1} 
Define 
\begin{equation*}
\begin{aligned}
         K & = \omega^{-1}\sum_{j=1}^n n^{-1/2} K_j = \omega^{-1}\sum_{j = 1}^n n^{-1/2}\mathcal{K}(j, \bm{X}) \\
        Y & = \omega^{-1}\sum_{j=1}^n n^{-1/2} Y_j
    \end{aligned}
\end{equation*}
where $Y_j  = \tilde\sigma_j\xi_j$, $\ \tilde\sigma_j^2 = \mathbb{E}[K_j^2 \mid \bm{X}^{(-j)}]$, and $\xi_j\stackrel{IID}{\sim}\mathcal{N}(0,1)$ is independent of everything else. Recall the variance $\omega^2 = \operatorname{Var}(K_j)$ may implicitly depend on $n$ but assumed to be greater than a constant for large enough $n$.

We will apply a version of the Portmanteau theorem to bound
\begin{equation*}
    |\mathbb{P}(K \leq x) - \mathbb{P}(Y \leq x)|.
\end{equation*}
Examining the proof of Lemma~2 in \cite{chin2022short} (or Theorem~12 in \cite{wasserman2014stein}), we know for any $\epsilon >0$, there exists a smooth indicator function $g = g_{\epsilon, x}$ that 1) is three times differentiable and 2) is bounded themselves and has bounded derivatives (note this bound is also dependent on $\epsilon$) that satisfies
\begin{equation*}
\begin{aligned}
   & |\mathbb{P}(K \leq x) - \mathbb{P}(Y \leq x)|\\
   & \leq |\mathbb{E}g(K) - \mathbb{E}g(Y)| + \epsilon
\end{aligned}
\end{equation*}


We are going to apply the Slepian's interpolation to bound $|\mathbb{E} g(K)-\mathbb{E} g(Y)|$. We will consider the following random variables as in the standard treatment:
\begin{equation}
    \begin{aligned}
        Z(t) & = \sqrt{t}K + \sqrt{1-t} Y \text{ for }t\in [0,1],\\
        Z_j(t) & = \omega^{-1}n^{-1/2}(\sqrt{t}K_j + \sqrt{1-t}Y_j) \Rightarrow Z(t) = \sum_{j=1}^n Z_j(t).
    \end{aligned}
\end{equation}

We also consider the following LOO version random variables:
\begin{equation}
    \begin{aligned}
        K^i & = \omega^{-1}\sum_{j\neq i} n^{-1/2} \mathcal{K}(j, \bm{X}^{i}) \text{ where }\bm{X}^{i}:= (\bm{X}\symbol{92}\{X_i\})\cup \{X_i^\prime\}, \\
        K_j^i & = \omega^{-1}n^{-1 / 2} \mathcal{K}\left(j, \bm{X}^{i}\right) \Rightarrow K^i = \sum_{j \neq i} K_j^i, \\
        Y^i & = \omega^{-1}\sum_{j\neq i}n^{-1/2} Y_j, \\
        Z^i(t) & = \sqrt{t}K^i + \sqrt{1-t} Y^i \text{ for }t\in[0,1].
    \end{aligned}
\end{equation}
As a remark, the quantity $K^i$ is in fact constructed by replacing $X_i$ with an IID copy $X_i'$. \emph{If we were to define} $K^i$ as $\omega^{-1}\sum_{j\neq i}n^{-1/2}K_j$, $K^i$ would have been dependent on $X_i$ through the second argument of the $\mathcal{K}$ mapping. We need to consider the data set $\bm{X}^{i}$ to completely eliminate $X_i$'s impact on the LOO version of $K$.

Now we proceed with our proof. Define $\Psi(t)=\Enb g(Z(t))$. Then, bounding $|\mathbb{E} g(K)-\mathbb{E} g(Y)|$ reduces to controlling $\int_{0}^{1} \Psi'(t) d t$. Its integrand yields the decomposition
\begin{equation} \label{eq: slepian integrand decomposition}
\begin{aligned}
\Psi'(t)
& =  \Enb g'(Z(t)) \sum_{i=1}^{n} Z_{i}'(t) = \sum_{i=1}^{n} \Enb Z_{i}'(t)\left[g'\left(Z^{i}(t)\right)+g'(Z(t))-g'\left(Z^{i}(t)\right)\right] \\
& \stackrel{(I)}{=} \sum_{i = 1}^n \E{ Z_i'(t) g'\left( Z^i(t)\right)} + \sum_{i=1}^{n} \E{Z_{i}'(t) \left(Z(t) - Z^i(t)\right) g^{\prime \prime}\left(Z^{i}(t)\right)} \\ 
& \hspace{4mm}+ \sum_{i=1}^n \E{Z'_i(t) \left(Z(t) - Z^i(t)\right)^{2} \int_{0}^{1}(1-u) g^{\prime \prime \prime}\left(Z^{i}(t)+u\{Z(t) - Z^i(t)\}\right) d u}\\
& =: \mathfrak{A} + \mathfrak{B} + \mathfrak{C}.
\end{aligned}
\end{equation}
In step $(I)$ we used the following Taylor expansion
\begin{equation}
    \begin{aligned}
f(x) & =f(a)+(x-a) f^{\prime}(a)+\int_a^x(x-t) f^{\prime \prime}(t) d t \\
& =f(a)+(x-a) f^{\prime}(a)+(x-a)^2 \int_0^1(1-u) f^{\prime \prime}(a+u(x-a)) d u.
\end{aligned}
\end{equation}
In our case, $f = g'$, $x = Z(t)$ and $a = Z^i(t)$. We are going to bound the three terms $\mathfrak{A}$, $\mathfrak{B}$ and $\mathfrak{C}$ separately. We will use the following explicit form of $Z'_i(t)$ to verify several properties later:
\begin{equation}
Z_i^{\prime}(t)=\frac{d Z_i(t)}{d t}=\frac{1}{2\omega \sqrt{n}}\left(\frac{K_i}{\sqrt{t}}-\frac{Y_i}{\sqrt{1-t}}\right).
\end{equation}
By construction, $Z^i(t)$ is independent of $X_i$ and $\xi_i$. We also have $\mathbb{E}\left[Z_i^{\prime}(t) \mid \boldsymbol{X}^{(-i)}\right]=0$. So
\begin{equation}
    \begin{aligned}
         \mathfrak{A} & = \mathbb{E}\left[Z_i^{\prime}(t) g^{\prime}\left(Z^i(t)\right)\right]\\
        & = \mathbb{E}\left[\mathbb{E}_{X_i,\xi_i}\left[Z_i^{\prime}(t) g^{\prime}\left(Z^i(t)\right)\right]\right]\\
        &\ \text{ the inner expectation is conditioned on everything except }X_i,\xi_i\\
        & = \mathbb{E}\left[g^{\prime}\left(Z^i(t)\right)\mathbb{E}_{X_i,\xi_i}\left[Z_i^{\prime}(t)\right]\right]\\
        & = \mathbb{E}\left[g^{\prime}\left(Z^i(t)\right)\cdot 0\right] = 0.
    \end{aligned}
\end{equation}

Now we bound $\mathfrak{B}$. In many simpler cases, $Z(t) - Z^i(t)$ is just $Z_i(t)$. Under our cross-validation case, we will have an extra term in $Z(t) - Z^i(t)$---which is ultimately due to our different definition of $K^i$ than the standard Slepian interpolation as discussed earlier. We denote the residual as $S^i(t)$:
\begin{equation}
    S^i(t) := Z(t) - Z^i(t) - Z_i(t).
\end{equation}
For the ease of notation, we will also need a quantity closely related to it:
\begin{equation}
    S^i = \omega \sqrt{n}(K - K^i ) - K_i = \sum_{j\ne i} \nabla_i K_j \implies S^i(t) = \sqrt{t} w^{-1} n^{-1/2} S^i.
\end{equation}
The $\mathfrak{B}$ term can be simplified as
\begin{equation}\label{eq: bound second term in Slepian}
\begin{aligned}
& \sum_{i=1}^n \mathbb{E}\left[Z_i^{\prime}(t)\left(Z_i(t)+S^i(t)\right) g^{\prime \prime}\left(Z^i(t)\right)\right] \\
& \stackrel{(II)}{=}\sum_{i=1}^n \mathbb{E}\left[Z_i^{\prime}(t) S^i(t) g^{\prime \prime}\left(Z^i(t)\right)\right] \\
& =\frac{1}{2 \omega^2 n} \sum_{i=1}^n \mathbb{E}\left[K_i S^i g^{\prime \prime}\left(Z^i(t)\right) - \frac{\sqrt{t}}{\sqrt{1-t}} Y_i S^i g^{\prime \prime}\left(Z^i(t)\right)\right] \\
& \hspace{-2mm}\stackrel{(III)}{=}\frac{1}{2 \omega^2 n} \sum_{i=1}^n \mathbb{E}\left[K_i S^i g^{\prime \prime}\left(Z^i(t)\right)\right]\\
& \lesssim \omega^{-2}\left\|g^{\prime \prime}\right\|_{\infty} \sup _i \mathbb{E}\left|K_i S^i\right|\\
& \leq \omega^{-2}\left\|g^{\prime \prime}\right\|_{\infty} \sup _i\left\|K_i\right\|_2\left\|S^i\right\|_2\\
& = \omega^{-1}\left\|g^{\prime \prime}\right\|_{\infty} \sup _i\left\|S^i\right\|_2.
\end{aligned}
\end{equation}
The details of step $(II)$ are presented in Lemma~\ref{lemma: handle cross term}. In step $(III)$, we used 
\begin{equation}
    \begin{aligned}
        & \mathbb{E}\left[Y_i S^i g^{\prime \prime}\left(Z^i(t)\right)\right]\\&=\mathbb{E}\left[\mathbb{E}_{\xi_i}\left[Y_i S^i g^{\prime \prime}\left(Z^i(t)\right)\right]\right]\\
        & =\mathbb{E}\left[S^i g^{\prime \prime}\left(Z^i(t)\right) \mathbb{E}_{\xi_i}\left[Y_i\right]\right]=0.
    \end{aligned}
\end{equation}
We also recall that $\omega^{2}$ was assumed to be bounded away from zero and $K_i$ takes value in a bounded interval, so the last line of \eqref{eq: bound second term in Slepian} is essentially $\|S^i\|_2$.

As for $\mathfrak{C}$ in \eqref{eq: slepian integrand decomposition}, we consider the following. Because the third derivative of $g$ is bounded, it suffices to derive an upper bound for $A_i := \E{Z'_i(t) (Z_i(t) + S^i(t))^2}$. Observe that 
$$
    \abs{A_i} \le \E{|Z_i'(t)| (Z_i(t) + S^i(t))^2} \lesssim \Enb |Z'_i (t)|Z^2_i(t) + \Enb |Z'_i(t)| \{S^{i} (t)\}^2. 
$$
Under the boundedness of $K_i$, one can show that both $\sup_i \norm{K_i}^3_3$ and $\sup_i \norm{Y_i}^3_3$ are bounded. It follows that
\begin{equation}\label{eq: third moment bound on K}
\begin{aligned}
\mathbb{E}\left|Z_i^{\prime}(t)\right| Z_i^2(t) & \leq \frac{1}{2 \omega^3 n^{3 / 2}} \mathbb{E}\left[\left(\frac{\left|K_i\right|}{\sqrt{t}}+\frac{\left|Y_i\right|}{\sqrt{1-t}}\right)\cdot\left(\sqrt{t} K_i+\sqrt{1-t} Y_i\right)^2\right] \\
& \lesssim \frac{1}{4 n^{3 / 2}}\left(\frac{1}{\sqrt{t}} \vee \frac{1}{\sqrt{1-t}}\right) \cdot \mathbb{E}\left(\left|K_i\right|+\left|Y_i\right|\right)\left(\left|K_i\right|+\left|Y_i\right|\right)^2 \\
& \leq \frac{2}{n^{3 / 2}}\left(\frac{1}{\sqrt{t}} \vee \frac{1}{\sqrt{1-t}}\right) \cdot \mathbb{E}\left[\max \left\{\left|K_i\right|^3,\left|Y_i\right|^3\right\}\right] \\
& \leq \frac{2}{n^{3 / 2}}\left(\frac{1}{\sqrt{t}} \vee \frac{1}{\sqrt{1-t}}\right) \cdot\left(\sup _i\left\|K_i\right\|_3^3+\sup _i\left\|Y_i\right\|_3^3\right)\\
& \lesssim n^{-3 / 2}.
\end{aligned}
\end{equation}
Similarly, we have
\begin{equation}\label{eq: where Ks2 appears}
\begin{aligned}
\mathbb{E}\left|Z_i^{\prime}(t)\right|\left\{S^i(t)\right\}^2 &\lesssim \frac{1}{\omega^3 n^{3 / 2}}\left(\frac{1}{\sqrt{t}} \vee \frac{1}{\sqrt{1-t}}\right) \cdot\left(\mathbb{E}\left|K_i\right|\left\{S^i\right\}^2+\mathbb{E}\left|Y_i\right|\left\{S^i\right\}^2\right) \\
& \lesssim \frac{1}{n^{3 / 2}}\left(\frac{1}{\sqrt{t}} \vee \frac{1}{\sqrt{1-t}}\right) \cdot\left(\mathbb{E}\left|K_i\right|\left\{S^i\right\}^2+\mathbb{E}\left|Y_i\right|\left\{S^i\right\}^2\right). 
\end{aligned}
\end{equation}
By the boundedness of $K_i$ again, we have $\Enb |K_i| \{S^i\}^2 \lesssim \Enb \{S^i\}^2$ and $\Enb|Y_i| \{S^i \}^2 =\Enb\abs{\xi_i} \cdot \E{\tilde{\sigma}_{i} \{S^i \}^2} \lesssim \Enb \{S^i \}^2$. To establish the last inequality we also used that the first moment of a folded normal distribution is finite. 
Therefore, for $\mathfrak{C}$ defined in~\eqref{eq: slepian integrand decomposition} we have
\begin{equation}
   \mathfrak{C}\lesssim \|g'''\|_{\infty}(n\cdot n^{-3/2} + n\cdot n^{-3/2}\sup_i\|S^i\|_2^2). 
\end{equation} 
Overall, we have shown
\begin{equation}
\begin{aligned}\label{eq: part 1 final bound}
 |\mathbb{P}(K \leq x)-\mathbb{P}(Y \leq x)|& \leq |\mathbb{E} g(K)-\mathbb{E} g(Y)| + \epsilon\\
 & \leq C_\epsilon(\omega^{-1} \sup_i\|S^i\|_2 + n^{-1/2} + n^{-1/2}\sup_i\|S^i\|_2^2) + \epsilon,
\end{aligned}
\end{equation}
where $C_\epsilon$ is a constant depending on $\epsilon$. Note that an explicit bound of $\|S^i\|_2$ based on stability is given in Lemma~\ref{lemma: Di}. Typically, the leading term in \eqref{eq: part 1 final bound} is $\omega^{-1}  \sup _i\left\|S^i\right\|_2$, which comes from the $\mathfrak{B}$ term. We have the following bound 
\begin{equation}
    \label{eq: stability bound}
    \begin{aligned}
       &\omega^{-1} \sup _i\left\|S^i\right\|_2 \\ 
       & \lesssim \left(\omega^{-2}\sum_{j \neq i}\left\|\nabla_i K_j\right\|_2^2+\omega^{-2}\sum_{j \neq i \neq k}\left\|\nabla_k \nabla_i K_j\right\|_2\left\|\nabla_j \nabla_i K_k\right\|_2\right)^{1 / 2}\\
& \lesssim \left(\sum_{j \neq i}\left\|\nabla_i K_j\right\|_2^2+ \sum_{j \neq i \neq k}\left\|\nabla_k \nabla_i K_j\right\|_2\left\|\nabla_j \nabla_i K_k\right\|_2\right)^{1 / 2}
    \end{aligned}
\end{equation}
\textbf{Part 2}
Now we need to bound $|\mathbb{P}(Y \leq x) - \Phi(x) |$ for all $x\in\mathbb{R}$ where $\Phi$ is the cumulative distribution function of the standard normal distribution. Let $W$ denote a standard normal random variable. We have
\begin{equation}
    \begin{aligned}
        & |\mathbb{P}(Y \leq x) - \mathbb{P}(W \leq x)|\\
        & \stackrel{(I)}{=} |\mathbb{P}(\tau W \leq x) - \mathbb{P}(W \leq x)|\\
        & \stackrel{(II)}{\leq} |\mathbb{E}[g(\tau W) - g(W)]| + \epsilon\\
        & \leq \|g'\|_\infty \mathbb{E}|\tau W - W|+ \epsilon.
    \end{aligned}
\end{equation}
The $\tau$ is step $(I)$ is a random variable such that $\tau^{2}=n^{-1} \omega^{-2}\sum_{i=1}^{n} \tilde{\sigma}_{i}^{2}$. In this step, we also used that $\tau^{-1}Y$ is standard normal. When conditioned on the data $\bm{X}$, $\tau^{-1}Y$ is standard normal, which implies $\tau^{-1} Y$ is standard normal marginally.
The $g = g_\epsilon$ in step $(II)$ is the smooth indicator function we used earlier. Now we analyze the $\mathbb{E}\left|\tau W - W\right|$ term:
\begin{equation}
    \begin{aligned}
        \mathbb{E}|\tau W-W|
        & \leq \|\tau - 1\|_2 \|W\|_2 \leq \|\tau^2 - 1\|_2  \\
        & \le \frac{1}{\omega^2 n} \sum_{i = 1}^n \norm{\tilde{\sigma}_i^2 - \omega^2}_2 \\
        &= \frac{1}{\omega^2 n} \sum_{i=1}^n \sqrt{\mathrm{Var}(\tilde{\sigma}_i^2)},
    \end{aligned}
\end{equation}
where the last step holds true because for any $i \in [n]$, we have $ \Enb[\tilde{\sigma}_i^2] = \E{\E{K_i^2 \mid \bm{X}^{(-j)}}} = \omega^2$. Now we bound $\mathrm{Var}(\tilde{\sigma}_i^2)$ for all $i \in [n]$ by the Efron-Stein's inequality.
\begin{equation}
    \begin{aligned}
\operatorname{Var}\left(\tilde{\sigma}_i^2\right) 
& = \operatorname{Var}(\mathbb{E}\left[K_i^2 \mid \bm{X}\right])\\
& = \operatorname{Var}(\mathbb{E}\left[\mathcal{K}^2(i,\bm{X}) \mid \bm{X}\right])\\
& \leq \frac{1}{2}\sum_{j \neq i} \mathbb{E}\left(\mathbb{E}\left[\mathcal{K}^2(i,\bm{X}) \mid \bm{X}\right] - \mathbb{E}\left[\mathcal{K}^2(i,\bm{X}^{j}) \mid \bm{X}^{j}\right]\right)^2\\
& = \frac{1}{2}\sum_{j \neq i} \mathbb{E}\left(\mathbb{E}\left[\mathcal{K}^2(i,\bm{X}) -\mathcal{K}^2(i,\bm{X}^{j}) \mid \bm{X}, X_j^\prime\right]\right)^2\\
& \stackrel{(I)}{\leq} \frac{1}{2}\sum_{j \neq i} \mathbb{E}\left[\left\{\mathcal{K}^2(i,\bm{X}) -\mathcal{K}^2(i,\bm{X}^{j})\right\}^2 \right]\\
& = \frac{1}{2}\sum_{j \neq i} \mathbb{E}\left[(\mathcal{K}(i,\bm{X}) + \mathcal{K}(i,\bm{X}^{j}))^2 (\mathcal{K}(i,\bm{X}) - \mathcal{K}(i,\bm{X}^{j}))^2  \right]\\
& \stackrel{(II)}{\lesssim} n \max_{j \neq i}\mathbb{E}\left[ (\mathcal{K}(i,\bm{X}) - \mathcal{K}(i,\bm{X}^{j}))^2  \right] = n\max_{j\neq i} \mathbb{E}\left(\nabla_j K_i\right)^2.
    \end{aligned}
\end{equation}
In step $(I)$ we used the Jensen's inequality and in step $(II)$ we applied the boundedness of $\mathcal{K}$-mapping. Therefore, we can conclude that
\begin{equation}
\begin{aligned}
    & |\mathbb{P}(Y \leq x)-\mathbb{P}(W \leq x) |\\
    \leq &  \left\|g^{\prime}\right\|_{\infty} \mathbb{E}|\tau W-W|+\epsilon\\
    \leq & C_\epsilon \sqrt{n \max_{i\neq j \in [n]} \mathbb{E}\left(\nabla_j K_i\right)^2} + \epsilon.
\end{aligned}
\end{equation}

\textbf{Conclusion} Combining Parts 1 and 2, we have: for any $x\in\mathbb{R}$ and any $\epsilon > 0$, there is a constant $C_\epsilon > 0$ that only depends on $\epsilon$ such that
\begin{equation}
    |\mathbb{P}(K \leq x)- \Phi(x)| \leq C_\epsilon ( n^{-1/2} + n^{1/2} \Delta_1 + n\Delta_2 + n^{1/2} \Delta_1^2+ n^{3/2} \Delta_2^2 ) + 2\epsilon,
\end{equation}
where $\Delta_1 = \sqrt{\max _{i \neq j \in[n]} \mathbb{E}\left(\nabla_j K_i\right)^2}$ and $\Delta_2=\sqrt{\max _{i \neq j \neq k \in[n]} \mathbb{E}\left(\nabla_k \nabla_i K_j\right)^2}$.


\end{proof}

\begin{remark}
As in the standard Slepian interpolation (e.g. \cite{wasserman2014stein}), matching the first and second moments of $K_j$ and $Y_j$ is what we actually needed (to cancel out certain terms in the proof). If we have chosen non-normal $Y_j$'s to achieve this, then in the \text{Part 2} of this proof, we need to engage with one extra central limit theorem to show $Y$ is approaching normal, which induces unnecessary steps.
\end{remark}

\subsection{Technical Lemmas for Theorem~\ref{th: weakly dependent CLT}}
\begin{lemma}\label{lemma: handle cross term} 
Following the same notation as in \eqref{eq: bound second term in Slepian}, we have 
    \begin{equation}
        \mathbb{E}\left[Z_i^{\prime}(t) Z_i(t) g^{\prime \prime}\left(Z^i(t)\right)\right] = 0.
    \end{equation}
\end{lemma}
\begin{proof}
    We use the definition of $Z_i(t)$ and $Z^i(t)$. A direct computation gives
    \begin{equation}
\begin{aligned}
\mathbb{E}_{X_i, \xi_i}\left[Z_i^{\prime}(t) Z_i(t)\right] & =\frac{1}{2\omega^2 n} \mathbb{E}_{X_i, \xi_i}\left[K_i^2-Y_i^2+\left(\frac{\sqrt{1-t}}{\sqrt{t}}-\frac{\sqrt{t}}{\sqrt{1-t}}\right) K_i Y_i\right] \\
& =\frac{1}{2 \omega^2 n}\left(\tilde{\sigma}_{i}^2-\tilde{\sigma}_{i}^2+\left(\frac{\sqrt{1-t}}{\sqrt{t}}-\frac{\sqrt{t}}{\sqrt{1-t}}\right)\mathbb{E}_{X_i, \xi_i}\left[K_i Y_i\right]\right) \\
& =\frac{1}{2\omega^2  n}\left(\frac{\sqrt{1-t}}{\sqrt{t}}-\frac{\sqrt{t}}{\sqrt{1-t}}\right) \mathbb{E}_{X_i, \xi_i}\left[K_i\right] \cdot \mathbb{E}_{X_i, \xi_i}\left[Y_i\right]=0,
\end{aligned}
\end{equation}
where $\mathbb{E}_{X,\xi}[\cdot]$ means taking expectation conditioned on everything except $X$ and $\xi$. This yields
\begin{equation}
    \begin{aligned}
        &\mathbb{E}\left[Z_i^{\prime}(t) Z_i(t) g^{\prime \prime}\left(Z^i(t)\right)\right]\\
        &=\mathbb{E}\left[\mathbb{E}_{X_i, \xi_i}\left[Z_i^{\prime}(t) Z_i(t) g^{\prime \prime}\left(Z^i(t)\right)\right]\right]\\
        &=\mathbb{E}\left[g^{\prime \prime}\left(Z^i(t)\right) \mathbb{E}_{X_i, \xi_i}\left[Z_i^{\prime}(t) Z^i(t)\right]\right]=0 .
    \end{aligned}
\end{equation}
\end{proof}

\begin{lemma} \label{lemma: Di}
Let $S^i$ be the random variable defined in the proof of Theorem~\ref{th: weakly dependent CLT}\, Then,
\begin{equation}
    \|S^i\|_2^2 \le\sum_{j \neq i} \left\|\nabla_i K_j\right\|_2^2+\sum_{j \neq i \neq k} \left\|\nabla_k \nabla_i K_j\right\|_2\left\|\nabla_j \nabla_i K_k\right\|_2.
\end{equation}
Moreover, $\left\|S^{i}\right\|_{2}=o(1)$ provided the $K_i$'s satisfying the stability conditions \eqref{eq: stability}.
\end{lemma}


\begin{proof}
The quantity $S^i$ was defined as
\begin{equation}
\begin{aligned}
S^i & =\omega \sqrt{n}\left(K-K^i\right)-K_i \\
& = \sum_{j \neq i}^n \left( \mathcal{K}\left(j, \bm{X}\right)-\mathcal{K}\left(j, \bm{X}^{i}\right)\right) \\
& = \sum_{j \neq i}^n \nabla_i K_j.
\end{aligned}
\end{equation}
A direct calculation gives
$$
\begin{aligned}
\Enb \{S^{i}\}^{2} & =\sum_{j,k \neq i} \Enb \nabla_{i} K_{j} \nabla_{i} K_{k} \\
& =\sum_{j \neq i } \Enb\left(\nabla_{i} K_{j}\right)^{2}+\sum_{j \neq i \neq k} \Enb \nabla_{i} K_{j} \nabla_{i} K_{k}.
\end{aligned}
$$
The first summation is $o(1)$ because we assumed $\sup _{j}\left\|\nabla_{i} K_{j}\right\|_{2}=o(n^{-1/2})$. The second summation would be $o(1)$ if we could show $\left|\mathbb{E} \nabla_{i} K_{j} \nabla_{i} K_{k}\right| = o\left(n^{-2}\right)$. Applying \Cref{lemma: nabla operation} twice, we obtain
\begin{align*}
\abs{\mathbb{E}\left(\nabla_{i} K_{j}\right)\left(\nabla_{i} K_{k}\right)}
&=\abs{\mathbb{E}\left(\nabla_{i} K_{j}\right)\left(\nabla_{j} \nabla_{i} K_{k}\right)}=\abs{\mathbb{E}\left(\nabla_{k} \nabla_{i} K_{j}\right)\left(\nabla_{j} \nabla_{i} K_{k}\right)} \\
& \le \left\|\nabla_{k} \nabla_{i} K_{j}\right\|_{2}\left\|\nabla_{j} \nabla_{i} K_{k}\right\|_{2}=o\left(n^{-2}\right)
\end{align*}
as desired.
\end{proof}

\begin{lemma}\label{lemma: general free nabla}
Let $A$, $B$, and $C$ be collections of random variables.  
For measurable functions $f$ and $g$, define $f(A)$ and $g(B)$.  
Assume
\begin{equation*}
\mathbb{E}\!\left[f(A)\,\middle|\,C\right]=0, 
\qquad\text{and}\qquad 
A \perp\!\!\!\perp B \mid C .
\end{equation*}
Then
\begin{equation*}
\mathbb{E}\!\left[f(A)\,g(B)\right]=0 .
\end{equation*}
\end{lemma}

\begin{proof}
\begin{equation*}
\begin{aligned}
\mathbb{E}\!\left[f(A)\,g(B)\right] 
&= \mathbb{E}\!\left[g(B)\,\mathbb{E}\!\left[f(A)\,\middle|\,B,C\right]\right] \\[4pt]
&= \mathbb{E}\!\left[g(B)\,\mathbb{E}\!\left[f(A)\,\middle|\,C\right]\right]  \\[4pt]
&= \mathbb{E}\!\left[g(B)\cdot 0\right] \\[4pt]
&= 0 .
\end{aligned}
\end{equation*}
\end{proof}

\begin{corollary}   
\label{lemma: nabla operation}
Let $i\neq j\neq  k \in [n]$ be a triplet of distinct elements. Let $f$ be a function of the data $\bm{X}$.
Then we have
\begin{equation*}
\mathbb{E}\left(\nabla_i K_j\right) f=\mathbb{E}\left(\nabla_i K_j\right)\left(\nabla_j f\right)
\end{equation*}
and
\begin{equation*}
\mathbb{E}\left(\nabla_k \nabla_i K_j\right) f=\mathbb{E}\left(\nabla_k \nabla_i K_j\right)\left(\nabla_j f\right) .
\end{equation*}
\end{corollary} 
\begin{proof}
For the first equality, it suffices to prove 
\[
\mathbb{E}\left(\nabla_i K_j\right) f\left(\boldsymbol{X}^j\right) = 0.
\]
Indeed,
\begin{equation*}
\mathbb{E}\left(\nabla_i K_j\right) f\left(\boldsymbol{X}^j\right)
= \mathbb{E}\left[ K_j\, f\left(\boldsymbol{X}^j\right) \right]
- \mathbb{E}\left[ K_j\left(\boldsymbol{X}^i\right)\, f\left(\boldsymbol{X}^j\right) \right].
\end{equation*}

We apply \Cref{lemma: general free nabla} with 
\begin{align*}
\text{First term:} \quad & A = \boldsymbol{X}, \quad B = \boldsymbol{X}^j, \quad C = \boldsymbol{X}^{(-j)}; \\
\text{Second term:} \quad & A = \boldsymbol{X}^i, \quad B = \boldsymbol{X}^j, \quad C = \boldsymbol{X}^i \setminus \{X_j\}.
\end{align*}

By the same argument, we can conclude the second equation.
\end{proof}


\section{Stability Properties of Exponential Weighting}
\label{app: stability}
In this section, we present some stability properties for the proposed exponential weighting scheme. They will eventually be employed to derive the central limit theorem of the test statistic $T_r$ in Algorithm~\ref{algorithm: exp weighting}. We first recall the notations in the main text. For any $i \in [n]$, $X'_i$ is an IID copy of $X_i$ and $\bm{X}^i$ is the sample obtained by replacing/perturbing the $i$-th sample $X_i$ with $X_i'$. Also, we define an operator $\nabla_i$ such that $\nabla_i f(\bm{X}) = f(\bm{X}) - f(\bm{X}^i)$ for any function $f$ with respect to the sample $\bm{X}$.

\begin{lemma}[First Order Stability] 
\label{lemma: first order stability}
    Let $r \in [p]$ be the dimension of interest and let $j \in [n]$ and $i \notin I_{v_j}$ be two sample indices. Define
    \begin{equation}
        K_{j,r} = \mathcal{K}_r(j, \bm{X}) := X_{j, r} - Q_{j, r}\QexcludingVj - \E{X_{j, r} - Q_{j, r}\QexcludingVj | \bm{X}^{(-v_j)}}.
    \end{equation} 
    If the dimensions of $X_1 - \mu$ are uniformly bounded by a constant $M > 0$ almost surely, then we have 
    \begin{equation}\label{eq: uniform first stability}
        \max_{i,j,r} \norm{\nabla_i K_{j,r}}_2 \leq C\lambda \varsigma M n^{-1} \leq C\lambda M^2 n^{-1} 
    \end{equation}
    for $n \geq 8\lambda M$, a universal constant $C > 0$ and $\varsigma^2 = \max_{r \in [p]} \mathrm{Var}(X_{1,r})$. Specifically, when $\lambda = o(\sqrt{n})$, we have $\max_{i,j,r}\norm{\nabla_i K_{j,r}}_2 = o(n^{-1/2})$.
\end{lemma}

\begin{proof}
According to the definition of $\nabla_i K_{j,r}$, one has
\begin{equation}
    \begin{aligned}
         &\Enb (\nabla_i K_{j,r})^2 \\
         &= \Enb(\mathcal{K}_r(j, \bm{X}) - \mathcal{K}_r(j, \bm{X}))^2 \\
        & = \Enb \left ( X_{j,r} - Q_{j,r}\QexcludingV - \Enb[X_{j,r} - Q_{j,r}\QexcludingV|\bm X^{(-v)}] - (X_{j,r} - Q_{j,r}\QexcludingVperturbi - \mathbb{E}[X_{j,r} - Q_{j,r}\QexcludingVperturbi|\bm X^{(-v),i}]) \right )^2\\
        & = \mathbb{E}\left(Q_{j,r}\QexcludingVperturbi-\mathbb{E}\left[Q_{j,r}\QexcludingVperturbi \mid \boldsymbol{X}^{(-v), i}\right]-Q_{j,r}\QexcludingV+\mathbb{E}\left[Q_{j,r}\QexcludingV \mid \boldsymbol{X}^{(-v)}\right]\right)^2.
    \end{aligned}
\end{equation}
The quantity $Q_{j,r}\QexcludingVperturbi$ represents the weighted competitor $Q_{j,r}\QexcludingV$ but whose weights are computed with $\hat{\mu}^{(-v),i}$, i.e., the out-of-fold mean with $X_i$ replaced by an IID copy $X_i'$. To simplify notation, 
we omit the superscript $(-v)$ for every exponential weighting $\hat w_{r,s}^{(-v)}$ and sample mean $\hat{\mu}^{(-v)}$ in the following. 


To derive an upper bound for $\Enb \{Q_{j,r}^i-\mathbb{E}[Q_{j,r}^i\mid \bm X^{(-v),i}] - (Q_{j,r} - \mathbb{E}[Q_{j,r}\mid \bm X^{(-v)}] )\}^2$, we observe that 
\begin{equation}
\begin{aligned}
& \left|Q_{j,r}^i-\mathbb{E}\left[Q_{j,r}^i \mid \boldsymbol{X}^{(-v), i}\right]-Q_{j,r}+\mathbb{E}\left[Q_{j,r} \mid \boldsymbol{X}^{(-v)}\right]\right| \\
& =\left|\sum_{s \ne r} \left(\hat{w}_{r,s}^i-\hat{w}_{r,s}\right)\left(X_{j, s}-\mu_s\right)\right| \leq \sum_{s\ne r}\left|\frac{\hat{w}_{r,s}^i}{\hat{w}_{r,s}}-1\right| \cdot \hat{w}_{r,s} \cdot\left|X_{j, s}-\mu_s\right| .
\end{aligned}
\end{equation}
In particular, we investigate the absolute difference between the ratio $\hat{w}_{r,s}^{ i}/\hat{w}_{r,s}$ and $1$. Define $\tilde n = n(1-1 / V)$. For any $s \ne r$, we have
\begin{equation}
    \label{eq: bound weight ratio}
    \begin{aligned}
    \hat{w}_{r,s}^{i}/\hat{w}_{r,s} & = \frac{\exp(-\lambda \hat{\mu}_s^{i})}{\sum_{t \ne r} \exp(-\lambda \hat{\mu}_{t}^{i})} \cdot \frac{\sum_{t \ne r} \exp(-\lambda \hat{\mu}_t)}{\exp(-\lambda \hat{\mu}_{s})}\\  
    &= \exp(-\lambda \tilde n^{-1} (X'_{i,s} - X_{i,s}))\frac{\sum_{t \ne r} \exp(-\lambda\hat{\mu}_t)}{\sum_{t \ne r} \exp(-\lambda\hat{\mu}^{i}_t)}\\  
    & = \exp(-\lambda \tilde n^{-1} (X'_{i,s} - X_{i,s})) \frac{\sum_{t \ne r} \exp(-\lambda\hat{\mu}^{i}_t)\exp(-\lambda \tilde n^{-1}  (X_{i,t} - X_{i,t}'))}{\sum_{t \ne r} \exp(-\lambda\hat{\mu}^{i}_t)}\\ 
    & \le \exp \left (2 \lambda \tilde n^{-1}  \max_{t \in [p]} \abs{X_{i, t}' - X_{i ,t}} \right ) \le \exp (4 \lambda \tilde n^{-1}  M). 
    \end{aligned}
\end{equation}
Then, the mean value theorem gives
\begin{equation}
    \hat{w}_{r,s}^{i}/\hat{w}_{r,s} - 1 \le 4 \lambda \tilde n^{-1}  M \exp (\xi)
\end{equation}
for some universal $\xi \in (0, 4 \lambda \tilde n^{-1}  M)$. Provided that $4 \lambda \tilde n^{-1}  M\le 1$, we further have 
$$
\hat{w}_{r,s}^{i}/\hat{w}_{r,s}- 1 \leq  4 
e\lambda\tilde n^{-1}  M.
$$ 
Similarly, one can obtain $\hat{w}_{r,s}^i / \hat{w}_{r,s}-1 \geq-4 \lambda \tilde{n}^{-1} M$.

It follows that
\begin{equation}  \label{eq: weight ratio from 1}
    \abs{\hat{w}_{r,s}^{i}/\hat{w}_{r,s}- 1} \leq 4e \lambda \tilde n^{-1}  M
\end{equation}
and
\begin{equation}\label{eq: bound Q difference}
\begin{aligned}
& \left|Q_{j,r}^i-\mathbb{E}\left[Q_{j,r}^i \mid \boldsymbol{X}^{(-v), i}\right]-Q_{j,r}+\mathbb{E}\left[Q_{j,r} \mid \boldsymbol{X}^{(-v)}\right]\right| \\
& \leq  4e\lambda\tilde n^{-1}  M \sum_{s\ne r} \hat{w}_{r,s}\left|X_{j, s}-\mu_s\right|. 
\end{aligned}
\end{equation}
It follows the Jensen's inequality that
\begin{equation}
\label{eq: jensen for stability}
\begin{aligned}
     &\Enb \left (Q_{j,r}^i-\mathbb{E}\left[Q_{j,r}^i \mid \boldsymbol{X}^{(-v), i}\right]-Q_{j,r}+\mathbb{E}\left[Q_{j,r} \mid \boldsymbol{X}^{(-v)}\right]\right)^2 \\
     & \le 16 e^2 \lambda^2 \tilde{n}^{-2} M^2 \Enb \left (\sum_{s \ne r} \hat{w}_{r,s} \abs{X_{j,s} - \mu_s} \right)^2 \\
    & = 16 e^2 \lambda^2 \tilde{n}^{-2} M^2 \Enb \left ( \E{ \sum_{s \ne r}\hat{w}_{r,s} \abs{X_{j,s}  - \mu_s} \mid \bm{X}^{(-v_j)}} \right)^2 \\
    & \le 16 e^2 \lambda^2 \tilde{n}^{-2} M^2 \E{\E{\sum_{s \ne r}\hat{w}_{r,s} (X_{j,s }- \mu_s)^2 \mid \bm{X}^{(-v_j)}}} \\
    & = 16 e^2 \lambda^2 \tilde{n}^{-2} M^2 \sum_{s \ne r} \E{\hat{w}_{r,s}} \mathrm{Var} (X_{j,s}) \\ 
    &= 16e^2\lambda^2 \varsigma^2 M^2 \tilde{n}^{-1}.  
\end{aligned}
\end{equation}
Since the bound does not depend on $i,j$ or $r$, we have the first uniform bound in \eqref{eq: uniform first stability}. The second uniform bound follows directly from the bound $\mathrm{Var}(X_{1,s}) \le M^2$ for all $s \in [p]$. 
\end{proof}

\begin{lemma}[Second Order Stability]
\label{lemma: second order stability} Let $r \in [p]$ be the dimension of interest, and let $j \in [n]$ and $i,k \notin I_{v_j}$ be some sample indices. Define

\begin{equation}
        K_{j,r} = \mathcal{K}_r(j, \bm{X}) := X_{j, r} - Q_{j,r}\QexcludingVj - \E{X_{j, r} - Q_{j,r}\QexcludingVj | \bm{X}^{(-v_j)}}.
    \end{equation} 
If the dimensions of $X_1 - \mu$ are uniformly bounded by a constant $M > 0$ almost surely, then we have 

\begin{equation}
    \max _{i, j, k, r}\left\|\nabla_i \nabla_k K_{j,r}\right\|_2 \le C\lambda^2 \varsigma M^2 n^{-2} \leq C\lambda^2 M^3 n^{-2}
    \end{equation}
for large enough $n$, a universal constant $C$ and $\varsigma^2 = \max_{r \in [p]} \mathrm{Var}(X_{1, r})$. In particular, when $\lambda = o(\sqrt{n})$, we have $\max_{i,j,k,r}\norm{\nabla_i \nabla_k K_{j,r}}_2 = o(n^{-1})$.

\end{lemma}

\begin{proof}
To simplify notation, we omit 
the superscript $(-v)$ for every exponential weighting $\hat{w}$ and sample mean $\hat{\mu}$. We also take $r=1$ and define $\tilde{n} = n(1 - 1/V)$. The bounds that we will establish are uniform over $i,j,k,r$.
    
By the definition of $\nabla_i \nabla_k K_{j,1}$ with $i,k \notin I_{v_j}$, one has
    \begin{equation}\label{eq: simplify second order}
        \begin{aligned}
            &\Enb(\nabla_{i}\nabla_{k} K_{j,1})^2 \\
            & = \Enb \left (\nabla_i[ \mathcal{K}_1 (X_j,\bm{X}^{(-v)})  - \mathcal{K}_1(X_j,\bm{X}^{(-v),k}) ] \right)^2\\
            & = \Enb \left ( \mathcal{K}_1(X_j,\bm{X}^{(-v)}) - \mathcal{K}_1(X_j,\bm{X}^{(-v),k}) - \mathcal{K}_1(X_j,\bm{X}^{(-v),i}) + \mathcal{K}_1(X_j,\bm{X}^{(-v),ik}) \right )^2\\
            & = \Enb \Bigg (Q_{j,1}^k - \E{Q_{j,1}^k \mid \bm{X}^{(-v),k}} - Q_{j,1} + \E{Q_{j,1} \mid \bm{X}^{(-v)}} \\
            & \hspace{1cm} + Q_{j,1}^{i} - \E{Q_{j,1}^{i} \mid \bm{X}^{(-v),i}} - Q_{j,1}^{ik} + \E{Q_{j,1}^{ik} \mid \bm{X}^{(-v), {ik}}} \Bigg )^2 \\
            & = \Enb \left (\sum_{s=2}^p (\hat{w}_{1,s}^k - \hat{w}_{1,s} + \hat{w}_{1,s}^{i} - \hat{w}_{1,s}^{ik}) (X_{j,s} - \mu_s) \right )^2.
        \end{aligned}
    \end{equation}
    The quantity $Q_{j, r}^i$ ($Q_{j, r}^{ik}$) represents the weighted competitor $Q_{j, r}$ whose weights are computed with $\hat{\mu}^{(-v),i}$ ($\hat{\mu}^{(-v), ik}$), i.e., the out-of-fold mean with $X_i$ replaced by $X_i'$ (with $X_i, X_k$ replaced by $X_i', X_k'$).  
    
    Observe that
    \begin{equation}\label{eq: disect Q}
        \begin{aligned}
            &\abs{ \sum_{s = 2}^p (\hat{w}_{1,s}^k - \hat{w}_{1,s} + \hat{w}_{1,s}^{i} - \hat{w}_{1,s}^{ik}) (X_{j,s} - \mu_s) } \\
            & \le \sum_{s = 2}^p \abs{\hat{w}_{1,s} - \hat{w}_{1,s}^{k} - \hat{w}_{1,s}^{i} + \hat{w}_{1,s}^{ik}} \cdot \abs{X_{j, s} - \mu_s} \\
            & = \sum_{s = 2}^p \abs{\left(\hat{w}_{1,s}\left(1-\frac{\hat{w}_{1,s}^{k}}{\hat{w}_{1,s}}\right) - \hat{w}_{1,s}^{i}\left(1-\frac{\hat{w}_{1,s}^{k}}{\hat{w}_{1,s}}\right)+ \hat{w}_{1,s}^{i}\left(1-\frac{\hat{w}_{1,s}^{k}}{\hat{w}_{1,s}}\right)-\hat{w}_{1,s}^{i}\left(1-\frac{\hat{w}_{1,s}^{ik}}{\hat{w}_{1,s}^{i}}\right)\right)} \cdot \abs{X_{j,s} - \mu_s}\\
            & =\sum_{s = 2}^p\abs{\left(\hat{w}_{1,s}-\hat{w}_{1,s}^{i}\right)\left(1-\frac{\hat{w}_{1,s}^{k}}{\hat{w}_{1,s}}\right)+\hat{w}_{1,s}^{i}\left(\frac{\hat{w}_{1,s}^{ik}}{\hat{w}_{1,s}^{i}}-\frac{\hat{w}_{1,s}^{k}}{\hat{w}_{1,s}}\right)} \cdot \abs{X_{j,s} - \mu_s}\\
            & \le \sum_{s = 2}^p \hat{w}_{1,s} \abs{\left(1-\frac{\hat{w}_{1,s}^{i}}{\hat{w}_{1,s}}\right)\left(1-\frac{\hat{w}_{1,s}^{k}}{\hat{w}_{1,s}}\right)} \cdot \abs{X_{j,s} - \mu_s} +\sum_{s = 2}^p \hat{w}_{1,s}^{i} \abs{\left(\frac{\hat{w}_{1,s}^{ik}}{\hat{w}_{1,s}^{i}}-\frac{\hat{w}_{1,s}^{k}}{\hat{w}_{1,s}}\right)} \cdot \abs{X_{j,s} - \mu_s}
        \end{aligned}
    \end{equation} 
One can follow the arguments in \eqref{eq: weight ratio from 1}~\eqref{eq: bound Q difference} and~\eqref{eq: jensen for stability} to bound the $L_2$ norm of the first summation in~\eqref{eq: disect Q}, up to a universal constant, by $\lambda^2 \varsigma M^2 \tilde{n}^{-2}$.

As for the second summation in~\eqref{eq: disect Q}, we investigate the absolute difference $\abs{\hat{w}_{1,s}^{ik}/\hat{w}_{1,s}^{i}-\hat{w}_{1,s}^{k}/\hat{w}_{1,s}}$ for each $s \in \{2, \ldots, p\}$.
Because
$$
    \frac{\exp \left(-\lambda \hat{\mu}_s^{ik}\right)}{\exp \left(-\lambda \hat{\mu}_s^i\right)} = \frac{\exp \left(-\lambda \hat{\mu}_s^{k}\right)}{\exp (-\lambda \hat{\mu}_s)} = \exp(-\lambda  \tilde n^{-1}(X_{k,s}^\prime - X_{k,s})),
$$
we have 
\begin{equation}\label{eq: weight for each Xk 2}
    \begin{aligned}
        \abs{\frac{\hat{w}_{1,s}^{ik}}{\hat{w}_{1,s}^{i}}-\frac{\hat{w}_{1,s}^{k}}{\hat{w}_{1,s}}} &= \exp \left(-\lambda \tilde n^{-1}\left(X_{k, s}'-X_{k,s}\right)\right) \abs{\left(\frac{\Xi^i}{\Xi^{ik}}-\frac{\Xi}{\Xi^k}\right)}\\
        & = \exp \left(-\lambda \tilde n^{-1}\left(X_{k, s}'-X_{k,s}\right)\right) \abs{\frac{\Xi}{\Xi^k}\left(\frac{\Xi^i\Xi^k}{\Xi^{ik}\Xi}-1\right)},
    \end{aligned}
\end{equation}
where $\Xi$ denotes the corresponding normalization constant for exponential weights $\{\hat{w}_{1,s}, s \in \{2, \ldots, p\}\}$. Similarly, $\Xi^i$, $\Xi^k$ and $\Xi^{ik}$ are for the one/two-sample perturbed weights, indicated by their respective superscripts. A direct calculation gives 
\begin{align*} \label{eq: bound ratio}
    \frac{\Xi^i \Xi^k}{\Xi^{ik} \Xi} &= \frac{\{\sum_{t=2}^p \exp(-\lambda \hat{\mu}_t^i)\}\{\sum_{t'=2}^p \exp(-\lambda \hat{\mu}_{t'}^k)\}}{\{\sum_{t=2}^p \exp(-\lambda \hat{\mu}_t^{ik})\}\{\sum_{t'=2}^p \exp(-\lambda \hat{\mu}_{t'})\}}\\
    & = \frac{\sum_{t=2}^p \exp(-\lambda (\hat{\mu}_t^i + \hat{\mu}_{t}^k)) + \sum_{2 \le t < t'} \left \{\exp(-\lambda (\hat{\mu}_t^i + \hat{\mu}_{t'}^k)) + \exp(-\lambda (\hat{\mu}_{t'}^i + \hat{\mu}_t^k)) \right\}}{\sum_{t=2}^p \exp(-\lambda (\hat{\mu}_t^{ik} + \hat{\mu}_t)) + \sum_{2 \le t < t'} \left\{ \exp(-\lambda (\hat{\mu}_t^{ik}+ \hat{\mu}_{t'})) + \exp(-\lambda (\hat{\mu}_{t'}^{ik} + \hat{\mu}_t)) \right\}}.
\end{align*}
To ease the notation, we write 
\begin{align*}
&E_{t,t}^{i,k} = \exp(-\lambda (\hat{\mu}_t^i + \hat{\mu}_t^k)), \hspace{3mm}E_{t,t'}^{i,k} = \exp(-\lambda (\hat{\mu}_t^i + \hat{\mu}_{t'}^k)), \hspace{3mm}E_{t',t}^{i,k} = \exp(-\lambda (\hat{\mu}_{t'}^i + \hat{\mu}_t^k)), \\
&E_{t,t}^{ik,\emptyset} = \exp (- \lambda(\hat{\mu}^{ik}_t + \hat{\mu}_t)), \hspace{3mm} E_{t,t'}^{ik,\emptyset} = \exp(-\lambda (\hat{\mu}_t^{ik}+ \hat{\mu}_{t'})),  \hspace{3mm} E_{t',t}^{ik,\emptyset} = \exp(-\lambda (\hat{\mu}_{t'}^{ik} + \hat{\mu}_t)).
\end{align*}
As $E_{t,t}^{i,k}=E_{t,t}^{ik, \emptyset}$, we obtain that
\begin{equation}\label{eq: bound on Xi}
\begin{aligned}
    \frac{\Xi^i \Xi^k}{\Xi^{ik} \Xi} - 1 &= \frac{\sum_{t=2}^p E_{t,t}^{i,k} + \sum_{2 \le t < t'} \left(E_{t,t'}^{i,k} + E_{t',t}^{i,k}\right)}{\sum_{t=2}^p E_{t,t}^{ik,\emptyset} + \sum_{2 \le t < t'} \left(E_{t,t'}^{ik,\emptyset}  + E_{t',t}^{ik,\emptyset}\right)} - 1\\
    & = \frac{\sum_{t=2}^p E_{t,t}^{ik,\emptyset} \left(E_{t,t}^{i,k}+E_{t,t}^{i,k}\right)\left(E_{t,t}^{ik, \emptyset}+E_{t,t}^{ik, \emptyset}\right)^{-1}}{\sum_{t=2}^p E_{t,t}^{ik,\emptyset} + \sum_{2 \le t < t'} \left(E_{t,t'}^{ik,\emptyset}  + E_{t',t}^{ik,\emptyset} \right)} \\
    &\hspace{4mm} + \frac{\sum_{2 \le t < t'}  (E_{t,t'}^{ik,\emptyset}  + E_{t',t}^{ik,\emptyset})(E_{t,t'}^{i,k} + E_{t',t}^{i,k})(E_{t,t'}^{ik,\emptyset}  + E_{t',t}^{ik,\emptyset})^{-1}}{\sum_{t=2}^p E_{t,t}^{ik,\emptyset} + \sum_{2 \le t < t'} \left(E_{t,t'}^{ik,\emptyset}  + E_{t',t}^{ik,\emptyset} \right)}  - 1\\
    & \leq \sup_{2\leq t \leq  t' \leq p} \left(E_{t,t'}^{i, k}+E_{t',t}^{i,k}\right)\left(E_{t,t'}^{ik, \emptyset}+E_{t',t}^{ik, \emptyset}\right)^{-1} - 1 \\
    & = \sup_{2 \le t \le t' \le p} \frac{E_{t,t'}^{i, k}+E_{t',t}^{i,k} - E_{t,t'}^{ik, \emptyset}-E_{t',t}^{ik, \emptyset}}{E_{t,t'}^{ik, \emptyset}+E_{t',t}^{ik, \emptyset}},
\end{aligned}
\end{equation}
and that $\frac{\Xi^i \Xi^k}{\Xi^{ik} \Xi} - 1 \ge \inf_{2 \le t \le t' \le p} \left( E_{t,t'}^{i, k}+E_{t',t}^{i,k} - E_{t,t'}^{ik, \emptyset}-E_{t',t}^{ik, \emptyset} \right)/\left ( E_{t,t'}^{ik, \emptyset}+E_{t',t}^{ik, \emptyset} \right)$ analogously.

Let $t, t' \in [p]$ be arbitrary. By the mean value theorem,
\begin{align*}
    &E_{t,t'}^{i, k} - E_{t,t'}^{ik, \emptyset} \\
    &= \exp \left(-\lambda\left(\hat{\mu}_t^i+\hat{\mu}_{t'}^k\right)\right) - \exp \left(-\lambda\left(\hat{\mu}_t^{ik}+\hat{\mu}_{t'}\right)\right) \\
    &= \exp(\xi_1)\lambda(\hat{\mu}_t^{ik}+\hat{\mu}_{t'} - \hat{\mu}_t^i-\hat{\mu}_{t'}^k)
\end{align*}
for some variable $\xi_1$ between $-\lambda\left(\hat{\mu}_t^i+\hat{\mu}_{t'}^k\right)$ and $-\lambda\left(\hat{\mu}_t^{ik}+\hat{\mu}_{t'}\right)$. Similarly, 
$$
E_{t',t}^{i,k}-E_{t',t}^{ik, \emptyset} = \exp(\xi_2)\lambda(\hat{\mu}_{t'}^{ik}+\hat{\mu}_t - \hat{\mu}_{t'}^i-\hat{\mu}_t^k)
$$
for some variable $\xi_2$ between $-\lambda\left(\hat{\mu}_{t'}^i+\hat{\mu}_t^k\right)$ and $-\lambda\left(\hat{\mu}_{t'}^{ik}+\hat{\mu}_t\right)$. 

Using this fact, we can achieve
\begin{equation}
    \begin{aligned}
        &E_{t,t'}^{i, k}+E_{t',t}^{i,k}-E_{t,t'}^{ik, \emptyset}-E_{t',t}^{ik, \emptyset} \\
        & = \exp \left(\xi_1\right) \lambda \left(\hat{\mu}_t^{ik}+\hat{\mu}_{t'}-\hat{\mu}_t^i-\hat{\mu}_{t'}^k\right) + \exp \left(\xi_2\right) \lambda\left(\hat{\mu}_{t'}^{ik}+\hat{\mu}_t-\hat{\mu}_{t'}^i-\hat{\mu}_t^k\right)\\
        & = \exp \left(\xi_1\right) \lambda \left(\hat{\mu}_t^{ik}+\hat{\mu}_{t'}-\hat{\mu}_t^i-\hat{\mu}_{t'}^k\right) + \exp \left(\xi_1\right) \lambda \left(\hat{\mu}_{t'}^{ik}+\hat{\mu}_t-\hat{\mu}_{t'}^i-\hat{\mu}_t^k\right)\\
        &\quad -\exp \left(\xi_1\right) \lambda \left(\hat{\mu}_{t'}^{ik}+\hat{\mu}_t-\hat{\mu}_{t'}^i-\hat{\mu}_t^k\right)+\exp \left(\xi_2\right) \lambda\left(\hat{\mu}_{t'}^{ik}+\hat{\mu}_t-\hat{\mu}_{t'}^i-\hat{\mu}_t^k\right)\\
        & = \{\exp(\xi_2) - \exp(\xi_1)\}\lambda \left(\hat{\mu}_{t'}^{ik}+\hat{\mu}_t-\hat{\mu}_{t'}^i-\hat{\mu}_t^k\right), 
    \end{aligned}
\end{equation}
where the last equality holds true because 
\begin{align*}
    \hat{\mu}^{ik}_{t} + \hat{\mu}_{t'} - \hat{\mu}^i_{t} - \hat{\mu}^k_{t'} &= \tilde{n}^{-1}(X'_{k,t} - X_{k,t} + X_{k,t'} - X_{k,t'}')  \\
    & = - (\hat{\mu}^{ik}_{t'} + \hat{\mu}_t - \hat{\mu}^i_{t'} - \hat{\mu}_t^k).
\end{align*}

Note that $\abs{\lambda \left(\hat{\mu}_{t'}^{ik}+\hat{\mu}_t-\hat{\mu}_{t'}^i-\hat{\mu}_t^k\right)} \le 4 \lambda  \tilde n^{-1} M$. By the mean value theorem again, the difference between the exponentials is $\exp \left(\xi_2\right)-\exp \left(\xi_1\right) = \exp(\xi_3) (\xi_2 - \xi_1)$ for some $\xi_3$ between $\xi_1$ and $\xi_2$. 
Particularly, the absolute difference between $\xi_2$ and $\xi_1$ is bounded by
\begin{align*}
\max \bigg \{ &\abs{\lambda\left(\hat{\mu}_{t'}^{ik}+\hat{\mu}_t-\hat{\mu}_t^i-\hat{\mu}_{t'}^k \right)}, \abs{\lambda\left(\hat{\mu}_{t'}^i+\hat{\mu}_t^k-\hat{\mu}_t^{ik}-\hat{\mu}_{t'}\right)}, \\ 
&\abs{\lambda\left (\hat{\mu}^{ik}_{t'} + \hat{\mu}_t - \hat{\mu}^{ik}_t - \hat{\mu}_{t'} \right)}, \abs{\lambda \left(\hat{\mu}^i_{t'} + \hat{\mu}^k_t - \hat{\mu}^i_t - \hat{\mu}^k_{t'} \right)} \bigg \}, 
\end{align*}
which is in turn bounded by $8 \lambda  \tilde n^{-1}M$. Furthermore, provided that $4 \lambda  \tilde n^{-1} M < 1$,
\begin{equation}
    \begin{aligned}
        &\frac{\exp(\xi_3)}{E_{t,t'}^{ik, \emptyset}+E_{t',t}^{ik, \emptyset}}\\
        & \le \frac{\exp \left ( \max\{-\lambda\left(\hat{\mu}_t^i+\hat{\mu}_{t'}^k\right), -\lambda\left(\hat{\mu}_t^{ik}+\hat{\mu}_{t'}\right), -\lambda\left(\hat{\mu}_{t'}^i+\hat{\mu}_t^k\right), -\lambda\left(\hat{\mu}_{t'}^{ik}+\hat{\mu}_t \right)\} \right)}{E_{t,t'}^{ik, \emptyset}+E_{t',t}^{ik, \emptyset}}\\
        &\le e,  
    \end{aligned}
\end{equation}
where one can obtain the last inequality by following the argument in~\eqref{eq: bound weight ratio} and discussing the values of 
\begin{align*}
    &\max\{-\lambda\left(\hat{\mu}_t^i+\hat{\mu}_{t'}^k\right), -\lambda\left(\hat{\mu}_t^{ik}+\hat{\mu}_{t'}\right), -\lambda\left(\hat{\mu}_{t'}^i+\hat{\mu}_t^k\right), -\lambda\left(\hat{\mu}_{t'}^{ik}+\hat{\mu}_t \right)\}, \\
    &\min\{-\lambda\left(\hat{\mu}_t^{ik}+\hat{\mu}_{t'}\right), -\lambda\left(\hat{\mu}_{t'}^{ik}+\hat{\mu}_t \right)\}.
\end{align*} 

Therefore, we have
$$
\sup_{2\le t \le t'\le p}\frac{E_{t,t'}^{i, k}+E_{t',t}^{i,k} - E_{t,t'}^{ik, \emptyset}-E_{t',t}^{ik, \emptyset}}{E_{t,t'}^{ik, \emptyset}+E_{t',t}^{ik, \emptyset}} \le 32 e\lambda^2 \tilde{n}^{-2}M^2.
$$
Similarly, one can obtain 
$$
    \inf_{2 \le t \le t' \le p} \frac{ E_{t,t'}^{i, k}+E_{t',t}^{i,k} - E_{t,t'}^{ik, \emptyset}-E_{t',t}^{ik, \emptyset}}{E_{t,t'}^{ik, \emptyset}+E_{t',t}^{ik, \emptyset} } \ge - \tilde{C} \lambda^2 \tilde{n}^{-2} M^2
$$
for some universal constant $\tilde{C} > 0$. 

Plugging these into \eqref{eq: bound on Xi}, we overall have $\abs{\frac{\Xi^i \Xi^k}{\Xi^{ik} \Xi} - 1} \leq  C\lambda^2 n^{-2} M^2$ for any $i, k \notin I_{v_j}$ with a universal constant $C = \max\{32e, \tilde{C}\}$. Going back to \eqref{eq: weight for each Xk 2}, it follows that when $4 \lambda \tilde n^{-1} M < 1$,
\begin{align*}
    \abs{\frac{w_{r,s}^{ik}}{w_{r,s}^{i}}-\frac{w_{r,s}^{k}}{w_{r,s}}} &= \exp \left(-\lambda\tilde{n}^{-1}\left(X_{k, s}'-X_{k,s}\right)\right) \abs{\frac{\Xi}{\Xi^k}\left(\frac{\Xi^i\Xi^k}{\Xi^{ik}\Xi}-1\right)} \\
    & \le \exp \left (4 \lambda \tilde{n}^{-1} M \right )  \cdot \abs{\frac{\Xi^i\Xi^k}{\Xi^{ik}\Xi}-1}\\
    &\le C\lambda^2 \tilde{n}^{-2} M^2, 
\end{align*}
where the first inequality holds true by the argument in~\eqref{eq: bound weight ratio}. This in turn makes the $L_2$ norm of the second summation in~\eqref{eq: disect Q} of order $O\left(\lambda^2 \varsigma M^2 \tilde{n}^{-2} \right)$ by Jensen's inequality, from which we can conclude the proof.
\end{proof}

\subsection{Stability under Extra Structure}

\begin{lemma}\label{lemma: ultra stable}
    Assume the same conditions as Lemma~\ref{lemma: first order stability}. We further assume that there exists an index $r^*\in[p]\symbol{92}\{r\}$ such that
    \begin{equation}
        \mu_{r^*} \leq \mu_s - \Delta_n,
    \end{equation}
    for all $s \in [p]\symbol{92}\{r,r^*\}$, and the gap $\Delta_n>0$ satisfies
    \begin{equation}\label{eq: definition of delta}
        \lim_{n\rightarrow\infty} \Delta_n \sqrt{n}/\log(pn) =\infty.
    \end{equation}
    Then for any choices of $\lambda = \lambda_n$, we always have $\max_{i,j}\left\|\nabla_i K_{j,r}\right\|_2=o\left(n^{-1 / 2}\right)$.
\end{lemma}

\begin{proof}
    To simplify notation, we omit the superscript $(-v)$ for every exponential weighting $\hat{w}$. Without loss of generality, we assume $r = 1$ and $r^* = 2$ to simplify the notation. We already established in the proof of Lemma~\ref{lemma: first order stability} that
    \begin{equation}
\mathbb{E}\left(\nabla_i K_{j,1}\right)^2 \leq \mathbb{E}\left[\left(\sum_{s=2}^p\left(\hat{w}_{1,s}^i-\hat{w}_{1,s}\right)\left(X_{j, s}-\mu_s\right)\right)^2\right] =:  \mathbb{E}[A^2].
\end{equation}
Define
\begin{equation}
\begin{array}{r}
\left|A_2\right|=\left|\hat{w}_{1,2}-1\right|\left|X_{j, 2}-\mu_2\right|, \\
\left|A_2^i\right|=\left|\hat{w}_{1,2}^i-1\right|\left|X_{j, 2}-\mu_2\right|, \\
\left|A_{+}\right|=\sum_{s=3}^p \hat{w}_{1,s}\left|X_{j, s}-\mu_s\right|, \\
\left|A_{+}^i\right|=\sum_{s=3}^p \hat{w}_{1,s}^i\left|X_{j, s}-\mu_s\right|.
\end{array}
\end{equation}
We have 
\begin{equation}
\begin{aligned}
|A| \leq & \left|A_2\right|+\left|A_2^i\right|+\left|A_{+}\right|+\left|A_{+}^i\right| \\
& \leq M\left(\left|\hat{w}_{1,2}-1\right|+\left|\hat{w}_{1,2}^i-1\right|+\sum_{s=3}^p \hat{w}_{1,s}+\sum_{s=3}^p \hat{w}_{1,s}^i\right) \\
& =2 M\left(\sum_{s=3}^p \hat{w}_{1,s}+\sum_{s=3}^p \hat{w}_{1,s}^i\right).
\end{aligned}
\end{equation}
Therefore
\begin{equation}
\mathbb{E}\left(\nabla_i K_{j,1}\right)^2 \leq 8 M^2 \mathbb{E}\left[\left(\sum_{s=3}^p \hat{w}_{1,s}\right)^2\right].
\end{equation}
Define an event 
\begin{equation}
\mathcal{E}_0=\left\{\hat{w}_{1,s} \leq \exp (-\lambda (\mu_s - \mu_2) / 2), \hspace{1mm}\forall s \geq 3\right\}.
\end{equation}
Then applying Lemma \ref{lemma: vanishing weights}, we have 
\begin{equation}
\mathbb{P}\left(\mathcal{E}_0^c\right) \leq 2 pe^{-n\Delta^2/(8M^4)}.
\end{equation}
Then we have 
\begin{equation}
\begin{aligned}         
\mathbb{E}\left(\nabla_i K_{j,1}\right)^2& \leq 8M^2\mathbb{E}\left[\left(\sum_{s=3}^p \hat{w}_{1,s}\right)^2 \mid \mathcal{E}_0\right] + 16M^2p e^{-n \Delta^2 /\left(8 M^4\right)}\\
& \leq 8M^2p^2e^{-\lambda\Delta} + 16 M^2 p e^{-n \Delta^2 /\left(8 M^4\right)}.
\end{aligned}
\end{equation} 
By the definition of $\Delta$ in~\eqref{eq: definition of delta}, $16 M^2 p e^{-n \Delta^2 /\left(8 M^4\right)} = o(n^{-1})$. For any $\lambda$ such that $\liminf_{n\rightarrow\infty} \lambda n^{-1/2} > 0$, we can verify that $8 M^2 p^2 e^{-\lambda \Delta}$ is also of order $o(n^{-1})$. The results in Lemma~\ref{lemma: first order stability} stated that for $\lambda = o(\sqrt{n})$, $\mathbb{E}\left(\nabla_i K_{j,1}\right)^2$ is of order $o(n^{-1})$. Combining these two pieces, we conclude that the stability result always holds regardless of the choice of $\lambda$.
\end{proof}

\section{Normality without Boundedness Assumptions}
\label{app: remove boundedness}
In this section, we extend our analysis of asymptotic normality to settings where $X_i$'s do not have bounded entries. To accommodate general light-tail distributions, we introduce a flexible class of random variables that take both sub-Gaussian and sub-exponential distributions as examples:

\begin{definition}
A random variable $X$ is said to follow a \emph{sub-Weibull} distribution with parameters $(\alpha, K)$, denoted as $(\alpha, K)$-sub-Weibull, if there exist constants $\alpha > 0$ and $K > 0$ such that for all $t \geq 0$,
\begin{equation*}
\mathbb{P}(|X| > t) \leq 2 \exp\left( -  t^\alpha/K \right).
\end{equation*}
An equivalent characterization is that there exists a constant $C = C(\alpha, K) > 0$ such that for all $q \geq 1$,
\begin{equation*}
\left( \mathbb{E}|X|^q \right)^{1/q} \leq C\, q^{1/\alpha}.
\end{equation*}
\end{definition}

\begin{proposition}\label{th: unbounded random center}
Let $X_i \in \mathbb{R}^p$ be a collection of IID random vectors. Suppose each coordinate $X_{1,t}$ for $t \in [p]$ is $(\alpha, K)$-sub-Weibull with parameter $\alpha \geq 1$. We further assume
    \begin{enumerate}
        \item The smallest eigenvalue of covariance matrix $\operatorname{Cov}(X_1)$ is bounded away from zero by a positive constant. The largest eigenvalue is bounded from above.
        \item The weighting parameter in Algorithm~\ref{algorithm: exp weighting} satisfies $\lambda=\lambda_n = o(\sqrt{n} / \log^{1/\alpha}p)$.
    \end{enumerate}
Then for any $x\in\mathbb{R}$:
    \begin{equation*}
       \lim _{n \rightarrow \infty} \max _{r \in [p]}\left|\mathbb{P}\left(\tilde{T}_r \leq x\right)-\Phi(x)\right|=0.
    \end{equation*}
\end{proposition}
\begin{proof}
The proof is identical to that of \Cref{th: weakly dependent CLT}, but specified to the case 
\begin{equation*}
K_{i}=\mathcal{K}_r(i, \boldsymbol{X}):=X_{i, r}-Q_{i, r}-\mathbb{E}\left[X_{i, r}-Q_{i, r} \mid \boldsymbol{X}^{\left(-v_i\right)}\right].
\end{equation*}
We previously relied on the assumption that $K_i$ is almost surely bounded to control its second and third moments in \Cref{eq: bound second term in Slepian} and \Cref{eq: third moment bound on K}. These bounds continue to hold under the weaker sub-Weibull assumptions. 

Under this formulation, and assuming the largest eigenvalue of ${\rm Cov}(X)$ is bounded from above, we can similarly show that $\tilde{\sigma}_i^2 := \mathbb{E}[K_i^2 \mid \boldsymbol{X}^{(-i)}]$ remains bounded almost surely. (One needs the summation of weights to be $1$ here.)

The part that needs the most significant modification is deriving the bound of $\mathbb{E}\left|K_i\right| \left(S^i\right)^2$ in \eqref{eq: where Ks2 appears}. We present the details in \Cref{lemma: bound KiSi2} and the technical lemmas it invokes. We bound this quantity by a combination of $\Delta_1$ and $\Delta_2$. Note that, in the bounded case $X_{i,s}$ case, we can show $\mathbb{E}\left|K_i\right|\left(S^i\right)^2$ is $o(1)$ when $\lambda = o(\sqrt{n})$. But $\mathbb{E}\left|K_i\right|\left(S^i\right)^2 = o(\sqrt{n})$ with $\lambda=o\left(\sqrt{n} \log ^{-1 / \alpha} p\right)$ is good enough to guarantee $\mathfrak{C}$ in the proof of \Cref{th: weakly dependent CLT} to be $o(1)$.

With the bounds on $\Delta_1$ and $\Delta_2$ in \Cref{lemma: moment bound K} and \Cref{lemma: moment bound second order K}, we can also derive $\|S^i\|_2 = o(1)$.
\end{proof}

\begin{lemma}\label{lemma: bound KiSi2}
Let $X_i \in \mathbb{R}^p$ be a collection of IID random vectors. Suppose each coordinate $X_{1,t}$ for $t \in [p]$ is $(\alpha, K)$-sub-Weibull with parameter $\alpha \geq 1$.

Define
\begin{equation*}
S^i := \sum_{j \neq i} \nabla_i K_j.
\end{equation*}

Assume that $\lambda$ satisfies
\begin{equation*}
\lambda < \frac{n}{72\left(16K \vee 16K \log (p)\right)^{1/\alpha}}.
\end{equation*}

Then
\begin{equation*}
\mathbb{E}\left|K_i\right| \left(S^i\right)^2 
\lesssim \lambda^2 n^{-1} \log^{2/\alpha} p 
+ \lambda^3 n^{-1} \log^{3/\alpha} p \left(1 + \lambda n^{-1} \log^{1/\alpha} p\right).
\end{equation*}

Specifically, when $\lambda = o\left(\sqrt{n} \log^{-1/\alpha} p\right)$,
\begin{equation*}
\mathbb{E}\left|K_i\right| \left(S^i\right)^2 = o(\sqrt{n}).
\end{equation*}
\end{lemma}

\begin{proof}
Expand the quantity of interest:
\begin{equation*}
\mathbb{E}\left|K_i\right|\left(S^i\right)^2 
= \mathbb{E}\left|K_i\right| \left( \sum_{j \neq i} \nabla_i K_j \right)^2
= \sum_{j \neq i} \mathbb{E}\left[\left|K_i\right| \left(\nabla_i K_j\right)^2\right] 
+ \sum_{\substack{j \neq i \\ k \neq i,\, k \neq j}} \mathbb{E}\left[\left|K_i\right| \nabla_i K_j \nabla_i K_k\right].
\end{equation*}

For the first term, we apply Hölder's inequality:
\begin{equation*}
\sum_{j \neq i} \mathbb{E}\left[\left|K_i\right| \left(\nabla_i K_j\right)^2\right]
\leq n \left\|K_i\right\|_3 \cdot \sup_j \left\|\nabla_i K_j\right\|_3^2.
\end{equation*}

Now consider the second term. 
\begin{equation*}
\begin{aligned}
\mathbb{E}\left[\left|K_i\right| \nabla_i K_j \nabla_i K_k\right] & \stackrel{(I)}{=} \mathbb{E}\left[\left(\nabla_i K_j\right) \nabla_j\left(\left|K_i\right| \nabla_i K_k\right)\right] \\
& =\mathbb{E}\left[\left(\nabla_j \nabla_i K_k\right)\left(\nabla_i K_j\right)\left|K_i\right|\right]+\mathbb{E}\left[\left(\nabla_i K_j\right)\left(\nabla_j\left|K_i\right|\right) \nabla_i K_k\left(\mathbf{X}^j\right)\right]
\end{aligned}
\end{equation*}
The first term above can be further simplified as
\begin{equation*}
\begin{aligned}
& \mathbb{E}\left[\left(\nabla_j \nabla_i K_k\right)\left(\nabla_i K_j\right)\left|K_i\right|\right]\\
\stackrel{(I I)}{=} & \mathbb{E}\left[\left(\nabla_j \nabla_i K_k\right) \nabla_k\left(\left|K_i\right| \nabla_i K_j\right)\right]\\
= & \mathbb{E}\left[\left(\nabla_j \nabla_i K_k\right)\left|K_i\right|\left(\nabla_k \nabla_i K_j\right)\right]+\mathbb{E}\left[\left(\nabla_j \nabla_i K_k\right)\left(\nabla_k\left|K_i\right|\right) \nabla_i K_j\left(\mathbf{X}^k\right)\right] .
\end{aligned}
\end{equation*}
In steps (I) and (II) we applied \Cref{lemma: nabla operation}.

We then apply Hölder’s inequality to each term:
\begin{align*}
\mathbb{E}\left[\left|K_i\right| \nabla_i K_j \nabla_i K_k\right] 
&\leq \sup_{i,j,k} \left\|\nabla_j \nabla_i K_k\right\|_3^2 \cdot \left\|K_i\right\|_3 \\
&\quad + \sup_k \left\|\nabla_k\left|K_i\right|\right\|_3 \cdot \sup_j \left\|\nabla_i K_j\right\|_3 \cdot \sup_{j,k} \left\|\nabla_j \nabla_i K_k\right\|_3 \\
&\quad + \sup_j \left\|\nabla_i K_j\right\|_3^2 \cdot \sup_j \left\|\nabla_j\left|K_i\right|\right\|_3.
\end{align*}

Combining and simplifying:
\begin{align*}
\mathbb{E}\left[\left|K_i\right| \nabla_i K_j \nabla_i K_k\right] 
&\leq \sup_{i,j,k} \left\|\nabla_j \nabla_i K_k\right\|_3^2 \cdot \left\|K_i\right\|_3 
+ \sup_{i,j} \left\|\nabla_i K_j\right\|_3^2 \cdot \sup_{i,j,k} \left\|\nabla_j \nabla_i K_k\right\|_3 \\
&\quad + \sup_{i,j} \left\|\nabla_i K_j\right\|_3^3.
\end{align*}

By applying \Cref{lemma: moment bound weighted sum}, \Cref{lemma: moment bound K}, and \Cref{lemma: moment bound second order K}, we obtain the overall bound:
\begin{equation*}
\mathbb{E}\left|K_i\right|\left(S^i\right)^2 
\lesssim C(\alpha, K) \left[ \lambda^2 n^{-1} \log^{2/\alpha} p 
+ \lambda^3 n^{-1} \log^{3/\alpha} p \left(1 + \lambda n^{-1} \log^{1/\alpha} p \right) \right].
\end{equation*}

Finally, when $\lambda = o\left(\sqrt{n} \log^{-1/\alpha} p\right)$, we have
\begin{equation*}
\mathbb{E}\left|K_i\right|\left(S^i\right)^2 = o(\sqrt{n}),
\end{equation*}
as desired.
\end{proof}

\begin{lemma}\label{lemma: moment bound weighted sum}
Let \( q \geq 1 \), and suppose each \( X_{j,s} \sim \textup{sub-Weibull}(\alpha, K) \). Then there exists a constant \( C(q, \alpha, K) \) such that:
\[
\left\| \sum_{s \neq r} \hat{w}_{r,s} |\tilde{X}_{j,s}| \right\|_q \leq \max_{s \in [p]} \left\| \tilde{X}_{j,s} \right\|_q.
\]
and
\[
\left\| K_{j,r} \right\|_q \leq 4\max_{s \in [p]} \left\| X_{j,s} \right\|_q,
\]

As a result, \( K_{j,r} \) is itself sub-Weibull as well.
\end{lemma}

\begin{proof}
By definition of $K_{j,r}$:
\[
K_{j,r} = \tilde{X}_{j,r} - \sum_{s \neq r} \hat{w}_{r,s} \tilde{X}_{j,s},
\]
where \( \tilde{X}_{j,s} := X_{j,s} - \mu_s \) is the centered version of the covariate.

Using the triangle inequality:
\[
\|K_{j,r}\|_q \leq \|\tilde{X}_{j,r}\|_q + \left\| \sum_{s \neq r} \hat{w}_{r,s} \tilde{X}_{j,s} \right\|_q.
\]

Now, applying Jensen’s inequality with respect to the weights (recall $\sum_{s\neq r} \hat w_{r,s} = 1$):
\begin{equation}\label{eq: weight jensen}
\mathbb{E}\left|\sum_{s \neq r} \hat{w}_{r, s} \tilde{X}_{j, s}\right|^q \leq \mathbb{E}\left(\sum_{s \neq r} \hat{w}_{r, s}\left|\tilde{X}_{j, s}\right|^q\right)=\mathbb{E}\left[\sum_{s \neq r} \hat{w}_{r, s}\left\|\tilde{X}_{j, s}\right\|_q^q\right] \leq \max _{s \neq r}\left\|\tilde{X}_{j, s}\right\|_q^q.
\end{equation}

A similar bound holds for $\left\|\sum_{s \neq r} \hat{w}_{r, s}\left|\tilde{X}_{j, s}\right|\right\|_q$ as
\begin{equation*}
    \mathbb{E}\left|\sum_{s \neq r} \hat{w}_{r, s} |\tilde{X}_{j, s}|\right|^q \leq \mathbb{E}\left(\sum_{s \neq r} \hat{w}_{r, s}\left|\tilde{X}_{j, s}\right|^q\right)
\end{equation*}
as well.

Taking \( q \)-th roots:
\[
\left\| \sum_{s \neq r} \hat{w}_{r,s} \tilde{X}_{j,s} \right\|_q 
\vee \left\| \sum_{s \neq r} \hat{w}_{r,s} |\tilde{X}_{j,s}| \right\|_q\leq \max_{s \neq r} \|\tilde{X}_{j,s}\|_q.
\]

Therefore,
\[
\|K_{j,r}\|_q 
\leq \|\tilde{X}_{j,r}\|_q + \max_{s \neq r} \|\tilde{X}_{j,s}\|_q 
\leq 2 \max_{s \in [p]} \|\tilde{X}_{j,s}\|_q 
\leq 4 \max_{s \in [p]} \|X_{j,s}\|_q,
\]
where we used that \( \|\tilde{X}_{j,s}\|_q \leq 2 \|X_{j,s}\|_q \) for sub-Weibull variables.

The moment bounds imply that \( K_{j,r} \) inherits sub-Weibull tail behavior.
\end{proof}

\begin{lemma}\label{lemma: moment bound K}
Let \( q \geq 1 \), and suppose \( X_{i,t} \sim \textup{sub-Weibull}(\alpha, K) \) with parameter \( \alpha \geq 1 \). Assume that \( \lambda \) satisfies:
\[
\lambda < \frac{n}{12q(4K \vee 2K \log p)^{1/\alpha}}.
\]
Then the stability term \( \nabla_i K_{j,r} \) satisfies the moment bound:
\[
\left\| \nabla_i K_{j,r} \right\|_q 
\leq C(q, \alpha, K) \cdot \lambda n^{-1} (\log p)^{1/\alpha},
\]
where \( C(q, \alpha, K) \) is a constant depending only on the sub-Weibull parameters and \( q \).
\end{lemma}
\begin{proof}
By definition,
\begin{equation}
    \nabla_i K_{j, r} 
    = Q_{j, r}^i - \mathbb{E}\left[Q_{j, r}^i \mid \boldsymbol{X}^{(-v), i}\right]
    - Q_{j, r} + \mathbb{E}\left[Q_{j, r} \mid \boldsymbol{X}^{(-v)}\right].
\end{equation}
This simplifies to:
\[
\nabla_i K_{j,r} = \sum_{s \neq r} \left(\hat{w}_{r, s}^i - \hat{w}_{r, s}\right) \left(X_{j, s} - \mu_s\right)
= \sum_{s \neq r} \hat{w}_{r, s} \left( \frac{\hat{w}_{r, s}^i}{\hat{w}_{r, s}} - 1 \right) (X_{j, s} - \mu_s).
\]

Using the bound from \eqref{eq: bound weight ratio}, we have:
\[
\frac{\hat{w}_{r, s}^i}{\hat{w}_{r, s}} \leq \exp\left(2 \lambda \tilde{n}^{-1} \max_{t \in [p]} |X_{i,t}' - X_{i,t}|\right).
\]
Applying a Taylor expansion yields:
\[
\frac{\hat{w}_{r, s}^i}{\hat{w}_{r, s}} - 1
\leq 2 \lambda \tilde{n}^{-1} \max_{t \in [p]} |X_{i,t}' - X_{i,t}| \cdot \exp\left(2 \lambda \tilde{n}^{-1} \max_{t \in [p]} |X_{i,t}' - X_{i,t}|\right).
\]
Define
\[
R := 2\lambda \tilde{n}^{-1}\max_{t \in [p]} |X_{i,t}' - X_{i,t}|,
\quad \text{so that} \quad 
\frac{\hat{w}_{r, s}^i}{\hat{w}_{r, s}} - 1 \leq R \exp(R).
\]

Therefore,
\[
|\nabla_i K_{j,r}| \leq R \exp(R) \sum_{s \neq r} \hat{w}_{r,s} |\tilde{X}_{j,s}|,
\]
where \( \tilde{X}_{j,s} := X_{j,s} - \mu_s \) denotes the centered covariate.

Applying Hölder’s inequality:
\[
\left\| \nabla_i K_{j,r} \right\|_q^q 
\leq \|R\|_{3q}^q \cdot \|\exp(R)\|_{3q}^q \cdot 
\left\| \sum_{s \neq r} \hat{w}_{r,s} |\tilde{X}_{j,s}| \right\|_{3q}^q.
\]

Finally, applying known moment bounds for each term, established in \Cref{lemma: maximum subW}, \Cref{lemma: finite mgf sub-W} and \Cref{lemma: moment bound weighted sum}, we conclude that:
\[
\left\| \nabla_i K_{j,r} \right\|_q 
\leq C(q, \alpha, K) \cdot \lambda n^{-1} (\log p)^{1/\alpha}.
\]
Note that when $X_{i,t}$ is sub-Weibull, $X_{i, t}^{\prime}-X_{i, t}$ is also sub-Weibull.
\end{proof}

\begin{lemma}\label{lemma: moment bound second order K}
Let \( q \geq 1 \), and suppose \( X_{i,t} \sim \textup{sub-Weibull}(\alpha, K) \) with parameter \( \alpha \geq 1 \). Assume that \( \lambda \) satisfies:
\[
\lambda < \frac{n}{24q(16K \vee 16K \log (p))^{1/\alpha}}.
\]
Then the second-order stability term \( \nabla_j \nabla_i K_{k,r}  \) satisfies the moment bound:
\[
\left\| \nabla_i \nabla_k K_{j,r} \right\|_q 
\leq C(q, \alpha, K) \cdot \lambda^2 n^{-2} (\log p)^{2/\alpha},
\]
where \( C(q, \alpha, K) \) is a constant depending only on the sub-Weibull parameters and \( q \).
\end{lemma}

\begin{proof}
Without loss of generality, we assume $r = 1$. Similar to Equation~\eqref{eq: simplify second order}, we can simplify the second-order stability of $K_{j,1}$ as
\begin{equation*}
\nabla_i \nabla_k K_{j, 1} = \sum_{s \neq 1} 
\left( 
\hat{w}_{1, s}^k - \hat{w}_{1, s} + \hat{w}_{1, s}^i - \hat{w}_{1, s}^{ik} 
\right)(X_{j, s} - \mu_s).
\end{equation*}

Rewriting the expression using $\tilde{X}_{j, s} := X_{j, s} - \mu_s$, we obtain
\begin{equation*}
\nabla_i \nabla_k K_{j, 1} 
= \sum_{s \neq 1} \left\{ 
\left( \hat{w}_{1, s} - \hat{w}_{1, s}^i \right) 
\left( \frac{\hat{w}_{1, s}^k}{\hat{w}_{1, s}} - 1 \right)
+ \hat{w}_{1, s}^i 
\left( 
\frac{\hat{w}_{1, s}^k}{\hat{w}_{1, s}} 
- \frac{\hat{w}_{1, s}^{ik}}{\hat{w}_{1, s}^i} 
\right)
\right\} \tilde{X}_{j, s}.
\end{equation*}

This can be further decomposed into two terms:
\begin{equation*}
\sum_{s \neq 1} 
\hat{w}_{1, s}^i 
\left( 
\frac{\hat{w}_{1, s}}{\hat{w}_{1, s}^i} - 1 
\right) 
\left( 
\frac{\hat{w}_{1, s}^k}{\hat{w}_{1, s}} - 1 
\right) 
\tilde{X}_{j, s}
+ \sum_{s \neq 1} 
\hat{w}_{1, s}^i 
\left( 
\frac{\hat{w}_{1, s}^k}{\hat{w}_{1, s}} 
- \frac{\hat{w}_{1, s}^{ik}}{\hat{w}_{1, s}^i} 
\right) 
\tilde{X}_{j, s}.
\end{equation*}

Therefore, the $\|\cdot\|_q$ norm of $|\nabla_i \nabla_k K_{j, 1}|$ can be bounded by the sum of the following two terms:
\begin{equation}\label{eq: second order split into two}
\begin{aligned}
& \left\|\sum_{s \neq 1} \hat{w}_{1, s}^i\left(\frac{\hat{w}_{1, s}}{\hat{w}_{1, s}^i}-1\right)\left(\frac{\hat{w}_{1, s}^k}{\hat{w}_{1, s}}-1\right) \tilde{X}_{j, s}\right\|_q, \\
& \left\|\sum_{s \neq 1} \hat{w}_{1, s}^i\left(\frac{\hat{w}_{1, s}^k}{\hat{w}_{1, s}}-\frac{\hat{w}_{1, s}^{i k}}{\hat{w}_{1, s}^i}\right) \tilde{X}_{j, s}\right\|_q .
\end{aligned}
\end{equation}

\textbf{Step 1.} We will use the same argument as in \Cref{lemma: moment bound K} to bound the first term. Define
\begin{equation*}
\mathcal{A} := \left| 
\sum_{s \neq 1} 
\hat{w}_{1, s}^i 
\left( 
\frac{\hat{w}_{1, s}}{\hat{w}_{1, s}^i} - 1 
\right) 
\left( 
\frac{\hat{w}_{1, s}^k}{\hat{w}_{1, s}} - 1 
\right) 
\tilde{X}_{j, s} 
\right|.
\end{equation*}

We then have
\begin{equation*}
\mathcal{A} \leq R_i R_k \exp(R_i) \exp(R_k) 
\left( \sum_{s \neq 1} \hat{w}_{1, s}^i |\tilde{X}_{j, s}| \right),
\end{equation*}
where
\begin{align*}
R_i &= 2 \lambda \tilde{n}^{-1} \max_{t \in [p]} \left| X_{i, t}' - X_{i, t} \right|, \\
R_k &= 2 \lambda \tilde{n}^{-1} \max_{t \in [p]} \left| X_{k, t}' - X_{k, t} \right|.
\end{align*}

Then we have 
\[
\|\mathcal{A}\|_q \leq C(q, \alpha, K) \lambda^2 n^{-2} (\log p)^{2/\alpha}
\]
by applying the five-term Hölder's inequality. The $\lambda n^{-1} (\log p)^{1/\alpha}$ terms come from $\|R_i\|_{5q}$ and $\|R_k\|_{5q}$.

\textbf{Step 2.} Now we bound the second term in \eqref{eq: second order split into two}:
\begin{equation*}
\mathcal{B} := \left| \sum_{s \neq 1} \hat{w}_{1, s}^i 
\left( 
\frac{\hat{w}_{1, s}^k}{\hat{w}_{1, s}} 
- 
\frac{\hat{w}_{1, s}^{ik}}{\hat{w}_{1, s}^i} 
\right) 
\tilde{X}_{j, s} 
\right|.
\end{equation*}

In the proof of \Cref{lemma: second order stability}, we established (using the same notation therein)
\begin{equation*}
\left| 
\frac{w_{r, s}^{ik}}{w_{r, s}^i} 
- 
\frac{w_{r, s}^k}{w_{r, s}} 
\right| 
\leq \exp(2 R_k) \cdot 
\left| 
\frac{\Xi^i \Xi^k}{\Xi^{ik} \Xi} - 1 
\right|.
\end{equation*}

Moreover,
\begin{equation*}
\left| 
\frac{\Xi^i \Xi^k}{\Xi^{ik} \Xi} - 1 
\right| 
\leq 
\sup_{2 \leq t \leq t' \leq p} 
\left| 
\frac{
E_{t, t'}^{i, k} + E_{t', t}^{i, k} 
- 
E_{t, t'}^{ik, \emptyset} 
- 
E_{t', t}^{ik, \emptyset}
}{
E_{t, t'}^{ik, \emptyset} + E_{t', t}^{ik, \emptyset}
} 
\right|
\end{equation*}
\begin{equation*}
= 
\sup_{2 \leq t \leq t' \leq p} 
\left| 
\frac{
\left\{ 
\exp(\xi_2) - \exp(\xi_1) 
\right\} 
\lambda 
\left( 
\hat{\mu}_{t'}^{ik} + \hat{\mu}_t 
- \hat{\mu}_{t'}^i - \hat{\mu}_t^k 
\right)
}{
E_{t, t'}^{ik, \emptyset} + E_{t', t}^{ik, \emptyset}
} 
\right|
\end{equation*}
\begin{equation*}
\leq 
\sup_{t, t'} |\xi_2 - \xi_1| \cdot 
\sup_{t, t'} 
\lambda \tilde{n}^{-1} 
\left| 
X_{k, t}' - X_{k, t} + X_{k, t'} - X_{k, t'}' 
\right|
\cdot 
\sup_{t, t'} 
\frac{\exp(\xi_3)}{E_{t, t'}^{ik, \emptyset} + E_{t', t}^{ik, \emptyset}}.
\end{equation*}

Therefore, $\|\mathcal{B}\|_q$ can be bounded by the product of the following terms:
\begin{equation*}
\begin{aligned}
\left\|\exp \left(R_k\right)\right\|_{6 q}^2 &\leq C \\
\left\|\sup _{t, t^{\prime}}\left(\xi_2-\xi_1\right)\right\|_{6 q} &\leq C(q, \alpha, K) \cdot \lambda n^{-1}(2 \log p)^{1 / \alpha} \\
\left\|\sup _{t, t^{\prime}} \lambda \tilde{n}^{-1}\left(X_{k, t}^{\prime}-X_{k, t}+X_{k, t^{\prime}}-X_{k, t^{\prime}}^{\prime}\right)\right\|_{6 q} & \leq C(q, \alpha, K) \cdot \lambda n^{-1}(2 \log p)^{1 / \alpha} \\
\left\|\sum_{s \neq 1} \hat{w}_{1, s}^i\left|\tilde{X}_{j, s}\right|\right\|_{6 q} & \leq C(\alpha, K)
\end{aligned}
\end{equation*}
and
\begin{equation}\label{eq: new to second order}
\left\|\sup _{2 \leq t \leq t^{\prime} \leq p} \frac{\exp \left(\xi_3\right)}{E_{t, t^{\prime}}^{i k, \emptyset}+E_{t^{\prime}, t}^{i k, \emptyset}}\right\|_{6 q} \leq C .
\end{equation}
We only need to elaborate on the last term; the rest follow directly from 
\Cref{lemma: moment bound weighted sum}, 
\Cref{lemma: finite mgf sub-W}, 
and 
\Cref{lemma: maximum subW}. 

Observe that
\begin{align*}
&\sup_{2 \leq t \leq t' \leq p} 
\frac{\exp(\xi_3)}{E_{t, t'}^{ik, \emptyset} + E_{t', t}^{ik, \emptyset}} \\
&\leq 
\sup_{t, t'} 
\frac{
\exp\left(
\max \left\{ 
-\lambda(\hat{\mu}_t^i + \hat{\mu}_{t'}^k), 
-\lambda(\hat{\mu}_t^{ik} + \hat{\mu}_{t'}), 
-\lambda(\hat{\mu}_{t'}^i + \hat{\mu}_t^k), 
-\lambda(\hat{\mu}_{t'}^{ik} + \hat{\mu}_t) 
\right\} 
\right)
}{
E_{t, t'}^{ik, \emptyset} + E_{t', t}^{ik, \emptyset}
} \\
&\leq 2 
+ 
\sup_{t, t'} 
\frac{
\exp\left( -\lambda(\hat{\mu}_t^i + \hat{\mu}_{t'}^k) \right) 
+ 
\exp\left( -\lambda(\hat{\mu}_{t'}^i + \hat{\mu}_t^k) \right)
}{
\exp\left( -\lambda(\hat{\mu}_t^{ik} + \hat{\mu}_{t'}) \right) 
+ 
\exp\left( -\lambda(\hat{\mu}_{t'}^{ik} + \hat{\mu}_t) \right)
} \\
&\leq 2 
+ 
\sup_{t, t'} 
\exp\left( 
-\lambda n^{-1} (X_{k, t} - X_{k, t}' + X_{k, t'}' - X_{k, t'}) 
\right) \\
&\quad 
+ 
\sup_{t, t'} 
\exp\left( 
-\lambda n^{-1} (X_{k, t'} - X_{k, t'}' + X_{k, t}' - X_{k, t}) 
\right).
\end{align*}

Then we apply \Cref{lemma: finite mgf sub-W} to conclude the bound in \eqref{eq: new to second order}, under the given condition on $\lambda$.
\end{proof}

\begin{lemma}\label{lemma: finite mgf sub-W}
Let \( R := 4\lambda n^{-1}\max_{s \in [p]} |X_{i,s}| \), and suppose \( X_{i,s} \sim \textup{sub-Weibull}(\alpha, K) \).  
Assume that \( \alpha \geq 1 \), and the parameters satisfy
\begin{equation*}
\lambda<\frac{n}{4 q(4 K \vee 2 K \log p)^{1 / \alpha}}.
\end{equation*}
Then \( \| \exp(R) \|_q^q < C \) for a universal constant $C$.
\end{lemma}
\begin{proof}
We use the tail integral representation for the moment generating function:
\[
\mathbb{E}[\exp(qR)] = \int_0^\infty \mathbb{P}(\exp(qR) \geq t) \, dt.
\]
Splitting the integral at \( t = e \), we obtain:
\[
\leq \exp(1) + \int_e^\infty \mathbb{P}(qR > \log t) \, dt.
\]
By the definition of \( R = \frac{4\lambda}{n} \max_{s \in [p]} |X_{i,s}| \), this becomes:
\[
= \exp(1) + \int_e^\infty \mathbb{P}\left( \max_{s \in [p]} |X_{i,s}| > \frac{n}{4\lambda q} \log t \right) dt.
\]
Applying the union bound and the sub-Weibull tail bound:
\[
\leq \exp(1) + \int_e^\infty p \cdot \mathbb{P}\left( X > \frac{n}{4\lambda q} \log t \right) dt
\leq \exp(1) + \int_e^\infty p \exp\left( -\frac{T^\alpha}{K} \right) dt,
\]
where \( T := \frac{n}{4\lambda q} \log t \).

By assumption, \( \lambda < \frac{n}{4q (2K \log p)^{1/\alpha}} \). Then,
\[
\frac{T^\alpha}{K} = \frac{n^\alpha}{(4\lambda q)^\alpha K} \log^\alpha t.
\]
Since \( \log^\alpha t \ge 1 \) for all \( t \ge e \), it suffices to show
\[
\frac{n^\alpha}{(4\lambda q)^\alpha K} > 2 \log p,
\]
which holds by the assumption. Thus,
\[
\frac{T^\alpha}{K} > 2 \log p \quad \text{for all } t \ge e.
\]
Thus,
\[
\mathbb{E}[\exp(qR)]
\leq \exp(1) + \int_e^\infty p \exp\left( -\frac{T^\alpha}{2K} \right) dt.
\]
Substituting \( T = \frac{n}{4\lambda q} \log t \), we have:
\[
\leq \exp(1) + \int_e^\infty p \exp\left( -\frac{n^\alpha \log^\alpha t}{(4\lambda q)^\alpha \cdot 2K} \right) dt.
\]
Let \( C := \frac{n^\alpha}{(4\lambda q)^\alpha \cdot 2K} \). When \( C > 2 \), i.e., 
\[
\lambda < \frac{n}{4q (4K)^{1/\alpha}},
\]
and \(\alpha \geq 1\), the integrand decays strictly faster than \( t^{-1} \), ensuring that the integral is finite. Hence,
\[
\mathbb{E}[\exp(qR)] < \infty \quad \Rightarrow \quad \| \exp(R) \|_q^q < \infty.
\]
\end{proof}

\begin{lemma}\label{lemma: maximum subW}
Let \( q \geq 1 \). Suppose each \( X_{i,t} \) is sub-Weibull\((\alpha, K)\). Then
\[
\left\| \max_{t \in [p]} |X_{i,t}| \right\|_q^q \leq C(q, \alpha, K)\,(\log p)^{q/\alpha}.
\]
Consequently, for \( R := 4\lambda n^{-1} \max_{t \in [p]} |X_{i,t}| \), we have
\[
\|R\|_q \leq C(q, \alpha, K)\, \lambda n^{-1} (\log p)^{1/\alpha}.
\]
\end{lemma}
\begin{proof}
We start by expressing the moment using a tail integral:
\[
\mathbb{E} \left[ \max_{t \in [p]} |X_{i,t}|^q \right]
= \int_0^\infty \mathbb{P} \left( \max_{t \in [p]} |X_{i,t}|^q > t \right) dt.
\]
Making the change of variable \( t = u^q \), we obtain:
\[
= q \int_0^\infty u^{q-1} \, \mathbb{P} \left( \max_{t \in [p]} |X_{i,t}| > u \right) du.
\]
Take \( T = (2K \log p)^{1/\alpha} \). Then,
\[
\mathbb{E} \left[ \max_{t \in [p]} |X_{i,t}|^q \right]
\leq q \int_0^T u^{q-1} \, du 
     + q \int_T^\infty u^{q-1} \, \mathbb{P} \left( \max_{t \in [p]} |X_{i,t}| > u \right) du.
\]
Using the sub-Weibull tail bound and the union bound,
\[
\leq T^q + 2q \int_T^\infty u^{q-1} \, p \exp\left( -\frac{u^\alpha}{K} \right) du.
\]
\[
= T^q + 2q \int_T^\infty u^{q-1} \exp\left( \log p - \frac{u^\alpha}{K} \right) du.
\]
Since \( \log p = T^\alpha / (2K) < u^\alpha / (2K) \) for all \( u > T \), we have:
\[
\mathbb{E} \left[ \max_{t \in [p]} |X_{i,t}|^q \right]
\leq T^q + 2q \int_T^\infty u^{q-1} \exp\left( -\frac{u^\alpha}{2K} \right) du.
\]
Change variables via \( v = u^\alpha/(2K) \), so that \( u = (2Kv)^{1/\alpha} \) and
\[
du = \frac{(2K)^{1/\alpha}}{\alpha} v^{1/\alpha - 1} dv.
\]
Then the integral becomes:
\[
= T^q + 2q \int_{T^\alpha / (2K)}^\infty 
(2K v)^{(q - 1)/\alpha} e^{-v} \cdot \frac{(2K)^{1/\alpha}}{\alpha} v^{1/\alpha - 1} dv.
\]
\[
= T^q + \frac{2q (2K)^{q/\alpha}}{\alpha} 
\int_{T^\alpha / (2K)}^\infty v^{q/\alpha - 1} e^{-v} dv.
\]
Using the incomplete Gamma function bound:
\[
\leq T^q + \frac{2q (2K)^{q/\alpha}}{\alpha} \Gamma\left( \frac{q}{\alpha} \right).
\]
Substituting back \( T = (2K \log p)^{1/\alpha} \), we obtain:
\[
\mathbb{E} \left[ \max_{t \in [p]} |X_{i,t}|^q \right]
\leq q (2K \log p)^{q/\alpha} + \frac{2q (2K)^{q/\alpha}}{\alpha} \Gamma\left( \frac{q}{\alpha} \right).
\]

\end{proof}

\section{Berry–Esseen Bound under Weak Dependence (adapted from \cite{shao2025berry})}
\label{sec: BE bound from Shao and Zhang} 

Theorem~\ref{th: weakly dependent CLT} establishes a Berry-Esseen bound for weakly dependent data, where the dependence is controlled by the notion of stabilities, with the proof completed using the Slepian interpolation. Other stability-based CLT results are available in the literature. For instance, \cite{chatterjee2008new} introduced, and \cite{shao2025berry} further developed, the \textit{generalized perturbative approach}. In this section, we aim to derive a qualitatively similar Berry-Esseen bound for our weak dependence CLT, compared to Theorem~\ref{th: weakly dependent CLT}, yet exploring the paradigm of the generalized perturbative approach. Here, we would not assume uniformly bounded $K_i$. Instead, a uniformly bounded fourth moment suffices.  

\begin{theorem}\label{th: weakly dependent CLT, shao and zhang}
Let $\bm{X} = \{X_i\in\mathcal{X}, i \in [n]\}$ be a collection of IID random vectors. We consider $K_i=\mathcal{K}(i ; \boldsymbol{X}) \in \mathbb{R}$ and assume $K_i$ has the same distribution for all $i\in[n]$. We further assume $\mathbb{E}[K_i\mid \bm{X}^{(-i)}] = 0$ where $\bm{X}^{(-i)}\coloneqq \bm{X}\symbol{92}\{X_i\}$, the variance $\omega_n^2 \coloneqq\operatorname{Var}\left(K_{1}\right)$ satisfies $\liminf_n \omega_n > 0$ and $\Enb K_1^4 < \infty$.

Then there exists a constant $C >0$ depending on $\Enb K_1^4$ such that
\begin{equation}\label{eq: explicit bound by stability}
\begin{aligned}
& \sup_{x \in \R} \left|\mathbb{P}\left(\frac{1}{\omega_n \sqrt{n}}\left(\sum_{i=1}^n K_i\right) \leq x\right)-\Phi(x)\right| \\
\leq & C_n \left (\frac{1}{\sqrt{n}} +  n^{\frac{1}{4}} \sqrt{\Delta_{1,(4)}} + \sqrt{n} \Delta_{1, (4)}  +  \sqrt{n} \Delta_{1, (2)} + n \Delta_{2,(2)} + n \Delta_{1, (2)}^2
    + n^2 \Delta_{2,(2)}^2 \right )
\end{aligned}
\end{equation}
with $C_n \lesssim 1/\psi_n^2 + 1/\min\{\omega_n, \psi_n\}$ and $|\omega_n - \psi_n| \le \sqrt{n} \cdot \max_{i \ne j} \norm{\nabla_i K_j}_2$. 
Here $\Phi(x)$ is the CDF of $\mathcal{N}(0,1)$ and for any $q \ge 1$,  
$$
    \Delta_{1, (q)} := \max_{j \ne i}\norm{\nabla_i \mathcal{K}(j, \bm{X})}_q  \text{\hspace{4mm}and\hspace{4mm}} \Delta_{2, (q)} = \max_{j \ne k \ne i} \norm{\nabla_k \nabla_i \mathcal{K}(j, \bm{X})}_q. 
$$    Specifically, if $\Delta_{1,(4)} = o(1/\sqrt{n})$ and $\Delta_{2,(2)} = o(1/n)$, we have $C_n = O(1)$ and 
\begin{equation*}
\frac{1}{\omega_n \sqrt{n}}\bigg(\sum_{i=1}^n K_i\bigg) \stackrel{d}{\longrightarrow} \mathcal{N}(0,1).
\end{equation*}
\end{theorem}

\begin{proof}
Define $W = \frac{1}{ \sqrt{n}}\sum_{i = 1}^n \mathcal{K}(i, \bm{X})$ and $\psi_n^2 = \mathrm{Var}(W) = \mathrm{Var}(K_1) + \frac{1}{n}\sum_{i \ne j} \mathbb{E}K_i K_j$. Theorem 2.1 in \cite{shao2025berry} gives the Berry-Esseen bound: 
\begin{equation}
\label{eq: BE bound shao/zhang 2025}
    d_{\mathrm{K}}(W) := \sup_{x \in \R} \abs{\mathbb{P}\left (\frac{W - \Enb W}{\sqrt{\mathrm{Var}(W)}} \le x \right ) - \Phi(x)} \le \frac{1}{\psi_n^2} \Enb \abs{\E{V\mid W} - \psi_n^2} + \frac{4}{\psi_n^2} \Enb \abs{\E{V^\ast \mid W}}, 
\end{equation}
where $\Phi$ stands for the cumulative distribution function of the standard normal, and the random variables $V$ and $V^\ast$ are defined by 
$$
    V = \frac{1}{2} \sum_{i = 1}^n (\nabla_i W) (\nabla_i W^{A_i}) \text{\hspace{3mm} and \hspace{3mm}} V^\ast = \frac{1}{2} \sum_{i = 1}^n (\nabla_i W) \abs{\nabla_i W^{A_i}}
$$
with the set $A_i := \{1, \ldots, i - 1\}$. For any subset $A \subseteq [n]$, $W^A$ denotes the statistic obtained from $W$ by replacing $\{X_j : j \in A\}$ with IID copies $\{X_j' : j \in A\}$. We proceed to bound $\Enb \abs{\E{V \mid W} - \psi_n^2}$ and $\Enb\abs{\E{V^\ast \mid W}}$.  

For any $i \in [n]$, 
$$
    \nabla_i W = \nabla_i \left ( \frac{1}{\sqrt{n}} \sum_{j = 1}^n \mathcal{K}(j, \bm{X}) \right ) = \frac{1}{\sqrt{n}} \nabla_i \mathcal{K}(i, \bm{X}) + \frac{1}{\sqrt{n}} \sum_{j \ne i} \nabla_i \mathcal{K}(j, \bm{X})
$$
so that 
\begin{equation}
\label{eq: W product expansion}
\begin{aligned}
    &(\nabla_i W) (\nabla_i W^{A_i}) \\
    & = \frac{1}{n} \nabla_i \mathcal{K}(i, \bm{X}) \nabla_i\mathcal{K}(i, \bm{X}^{A_i}) + \frac{1}{n} \nabla_i \mathcal{K}(i, \bm{X}) \sum_{j\ne i} \nabla_i \mathcal{K}(j, \bm{X}^{A_i}) \\
    &\hspace{5mm} + \frac{1}{n}\nabla_i \mathcal{K}(i, \bm{X}^{A_i}) \sum_{j\ne i} \nabla_i \mathcal{K}(j, \bm{X}) + \frac{1}{n} \left(\sum_{j \ne i} \nabla_i \mathcal{K}(j, \bm{X}) \right) \left(\sum_{j \ne i} \nabla_i \mathcal{K}(j, \bm{X}^{A_i}) \right). 
\end{aligned}
\end{equation}
Lemma 4.1 in~\cite{shao2025berry} implies 
$$
    \mathrm{Var}(W) = \psi_n^2 = \E{V} = \frac{1}{2} \sum_{i = 1}^n \E{(\nabla_iW) (\nabla_iW^{A_i})}, 
$$
so we can obtain from the expansion~\eqref{eq: W product expansion} that  
\begin{equation}
\label{eq: bound for conditional expectation regarding V}
\begin{aligned}
    & \Enb\abs{\E{V \mid W} - \psi_n^2} \\
    &= \Enb \hspace{0.5mm} \Big | \hspace{0.5mm} \frac{1}{2} \sum_{i = 1}^n  \E{(\nabla_i W)(\nabla_i W^{A_i}) \mid W} - \frac{1}{2} \sum_{i = 1}^n  \E{(\nabla_i W) (\nabla_i W^{A_i})} \Big | \\
    & \le \Enb \abs{ \frac{1}{2n} \sum_{i = 1}^n \E{\nabla_i \mathcal{K}(i, \bm{X}) \nabla_i \mathcal{K}(i, \bm{X}^{A_i})\mid W} - \frac{1}{2n} \sum_{i = 1}^n \E{\nabla_i \mathcal{K}(i, \bm{X}) \nabla_i \mathcal{K}(i, \bm{X}^{A_i})}} \\
    & \hspace{5mm} + \frac{1}{n} \sum_{i = 1}^n \Enb \abs{\nabla_i \mathcal{K}(i, \bm{X}) \sum_{j \ne i} \nabla_i \mathcal{K}(j, \bm{X}^{A_i})} + \frac{1}{n} \sum_{i = 1}^n \Enb \abs{\nabla_i \mathcal{K}(i, \bm{X}^{A_i}) \sum_{j \ne i} \nabla_i \mathcal{K}(j, \bm{X})} \\
    & \hspace{5mm} + \frac{1}{n} \sum_{i = 1}^n \Enb \abs{\left(\sum_{j \ne i}\nabla_i \mathcal{K}(j, \bm{X}) \right) \left(\sum_{j \ne i}\nabla_i \mathcal{K}(j, \bm{X}^{A_i}) \right)} \\
    &:= \mathcal{T}_1 + \mathcal{T}_2 + \mathcal{T}_3 + \mathcal{T}_4. 
\end{aligned}
\end{equation}
We now bound each term respectively. 
Begin with term $\mathcal{T}_1$. Define 
$$
\mathscr{K} := \frac{1}{2n} \sum_{i = 1}^n \nabla_i \mathcal{K}(i, \bm{X}) \nabla_i \mathcal{K}(i, \bm{X}^{A_i}).
$$
Observe that 
$\mathcal{T}_1 \le \sqrt{\mathrm{Var}\left (\E{\mathscr{K} \mid W} \right )} \le \sqrt{\mathrm{Var}(\mathscr{K})}$, 
where the last variance is equal to 
\begin{align*}
    &\frac{1}{4n^2} \sum_{i = 1}^n \Enb \left (
    \nabla_i \mathcal{K}(i, \bm{X}) \nabla_i \mathcal{K}(i, \bm{X}^{A_i}) - \E{\nabla_i\mathcal{K}(i, \bm{X} ) \nabla_i \mathcal{K}(i, \bm{X}^{A_i})}\right )^2 \\
    & + \frac{1}{4n^2} \sum_{i \ne j} \Enb \left \{ \nabla_i \mathcal{K}(i, \bm{X}) \nabla_i \mathcal{K}(i, \bm{X}^{A_i}) - \E{\nabla_i \mathcal{K}(i, \bm{X}) \nabla_i \mathcal{K}(i, \bm{X}^{A_i})} \right\} \\
    & \hspace{2.5cm} \times \left\{ \nabla_j \mathcal{K}(j, \bm{X}) \nabla_j \mathcal{K}(j, \bm{X}^{A_j}) - \E{\nabla_j \mathcal{K}(j, \bm{X}) \nabla_j \mathcal{K}(j, \bm{X}^{A_j})} \right\} \\
    & := S_1 + S_2. 
\end{align*}
There exists $C_1 > 0$ only depending on the fourth moment such that
\begin{equation}
\label{eq: vanishing conditional variance 1}
    S_1 \le \frac{C_1}{n} = o(1). 
\end{equation}
The second sum $S_2$ requires some careful treatments. Define $\bm{X}_{A}' = \{X_k', k \in A\}$ for any $A \subseteq [n]$ and let $X_i{''}, X_j{''}$ be two extra IID samples (independent of $X_i'$ and $X_{j}'$ as well). The sum $S_2$ is equal to 
$$
\frac{1}{4n^2} \sum_{i \ne j} \left\{  \E{U_i U_j} + \E{U_i V_j} + \E{V_i U_j} + \E{V_i V_j} \right \}, 
$$
where for any $i \in [n]$, we define
\begin{align*}
    & U_i = \nabla_i \mathcal{K}(i, \bm{X}) \nabla_i \mathcal{K}(i, \bm{X}^{A_i}) - \E{\nabla_i \mathcal{K}(i, \bm{X}) \nabla_i \mathcal{K}(i, \bm{X}^{A_i}) \mid \bm{X}^{(-i)}, \bm{X}_{A_i}'}; \\
    & V_i = \E{\nabla_i \mathcal{K}(i, \bm{X}) \nabla_i \mathcal{K}(i, \bm{X}^{A_i}) \mid \bm{X}^{(-i)}, \bm{X}_{A_i}'} - \E{\nabla_i \mathcal{K}(i, \bm{X}) \nabla_i \mathcal{K}(i, \bm{X}^{A_i})}. 
\end{align*}
In the sequel, additional samples $X_i^{\prime\prime}$ and $X_j^{\prime\prime}$ are used for $X_i$'s and $X_j$'s perturbations whenever $(i),(j)$ appear as superscripts on variables and as subscripts on the $\nabla$ operator. Because
\begin{align*}
    &\E{U_j \mid \bm{X}^{(-j)}, \bm{X}_{A_i \cup A_j}', X_i', X_j{''}} = \E{U_j \mid \bm{X}^{(-j)}, \bm{X}_{A_j}'} = 0; \\ 
   &\E{U_j^{(i)} \mid \bm{X}^{(-j)}, \bm{X}_{A_i\cup A_j}', X_i', X_i{''}, X_j{''}} = \E{U_j^{(i)} \mid \bm{X}^{(-j),(-i)}, \bm{X}_{A_j}', X_i{''}} = 0
\end{align*}
for any $i \ne j$, we have 
\begin{itemize}
    \item $\E{U_i^{(j)} U_j} = \E{U_i^{(j)} \cdot \E{U_j \mid \bm{X}^{(-j)}, \bm{X}'_{A_i \cup A_j}, X_i', X_j{''}}} = 0$;
    \item  $\E{U_i U_j^{(i)}} = 0$ by symmetry; 
    \item $\E{U_i^{(j)}U_j^{(i)}} = \E{U_i^{(j)} \cdot \E{U_j^{(i)} \mid \bm{X}^{(-j)}, \bm{X}_{A_i \cup A_j}', X_i', X_i{''}, X_j{''}}} = 0.$
\end{itemize}
Accordingly, $\E{U_i U_j} = \E{(\nabla_{(j)} U_i) (\nabla_{ (i)} U_j)} \le \norm{\nabla_{(i)} U_j}_2^2 = \norm{\nabla_i U_j}_2^2$ by Lemma~\ref{lemma: general free nabla} and the Cauchy-Schwarz inequality. Next, 
\begin{align*}
    \norm{\nabla_i U_j}_2 &\le 2 \norm{\nabla_i \left \{\nabla_j \mathcal{K}(j, \bm{X}) \nabla_j\mathcal{K}(j, \bm{X}^{A_j})\right\}}_2 \\
    & = 2 \norm{\left \{ \nabla_i \nabla_j \mathcal{K}(j, \bm{X}) \right\} \nabla_j \mathcal{K}(j, \bm{X}^{A_j}) + \nabla_j \mathcal{K}(j, \bm{X}^i) \left\{\nabla_i \nabla_j \mathcal{K}(j, \bm{X}^{A_j}) \right\} }_2 \\
    & \le 4 \cdot \norm{\nabla_i \mathcal{K}(j, \bm{X})}_4 \cdot  \norm{\nabla_j \mathcal{K}(j, \bm{X}^{A_j})}_4 + 4 \cdot \norm{\nabla_j \mathcal{K}(j, \bm{X}^i)}_4 \cdot  \norm{\nabla_i \mathcal{K}(j, \bm{X}^{A_j})}_4 \\
    &\le  8 \cdot \norm{\nabla_i \mathcal{K}(j, \bm{X})}_4 \cdot  \norm{\nabla_j \mathcal{K}(j, \bm{X}^{A_j})}_4 \\
    & \le C_2 \norm{\nabla_i \mathcal{K}(j, \bm{X})}_4
\end{align*}
for some $C_2 > 0$ only depending on $\Enb\mathcal{K}^4(1, \bm{X})$. Therefore, 
\begin{equation}
\label{eq: vanishing conditional variance 2}
    \frac{1}{4n^2} \sum_{i \ne j} \E{U_i U_j} \le \frac{C_2^2}{4} \norm{\nabla_i \mathcal{K}(j, \bm{X})}_4^2. 
\end{equation}
Moreover, we know $\norm{U_i}_2 = O(1)$, where the bounding constant $C_3$ can be specified in terms of the fourth moment $\Enb \mathcal{K}^4(1, \bm{X})$ up to a universal constant. \\
Last, we analyze $\norm{V_i}_2^2$ for any fixed $i$. To ease our notation, consider the relabeling $X_{n + \ell} = X_\ell'$ for any $\ell \in \{1, \ldots, i - 1\}$ in the following. 
Note that
\begin{align*}
    V_i & = \sum_{k \ne i} \left \{ \E{\nabla_i \mathcal{K}(i, \bm{X}) \nabla_i \mathcal{K}(i, \bm{X}^{A_i}) \mid \bm{X}_{1:k}^{(-i)}} - \E{\nabla_i \mathcal{K}(i, \bm{X}) \nabla_i \mathcal{K}(i, \bm{X}^{A_i}) \mid \bm{X}_{1:(k-1)}^{(-i)}} \right\} \\
    &= \sum_{k \ne i} \E{\nabla_{k} \left\{\nabla_i \mathcal{K}(i, \bm{X}) \nabla_i \mathcal{K}(i, \bm{X}^{A_i}) \right \} \mid \bm{X}_{1:k}^{(-i)}}, 
\end{align*}
where $\bm{X}_{1:k}^{(-i)} = \{X_1, \ldots, X_k\} \!\setminus\!\{X_i\}$ for any $k \in \{1, \ldots, n + i - 1\} \!\setminus\!\{i\}$, and $\bm{X}_{1:0}^{(-i)}$ is the trivial $\sigma$-field such that $\Enb[A \mid \bm{X}_{1:0}^{(-i)}] = \Enb[A]$ for any random variable $A$. The sequence 
$$
    \left \{\E{\nabla_k \left\{\nabla_i \mathcal{K}(i, \bm{X}) \nabla_i \mathcal{K}(i, \bm{X}^{A_i}) \right \} \mid \bm{X}_{1:k}^{(-i)}} \right \}_{k \ne i}
$$
forms a martingale increment sequence so that
\begin{align*}
     \norm{V_i}_2^2 & = \sum_{k \ne i} \norm{\E{\nabla_k \left\{\nabla_i \mathcal{K}(i, \bm{X}) \nabla_i \mathcal{K}(i, \bm{X}^{A_i}) \right \} \mid \bm{X}_{1:k}^{(-i)}}}_2^2 \\
    & \le \sum_{k \ne i} \norm{\nabla_k \left \{ \nabla_i \mathcal{K}(i, \bm{X}) \nabla_i \mathcal{K}(i, \bm{X}^{A_i})\right\}}_2^2 \\
    & = \sum_{k \ne i} \norm{\left\{\nabla_k\nabla_i \mathcal{K}(i, \bm{X})\right\} \nabla_i \mathcal{K}(i, \bm{X}^{A_i}) + \nabla_i \mathcal{K}(i, \bm{X}^k) \left\{\nabla_k \nabla_i \mathcal{K}(i, \bm{X}^{A_i})\right\}}_2^2 \\
    & \le \sum_{k \ne i} \left( \norm{\nabla_k \nabla_i \mathcal{K}(i, \bm{X})}_4 \cdot \norm{\nabla_i \mathcal{K}(i, \bm{X}^{A_i})}_4 + \norm{\nabla_i \mathcal{K}(i, \bm{X}^k)}_4 \cdot \norm{\nabla_k\nabla_i \mathcal{K}(i, \bm{X}^{A_i})}_4 \right)^2 \\
    & \le 16 \sum_{k \ne i}
    \norm{\nabla_k \mathcal{K}(i, \bm{X})}_4^2 \cdot \norm{\nabla_i \mathcal{K}(i, \bm{X})}_4^2 \\
    & \le 16 C_4 \sum_{k \ne i}
    \norm{\nabla_k \mathcal{K}(i, \bm{X})}_4^2 \le 32 C_4 \cdot n \cdot \max_{k \ne i} \norm{\nabla_k \mathcal{K}(i, \bm{X})}_4^2 
\end{align*}
for some $C_4 > 0$ only depending on the fourth moment $\Enb \mathcal{K}^4(1, \bm{X})$. Therefore, 
\begin{equation}
\label{eq: vanishing conditional variance 3}
\begin{aligned}
    & \frac{1}{4n^2} \sum_{i \ne j} \left \{ \E{U_i V_j} + \E{V_i U_j} + \E{V_i V_j} \right\} \\
    & \le \frac{1}{4} \left \{ 8 \sqrt{2} C_3 \sqrt{C_4} \sqrt{n} \max_{i \ne j} \norm{\nabla_j \mathcal{K}(i, \bm{X})}_4 + 32 C_4 \cdot n \cdot \max_{i \ne j} \norm{\nabla_j \mathcal{K}(i, \bm{X})}_4^2 \right\}. 
\end{aligned}
\end{equation}

As for the rest, i.e., terms $\mathcal{T}_2, \mathcal{T}_3, \mathcal{T}_4$, we note that for any fixed $i \in [n]$,
\begin{equation}
\label{eq: residual terms}
\begin{aligned}
    & \Enb \abs{\nabla_i \mathcal{K}(i, \bm{X}) \sum_{j \ne i} \nabla_i \mathcal{K}(j, \bm{X}^{A_i})} \le \norm{\nabla_i \mathcal{K}(i, \bm{X})}_2 \cdot \norm{\sum_{j \ne i} \nabla_i \mathcal{K}(j, \bm{X}^{A_i})}_2; \\
    & \Enb \abs{\nabla_i \mathcal{K}(i, \bm{X}^{A_i}) \sum_{j \ne i} \nabla_i \mathcal{K}(j, \bm{X})} \le \norm{\nabla_i \mathcal{K}(i, \bm{X}^{A_i})}_2 \cdot \norm{\sum_{j \ne i} \nabla_i \mathcal{K}(j, \bm{X})}_2; \\
    & \Enb \abs{\left(\sum_{j \ne i}\nabla_i \mathcal{K}(j, \bm{X}) \right) \left(\sum_{j \ne i}\nabla_i \mathcal{K}(j, \bm{X}^{A_i}) \right)} \le \norm{\sum_{j \ne i} \nabla_i \mathcal{K}(j, \bm{X})}_2 \cdot \norm{\sum_{j \ne i} \nabla_i \mathcal{K}(j, \bm{X}^{A_i})}_2
\end{aligned}
\end{equation}
by the Cauchy-Schwarz inequality. The quantity $\norm{\nabla_i \mathcal{K}(i, \bm{X})}_2$ is of order $O(\omega_n) = O(\norm{\mathcal{K}(i,\bm{X})}_4)= O(1)$, so it reduces to prove $\norm{\sum_{j \ne i} \nabla_i \mathcal{K}(j, \bm{X})}_2 = o(1)$ to have $\mathcal{T}_2 + \mathcal{T}_3 + \mathcal{T}_4$ vanish. Indeed, we have from Lemma~\ref{lemma: Di} that 
\begin{equation}
\label{eq: free nabla trick}
\begin{aligned}
    \norm{\sum_{j \ne i} \nabla_i \mathcal{K}(j, \bm{X})}_2^2
    &= \sum_{j \ne i} \Enb(\nabla_i \mathcal{K}(j, \bm{X}))^2 + \sum_{j \ne k \ne i} \E{\nabla_i \mathcal{K}(j, \bm{X})\nabla_i \mathcal{K}(k, \bm{X})} \\
    &= \sum_{j \ne i} \Enb(\nabla_i \mathcal{K}(j, \bm{X}))^2 + \sum_{j \ne k \ne i} \E{\nabla_k \nabla_i \mathcal{K}(j, \bm{X}) \nabla_j \nabla_i \mathcal{K}(k, \bm{X})} \\
    & \le n \cdot \max_{j \ne i} \norm{\nabla_i \mathcal{K}(j, \bm{X})}_2^2 + n^2 \cdot \max_{j \ne k \ne i} \norm{\nabla_k \nabla_i \mathcal{K}(j, \bm{X})}_2^2. 
\end{aligned}
\end{equation}
Bounding the term \(\Enb \abs{\E{V^\ast \mid W}}\) in~\eqref{eq: BE bound shao/zhang 2025} can be carried out using similar arguments. Expanding $\nabla_i W$ and applying the triangle inequality yield
\begin{equation*}
\begin{aligned}
    & \Enb \abs{\E{V^\ast \mid W}} \le \Enb \abs{V^\ast} \\
    & \le \Enb \abs{\frac{1}{2n} \sum_{i = 1}^n \nabla_i \mathcal{K}(i, \bm{X}) \abs{\nabla_i W^{A_i}}} + \frac{1}{2n} \sum_{i = 1}^n \Enb \abs{\nabla_i \mathcal{K}(i, \bm{X}^{A_i}) \sum_{j \ne i} \nabla_i \mathcal{K}(j, \bm{X})} \\
    & \hspace{5mm} + \frac{1}{2n} \sum_{i = 1}^n \Enb \abs{\left(\sum_{j \ne i}\nabla_i \mathcal{K}(j, \bm{X}) \right) \left(\sum_{j \ne i}\nabla_i \mathcal{K}(j, \bm{X}^{A_i}) \right)}. 
\end{aligned}
\end{equation*}
By Lemma~\ref{lemma: absolute value expansion}, we know 
$$
    \abs{\nabla_i W^{A_i}} = \abs{\nabla_i \mathcal{K}(i, \bm{X}^{A_i})} + \operatorname{sgn}(\nabla_i \mathcal{K}(i, \bm{X}^{A_i})) \sum_{j \ne i} \nabla_i \mathcal{K}(j, \bm{X}^{A_i}) + \mathfrak{R}
$$
with the residual $\mathfrak{R}$ satisfying $0 \le \mathfrak{R} \le 2 \abs{\sum_{j \ne i} \nabla_i \mathcal{K}(j, \bm{X}^{A_i})}$. 
Therefore,
\begin{equation*}
\begin{aligned}
    &\Enb \abs{\frac{1}{2n} \sum_{i = 1}^n \nabla_i \mathcal{K}(i, \bm{X}) \abs{\nabla_i W^{A_i}}} \\
    & \le \Enb \abs{\frac{1}{2n} \sum_{i = 1}^n \nabla_i \mathcal{K}(i, \bm{X}) \abs{\nabla_i \mathcal{K}(i, \bm{X}^{A_i})}} + \frac{3}{2n} \sum_{i = 1}^n \Enb \abs{\nabla_i \mathcal{K}(i, \bm{X}) \sum_{j \ne i} \nabla_i \mathcal{K}(j, \bm{X}^{A_i})}
\end{aligned}
\end{equation*}
which in turn gives 
\begin{equation}
\label{eq: bound for conditional expectation regarding V*}
    \Enb\abs{\E{V^\ast \mid W}} \lesssim \mathcal{T}_5 + \mathcal{T}_2 + \mathcal{T}_3 + \mathcal{T}_4, \hspace{3mm} \mathcal{T}_5 := \Enb \abs{\frac{1}{2n} \sum_{i = 1}^n \nabla_i \mathcal{K}(i, \bm{X}) \abs{\nabla_i \mathcal{K}(i, \bm{X}^{A_i})}}.  
\end{equation}
It remains to bound $\mathcal{T}_5$. Define $\mathscr{M} = \frac{1}{2n} \sum_{i = 1}^n M_i, \hspace{1mm} M_i := \nabla_i\mathcal{K}(i, \bm{X}) \abs{\nabla_i \mathcal{K}(i, \bm{X}^{A_i})}$. We know $\E{M_i} = 0$ and $\mathcal{T}_5 \le \sqrt{\Enb\mathscr{M}^2}$,
where the inner expectation $\Enb \mathscr{M}^2$ is equal to 
\begin{align*}
    &\frac{1}{4n^2} \sum_{i = 1}^n \Enb \left( \nabla_i \mathcal{K}(i, \bm{X}) \abs{\nabla_i \mathcal{K}(i, \bm{X}^{A_i})} \right)^2 \\
    & + \frac{1}{4n^2} \sum_{i \ne j} \E{\nabla_i \mathcal{K}(i, \bm{X}) \abs{\nabla_i \mathcal{K}(i, \bm{X}^{A_i})} \nabla_j \mathcal{K}(j, \bm{X}) \abs{\nabla_j \mathcal{K}(j, \bm{X}^{A_j})}} \\
    & := \tilde{S}_1 + \tilde{S}_2. 
\end{align*}
There exists $C_5 > 0$ only depending on the fourth moment such that
\begin{equation}
\label{eq: vanishing conditional variance 4}
    \tilde{S}_1 \le \frac{C_5}{n}. 
\end{equation}
We apply the previous free-nabla and martingale increment techniques for $\tilde{S}_2$. A similar decomposition gives
\begin{align*}
    \tilde{S}_2 = \frac{1}{4n^2} \sum_{i \ne j} \left \{\E{\tilde{U}_i \tilde{U}_j} + \E{\tilde{U}_i \tilde{V}_j} + \E{\tilde{V}_i \tilde{U}_j} + \E{\tilde{V}_i \tilde{V}_j} \right \}, 
\end{align*}
where for any $i \in [n]$,
\begin{align*}
    \tilde{U}_i = M_i - \E{M_i \mid \bm{X}^{(-i)}, \bm{X}_{A_i}'} \hspace{3mm} \text{and}\hspace{3mm} \tilde{V}_i = \E{M_i \mid \bm{X}^{(-i)}, \bm{X}_{A_i}'}- \E{M_i} =  \E{M_i \mid \bm{X}^{(-i)}, \bm{X}_{A_i}'}. 
\end{align*}
Lemma C.4 and the Cauchy-Schwarz inequality again yields $\E{\tilde{U}_i \tilde{U}_j} = \E{(\nabla_{(j)} \tilde{U}_i) (\nabla_{(i)} \tilde{U}_j)} \le \norm{\nabla_{(i)} \tilde{U}_j}_2^2 = \norm{\nabla_i \tilde{U}_j}_2^2. 
$
In particular, 
\begin{align*}
    \norm{\nabla_i \tilde{U}_j}_2 &\le\norm{\nabla_i \left \{\nabla_j \mathcal{K}(j, \bm{X}) \abs{\nabla_j\mathcal{K}(j, \bm{X}^{A_j})}\right\}}_2 \\
    & = \norm{\left \{ \nabla_i \nabla_j \mathcal{K}(j, \bm{X}) \right\} \abs{\nabla_j \mathcal{K}(j, \bm{X}^{A_j})} + \nabla_j \mathcal{K}(j, \bm{X}^i) \left\{\nabla_i \abs{\nabla_j \mathcal{K}(j, \bm{X}^{A_j})} \right\} }_2 \\
    & \le 2 \cdot \norm{\nabla_i \mathcal{K}(j, \bm{X})}_4 \cdot  \norm{\nabla_j \mathcal{K}(j, \bm{X}^{A_j})}_4 + 2 \cdot \norm{\nabla_j \mathcal{K}(j, \bm{X}^i)}_4 \cdot  \norm{\nabla_i \abs{\mathcal{K}(j, \bm{X}^{A_j})}}_4 \\
    &\le  4 \cdot \norm{\nabla_i \mathcal{K}(j, \bm{X})}_4 \cdot  \norm{\nabla_j \mathcal{K}(j, \bm{X}^{A_j})}_4 \\
    & \le C_6 \norm{\nabla_i \mathcal{K}(j, \bm{X})}_4
\end{align*}
for some $C_6 > 0$ only depending on $\Enb\mathcal{K}^4(1, \bm{X})$. Therefore, 
\begin{equation}
\label{eq: vanishing conditional variance 5}
    \frac{1}{4n^2} \sum_{i \ne j} \E{\tilde{U}_i \tilde{U}_j} \le \frac{C_6^2}{4} \norm{\nabla_i \mathcal{K}(j, \bm{X})}_4^2. 
\end{equation}
Plus, we know $\norm{\tilde{U}_i}_2 = O(1)$, where the bounding constant $C_7$ can be specified in terms of the fourth moment $\Enb \mathcal{K}^4(1, \bm{X})$ up to a universal constant. Last, 
we have for any fixed $i$
\begin{align*}
    \tilde{V}_i &= \sum_{k \ne i} \left \{ \E{M_i \mid \bm{X}_{1:k}^{(-i)}} - \E{M_i \mid \bm{X}_{1:(k-1)}^{(-i)}} \right\} = \sum_{k \ne i} \E{\nabla_k M_i \mid \bm{X}_{1:k}^{(-i)}}
\end{align*}
with $\E{M_i \mid \emptyset} = \E{M_i} = 0$. The sequence $\{\E{\nabla_k M_i} \mid \bm{X}_{1:k}^{(-i)}\}_{k \ne i}$ forms a martingale increment sequence so that 
\begin{equation*}
\begin{aligned}
 \norm{\tilde{V}_i}_2^2 &= \sum_{k \ne i} \norm{\E{\nabla_k M_i \mid \bm{X}_{1:k}^{(-i)}}}_2^2 \le \sum_{k \ne i} \norm{\nabla_k \left \{\nabla_i \mathcal{K}(i, \bm{X}) \abs{\nabla_i \mathcal{K}(i, \bm{X}^{A_i})} \right\}}_2^2 \\
 & \le \sum_{k \ne i} \left( \norm{\nabla_k \nabla_i \mathcal{K}(i, \bm{X}) \abs{\nabla_i \mathcal{K}(i, \bm{X}^{A_i})}}_2 + \norm{\nabla_i\mathcal{K}(i, \bm{X}^k) \nabla_k \abs{\nabla_i \mathcal{K}(i, \bm{X}^{A_i})}}_2 \right)^2 \\
 & \lesssim \sum_{k \ne i} \norm{\nabla_k \mathcal{K}(i, \bm{X})}_4^2 \cdot \norm{\nabla_i \mathcal{K}(i, \bm{X})}_4^2 \le C_8 \cdot n \cdot \max_{k \ne i} \norm{\nabla_k \mathcal{K}(i, \bm{X})}_4^2.  
\end{aligned}
\end{equation*}
for some $C_8 > 0$ only depending on $\Enb \mathcal{K}^4(i, \bm{X})$. 
Therefore, 
\begin{equation}
\label{eq: vanishing conditional variance 6}
\begin{aligned}
    &\frac{1}{4n^2} \sum_{i \ne j} \left \{ \E{\tilde{U}_i \tilde{V}_j} + \E{\tilde{V}_i \tilde{U}_j} + \E{\tilde{V}_i \tilde{V}_j} \right \} \\
    & \lesssim C_7 \sqrt{C_8} \cdot \sqrt{n} \cdot  \max_{i \ne j} \norm{\nabla_j \mathcal{K}(i, \bm{X})}_4 + C_8 \cdot n \cdot \max_{i \ne j} \norm{\nabla_j \mathcal{K}(i, \bm{X})}_4^2. 
\end{aligned}
\end{equation}

Collecting all the bounds from~\eqref{eq: BE bound shao/zhang 2025} to~\eqref{eq: vanishing conditional variance 6}, one can conclude that there exists a universal constant $C > 0$ only depending on the fourth moment $\Enb\mathcal{K}^4(1, \bm{X})$ such that 
$$
    d_{\mathrm{K}}(W) \le \frac{C}{\psi_n^2} \left (\frac{1}{\sqrt{n}} +  n^{\frac{1}{4}} \sqrt{\Delta_{1,(4)}} + \sqrt{n} \Delta_{1, (4)} +  \sqrt{n} \Delta_{1, (2)} + n \Delta_{2,(2)} + n \Delta_{1, (2)}^2
    + n^2 \Delta_{2,(2)}^2 \right )
$$
with 
$$
    \Delta_{1, (q)} := \max_{j \ne i}\norm{\nabla_i \mathcal{K}(j, \bm{X})}_q  \text{\hspace{4mm}and\hspace{4mm}} \Delta_{2, (q)} = \max_{j \ne k \ne i} \norm{\nabla_k \nabla_i \mathcal{K}(j, \bm{X})}_q. 
$$
It then follows from the triangular inequality and Lemma 3-(i) in~\cite{korolev2017bounds} that 
\begin{align*}
    &\sup_{x \in \R} \left|\mathbb{P}\left(\frac{1}{\omega_n \sqrt{n}}\left(\sum_{i=1}^n K_i\right) \leq x\right)-\Phi(x)\right| \\
    &\le d_{\mathrm{K}}(W) + \sup_{x \in \R} \abs{\Phi \left(\frac{\psi_n}{\omega_n}x \right) - \Phi(x)} \\
    &\le d_{\mathrm{K}}(W) + \max\left\{ \frac{\psi_n}{\omega_n}, \frac{\omega_n}{\psi_n}\right\} - 1.  
\end{align*}
By Lemma~\ref{lemma: general free nabla} and the Cauchy-Schwarz inequality, we have \begin{align*}
    &\psi_n^2 = \omega_n^2 + \frac{1}{n} \sum_{i \ne j} \Enb{\nabla_jK_i \nabla_i K_j} \in \left [\omega_n^2 - n \cdot \max_{i \ne j} \norm{\nabla_iK_j}_2^2, \omega_n^2 + n \cdot \max_{i \ne j} \norm{\nabla_iK_j}_2^2 \right ] \\
    & \implies
    \abs{\omega_n -  \psi_n} \le \sqrt{n} \cdot \max_{i \ne j} \norm{\nabla_i K_j}_2
\end{align*}
which in turn yields 
\begin{align*}
    \max \left\{\psi_n/\omega_n, \omega_n/\psi_n\right\} - 1 \le \sqrt{n} \cdot \max_{i \ne j} \norm{\nabla_i K_j}_2/\min\{\omega_n, \psi_n\}. 
\end{align*} Overall, we have \begin{align*}
    &\sup_{x \in \R} \left|\mathbb{P}\left(\frac{1}{\omega_n \sqrt{n}}\left(\sum_{i=1}^n K_i\right) \leq x\right)-\Phi(x)\right| \\
    & \le C_n  \left (\frac{1}{\sqrt{n}} +  n^{\frac{1}{4}} \sqrt{\Delta_{1,(4)}} + \sqrt{n} \Delta_{1, (4)} +  \sqrt{n} \Delta_{1, (2)} + n \Delta_{2,(2)} + n \Delta_{1, (2)}^2
    + n^2 \Delta_{2,(2)}^2 \right )
\end{align*}
with $C_n \lesssim 1/\psi_n^2 + 1/\min\{\omega_n, \psi_n\}$. 
\end{proof}
\begin{lemma}[First-order expansion of the absolute value]
\label{lemma: absolute value expansion}
For all $a,b \in \mathbb{R}$, we have
\[
|a+b|
= |a| + \operatorname{sgn}(a)\, b + R(a,b),
\]
where the remainder term satisfies $0 \le R(a,b) \le 2|b|\,\mathbf{1}_{\{|b|\ge |a|\}}$ and $\operatorname{sgn}(0) = 0$ by convention.
\end{lemma}

\begin{proof}
Define $R(a,b) := |a+b| - |a| - \operatorname{sgn}(a)\, b$.
Since $|\cdot|$ is convex, for any $a,b\in\mathbb{R}$ we have $|a+b| \ge |a| + g\, b \text{ for any } g \in \partial |a|$.

If $a \neq 0$, the subdifferential is $\partial|a|=\{\operatorname{sgn}(a)\}$, and hence $|a+b| \ge |a| + \operatorname{sgn}(a)\, b$,
which implies $R(a,b)\ge 0$. Otherwise, if $a=0$, we have $\operatorname{sgn}(0)=0$ and $R(0,b)=|b|\ge 0$ trivially.
Thus $R(a,b)\ge 0$ for all $a,b \in \R$.

As for its upper bound, we distinguish cases and analyze them respectively. First consider the case $|b| < |a|$.
Then $a$ and $a+b$ must have the same sign, so $|a+b| = \operatorname{sgn}(a) (a + b) = |a| + \operatorname{sgn}(a)\, b$, yielding $R(a,b)=0$. Now suppose that $|b| \ge |a|$. Since $|a+b| \le |a| + |b|$,
one can obtain
\[
R(a,b)
= |a+b| - |a| - \operatorname{sgn}(a)\,b
\le |b| - \operatorname{sgn}(a)\,b
\le 2|b|.
\]
Combining the two cases gives $R(a,b) \le 2|b|\,\mathbf{1}_{\{|b|\ge |a|\}}$
as claimed.
\end{proof}

\begin{remark}[comparison between~\cite{chatterjee2008new} and~\cite{shao2025berry}]
The result~\eqref{eq: BE bound shao/zhang 2025} is a Berry-Esseen bound, whereas Theorem~2.2 of~\cite{chatterjee2008new}: 
\begin{equation}
\begin{aligned}
\label{eq: wasserstein distance bound by chatterjee}
    d_{\mathcal{W}}(W) &\le \frac{1}{\omega_n^2} \Enb \abs{\E{V \mid W} - \omega_n^2} + \frac{1}{2\omega_n^3} \sum_{j = 1}^n \Enb \abs{\nabla_j W}^3 \\
    &\le \frac{1}{\omega_n^2} \sqrt{\mathrm{Var}(T|W)} + \frac{1}{2\omega_n^3} \sum_{j = 1}^n \Enb \abs{\nabla_j W}^3
\end{aligned}
\end{equation}
establishes a bound in the Wasserstein distance. Here $d_{\mathcal{W}}(W)$ denotes the Wasserstein distance between the statistic $W$ and the standard normal, and the random variable $T$ is given by 
$$
    T = \frac{1}{2} \sum_{A \subsetneq [n]} \frac{1}{C_A} \sum_{j \notin A}(\nabla_j W) (\nabla_j W^A), \hspace{2mm} C_A = \binom{n}{\abs{A}} (n - \abs{A}). 
$$
By comparison, $T$ in~\eqref{eq: wasserstein distance bound by chatterjee} is replaced by the random variable $V$ in~\eqref{eq: BE bound shao/zhang 2025} that avoids summation over a power set. 
Although \cite{shao2025berry} emphasizes this replacement as a substantial improvement, it does not play a central role in our context, i.e., establishing a stability-based asymptotic normality for $\mathcal{K}$-mapping. All the bounds and relevant expansions regarding $\E{V\mid W}$ we have trivially apply to $\E{T\mid W}$. 
More notably, the second term in~\eqref{eq: BE bound shao/zhang 2025} eliminates the need of any absolute third moment in the bound, compared to~\eqref{eq: wasserstein distance bound by chatterjee}. 
This constitutes the most consequential refinement for our purpose. It allows us to take advantage of Lemma~\ref{lemma: general free nabla} to justify that the second term (thereby $d_{\mathrm{K}}(W)$ overall) vanishes asymptotically, using only the first- and second-order stabilities. 
\end{remark}

\section{Proofs of \Cref{th: random center} and \Cref{cor: coverage}}\label{app:thm3.1_cor3.4}
The proof of \Cref{th: random center} can be obtained by combining \Cref{th: weakly dependent CLT} and the stability of $K_{i,r} = X_{i, r}-Q_{i, r}\QexcludingVi-d_{i,r}\dexcludingVi$.
\begin{proof}[Proof of \Cref{th: random center}]
    For a fixed index $r$, let $K_{i,r} = X_{i, r}-Q_{i, r}-d_{i,r}$. By definition, $\mathbb{E}\left[K_{i,r} \mid \bm{X}^{(-i)}\right]$ $=0$. We also have the boundedness of $K_{i,r}$ from the assumptions on $X_{i}$. Assuming $\operatorname{Cov}(X_1)$ is positive definite ensures ${\rm Var}(K_{1,r}) > 0$ uniformly over $r\in[p]$ and $i \in [n]$. 
Lemma~\ref{lemma: first order stability} and Lemma~\ref{lemma: second order stability} ensure that when $\lambda=o(\sqrt{n})$, we have $(\sqrt{n}\Delta_1)\wedge (n\Delta_2)=o(1)$ uniformly over $r$.  The claimed result then follows from \Cref{th: weakly dependent CLT}.
\end{proof}

\begin{proof}[Proof of \Cref{cor: coverage}]
Given any $r\in\Theta$, we have
$d_{i,r}\le 0$ almost surely for all $i$ as explained in \Cref{rem:random_centering}.
By consistency of $\hat\sigma_r$ and its non-degeneracy, which is guaranteed by the minimal eigenvalue condition assumed in the theorem, we have  $\hat\sigma_r/\sigma_r\rightarrow 1$ in probability. Thus, by asymptotic normality of $\tilde T_r$ and the Slutsky's theorem,
\begin{align*}
    \mathbb{P}(r\not\in\widehat{C})=&\mathbb P\left(\frac{1}{\sqrt{n}\hat\sigma_r}\sum_{i=1}^n(X_{i,r}-Q_{i,r}\QexcludingVi)\ge z_{1-\alpha}\right)\\
 \le &\mathbb P\left(\frac{1}{\sqrt{n}\hat\sigma_r}\sum_{i=1}^n(X_{i,r}-Q_{i,r}\QexcludingVi-d_{i,r}\dexcludingVi)\ge z_{1-\alpha}\right)\\
 = & \mathbb P\left(\frac{\sigma_r}{\hat\sigma_r}\tilde T_r\ge z_{1-\alpha}\right)=\alpha+o(1)\,. \qedhere
\end{align*}
\end{proof}

\section{Proofs for Bias and Power Results}

\subsection{Proof of Theorem~\ref{thm:bias}}\label{app: proof of bias}
\begin{proof}
    Let $\mathcal T_1=\{t:\mu_t-\mu_{s_r}\ge C\log (p+n)/\lambda_n\}$ and
    $\mathcal T_0=[p]\backslash (\mathcal T_1\cup\{r\})$ for some $C > 2$.

    Let $\mathcal E$ be the event that
    $|\hat\mu_t^{(-i)}-\mu_t|\le c\sqrt{\log(p+n)}/\sqrt{n}$ for all $(t,i)\in[p]\times[n]$.
    By choosing $c$ large enough, we have
    $\mathbb P(\mathcal E)\ge 1-(p+n)^{-1}$. That is, $\mathbb P(\mathcal E^c) \le (p + n)^{-1 }$. 

    By assumption $\lambda_n=o(\sqrt{n})$. So when $n$ is large enough we have $c\sqrt{\log(p+n)}/\sqrt{n}\le (C/4)\log(p+n)/\lambda_n$.
Thus on the event $\mathcal E$, we have 
    $\hat w_{r,t}^{(-i)}\le e^{-\lambda_n\hat\mu_t^{(-i)}}/e^{-\lambda_n\hat\mu_{s_r}^{(-i)}}\le e^{-\lambda_n (C/2) \frac{\log(p+n)}{\lambda_n}}=(p+n)^{-(C/2)}$
for all $t\in\mathcal T_1$. 

Now on the event $\mathcal E$,
\begin{align*}
    \theta_r-\hat \theta_r \leq & \sum_{t\in\mathcal T_1} (p+n)^{-C/2}M + C\frac{\log(p+n)}{\lambda_n} \frac{1}{n}\sum_{i=1}^n\sum_{t\in\mathcal T_0}\hat w_{r,t}^{(-i)}\\
    \le & p(p+n)^{-C/2}M + C\frac{\log(p+n)}{\lambda_n}\\
    \leq & C(p+n)^{-1} + C\log (p+n)/\lambda_n \,.
\end{align*}
Finally, we have
\begin{align*}
    \theta_r-\hat \theta_r=&(\theta_r-\hat \theta_r)\mathds{1}(\mathcal E)+(\theta_r-\hat\theta_r)\mathds{1}(\mathcal E^c)\\
    \leq & C(p+n)^{-1}+C \frac{\log (p+n)}{\lambda_n} +O_P(\mathbb P(\mathcal E^c))\\
    = & C\log(p+n)/\lambda_n + O_P((p+n)^{-1})\,. \qedhere
\end{align*}
\end{proof}

\subsection{Proof of \Cref{cor: improved bias}}

\begin{proof}
    Define
\begin{equation*}
\mathcal{T}_1=\left\{t: \mu_t-\mu_{s_r} \geq m\right\}
\end{equation*}
and 
\begin{equation*}
\mathcal{T}_0=\{t: \mu_t = \mu_{s_r}\}.
\end{equation*}
It is easy to verify that $\mathcal T_1 \cup \mathcal T_0 = [p]\backslash\{r\}$. So
\begin{equation}\label{eq: split theta}
   \theta_r - \hat\theta_r = n^{-1}\sum_{i=1}^n \sum_{t \in \mathcal T_1} \hat w_{r,t}^{(-i)} (\mu_t - \mu_{s_r}). 
\end{equation}
Let $c>\sqrt{6}M$ be a  fixed constant, and define $\mathcal{E}$ be the event that $$\left|\hat{\mu}_t^{(-i)}-\mu_t\right| \leq c \sqrt{\log (p+n)} / \sqrt{n}$$ for all $(t, i) \in[p] \times[n]$. We have
\begin{equation*}
\mathbb{P}\left(\mathcal{E}^c\right) \leq 2 p n(p+n)^{- c^2/(4M^2)} = o(n^{-1/2}).
\end{equation*}

On the event $\mathcal E$, for $t \in \mathcal T_1$
$$
\hat w_{r,t}^{(-i)} \leq \exp(-\lambda\{|\mu_t - \mu_{s_r}| - 2c\sqrt{\log (p+n)} / \sqrt{n}\})
$$
$$
\leq \exp(-\lambda\{m - 2 c \sqrt{\log (p+n)} / \sqrt{n}\})
$$
For any $n \geq \frac{16 c^2}{m^2} \log (p+n)$,
$$
\hat{w}_{r, t}^{(-i)} \leq \exp \left(-\lambda \cdot \frac{m}{2}\right).
$$
We continue \Cref{eq: split theta}, on event $\mathcal E$:
$$
\theta_r-\hat{\theta}_r \leq pM\exp(-\lambda m /2)
$$
Finally, we have
$$
\theta_r-\hat{\theta}_r=\left(\theta_r-\hat{\theta}_r\right) 1(\mathcal{E})+\left(\theta_r-\hat{\theta}_r\right) 1\left(\mathcal{E}^c\right)
$$
$$
\leq p M \exp (-\lambda m / 2) + o_P(n^{-1/2}).
$$
When $\lambda = \omega(\log n)$, 
$$
\theta_r - \hat \theta_r \leq o_p(n^{-1/2}).
$$
\end{proof}

\subsection{Proof of Theorem~\ref{thm:power_new}}\label{app: new adaptive}
\begin{proof}
Without loss of generality, we assume $r=1$, and $\tilde\mu=\mu_2\le\mu_3\le...\le \mu_p$. 
Under the assumption of the theorem, $\sup_{i,r}|X_{i,r}|$ is bounded by a constant. For simplicity, we assume that $\sup_{i,r} \abs{X_{i,r}} \le 1$ so that $X_{i,r}$ is $1$-sub-Gaussian for all $i \in [n], r \in [p]$. Thus, as $\hat{\sigma}_1$ remains bounded for all $n$, it suffices to show that with high probability, 
$$
    \frac{1}{\sqrt{n}} \sum_{v=1}^V \sum_{i \in I_v}d_{i,r} := \frac{\sqrt{n}}{V}\sum_{v=1}^V d^{(-v)}\rightarrow\infty,
$$ 
where $d^{(-v)}=\sum_{s\ge 2}\hat w_{1,s}^{(-v)}(\mu_1-\mu_s)$.

Define the indices and events
\begin{align*}
    s_1=&\max\{s\ge 2:\mu_s\le \tilde\mu+\alpha_n/\lambda_n\}\,,\\
    s_2=&\max\{s\ge 2:\mu_s\le \tilde \mu+\beta_n/\lambda_n\}\,,\\
    \mathcal E_0=&\{\hat w_{1,s}^{(-v)}\le e^{-\lambda_n(\mu_s-\tilde\mu)/2}:~\forall~s>s_2\,,~v\in[V]\}\,,\\
    \mathcal E_v=&\left\{\sup_{s_1<s\le s_2}\hat w_{1,s}^{(-v)}\le e^{-\frac{\alpha_n}{2}} \right\}.
\end{align*}
Under the theorem assumption, $X_{i,r}$ is $1$-sub-Gaussian.  According to Lemma \ref{lemma: vanishing weights} and union bound, we have
$$\mathbb P(\mathcal E_0^c)\le 2pVe^{-\frac{n\beta_n^2}{8 \lambda_n^2}}=o(1)
$$
and
$$
\mathbb P(\mathcal E_v^c)\le 2 (s_2-s_1) e^{-\frac{n\alpha_n^2}{8 \lambda_n^2}}=o(1)\,. 
$$

\textbf{Case 1.} $\mathbb C = \mathbb{C}(r)>0$. In this case $s_1 < s_2$ and $\alpha_n<\beta_n$.
\begin{align}
d^{(-v)}=&\sum_{s\ge 2}\hat w_{1,s}^{(-v)}(\mu_1-\mu_s)\nonumber\\
=&\sum_{2\le s\le s_1}\hat w_{1,s}^{(-v)}(\mu_1-\mu_s)+\sum_{s_1<s\le s_2}\hat w_{1,s}^{(-v)}(\mu_1-\mu_s)+\sum_{s>s_2}\hat w_{1,s}^{(-v)} (\mu_1-\mu_s)\nonumber\\
\ge&\frac{\alpha_n}{\lambda_n}\sum_{2\le s\le s_1}\hat w_{1,s}^{(-v)} -\frac{\beta_n}{\lambda_n}\sum_{s_1<s\le s_2}\hat w_{1,s}^{(-v)} +\sum_{s>s_2}\hat w_{1,s}^{(-v)}(\mu_1-\mu_s)\,.\label{eq:case1_init}
\end{align}

Introduce the following quantities
\begin{align*}
    d_0^{(-v)}=&\sum_{2\le s\le s_1}\hat w_{1,s}^{(-v)}(\mu_1-\mu_s),\\
    d_1^{(-v)}=&\sum_{s_1< s\le s_2}\hat w_{1,s}^{(-v)}(\mu_1-\mu_s),\\
    d_2^{(-v)}=&\sum_{s> s_2}\hat w_{1,s}^{(-v)}(\mu_1-\mu_s),\\
    W_1^{(-v)}=&\sum_{s_1< s\le s_2}\hat w_{1,s}^{(-v)},\\
    W_2^{(-v)}=&\sum_{s> s_2}\hat w_{1,s}^{(-v)}.
\end{align*}
    
Then \eqref{eq:case1_init} can be written as
\begin{align}
    d^{(-v)} =& d_0^{(-v)}+d_1^{(-v)}+d_2^{(-v)}\nonumber\\
    \ge &(1-W_1^{(-v)}-W_2^{(-v)})\frac{\alpha_n}{\lambda_n}-W_1^{(-v)}\frac{\beta_n}{\lambda_n}+d_2^{(-v)}\,.\label{eq:case1}
\end{align}
 
On the event $\mathcal E_0\cap\mathcal E_v$ we have
\begin{align}\label{eq:case1_Ia}
   W_1^{(-v)}=& \sum_{s_1<s\le s_2}\hat w_{1,s}^{(-v)}\le (s_2-s_1)e^{-\alpha_n/2}\le 1/4
   \end{align}
   and
   \begin{align}
   W_2^{(-v)}=&\sum_{s> s_2}\hat w_{1,s}^{(-v)}
   \le  (p-s_2)e^{-\beta_n/2}
   \le  1/4 ,\label{eq:case1_Ib}
\end{align}
where the inequalities hold true since by assumption $\alpha_n\ge 2(\log 4+\log \mathbb C)$ and
$\beta_n\ge 2(\log 4+\log p)$.

On event $\mathcal E_v$ we have
\begin{align}
-\frac{\beta_n}{\lambda_n}W_1^{(-v)}\ge & -\frac{\beta_n}{\lambda_n}(s_2-s_1)e^{-\alpha_n/2} \nonumber\\
\ge&-\frac{\alpha_n}{8\lambda_n}\,,\label{eq:case1_II}
\end{align}
where the last inequality holds true whenever $\alpha_n\ge 2(\log 4+\log\beta_n+\log(s_2-s_1))$, which is guaranteed by assumption when $(n,p)$ are large enough.

On $\mathcal E_0$,
\begin{align}
  d_2^{(-v)} = &  \sum_{s>s_2}\hat w_{1,s}^{(-v)}(\mu_1-\mu_s)\nonumber\\
  \ge & -\sum_{s>s_2}e^{-\frac{\lambda_n}{2}(\mu_s-\tilde\mu)}(\mu_1-\mu_s)_-\nonumber\\
  \ge &-e^{-\frac{\beta_n}{3}}\sum_{s>s_2}e^{-\frac{\lambda_n}{6}(\mu_s-\tilde\mu)}(\mu_1-\mu_s)_-\nonumber\\
    \ge &-e^{-\frac{\beta_n}{3}}\sum_{s>s_2}e^{-\frac{\lambda_n}{6}(\mu_s-\mu_1)}(\mu_s-\mu_1)_+\nonumber\\
  \ge &-e^{-\frac{\beta_n}{3}}\cdot p\cdot \sup_{x\ge 0} e^{-\frac{\lambda_n}{6}x}x\nonumber\\
\ge & -\frac{6}{e}\frac{1}{\lambda_n}pe^{-\beta_n/3}\nonumber\\
\ge &-\frac{\alpha_n}{8\lambda_n}\label{eq:case1_III}
\end{align}
where the last inequality holds whenever $\beta_n\ge 3[\log(24/e)+\log p]$ since $\alpha_n\ge 2$ under the theorem assumption.

Now plugging in \eqref{eq:case1_Ia}, \eqref{eq:case1_Ib}, \eqref{eq:case1_II}, and \eqref{eq:case1_III} into \eqref{eq:case1}, we obtain, under the event $\mathcal E_0\cap\mathcal E_v$
$$
d^{(-v)}\ge \frac{\alpha_n}{4\lambda_n}
$$
for sufficiently large $(n,p)$. 

This completes the proof for the case of $V=O(1)$ since we can use a union bound to show that $\bigcap_{v=1}^V\mathcal E_v$ has probability tending to one.

In the case of diverging $V$, such as $V=n$,
let $\epsilon_0=\mathbb P(\mathcal E_0^c)$.

We need to show that $\bar d = \frac{1}{V}\sum_{v=1}^V d^{(-v)}\ge c \frac{\alpha_n}{\lambda_n}$ for some positive constant $c$.

Define 
\begin{align*}
   \bar W_1=&\frac{1}{V}\sum_{v=1}^V W_1^{(-v)}\,.
\end{align*}

We have shown that on $\mathcal E_0$, $d_2^{(-v)}\ge -\frac{\alpha_n}{8\lambda_n}$ and $W_2^{(-v)}\le 1/4$. As a result, on $\mathcal E_0$ we have
\begin{align*}
d^{(-v)}\ge&(1-W_1^{(-v)}-W_2^{(-v)})\frac{\alpha_n}{\lambda_n}-\frac{\beta_n}{\lambda_n}W_1^{(-v)}-\frac{\alpha_n}{8\lambda_n}\\
\ge &(5/8-W_1^{(-v)})\frac{\alpha_n}{\lambda_n}-\frac{\beta_n}{\lambda_n}W_1^{(-v)}\,
\end{align*}
and hence
$$
\bar d\ge (5/8-\bar W_1)\frac{\alpha_n}{\lambda_n}-\frac{\beta_n}{\lambda_n}\bar W_1\,.
$$
Then on the event $\mathcal E_0\cap\{\bar W_1\le \alpha_n/(4\beta_n)\}$ we have (recall that since $s_1<s_2$ we must have $\alpha_n<\beta_n$ so $\alpha_n/(4\beta_n)<1/4$)
$$
\bar d \ge (5/8-1/4)\frac{\alpha_n}{\lambda_n}-\frac{\beta_n}{\lambda_n}\frac{\alpha_n}{4\beta_n}=\frac{\alpha_n}{8\lambda_n}.
$$
Therefore
\begin{align*}
    \mathbb P(\bar d\ge \alpha_n/(8\lambda_n))\ge & 1-\epsilon_0-\mathbb P(\bar W_1>\alpha_n/(4\beta_n))\\
    \ge &1-\epsilon_0-\frac{\mathbb E \bar W_1}{\alpha_n/(4\beta_n)}\\
    = &1-\epsilon_0-\frac{\mathbb E W_1^{(-v)}}{\alpha_n/(4\beta_n)}\\
    \ge &1-\epsilon_0-3(s_2-s_1)e^{-\alpha_n/2}\frac{4\beta_n}{\alpha_n}\\
    =&1-o(1)\,,
\end{align*}
where the last inequality uses $\alpha_n\ge 2$ and that (for large enough $n,p$)
\begin{align*}
    \mathbb E W_1^{(-v)} \le & \mathbb P(\mathcal E_v)(s_2-s_1)e^{-\frac{\alpha_n}{2}}+\mathbb P(\mathcal E_v^c)\\
    \le & (s_2-s_1)e^{-\alpha_n/2}+2(s_2-s_1)e^{-\frac{n\alpha_n^2}{8 \lambda_n^2}}\\
    \le & 3(s_2-s_1)e^{-\alpha_n/2}\,.
\end{align*}

\textbf{Case 2.} $\mathbb C=0$. In this case we have, on event $\mathcal E_0$,
$$
d^{(-v)}=d_0^{(-v)}+d_2^{(-v)}\ge (1-W_2^{(-v)})\frac{\alpha_n}{\lambda_n}-\frac{\alpha_n}{8\lambda_n}\ge \frac{5\alpha_n}{8\lambda_n}\,.
$$
So for $\sqrt{n}\bar d\rightarrow\infty$ with high probability it suffices to have $\sqrt{n}\alpha_n/\lambda_n\rightarrow\infty$ which is equivalent to $\sqrt{n}(\mu_1-\tilde\mu)\rightarrow \infty$.
\end{proof}

\begin{lemma} \label{lemma: vanishing weights}
    Suppose that every dimension of the sample vector $X$ is $\delta$-sub-Gaussian for some constant $\delta > 0$. Define $\mu = \mathbb{E} X$. Let $f \in [p]$ be a fixed index and $\mathcal{S} \subseteq [p]$ be an index set such that for all $s \in \mathcal{S}, \mu_s - \mu_f \ge \gamma\ge 0$. Then, we have 
    $$
        \mathbb{P} (\hat{w}^{(-v)}_{r,s} > \exp ( - \lambda_n (\mu_s - \mu_f)/2) \text{ for some } s \in \mathcal{S}) \le 2 \exp \left(- \frac{n \gamma^2}{8\delta^2} + \log \abs{\mathcal{S}}\right),
    $$
    where $\abs{\mathcal{S}}$ denotes the cardinality of $\mathcal{S}$. 
\end{lemma}

\begin{proof}
    We apply the sub-Gaussian tail bound of the sample mean (e.g., Theorem 2.6.2 in~\cite{vershynin2018high}), and directly obtain
\begin{align*}
    \mathbb{P} &\left(\hat{w}^{(-v)}_{r,s} > \exp(-\lambda_n (\mu_s - \mu_f)/2) \text{ for some } s \in \mathcal{S} \right) \\
    &\le \sum_{s \in \mathcal{S}} \mathbb{P} \left (\hat{w}^{(-v)}_{r,s} > \exp(-\lambda_n (\mu_s - \mu_{f})/2)\right) \\
    &= \sum_{s \in \mathcal{S}} \mathbb{P} \left( \frac{\exp (-\lambda_n \hat{\mu}^{(-v)}_s)}{\sum_{t \ne r} \exp (-\lambda_n \hat{\mu}^{(-v)}_t) } > \exp(-\lambda_n (\mu_s - \mu_f)/2) \right) \\
    & \le \sum_{s \in \mathcal{S}} \mathbb{P} \left ( \exp ( -\lambda_n (\hat{\mu}^{(-v)}_s - \hat{\mu}^{(-v)}_{f}) ) > \exp(-\lambda_n (\mu_s - \mu_{f})/2) \right) \\
    & = \sum_{s \in \mathcal{S}} \mathbb{P} \left(\hat{\mu}_s^{(-v)} - \hat{\mu}^{(-v)}_{f} < (\mu_s - \mu_{f})/2 \right) \\
    & \le \sum_{s  \in \mathcal{S}} \mathbb{P} \left( \abs{\hat{\mu}^{(-v)}_s - \hat{\mu}^{(-v)}_{f} - (\mu_s - \mu_{f})} > (\mu_s - \mu_{f})/2 \right) \\
    & \le \sum_{s \in \mathcal{S}} 2\exp \left ( - (1-1/V) n (\mu_s - \mu_{f})^2/(8\delta^2) \right) \\
    & \le 2 \exp \left (-n \gamma^2/(8\delta^2) + \log \abs{\mathcal{S}} \right),
\end{align*}
where the second last inequality follows from sub-Gaussian concentration and the last inequality follows from $V\ge 2$. 
\end{proof}


\section{Initial candidate for data-driven weighting parameter tuning} \label{appendix: initial lambda_0}
To determine the largest $\lambda$ that sustains the desired coverage, we proposed an iterative algorithm in Section~\ref{section: data-driven}. Here we provide details about determining the initial value $\lambda_0$ in the algorithm. We set
\begin{equation*}
\lambda_0 = \frac{\sqrt{n}}
{2.5 \cdot \mathrm{sd}(X_{i,r} - X_{i, \hat{s}})},
\end{equation*}
where the estimated index $\hat s = \hat{s}^{(-v_i)}_r = \argmin_{s \ne r} \sum_{j \notin I_{v_i}} X_{j,s}$ is a generalization of the LOO definition in~\eqref{eq:loo-arg-min} to any fold number. The quantity $\mathrm{sd}(X_{i,r} - X_{i, \hat{s}})$ denotes the sample standard deviation of $\{X_{i,r} - X_{i, \hat{s}}, i\in [n]\}$. This initial value $\lambda_0$ is motivated by the theoretical analysis in Lemma~\ref{lemma: first order stability}. In the last step of its proof~\eqref{eq: bound Q difference}, we essentially seek to have the bound $C_{V} \mathbb{M}_j \lambda n^{-1} \sum_{s \ne r} \hat{w}_{r,s}\left|X_{i,r} - X_{i, s}-(\mu_r - \mu_s)\right| \ll n^{-1/2}$, where $C_V$ is a constant only dependent of the number of folds $V$ and the quantity $\mathbb{M}_j$ is given by 
$$
\mathbb{M}_j = \max_{s \in [p]} \abs{X'_{j, r} - X'_{j,s} - (X_{j, r} - X_{j,s})}
$$
for some $j \notin I_{v_i}$. Because we standardize the difference $X_{i,r} - X_{i,s}$ for all $s \ne r$ before choosing the initial candidate $\lambda_0$, one can expect that the deviation $\abs{X_{j,r} - X_{j,s}}$ is essentially bounded by a constant. We thus regard the quantity $\mathbb{M}_j$ as a fixed value for all $j \notin I_{v_i}$. 

Suppose that the sample size $n$ is sufficiently large. The exponential weightings would nearly recover the argmin, so the summation over $s$ is close to the absolute deviation $\abs{X_{i,r} - X_{i, \hat{s}} - (\mu_r - \mu_{\hat{s}})}$. Typically, one may expect it to be roughly bounded by a constant multiplying the standard deviation of $X_{i,r} - X_{i, \hat{s}}$ with high probability. This intuition ultimately leads to the choice of the given $\lambda_0$, where the conservative constant $2.5$ was selected empirically across a variety of simulation setups to ensure that
$\lambda_0$ itself can maintain asymptotic normality. 

To point out, when the estimated index $\hat{s}$ does not vary with $i$, the standard deviation $\mathrm{sd}(X_{i,r} - X_{i, \hat{s}})$ would be exactly $1$ as a result of our pre-processing. Alternatively, if the estimated index $\hat{s}$ varies with $i$, we may still expect the standard deviation to be approximately $1$ as the varying $\hat{s}$ phenomenon occurs when the (scaled) mean differences are sufficiently close.
In practice, there is a tiny chance to obtain $\mathrm{sd}(X_{i,r} - X_{i, \hat{s}}) = 0$ (this can happen with binary data that has highly correlated dimensions due to fold-splitting). In the case, we reset it to $1$ to avoid dividing by zero. 

\section{Heuristic Simultaneous Inference}
\label{sec: heuristic simultaneous inference}

\begin{figure}[!htbp]
  \centering  \includegraphics[width=0.99\textwidth]{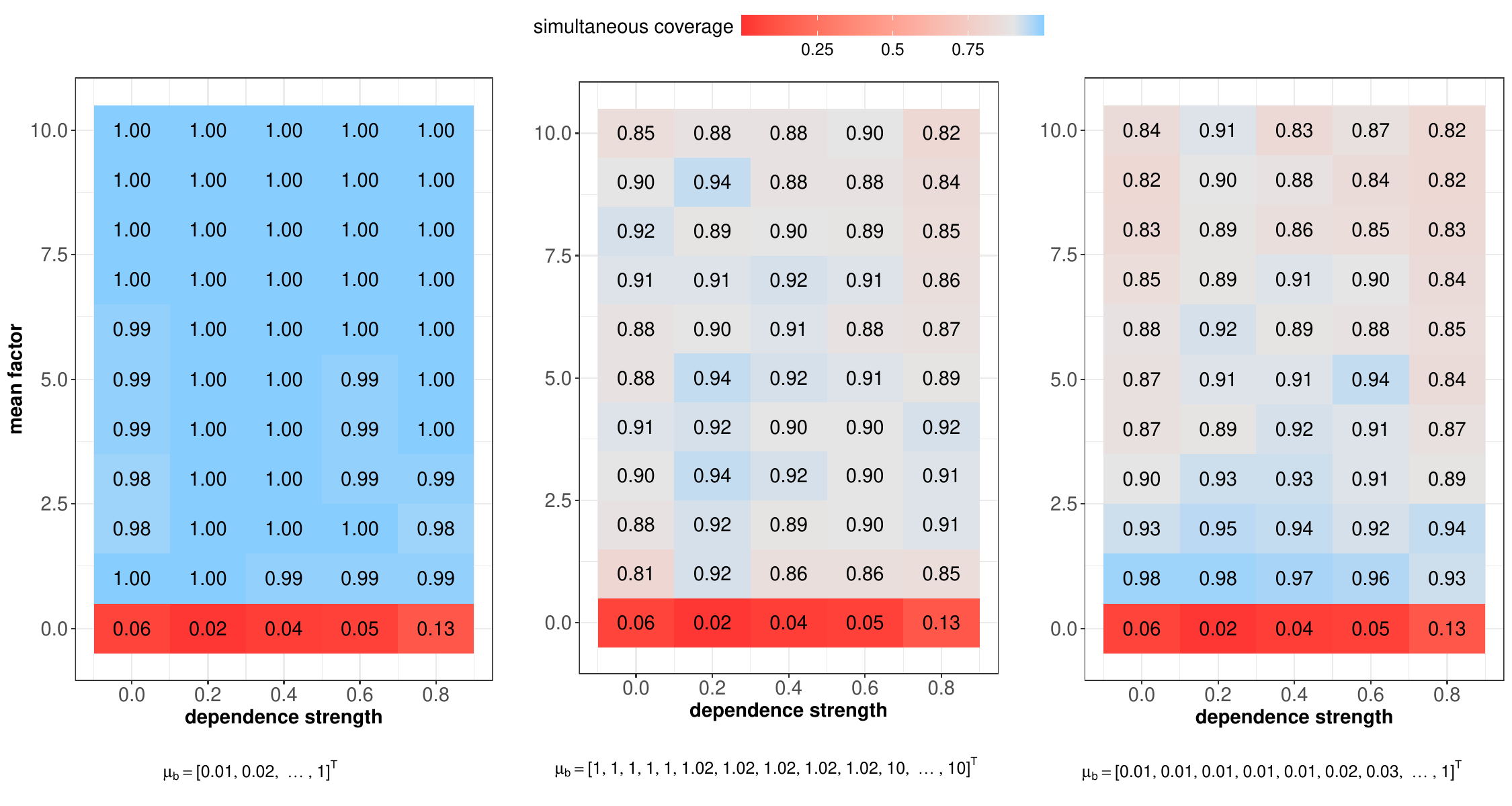}
  \caption{The proposal cannot always achieve simultaneous coverage under three base mean landscapes. The simultaneous coverage guarantee is evaluated under varying settings of the mean factor (signal strength) $f$ and dependence strength $\varrho$. For each configuration $(\mu_b, f, \varrho)$, we perform 100 simulation repetitions with a sample size of 1000.}
  \label{fig: fail simultaneous coverage}
\end{figure}

This work focuses on developing an algorithm that satisfies the marginal coverage~\eqref{eq: marginal validity}. It is not expected to achieve simultaneous coverage $\mathbb{P}(\Theta \subset \widehat{C}) > 1-\alpha$, demonstrated in \Cref{fig: fail simultaneous coverage}. When the mean factor is $0$, $\Theta = [p]$ and simultaneous coverage falls well below the nominal level. In the left subplot, with a positive mean factor, $\Theta = \{1\}$ so marginal and simultaneous coverage coincide, giving us valid results. In other cases, where $\Theta$ includes 5 elements, the method under-covers.

\begin{figure}[t!]
  \centering
  \includegraphics[width=0.99\textwidth]{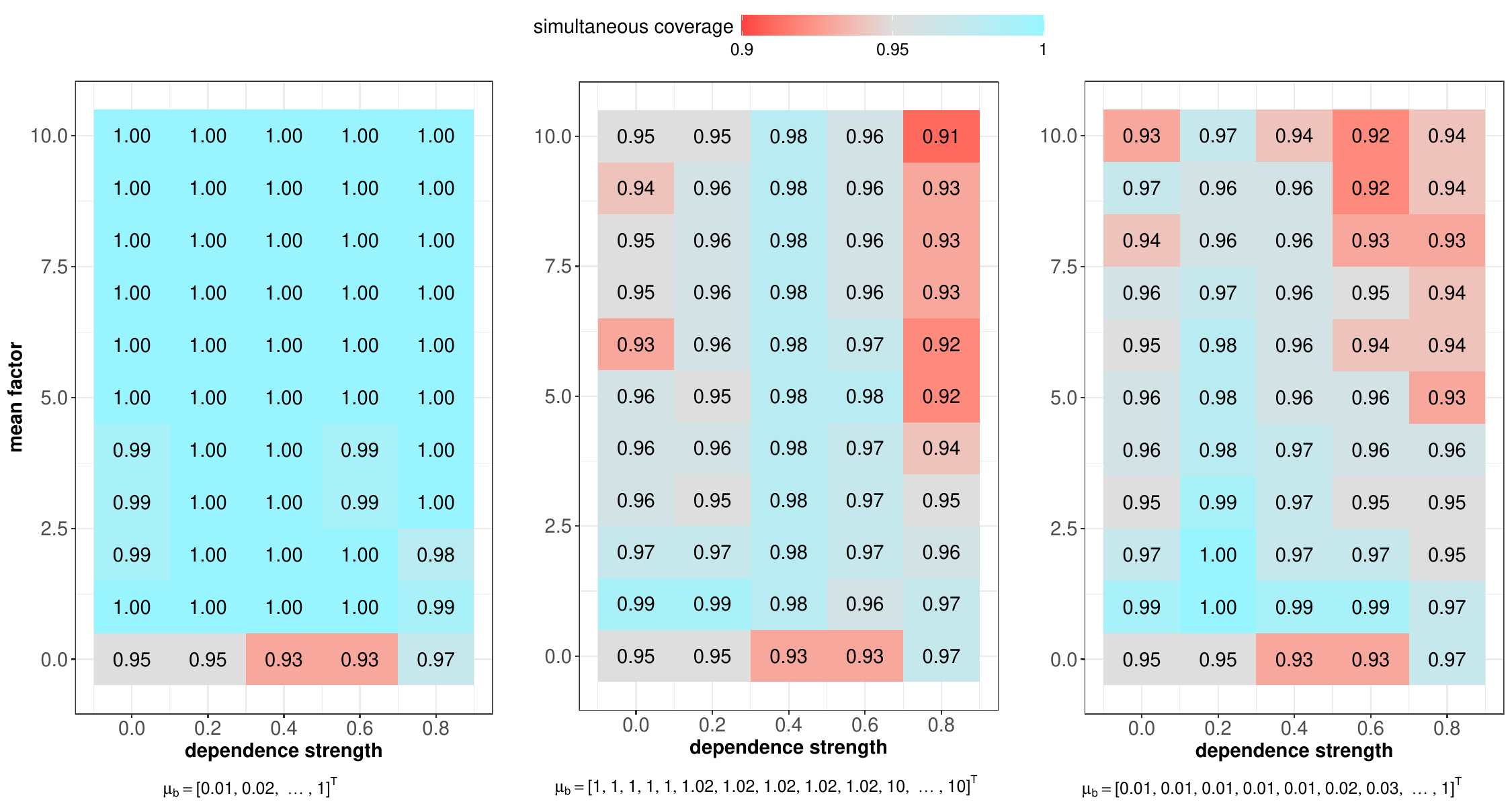}
  \caption{Simultaneous coverage of the heuristic two-step approach under three base mean landscapes. The simultaneous coverage guarantee is evaluated under varying settings of the mean factor (signal strength) $f$ and dependence strength $\varrho$. For each configuration $(\mu_b, f, \varrho)$, we perform 100 simulation repetitions with a sample size of 1000.}
  \label{fig: simultaneous coverage}
\end{figure}

In contrast to the approach of~\cite{mogstad2024inference}, our proposed LOO algorithm does not readily lend itself to extensions that enjoy theoretically guaranteed simultaneous coverage. With this in mind, we limit our current focus to exploring a heuristic adjustment---while acknowledging its lack of theoretical guarantee---with a hope that it serves as a practical starting point to inspire future developments toward procedures ensuring simultaneous coverage.

\subsection{Methods}
A straightforward approach to achieving simultaneous coverage in theory is to apply Bonferroni correction or Holm’s procedure directly to our established LOO algorithm. However, this method remains unrefined, as it overlooks the underlying mean landscape and the dependence among the test statistics $T_r$. Indeed, when the number of dimensions $p$ is large, it mostly probably would exhibit very limited power. In view of this, we suggest a two-step procedure:
\begin{enumerate}[label=(\arabic*.)]
    \item Run Algorithm~\ref{algorithm: exp weighting} with the critical value $z_{1-\alpha}$ for the given data, and compute the cardinality of the resulting confidence set $\widehat{C}$. Let  $\hat{\mathfrak{N}}_\alpha = \min\{\lceil\frac{|\widehat{C}|}{(1- \alpha)} \rceil, p\}$. 
    \item Run Algorithm~\ref{algorithm: exp weighting} with the adjusted critical value $z_{1 - \alpha/\hat{\mathfrak{N}}_\alpha}$ to output the final confidence $\widehat{C}_{u}$. 
\end{enumerate}
The first step functions as a pre-screening stage to approximate the cardinality of the true argmin set $\Theta$. According to our theoretical guarantee for marginal coverage, the mean of $\hat{\mathfrak{N}}_\alpha$
provides an upper bound on 
$|\Theta|$ in expectation, although this bound does not hold with high probability in general. 

\subsection{Simultaneous Coverage}
Numerically, we evaluate the simultaneous coverage of the two-step procedure under the settings described in Section~\ref{sec: method comparison setup}, as well as on a new mean landscape that is flat across the first five dimensions and exhibits a gradual increasing pattern thereafter. As illustrated in Figure~\ref{fig: simultaneous coverage}, the proposed two-step approach generally attains simultaneous coverage at the significance level $\alpha = 0.05$. More specifically, under any strictly increasing mean landscape, simultaneous coverage coincides with marginal coverage, so it is as expected that simultaneous coverage is maintained across all such scenarios depicted in the left plot. When the argmin is non-unique, the two-step approach typically attains simultaneous coverage although it can be somewhat liberal in scenarios where both the signal strength (mean differences) and the dependence among dimensions are relatively strong.

\subsection{Finite-Sample Power}
\begin{figure}[t!]
  \centering
\includegraphics[width=\textwidth]{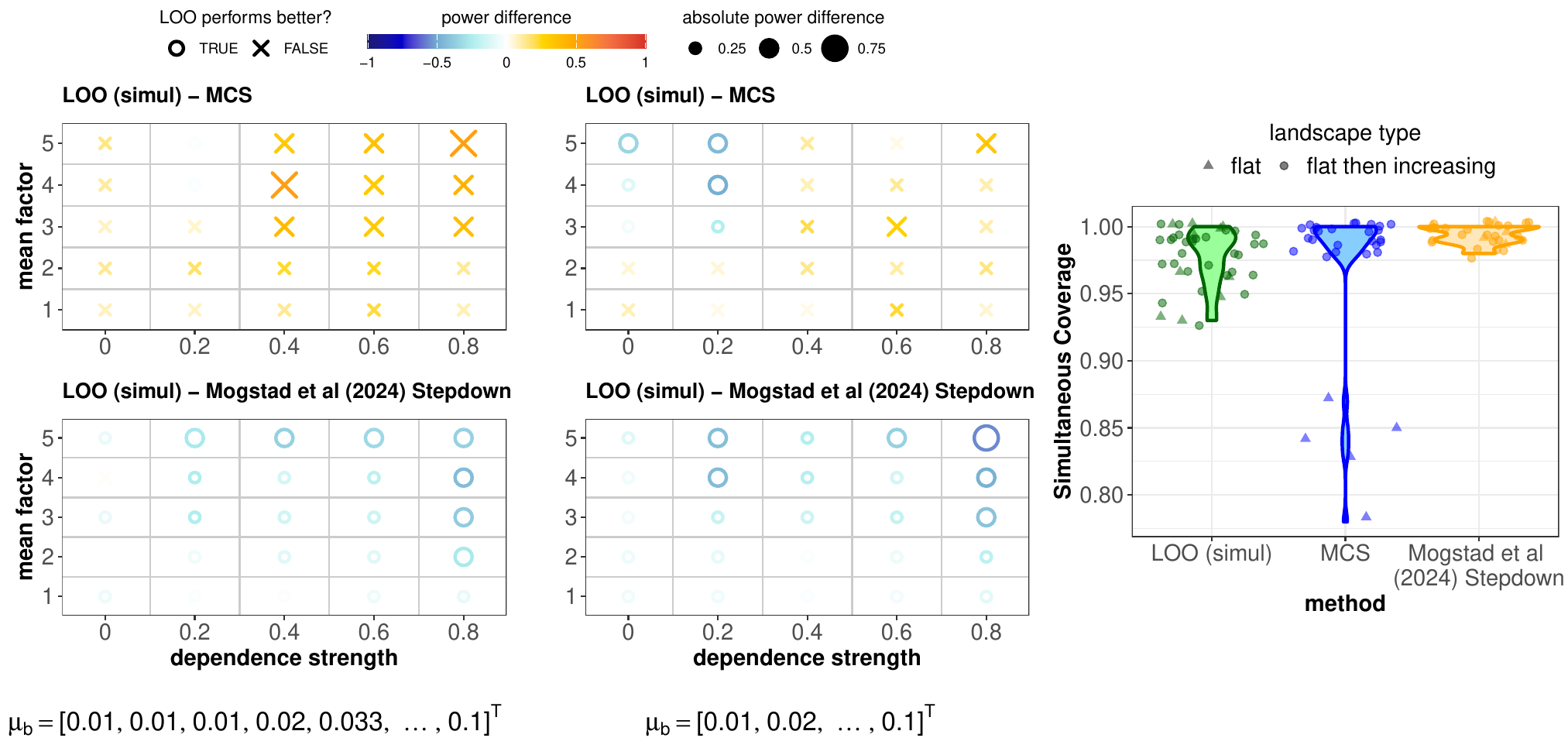}
  \caption{Performance comparison among methods for simultaneous coverage: statistical power under ``increasing" (leftmost) and ``flat-then-increasing" (middle) landscapes and corresponding empirical simultaneous coverages (rightmost). Each cell in the heatmaps corresponds to a different simulation setting. The x-axis corresponds to different dependency strength $\varrho$, and on the y-axis, signal strength $f$ is varied. The color in each cell illustrates the difference in the average number of false negatives between the proposed heuristic two-step LOO procedure and one literature method. A more negative value indicates a greater advantage of the proposed method over its competitor in rejecting sub-optimal dimensions. The rightmost violin plot displays the empirical simultaneous coverage achieved by the evaluated methods across various simulation settings under the ``flat'' and ``flat-then-increasing'' landscape scenarios. Every simulation result is conducted over $100$ repetitions with a sample size $n = 100$.}
  \label{fig: simultaneous power comparison}
\end{figure}

To evaluate the power of the proposed two-step approach, we compare its performance with the simultaneous coverage methods developed by~\cite{hansen2011model} and~\cite{mogstad2024inference}. For ease of reference, we refer to the two methods as \textit{MCS} and \textit{Stepdown}, respectively. Although one can naively invert simultaneous confidence intervals for ranks to obtain an argmin confidence set with simultaneous coverage,~\cite{mogstad2024inference} specifically developed a stepdown procedure tailored to the simultaneous argmin inference problem. Our comparison focuses on this stepdown method. Given the computational intensity of the MCS procedure, we restrict our comparison to a few representative cases. In particular, we focus on settings with a moderate dimension of $p = 10$ and a sample size of $n = 100$. As shown in the left two plots of Figure~\ref{fig: simultaneous power comparison}, the MCS approach generally demonstrates the highest power (excluding the most number of sub-optimal dimensions). However, we note that this method implicitly assumes a fixed cardinality for $\Theta$ to guarantee theoretical asymptotic coverage. Even under this assumption, achieving such coverage may require applying an exponentially large union bound (see Appendix A.10 in~\cite{kim2025locally} for a discussion), which can hinder the method’s ability to maintain simultaneous coverage in finite samples. Our numerical experiments indicate that its simultaneous coverage may drop to approximately $0.85$---below the nominal level of $1 - \alpha = 0.95$---in scenarios where the underlying mean landscape is flat, as shown in the rightmost plot of Figure~\ref{fig: simultaneous power comparison}. Conversely, the stepdown procedure by~\cite{mogstad2024inference} achieves simultaneous coverage in our numerical results, consistent with its theoretical guarantee, but this comes at the cost of reduced power. 

\begin{figure}[t!]
  \centering
\includegraphics[width=\textwidth]{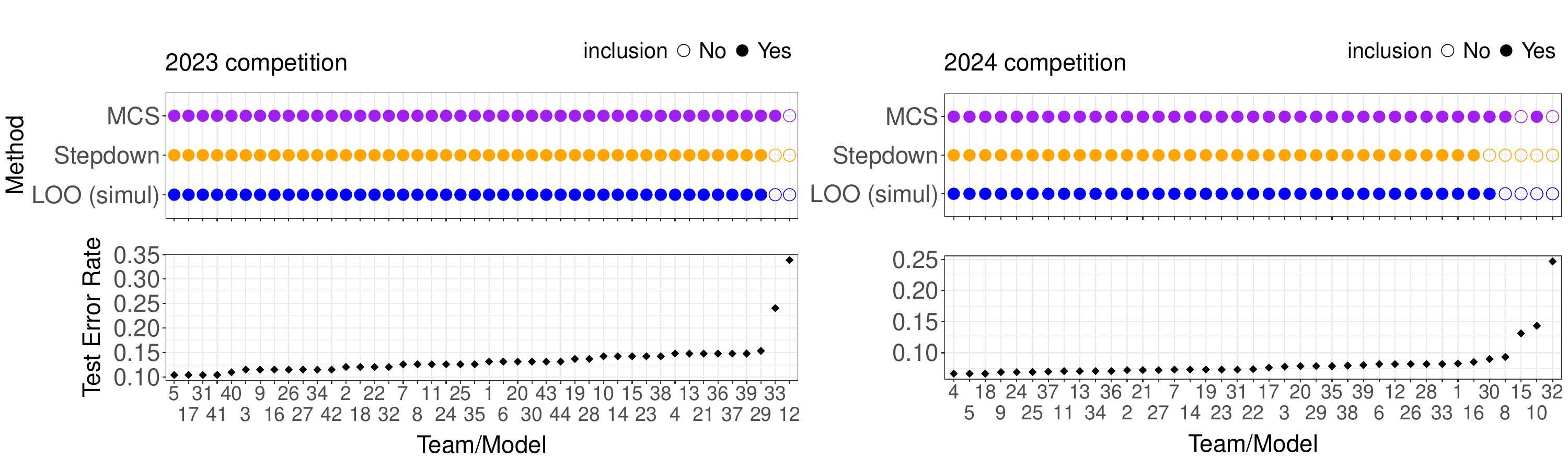}
\vspace{-20pt}
\caption{Confidence sets with real data (simultaneous coverage methods). We compare the proposed heuristic two-step LOO algorithm (\textit{LOO (simul)}), the stepdown procedure by~\cite{mogstad2024inference} (\textit{Stepdown}), and the model confidence set (\textit{MCS}) \citep{hansen2011model} over the test results in 2023 (left) and 2024 (right) classification competitions.}
  \label{fig: real data simultaneous inference}
\end{figure}

Here we also investigate their performance over the two real data sets in Section~\ref{section: real_data}. Following the same protocol as before, we repeat each method $100$ times to capture the variability due to inherent randomness. Figure~\ref{fig: real data simultaneous inference} presents a comparison for the confidence sets for one realization. In practice, the simultaneous coverage methods tend to eliminate only the most clearly suboptimal models. For instance, in the 2023 classification competition, the average sizes of the resulting confidence sets were $42$, $42$, and $43$ for the LOO (simul), Stepdown, and MCS methods respectively. Similarly, for the 2024 competition, these methods produced average cardinalities of $33.71, 34.01$ and $37$.

\begin{figure}[t!]
  \centering
\includegraphics[width=0.75\textwidth]{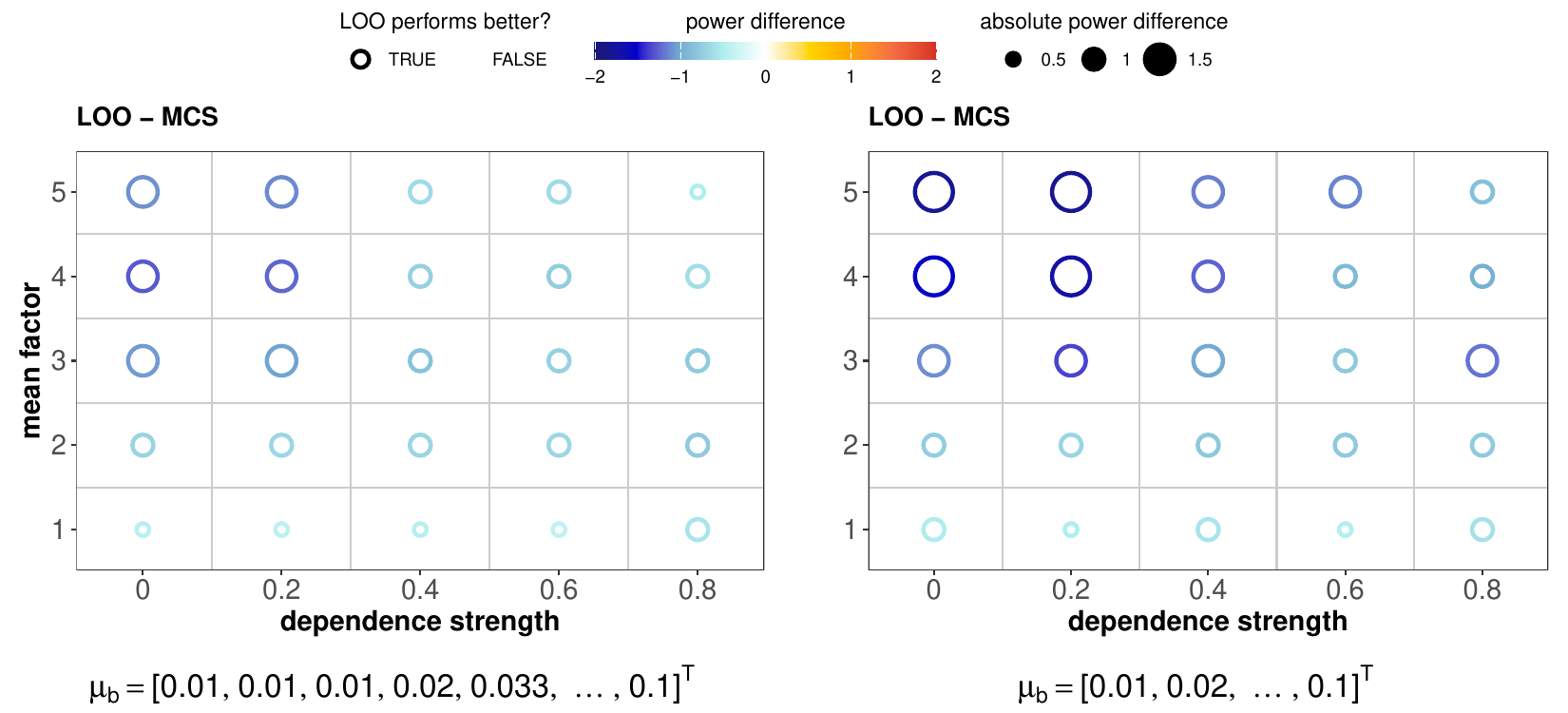}
  \caption{Performance comparison between the LOO algorithm (marginal coverage) and the MCS (\cite{hansen2011model}) procedure (simultaneous coverage): statistical power under ``increasing" (left) and ``flat-then-increasing" (right) landscapes. Each cell in the heatmaps corresponds to a different simulation setting. The x-axis corresponds to different dependency strength $\varrho$, and in the y-axis, signal strength $f$ is varied. The color in each cell illustrates the difference in the average number of false negatives between the proposed marginal LOO method and the MCS procedure. A more negative value indicates a greater advantage of the proposed LOO method in rejecting sub-optimal dimensions. Every simulation result is conducted over $100$ repetitions with a sample size $n = 100$.}
  \label{fig: LOO vs MCS}
\end{figure}

Such loss in statistical power is also evident in our simulation studies. This is illustrated in Figure~\ref{fig: LOO vs MCS}, where we compare the performance of the marginal LOO method against the MCS procedure. Since the LOO method is designed to guarantee only marginal coverage, it conceptually allows for more aggressive exclusion of sub-optimal dimensions, which often translates to higher statistical power in practice. In contrast, the MCS procedure enforces simultaneous coverage, which provides stronger inferential guarantees but can lead to very conservative selections. Ultimately, the choice between prioritizing simultaneous coverage or maximizing statistical power hinges on the specific goals of the application at hand as well as the broader scientific and practical incentives driving the analysis.

\section{Validity violations of the method by~\cite{futschik1995confidence}}\label{app:futschik.validity.violation}

In \cite{futschik1995confidence}, the authors require the true variance of $X_{1,s}$, $\sigma_s^2$, to be known to ensure the validity of their method. This is rarely the case in practice, so one may intend
to replace $\sigma_s^2$ by its estimate $\hat \sigma_s^2$. In this section, we illustrate that replacing the true $\sigma_s$ in the statistic~\eqref{eq: selection rule Gupta} with its sample estimate $\hat{\sigma}_s$ generally leads to validity violations; formally, using $\hat{\sigma}_s$ instead of $\sigma_s$ typically yields $\mathbb{P}(r \in \widehat{C}_1 \cap \widehat{C}_2) < 1 - \alpha$ for $r \in \Theta$ at the significance level of $\alpha$. In fact, because the method by~\cite{futschik1995confidence} is a two-step variant of the selection rule by~\cite{gupta1965some}, it suffices to show that the latter does not adapt to sample standard deviation $\hat{\sigma}_s$. 

We will illustrate the difficulty in a simple $p = 2$ case. Let $X_1, \ldots, X_n \in \R^2$ be IID samples such that $\Enb X_1 := \mu = [\upmu, \upmu]^\top$ for some $\upmu \in \R$, $X_{1,1}$ is independent of $X_{1,2}$, and $\max \{ \Enb X_{1,1}^4, \Enb X_{1,2}^4\} < \infty$. Let $\sigma^2 = \mathrm{Var}(X_{1,1}) = \mathrm{Var}(X_{1,2})$. Suppose that we want to test the first dimension using the selection rule in~\eqref{eq: selection rule Gupta}. In the case, the correct statistic that one should use reduces to the difference $T^{Gupta}_1 = \sqrt{n} \left (\frac{\hat{\mu}_1}{\sigma} - \frac{\hat{\mu}_2}{\sigma} \right)$ and the quantile $q_{(1 - \alpha), 2}$ is simply the $\alpha$ upper quantile of $N(0,2)$. At the level of $\alpha = 0.05$, the quantile $q_{(1-\alpha), 2}$ is approximately equal to $2.33$. If $\sigma$ is known, we know from the pairwise independence that $\sqrt{n} ( \frac{\hat{\mu}_1 - \upmu}{\sigma})$ and $\sqrt{n} (\frac{\hat{\mu}_2 - \upmu}{\sigma})$ are asymptotically distributed as two independent $N(0,1)$'s. Thus, the asymptotic distribution of their difference $T^{Gupta}_1$ is asymptotically $N(0, 2)$, where we can directly conclude the validity. 

However, if we instead considered the statistic $G = \sqrt{n} \left (\frac{\hat{\mu}_1}{\hat{\sigma}_1} - \frac{\hat{\mu}_2}{\hat{\sigma}_2} \right)$, where $\hat{\sigma}_1, \hat{\sigma}_2$ are the sample standard deviations for $\sigma$ using the samples of the first and second dimensions respectively, the validity would be easily violated. As counter-intuitive it may sound, the violation can be theoretically justified. For now, we assume that neither $X_{1,1}$ nor $X_{1,2}$ is a linear transformation of a Bernoulli random variable. Note that $G$ is the difference of two independent non-central t statistics. By Theorem 2.1 (ii) in~\cite{Bentkus2007}, they admit the stochastic convergences
\begin{equation*} \label{eq: non-standard normal asymptotics}
    \sqrt{n} \left (\frac{\hat{\mu}_1}{\hat{\sigma}_1} - \frac{\upmu}{\sigma} \right) \cind Z_1 \text{ and } \sqrt{n} \left (\frac{\hat{\mu}_2}{\hat{\sigma}_2} - \frac{\upmu}{\sigma} \right) \cind Z_2, \\
\end{equation*}
where $Z_1 \sim N(0, \tau^2_1)$ is independent of $Z_2 \sim N(0, \tau^2_2)$ with $\tau^2_r = 1 - \frac{M_{3,r} \upmu}{\sigma} + \frac{(M_{4,r} - 1)\upmu^2}{4 \sigma^2}$ for $r \in \{ 1, 2\}$. The constant $M_{k,r}$ denotes the scaled central moment $\Enb(X_{1,r} - \upmu)^k/\sigma^k$, $k \in \{3, 4\}$. By the pairwise independence and the continuous mapping theorem, we thus end up with $G \cind N(0, \upsilon^2)$, where
$$
\upsilon^2 = \tau^2_1 + \tau^2_2 = 2 - \frac{(M_{3,1} + M_{3,2})\upmu}{\sigma} + \frac{(M_{4,1} + M_{4,2} - 2)\upmu^2}{4\sigma^2}.
$$
The variance is not equal to $2$ in general (unless $\upmu = 0$ for instance), so the validity no longer holds. In fact, if at least one of $X_{1,1}$ and $X_{1,2}$ were a linear transformation of a Bernoulli random variable, the violation could be even more apparent because the asymptotic distribution of $G$ would not be normally distributed, shown by Theorem 2.1 (i) in~\cite{Bentkus2007}. 

\begin{figure}[!tbp]
    \centering
    \includegraphics[width =\textwidth]{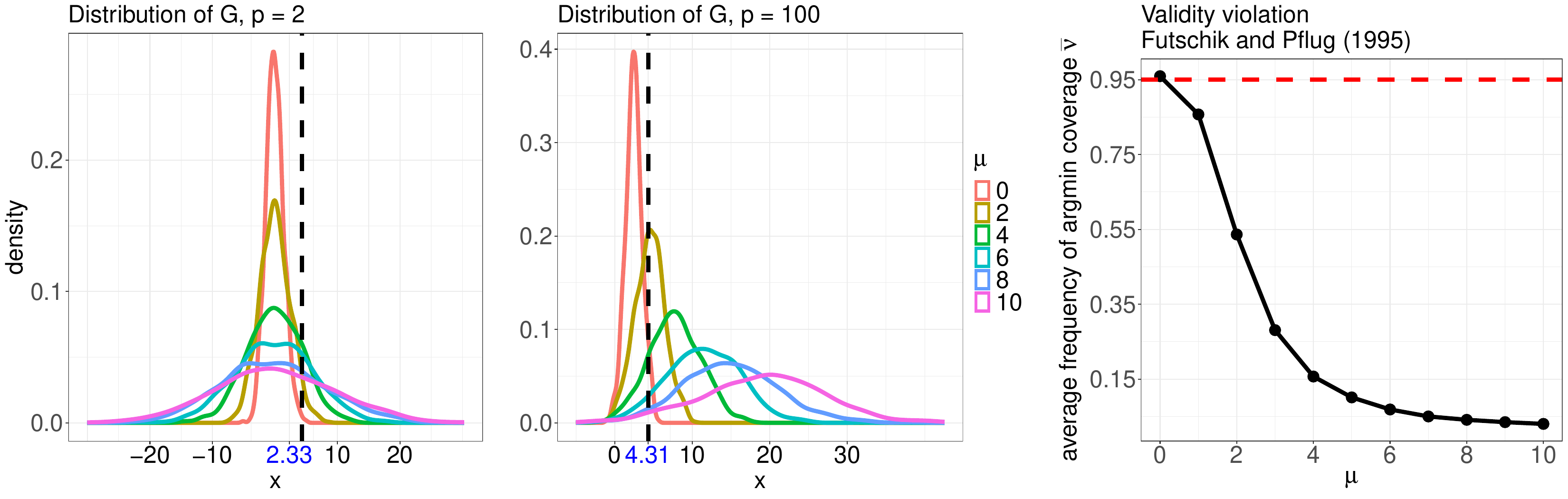}
    \caption{Validity violations of the selection rule by~\cite{gupta1965some} and the two-step procedure by~\cite{futschik1995confidence} in the presence of sample standard deviation $\hat{\sigma}_s$. The left and middle plots display the distributions of their test statistic $G$ when replacing with sample standard deviations $\hat{\sigma}_s$, with dimensions $p = 2, 100$. The densities are based on $1000$ multivariate normal samples with a flat mean landscape $\mu = [\upmu, \ldots, \upmu]^\top, \upmu \in \{0, 2, \ldots, 10\}$ and identity covariance matrix over $1000$ repetitions. The quantile $q_{(1-\alpha), p}$ at the level of $\alpha = 0.05$ is marked in blue on the x-axis. The proportion of $G$ exceeding $q_{(1-\alpha), p}$ approximates the expected type I error (false positive rate) in the selection rule by~\cite{gupta1965some} using sample standard deviations $\hat{\sigma}_s$. In the right plot, we show how the validity violations of the two-step procedure by~\cite{futschik1995confidence} escalates with increasing $\upmu$, using the same simulation setup in the middle plot.}
    \label{fig: gupta validity violation}
\end{figure}

\subsection{Numerical Verification}
We now support the previous discussion with a simple numerical example. Let $X_{1,1}, X_{1,2} \overset{\mathrm{IID}}{\sim} N(\upmu, 1)$. It follows that $M_{3,1} = M_{3,2} = 0$, $M_{4,1} = M_{4.2} = 3\sigma^3/\sigma^3 = 3$ and $\upsilon^2 = 2 + \upmu^2$. In the left subplot of Figure~\ref{fig: gupta validity violation}, the distribution of their resulting $G$ is simulated with sample size $n =1000$ and different values of $\upmu \in \{0, 2, \ldots, 10\}$. The distributions are based on $1000$ repetitions. It is evident that the variance of $G$ grows along with an increase in $\upmu$. This would cause the selection rule using $G$ to undercover any argmin dimension $r \in \Theta$ since a greater dispersion would result in a higher-than-expected proportion of $G$ realizations exceeding $q_{(1-\alpha),2} \approx 2.33$.

One can expect that such validity violations persist when considering any  fixed dimension $p > 2$ under a flat true mean landscape, i.e., $\mu = [\upmu, \ldots \upmu]^\top$. In the case, the relevant statistic $G(p)$ would be
$$
    G(p) := \sqrt{n} \left ( \frac{\hat{\mu}_r}{\hat{\sigma}_r} - \min_{s \ne r} \frac{\hat{\mu}_s}{\hat{\sigma}_s}  \right ) = \sqrt{n} \left (\frac{\hat{\mu}_r}{\hat{\sigma}_r} - \frac{\upmu}{\sigma} \right) - \min_{s \ne r} \left \{ \sqrt{n} \left (\frac{\hat{\mu}_s}{\hat{\sigma}_s} - \frac{\upmu}{\sigma} \right) \right \}.
$$
As the key issue lies in the scaled variance for the asymptotic distributions of non-central t statistics centered at $\upmu/\sigma$, the statistic $G_r$, in general, does not follow the distribution of $G^\ast(p) := \varepsilon_r - \min_{ s \ne r} \varepsilon_s$ with $\varepsilon_s \overset{\mathrm{IID}}{\sim} N(0, 1)$. In the middle subplot of Figure~\ref{fig: gupta validity violation}, we simulate samples from $N(\mu, I_{100})$ with $\mu = [\upmu, \ldots, \upmu]^\top \in \R^{100}$ and present the distributions of their resulting $G$ for different values of $\upmu$. This is based on a sample size of $n = 1000$ and $1000$ repetitions as earlier. More severe validity violations are observed: the $95\%$ percentile of $G^\ast(100)$ is roughly $4.31$, but it can be lower than the median of $G(100)$ for certain $\upmu$. 

Finally, note that the issue carries over to the two-step procedure by~\cite{futschik1995confidence}. In the right subplot of Figure~\ref{fig: gupta validity violation}, we consider the same simulation setup for $p = 100$, and compute the \textit{average frequency of argmin coverage} $\overline{\nu}$ (see Appendix~\ref{app: validity simulation}) for the method. As replicated in the numerical experiments, the coverage frequency $\overline{\nu}$ is predicted to fall below $0.95$ at the level of $\alpha = 0.05$ whenever $\upmu \ne 0$.

\section{Validity of Methods} \label{app: validity simulation}
We compared the power of several methods in $100$ settings in Section~\ref{sec: method comparison setup}. To assess the validity, we also include five flat landscape scenarios where $f = 0$ and $\abs{\Theta} = [p]$. In this case there is no difference between ``increasing" and ``3-tier" base mean $\mu_b$. 
We thus end up with 
105 distinct settings. 

\begin{figure}[!tb]
    \centering
    \includegraphics[width =0.85\textwidth]{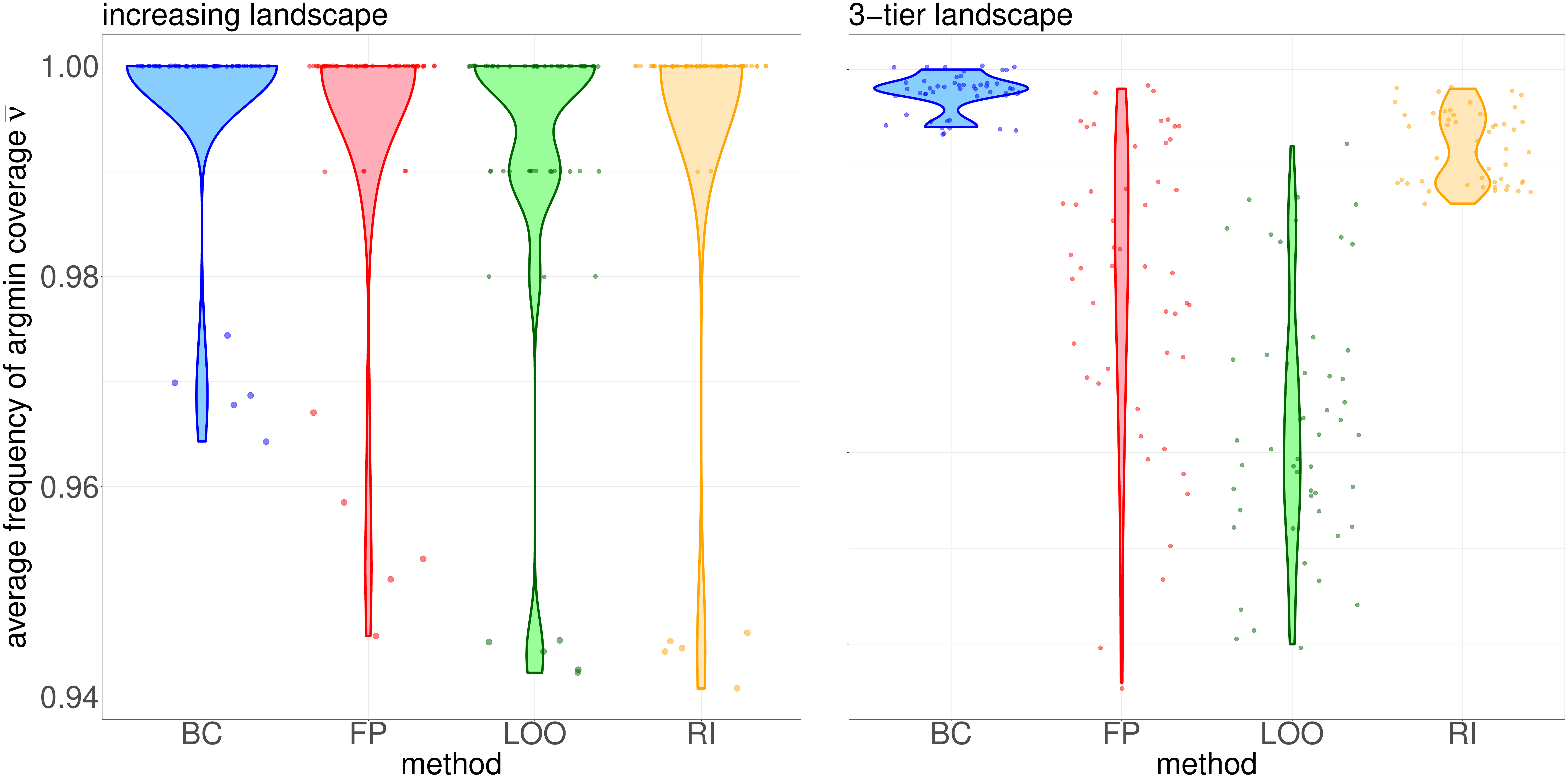}
    \caption{Method comparison, ``increasing" landscape (left) and ``3-tier" (right) landscapes. Comparison between the proposed 2-fold method and three other methods. 
    Each violin plot presents the distribution of the average frequencies of argmin coverage $\overline{\nu}$ over different simulation settings, and each jittered point corresponds to a value of $\overline{\nu}$ on the y-axis. The x-axis specifies the compared method: the proposed LOO method (LOO), the Bonferroni correction (BC), the method by~\cite{futschik1995confidence} (FP) and the rank inference method by~\cite{mogstad2024inference} (RI). A violin plot that deviates more from the nominal frequency $0.95$ indicates a more conservative method.}
    \label{fig: method comparison coverage}
\end{figure}

The same four methods---the proposed LOO method (LOO), the Bonferroni correction (BC), the selection rule by~\cite{futschik1995confidence} (FP) and the rank inference method by~\cite{mogstad2024inference} (RI)---are compared with respect to the \textit{average frequency of argmin coverage} $\overline{\nu}$. Formally, for all $r \in \Theta$ we record its frequency $\nu_r$ of inclusion in confidence sets over $100$ repetitions, and then $\overline{\nu}$ is computed by the average of $\nu_r$ over $r \in \Theta$. At the significance level of $\alpha = 0.05$, a method satisfying the validity~\eqref{eq: marginal validity} would result in $\overline{\nu}$ that is (approximately) greater than or equal to $0.95$, regardless of the mean factor $f$ and the dependency strength $\varrho$. In Figure~\ref{fig: method comparison coverage}, we display the distributions of $\overline{\nu}$ over the two types of mean landscapes using violin plots, with the specific values of $\overline{\nu}$ highlighted by jittered points.  We again present the results with dimension $p = 100$ and $n = 1000$.

In the left subplot, the true mean has the landscape $\mu = f \times \mu_b$ with $\mu_b = [0.01, 0.02, \ldots, 1]^\top$. So long as $f \ne 0$, the set of argmin would be a singleton, and the gap between the sub-optimal and optimal dimensions can lead to a negative mean shift in a procedure's test statistic when testing the latter. This inevitably makes the procedure conservative. Indeed, most values of $\overline{\nu}$ are close to $1$ for all the methods. However, when $f = 0$, the true mean is the flat landscape $\mu = [0, \ldots, 0]^\top$. In this scenario, we would expect the corresponding $\overline{\nu}$ to be near the nominal frequency 0.95. If it deviates significantly from the value, it suggests that the procedure is inherently conservative. The magnified jittered points in the left subplot represent the values of $\overline{\nu}$ under the flat landscape. We see all the compared methods produce values of $\overline{\nu}$ closer to 0.95, stressing the validity of each method. Notably, the values of $\overline{\nu}$ given by the Bonferroni correction appear inflated to some degree. This aligns with the common intuition that this multiple testing procedure is conservative because it relies on a union bound in its theoretical analysis.

In the right subplot, the violin plots are based on the 3-tier landscape $\mu = f \times \mu_b$ with $\mu_b = [1, 1, 1, 1, 1, 1.02, 1.02, $ $1.02, 1.02, 1.02, 10, \ldots, 10]^\top$. Despite the ties in the best tier, both the Bonferroni correction and the rank inference approach by~\cite{mogstad2024inference} appear to be quite conservative. This arises due to the existence of near ties in the second best tier. The true means of these sub-optimal dimensions are small enough that they are indistinguishable from the optimal dimensions using the methods. However, they are also sufficiently large to cause a negative mean shift in their test statistics, which leads to a conservative performance as explained earlier. By contrast, the proposed method leverages exponential weighting, as demonstrated in Theorem~\ref{thm:power_new}, to effectively detect the small signal level between the two tiers. It achieves this by eliminating the clearly inferior dimensions in the third tier. As a result, the proposed method essentially simplifies the mean landscape to the flat one. Indeed, the distribution of its resulting $\overline{\nu}$ gets closer to the nominal frequency $0.95$, which reveals a less conservative performance. As for the method by~\cite{futschik1995confidence}, it similarly shows a less conservative nature, but this is attributed to the screening-like step in its two-step selection procedure. In a way, the proposed method captures the benefit of their approach in a single step, showcasing the elegance of exponential weighting. 




\section{Real Data Applications under Differential Privacy}
\label{sec: differential private data}

The real data applications in Section~\ref{section: real_data} were performed on two binary loss matrices that evaluate student-trained classifiers over two datasets from \textit{Kaggle.com}. To ensure reproducibility  and facilitate research sharing while protecting student privacy, we anonymize the loss matrices and apply a differential privacy technique to generate differentially private versions. The differential privacy mechanism that we employ is the well-known \textit{binary randomized response} (see Section 2.3 in~\cite{dwork2014algorithmic}). Specifically, each binary loss value is flipped with probability $\pi \approx 0.047$. The choice of $\pi$ ensures the mechanism to achieve $3$-local differential privacy. The resulting differentially private loss matrices are publicly released on~\url{https://github.com/xu3cl4/argminCS}, and here we present the corresponding argmin confidence sets over these privatized data. 

\subsection{Analysis Results}

\begin{figure}[!tb]
    \centering
    \includegraphics[width =\textwidth]{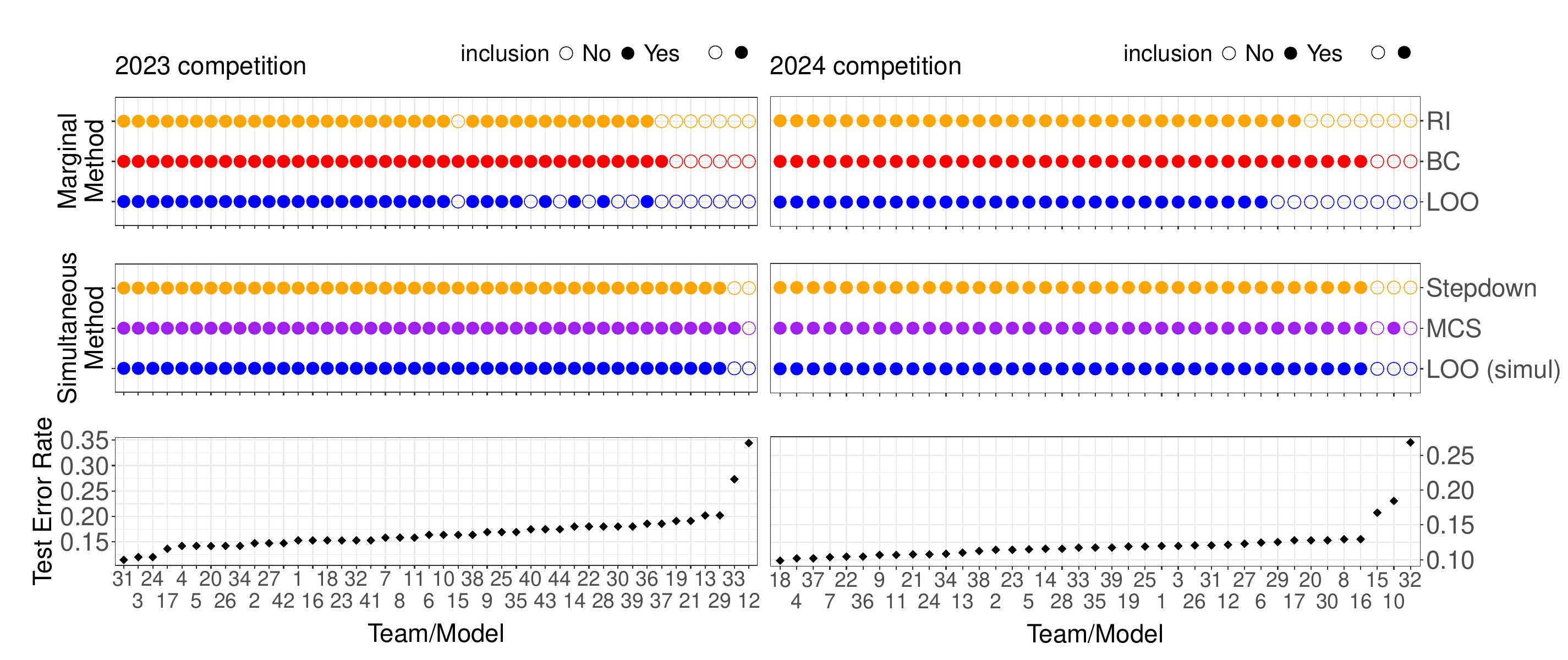}
    \vspace{-7mm}
    \caption{Confidence sets with differentially private data (both marginal and simultaneous coverage methods). We compare the proposed LOO algorithm (LOO), Bonferoni correction (BC), rank inference approach by~\cite{mogstad2024inference} (RI), the proposed heuristic two-step LOO algorithm (LOO (simul), Section~\ref{sec: heuristic simultaneous inference}), the stepdown procedure by~\cite{mogstad2024inference} (Stepdown), and the model confidence set (MCS) \citep{hansen2011model} over the differentially private test results in 2023 (left) and 2024 (right) classification competitions. The former three methods ensure marginal coverage, while the latter three aim for simultaneous coverage.}
    \label{fig: dp data applications}
\end{figure}

Figure~\ref{fig: dp data applications} presents the test error rates derived from the differentially private loss matrices. Due to the randomized flipping, the overall landscape of the test error rates, in terms of rate values and ordering, differs noticeably from the original in Figure~\ref{fig: data.application}. Following the framework in Section~\ref{section: real_data}, we construct $100$ argmin confidence sets for each method, respectively, to account for their internal randomness. Among the methods that ensure marginal coverage, the average cardinality of the confidence sets are $31.17$ (LOO), $38$ (BC, \textit{Bonferroni correction}) and $35.7$ (RI,~\cite{mogstad2024inference}) for the differentially private 2023 competition results, where the superiority of our proposed method is preserved. Moreover, the simultaneous coverage methods yield the average cardinalities of $41.79$ (LOO simul), $43$ (MCS,~\cite{hansen2011model}) and $42$ (Stepdown,~\cite{mogstad2024inference}). As before, simultaneous coverage methods often trade off power to favor its stronger coverage guarantee. Similar insights can be drawn from the argmin confidence sets for the differentially private 2024 competition results. Specifically, the marginal coverage methods have the average cardinalities of $30.02$ (LOO), $36$ (BC) and $31.98$ (RI), and the simultaneous coverage methods yield the average cardinalities of $36$ (LOO simul), $37$ (MCS) and $36$ (Stepdown). Figure~\ref{fig: dp data applications} presents one realization of the methods' resulting confidence sets for visualization. 

\section{Regarding reference \cite{kamath2015bounds}}\label{app: replicate results}
This unpublished technical report contains a straightforward treatment of the expectation of maximal statistics. It has received tens of citations at the time of drafting this manuscript. We found it intriguing to include results with sharp constants---as the author mentioned, this is a canonical question without clearly spelled-out answers in the literature. For our readers' convenience, we replicate the result and its proof below.

\begin{theorem}(Theorem~1 in \cite{kamath2015bounds})
Let $Y=\max _{1 \leq i \leq n} X_i$, where $X_i \sim \mathcal{N}\left(0, \sigma^2\right)$ are IID random variables. Then
$$
\frac{1}{\sqrt{\pi \log 2}} \sigma \sqrt{\log n} \leq \mathbb{E}[Y] \leq \sqrt{2} \sigma \sqrt{\log n} .
$$
\end{theorem}
\begin{proof}
    We comment that both constants which multiply $\sigma \sqrt{\log n}$ are tight. Indeed, as $n \rightarrow \infty, \mathbb{E}[Y] / \sqrt{\log n}$ converges to $\sqrt{2} \sigma$. On the other hand, by explicit calculations, one can verify the cases $n=1$ and $2$, for which $\mathbb{E}[Y]$ is 0 and $\sigma / \sqrt{\pi}$, respectively. In the former case, the inequality trivially holds for any multiplying constant, and in the latter, our inequality is tight.

    First, we show $\mathbb{E}[Y] \leq \sigma \sqrt{2} \sqrt{\log n}$. This result and method are folklore, but we include them here for completeness.
\begin{equation*}
\begin{aligned}
\exp (t \mathbb{E}[Y]) & \leq \mathbb{E}[\exp (t Y)] \\
& =\mathbb{E}\left[\max_i \exp \left(t X_i\right)\right] \\
& \leq \sum_{i=1}^n \mathbb{E}\left[\exp \left(t X_i\right)\right] \\
& =n \exp \left(t^2 \sigma^2 / 2\right)
\end{aligned}
\end{equation*}

The first inequality is Jensen's inequality, the second is the union bound, and the final equality follows from the definition of the moment generating function.
Taking the logarithm of both sides of this inequality, we get
\begin{equation*}
\mathbb{E}[Y] \leq \frac{\log n}{t}+\frac{t \sigma^2}{2}.
\end{equation*}

This can be minimized by setting $t=\sqrt{2 \log n}/\sigma$, which gives us the desired result
\begin{equation*}
\mathbb{E}[Y] \leq \sigma \sqrt{2} \sqrt{\log n}.
\end{equation*}

Next, we show the more difficult direction, the lower bound. We have already established that it holds for $n=1$ and 2. It can be verified for $n=3$ to 2834 using the Python 3 code provided in Section A (not replicated in this manuscript). Thus, for the remainder of the proof, we assume $n \geq 2835$.
Note that we have the following crude bound, which uses the Chernoff bound and the lower bound on $n$ :
\begin{equation*}
\begin{aligned}
\mathbb{E}[Y] & \geq \mathbb{P}\left(\left|\left\{i: X_i \geq 0\right\}\right| \geq\lceil n / 3\rceil\right) \cdot \mathbb{E}\left[\max _{1 \leq i \leq\lceil n / 3\rceil}\left|X_i\right|\right]+\\
& \quad \mathbb{P}\left(\left|\left\{i: X_i \geq 0\right\}\right|<\lceil n / 3\rceil\right) \cdot \mathbb{E}\left[-\left|X_i\right|\right] \\
& \geq 0.999 \cdot \mathbb{E}\left[\max _{1 \leq i \leq\lceil n / 3\rceil}\left|X_i\right|\right]-0.001 \sigma \cdot \sqrt{\frac{2}{\pi}} .
\end{aligned}
\end{equation*}


The second inequality uses the expected value of the half-normal distribution.

It remains to lower bound $\mathbb{E}\left[\max _{1 \leq i \leq k}\left|X_i\right|\right]$. We will show that
\begin{equation*}
\mathbb{P}\left(\left|X_i\right| \geq \sigma \sqrt{\log n}\right) \geq \frac{9}{n}.
\end{equation*}

This will imply the following lower bound:
\begin{equation*}
\begin{aligned}
\mathbb{E}\left[\max _{1 \leq i \leq\lceil n / 3\rceil}\left|X_i\right|\right] & \geq \sigma \sqrt{\log n} \cdot \mathbb{P}\left(\exists i:\left|X_i\right| \geq \sigma \sqrt{\log n}\right) \\
& \geq \sigma \sqrt{\log n} \cdot\left(1-\left(1-\frac{9}{n}\right)^{\lceil n / 3\rceil}\right) \\
& \geq\left(1-\frac{1}{e^2}\right) \sigma \sqrt{\log n}.
\end{aligned}
\end{equation*}

We compute the CDF of $\left|X_i\right|$ at the point $\sigma \sqrt{\log n}$.
\begin{equation*}
\begin{aligned}
\mathbb{P}\left(\left|X_i\right| \geq \sigma \sqrt{\log n}\right) & =1-\operatorname{erf}\left(\frac{\sqrt{\log n}}{\sqrt{2}}\right) \\
& \geq 1-\sqrt{1-\exp \left(-\frac{2}{\pi} \log n\right)} \\
& =1-\sqrt{1-n^{-\frac{2}{\pi}}}
\end{aligned}
\end{equation*}
where the first equality is based on the CDF of the half-normal distribution and the inequality is from the bound on the error function, $\operatorname{erf}(x) \leq \sqrt{1-\exp \left(-\frac{4}{\pi} x^2\right)}$. We require this value to be at least $\frac{9}{n}$ :
\begin{equation*}
\begin{aligned}
1-\sqrt{1-n^{-\frac{2}{\pi}}} & \geq \frac{9}{n} \\
1-\frac{9}{n} & \geq \sqrt{1-n^{-\frac{2}{\pi}}} \\
1-\frac{18}{n}+\frac{81}{n^2} & \geq 1-\frac{1}{n^{\frac{2}{\pi}}} \\
n^{2-\frac{2}{\pi}} & \geq 18 n-81 \\
\left(2-\frac{2}{\pi}\right) \log n & \geq \log (18 n-81) \\
\left(2-\frac{2}{\pi}\right) \frac{\log n}{\log (18 n-81)} & \geq 1
\end{aligned}
\end{equation*}

This inequality holds for all $n \geq 2835$, as desired. Putting these inequalities together, we have

\begin{equation*}
\mathbb{E}[Y] \geq 0.999\left(1-\frac{1}{e^2}\right) \sigma \sqrt{\log n}-0.001 \sigma \cdot \sqrt{\frac{2}{\pi}} \geq \frac{1}{\sqrt{\pi \log 2}} \sigma \sqrt{\log n}
\end{equation*}

where the second inequality holds for any integer $n>1$.
\end{proof}

\end{document}